\documentclass[11pt,letterpaper]{article}
\usepackage[letterpaper,margin=1in]{geometry}
\usepackage[T1]{fontenc}
\usepackage{lmodern}
\usepackage{amsmath,amssymb,amsthm,mathtools}
\usepackage{booktabs,array,enumitem,graphicx,xcolor,setspace,float}
\usepackage[labelfont=bf,labelsep=quad,font=small,skip=4pt]{caption}
\usepackage[round,authoryear]{natbib}
\setcitestyle{comma,aysep={}}%
\usepackage{tikz}
\usetikzlibrary{arrows.meta,positioning,shapes.geometric,fit,calc}
\usepackage[hidelinks]{hyperref}

\hypersetup{pdfauthor={Jiajun Qu, Yimo Yan, Yong-Hong KUO},pdftitle={Lattice-Cut Decomposition for Multistage Stochastic Integer Programs}}

\theoremstyle{plain}
\newtheorem{theorem}{Theorem}
\newtheorem{proposition}{Proposition}
\newtheorem{lemma}{Lemma}
\newtheorem{corollary}{Corollary}
\theoremstyle{definition}
\newtheorem{assumption}{Assumption}
\newtheorem{definition}{Definition}
\newtheorem{example}{Example}
\newtheorem{remark}{Remark}

\floatstyle{ruled}
\newfloat{algorithm}{tbp}{loa}
\floatname{algorithm}{Algorithm}

\newcommand{\Z}{\mathbb{Z}}
\newcommand{\R}{\mathbb{R}}
\newcommand{\LCD}{LCD}
\newcommand{\LCDiso}{\mbox{LCD-Iso}}

\newcommand{\figurenotes}[1]{\par\vspace{2pt}{\footnotesize\raggedright\emph{Notes.} #1\par}}
\newcommand{\figurenote}[1]{\par\vspace{2pt}{\footnotesize\raggedright\emph{Note.} #1\par}}
\newcommand{\coloronline}{(Color online) }
\newcounter{pcenv}
\newcounter{pcline}
\renewcommand{\thepcline}{\arabic{pcline}}

\newenvironment{pseudocode}{\stepcounter{pcenv}\setcounter{pcline}{0}\par\vspace{2pt}}{\par}
\newcommand{\pcl}[2][0]{\refstepcounter{pcline}\par\noindent
  \hangindent=\dimexpr3.5em+1.3em*#1\relax\hangafter=1\relax
  \makebox[2.2em][r]{\thepcline.\hspace{0.5em}}\hspace*{\dimexpr1.3em*#1\relax}#2\par}
\newcommand{\kw}[1]{\textbf{#1}}
\newif\ifboldbest \boldbestfalse
\newcommand{\tablebold}[1]{\ifboldbest\textbf{#1}\else#1\fi}
\newcommand{\boldnote}[1]{\ifboldbest#1\fi}



\begin{document}

\begingroup
\renewcommand{\thefootnote}{\fnsymbol{footnote}}
\begin{center}
{\LARGE\bfseries Lattice-Cut Decomposition for\\[1pt] Multistage Stochastic Integer Programs\par}
\vspace{10pt}
{\large Jiajun Qu\textsuperscript{1}\qquad
Yimo Yan\textsuperscript{2,}\footnote{Corresponding author.}\qquad
Yong-Hong KUO\textsuperscript{1}\par}
\vspace{6pt}
{\small \textsuperscript{1}Data and Systems Engineering, The University of Hong Kong\\
\textsuperscript{2}Department of Industrial Engineering, American University of Sharjah\\[3pt]
Contact: \href{mailto:jiajunqu25@connect.hku.hk}{\nolinkurl{jiajunqu25@connect.hku.hk}} (JQ);
\href{mailto:yyan@aus.edu}{\nolinkurl{yyan@aus.edu}} (YY);
\href{mailto:yhkuo@hku.hk}{\nolinkurl{yhkuo@hku.hk}} (YHK)\par}
\end{center}
\endgroup
\setcounter{footnote}{0}

\begin{singlespace}\small
\noindent\textbf{Abstract.}
We study cutting-plane methods for multistage stochastic integer programs in which orders, shipments, or crews come in batches. The cost-to-go function oscillates with the remainder of the state modulo the batch size, and no single affine cut can capture it. Existing remedies express every integer state in binary digits or attach a nonconvex correction to every cut. We propose lattice cuts, which have one slope and one intercept per remainder class, reproduce the periodic part of this function exactly, and add no binary variable per cut. The classes are computed from the constraint matrices, and closing cuts remove the nonconvexity that remains inside a class. We characterize the best bound of the two cut families by a linear program and prove finite convergence of the resulting algorithm, Lattice-Cut Decomposition. Computational results demonstrate that lattice cuts remove nearly all of the loss of affine cuts in exact tests and that, on horizons of up to 48 stages, the same algorithm certifies bounds significantly closer to the optimum with lattice cuts than without them on case-pack and co-production models. The results also show where the lattice does not help.

\end{singlespace}

\section{Introduction}\label{sec:intro}

We study cutting-plane methods for multistage stochastic integer programs (MSIPs) with integer state variables. Such models arise in planning problems that require a sequence of integer decisions under uncertainty: a retailer replenishes stock in case packs before demand is known \citep{wensing2018}, a hospital hires nurses in training cohorts, a plant runs recipes that yield products in fixed proportions \citep[see][for a co-production process with random yields]{bitran1994}, and a bike-sharing operator moves bikes by the truckload \citep[see][for models of bicycle redistribution]{shu2013}. Stochastic programs with integer decisions have also been applied, e.g., to the control of microgrids \citep{camisa2023} and to task offloading in fog computing \citep{zhang2021}. Because the scenario tree grows exponentially with the number of stages, MSIPs are solved by decomposition methods in which the problem is split by stage and the expected future cost of each stage, called the cost-to-go function, is approximated from below by cuts collected along sampled scenario paths \citep{pereira1991}. If all decisions are continuous, then the cost-to-go function is convex and affine cuts reproduce it exactly. With integer states, however, the cost-to-go function is not convex \citep{blair1982}, and affine cuts leave a gap that cannot be removed by adding more cuts.

Two lines of work have been proposed to close this gap. In the first, every integer state is represented by its binary digits \citep{zou2019}. Each state is then a vertex of a hypercube, at which Lagrangian cuts are tight. However, the dimension of the state is multiplied by the number of digits. In the second, the original state is kept and a nonconvex term, such as a reverse norm or a rectified linear function centered at the state at which the cut was generated, is added to each cut \citep{ahmed2022,zhang2022,deng2024,bansal2026}. Such a cut is tight at its center and loosens with the distance from it. Although these cuts avoid the binary expansion, the nonconvex term of each cut is represented in the stage model with auxiliary binary variables and, because the term is anchored at one state, tightness at another state requires another cut.

\begin{figure}[t]
\centering
\includegraphics[width=\textwidth]{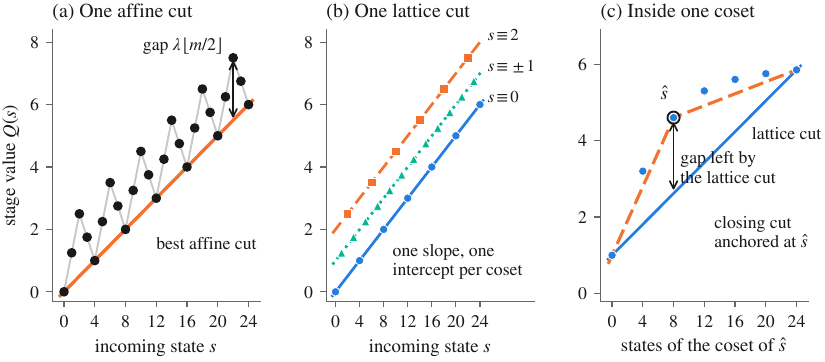}
\caption{\coloronline A Periodic Stage Value Function with the Best Affine Cut (a), a Lattice Cut (b), and a Closing Cut Inside One Remainder Class (c)}\label{fig:problem}
\figurenotes{Panels (a) and (b) show the stage value of the batch model of Example~\ref{ex:batch} with batch size $m=4$ and unit mismatch cost $\lambda=1$. In panel (a), the best affine cut underestimates the states with remainder 2, that is, every fourth state, by $\lambda\lfloor m/2\rfloor=2$, and the states with remainders 1 and 3 by $\lambda=1$. In panel (b), the lattice cut has one slope and one intercept for each remainder class and is exact. Panel (c) shows stylized values of the states of one class under a capacity limit, where the lattice cut leaves a gap at $\hat s$ that a closing cut anchored at $\hat s$ removes.}
\end{figure}

In this paper, we start from the observation that in many operational models the nonconvexity has a simple and regular source, namely, batches. When stock is ordered in packs of $m$ units, as in the batch-ordering model of \citet{veinott1965}, the cost of a stage depends on the remainder of the stock level after division by $m$, because a level that is a multiple of $m$ can be served without leftover, whereas a level halfway between two multiples cannot. The cost-to-go function is then the sum of a trend and a pattern that repeats with period $m$. Figure~\ref{fig:problem}(a) shows such a function together with the best affine cut, which passes through the states with remainder zero and underestimates the states in between by an amount that no choice of slope can reduce. Figure~\ref{fig:problem}(b) shows the way out: states with the same remainder lie on a common line, and the lines of different remainders are parallel.

We turn this observation into a family of cuts. A \emph{lattice cut} has one slope and a separate intercept for each remainder class. In several dimensions, the remainder classes are the cosets of a lattice, that is, of the set of state changes that leave all batch remainders unchanged, which gives the cut its name. A lattice cut is global in the sense that it is valid and, in the situation of Figure~\ref{fig:problem}(b), tight on entire cosets, and not only at the state at which it was generated. A lattice cut is also inexpensive to carry. Because the stage model identifies the coset of the next state with one set of indicator variables that all cuts share, a new cut adds a row but no variable. Batches are not the only source of nonconvexity, however, because bounds, capacities, and fixed charges make the cost-to-go function nonconvex within a coset (Figure~\ref{fig:problem}(c)). For this remaining part, we use \emph{closing cuts}, which subtract a penalty for the distance from an anchor state. When the periodic part accounts for most of the nonconvexity, closing cuts are needed at few states and with small penalties, because the lattice cut has already removed this part. We refer to the resulting algorithm as Lattice-Cut Decomposition (\LCD).

\subsection{Related Literature}\label{sec:literature}

Our work builds on three streams of research: cutting planes for stochastic programs with integer recourse, decomposition methods for multistage problems with integer states, and periodic structure in integer programming.

Cuts for stochastic programs originate in the partitioning procedure of \citet{benders1962} and in the L-shaped method of \citet{vanslyke1969}, which \citet{birge1985} extends to multistage problems in nested form. These methods rely on a convex recourse function. With integer recourse, however, the value function loses convexity \citep{blair1982}, and three approaches have been proposed for two-stage problems \citep[see][for an introduction]{kucukyavuz2017}. The first approach uses cuts that are valid for the nonconvex recourse function itself. \citet{laporte1993} propose integer L-shaped cuts, which are tight at binary first-stage points, and \citet{angulo2016} improve this method by alternating between linear and mixed-integer subproblems. \citet{caroe1998} generalize the L-shaped method to integer recourse with nonlinear feasibility and optimality cuts that are obtained from general duality theory. A limitation of this approach is that integer L-shaped cuts require binary first-stage variables, whereas the cuts of \citet{caroe1998} lead to nonlinear master problems. The second approach strengthens the linear relaxation with cutting planes from integer programming. \citet{sen2005} apply disjunctive programming to convexify the second-stage problems, and \citet{gade2014} and \citet{zhangkucuk2014} use Gomory cuts that are parameterized by the first-stage variables. For problems with integer first-stage variables and continuous recourse, \citet{bodur2017} compare split cuts that are derived in the extended variable space with split cuts that are derived directly in the Benders reformulation, and \citet{bodurluedtke2017} use mixed-integer rounding to improve the Benders cuts of a staffing and scheduling model. The third approach is based on the Lagrangian dual in which the nonanticipativity constraints are relaxed. \citet{caroe1999} introduce dual decomposition, in which this dual provides the bounds of a branch-and-bound algorithm, and \citet{lubin2013} revisit this algorithm with the aim of its parallelization. \citet{rahmaniani2020} combine dual decomposition with Benders decomposition, and \citet{chen2022} investigate techniques that generate Lagrangian cuts quickly. A limitation of this approach is that, for general integer first-stage variables, the Lagrangian bound is the convex envelope of the value function \citep{geoffrion1974}, and the gap between the envelope and the function itself is precisely what a periodic value function makes large.

For multistage problems, one approach works directly with the scenario tree. \citet{guan2009} generate cutting planes for the whole tree by combining inequalities that are valid for the individual scenarios, and \citet{lulli2004} develop a branch-and-price algorithm and validate it on stochastic batch sizing, which is the application in which periodicity is most visible. A limitation of the scenario-tree approach is that the size of the tree grows exponentially with the number of stages. \citet{daryalal2024} do not rely on a scenario tree. They restrict the multipliers of a Lagrangian dual of the problem to follow decision rules, which yields bounds on the optimal value and primal policies.

A second approach decomposes the problem by stage. Stochastic dual dynamic programming (SDDP), proposed by \citet{pereira1991}, approximates the cost-to-go functions by cuts that are collected along sampled scenario paths. \citet{philpott2008} prove the almost-sure convergence of a broad class of such sampling algorithms for multistage stochastic linear programs, and \citet{shapiro2011} analyzes the statistical properties and the convergence of SDDP under stagewise independence. \citet{zou2019} extend SDDP to MSIPs with binary state variables. In their method, which is called stochastic dual dynamic integer programming (SDDiP), Lagrangian cuts are tight at binary states, and general integer states are handled by binary expansion. Several authors build on SDDiP. \citet{quezada2022} extend it for stochastic uncapacitated lot sizing by a partial decomposition of the scenario tree and by cutting planes that strengthen the linear relaxations of the subproblems. \citet{bansal2024} compare Benders, integer L-shaped, and Lagrangian cuts for MSIPs and show that integer L-shaped cuts belong to the class of Lagrangian cuts. \citet{fullner2025} study Lipschitz regularization, and \citet{fullner2026} generate Lagrangian cuts by solving normalized dual problems. However, for general integer states, SDDiP relies on the binary expansion, under which the state vector grows with the number of digits.

Subsequent work retains the original states and strengthens the cut instead. \citet{philpott2020} approximate monotonic cost-to-go functions by step functions. \citet{ahmed2022} propose stochastic Lipschitz dynamic programming, in which reverse-norm cuts and augmented Lagrangian cuts approximate the cost-to-go functions from below. The tightness of augmented Lagrangian cuts rests on the result of \citet{feizollahi2017} that the augmented Lagrangian dual of a mixed-integer program closes the duality gap with a finite penalty when the augmenting function is a norm. \citet{zhang2022} regularize the value functions and derive cuts from generalized conjugacy, and \citet{fullner2022} refine a binary approximation of the states dynamically. \citet{deng2024} introduce rectified linear cuts and show that, for bounded integer variables, these cuts can be represented with a number of additional binary variables that depends on the ranges of the variables. \citet{bansal2026,bansal2026b} normalize rectified linear cuts and combine them with linear programming cuts in a convergent algorithm. In all of these methods, the nonconvex part of a cut is anchored at the state at which the cut was generated. Although the cuts are valid on the whole state space, none of them contains a term indexed by the remainder class of the state, and therefore tightness at one state need not carry over to the states that share its remainders. Our work differs from the above in the following respects: we index the nonconvex part of a cut by the coset of the state, so that one cut can be tight on entire cosets; we identify the coset with one set of indicator variables that all cuts of a stage share, so that a new cut adds a row but no variable; and we use anchored penalty terms only for the nonconvexity that remains inside a coset.

Our work is also related to research on the structure of value functions of integer programs. The observation that these value functions repeat along a lattice is classical in integer programming \citep[see][for the theory of the value function]{blair1982}, and it has been used in stochastic programming before. \citet{trapp2013} characterize the value function of a pure integer program by its level sets and use this characterization to solve two-stage stochastic programs through a value-function reformulation. For two-stage models with integer recourse, \citet{romeijnders2016tv} bound the expectation of periodic functions through the total variation of the density of the random right-hand side, and \citet{romeijnders2016convex} derive from these bounds a convex approximation of mixed-integer recourse models with a uniform error bound. \citet{vanderlaan2024} give a Benders decomposition for such models whose cuts recover the convex envelope of the expected recourse function, and \citet{romeijnders2024} extend these scaled cuts to multistage problems, for which the first-stage lower bound converges to the convex envelope of the first-stage expected cost-to-go function. In this literature, the periodic part of the recourse function is averaged out by a continuous distribution or removed by convexification, whereas we keep it and represent it exactly, coset by coset. Most recently, \citet{er2025} give a common constraint matrix for the integer hull of a polyhedron with a varying right-hand side such that, within each residue class of a suitable lattice, the right-hand side of the integer-hull description is an affine function of the original right-hand side. Although this result motivates our approach, its modulus is universal for the constraint matrix and can be very large. We instead take the lattice from the batch structure of the model through the Smith normal form of its balance matrices \citep{newman1972}, which yields few cosets, and we leave what the lattice does not capture to closing cuts. To the best of our knowledge, cuts indexed by the cosets of a lattice have not previously been used in stagewise decomposition.

\subsection{Our Contributions}\label{sec:contributions}

The main contributions of our work are summarized as follows.
\begin{enumerate}[leftmargin=*,itemsep=2pt,topsep=3pt]
\item We introduce lattice cuts and show how to choose the lattice. We compute, from the constraint matrices of the stages and by Smith normal form, the coarsest family of lattices that is consistent with the state transitions. We also show that, for a class of batch-balance models, the computed lattice is the period lattice of the periodic part of the cost-to-go function, so that no coarser lattice can represent it.
\item We characterize the bound that the two cut families can achieve. We show that, for a given lattice and a given cap on the penalties of the closing cuts, the best bound obtainable at the root is the value of a linear program. This result identifies the smallest cap at which the bound equals the optimal value, and it splits the gap of any run into a representation loss, which depends only on the lattice and the cap, and a finite-pool loss, which vanishes as cuts accumulate.
\item We make our method practical. In our formulation, all closing cuts of a stage share one binary variable per state threshold, so that the number of binary variables in a stage model is bounded by the ranges of the state variables and does not grow with the number of cuts. We also prove that \LCD\ converges to the optimal value in finitely many iterations with probability one.
\item We conduct a computational study on 15 problem families, in which we compare \LCD\ with stochastic dual dynamic integer programming and with a recent method based on rectified linear cuts on exact reference values, on horizons of up to 48 stages, and under equal time limits. We isolate the effect of the lattice by running our algorithm with and without it, and we report the cases in which the lattice does not help.
\end{enumerate}

The remainder of this paper is organized as follows. In Section~\ref{sec:model}, we state the problem and show why affine cuts leave a gap. In Section~\ref{sec:lattice-cuts}, we introduce lattice cuts and discuss the choice of the lattice. In Section~\ref{sec:cuts}, we add closing cuts and characterize the attainable bound. In Section~\ref{sec:algorithm}, we present our algorithm. In Section~\ref{sec:experiments}, we report the results of our computational study, and we conclude in Section~\ref{sec:conclusion}. Proofs that are not given in the main text are provided in Appendix~\ref{ec:proofs}, and the remaining appendices contain additional results and the details of the computational study.
\section{Problem Setting}\label{sec:model}

In this section, we state the model and the cuts that we consider (Section~\ref{sec:model-cuts}), and we show that, under a batch balance, every affine cut loses a fixed amount at states halfway between two batch multiples (Section~\ref{sec:sep}).

\subsection{Model and Cuts}\label{sec:model-cuts}

Let $T$ denote the number of stages. At the start of stage $t$, the system is in an integer state $s\in D_t\subseteq[\ell_t,u_t]\cap\Z^n$, where $\ell_t\le u_t$ are finite integer bounds, and the initial state $s_1$ is given. An outcome $w$ is then observed from a finite set $\Omega_t$ with probability $p_{tw}>0$, after which the decision maker chooses an action $a$, which may have integer and continuous components, and the next state $y\in D_{t+1}$. The pair $(a,y)$ must lie in a set $X_{tw}(s)$ that is described by mixed-integer linear constraints, and it incurs a linear cost $c_{tw}(s,a,y)$. After stage $T$, a terminal cost $V_{T+1}(y)$ is paid. Throughout this paper, we make the following assumption.

\begin{assumption}\label{assump:msip}
(a)~Each set $D_t$, $t=1,\ldots,T+1$, is finite and nonempty. (b)~The outcomes of different stages are independent. (c)~For all $t$, $w\in\Omega_t$, and $s\in D_t$, the set $X_{tw}(s)$ is nonempty, the terminal cost $V_{T+1}$ has a mixed-integer linear epigraph on $D_{T+1}$, and every stage problem attains a finite minimum.
\end{assumption}

Assumption~\ref{assump:msip}(a) reflects the finite bounds on the integer states, Assumption~\ref{assump:msip}(b) states that the outcomes are stagewise independent, and Assumption~\ref{assump:msip}(c) ensures that every stage problem is feasible and has a finite optimal value at every state of its stage. Stagewise independence is standard for sampling-based decomposition methods \citep{pereira1991,zou2019}. In some cases, dependence across stages can be dealt with by adding state variables to the model \citep{shapiro2011}. Under Assumption~\ref{assump:msip}, the cost-to-go functions $V_t$ satisfy, for $t=T,\ldots,1$,
\begin{equation}\label{eq:bellman}
 Q_{tw}(s)=\min_{(a,y)\in X_{tw}(s)}\bigl\{c_{tw}(s,a,y)+V_{t+1}(y)\bigr\},\qquad
 V_t(s)=\sum_{w\in\Omega_t}p_{tw}\,Q_{tw}(s),
\end{equation}
and the optimal value of the problem is $V^\star=V_1(s_1)$. Because a decomposition method does not know $V_{t+1}$, it works with a lower approximation $\underline V_{t+1}\le V_{t+1}$ that is built from cuts, and it solves the \emph{stage problem}
\begin{equation}\label{eq:current}
 F_{tw}(s)=\min_{(a,y)\in X_{tw}(s)}\bigl\{c_{tw}(s,a,y)+\underline V_{t+1}(y)\bigr\}\ \le\ Q_{tw}(s).
\end{equation}
A \emph{cut} for stage $t$ and outcome $w$ is a function $h$ with $h\le F_{tw}$ on $D_t$, and it is \emph{tight} at a state $\hat s$ if $h(\hat s)=F_{tw}(\hat s)$. A tight cut passes the exact value of $\hat s$ to stage $t-1$, where it enters the approximation $\underline V_t$. The classical choice is an affine cut $h(s)=\pi^\top s+\beta$. Throughout this paper, we keep the state in its original integer coordinates, that is, we do not expand it into binary digits.

\subsection{The Gap Left by Affine Cuts}\label{sec:sep}

We next show that affine cuts necessarily leave a gap when a stage contains a batch balance. Consider a stage in which a quantity $Cs$, which is an integer linear function of the state, has to be settled in batches of $m\ge 2$ units. The number of batches $z$ is an integer, and a surplus $e^+$ or a shortage $e^-$ costs $\lambda>0$ per unit, that is,
\[
 mz+e^+-e^-=Cs,\qquad z\in\Z,\quad e^+,e^-\in\Z_+ .
\]
Choosing $z$ optimally leaves a mismatch cost of $\lambda\operatorname{dist}(Cs,m\Z)$, where $\operatorname{dist}(r,m\Z)$ denotes the distance from $r$ to the nearest multiple of $m$. Suppose that the remaining actions $a\in\mathcal U$ cost $k(a)+\ell^\top s$ and lead to the next state $y=As+Ga+d_0$, which is valued by a function $W$. The stage value is then
\begin{equation}\label{eq:modfamily}
 F(s)=\ell^\top s+\lambda\operatorname{dist}(Cs,m\Z)
       +\min_{a\in\mathcal U}\bigl\{k(a)+W(As+Ga+d_0)\bigr\}.
\end{equation}
The following proposition shows that every affine cut for $F$ loses a fixed amount at a state that lies halfway between two batch multiples.

\begin{proposition}\label{prop:loss}
Let $d\in\Z^n$ satisfy $Ad=0$ and $Cd\equiv 1\pmod m$, let $q=\lfloor m/2\rfloor$, and let $s^0$ be a state with $Cs^0\equiv 0\pmod m$. If $s^0$, $s^0+qd$, and $s^0+md$ belong to $D_t$, then every affine function $u\le F$ on $D_t$ satisfies $F(s^0+qd)-u(s^0+qd)\ge\lambda q$.
\end{proposition}
\begin{proof}
Because $Ad=0$, the three states have the same set of reachable next states, and therefore the minimum in \eqref{eq:modfamily} takes the same value $M$ at all three. The mismatch cost is zero at $s^0$ and $s^0+md$ and equals $\lambda q$ at $s^0+qd$. Because $u$ is affine, we have $u(s^0+qd)=(1-q/m)\,u(s^0)+(q/m)\,u(s^0+md)\le(1-q/m)F(s^0)+(q/m)F(s^0+md)=\ell^\top(s^0+qd)+M=F(s^0+qd)-\lambda q$.
\end{proof}

By Proposition~\ref{prop:loss}, the loss $\lambda q$ does not depend on the slope of the cut and is not reduced by adding affine cuts. The loss is caused by averaging over states with different remainders. If only states with the same remainder as $s^0+qd$ are averaged, then the mismatch cost is the same constant $\lambda q$ at all of them, and it can be added to the intercept. Proposition~\ref{prop:loss} requires a direction $d$ with $Ad=0$, which excludes models in which the stock carries over to the next stage. In Appendix~\ref{app:add-theory}, we verify numerically on a Case-Pack instance that affine cuts lose a fixed amount at states between two batch multiples in such a model as well. We illustrate the observation of this section in the following example.

\begin{example}[Batch Model]\label{ex:batch}
Let the state be a single integer $s$, let $C=1$, and let the value be $F(s)=cs+\lambda\operatorname{dist}(s,m\Z)$, which is \eqref{eq:modfamily} with $W=0$ and no further actions. Figure~\ref{fig:problem}(a) plots $F$ for $m=4$. The affine function $cs$ is the best affine cut, and it underestimates the states with remainder $2$ by $2\lambda$. The states with remainder $r$ lie on the line $cs+\lambda\min\{r,m-r\}$. The function that takes the value of this line at every state with remainder $r$ has one slope, $c$, and $m$ intercepts, and it equals $F$ at every state (Figure~\ref{fig:problem}(b)).
\end{example}

\section{Lattice Cuts}\label{sec:lattice-cuts}

In this section, we first define lattice cuts and describe how they are computed (Section~\ref{sec:lattice-def}). We then show how the lattice is obtained from the constraint matrices of the stages (Section~\ref{sec:lattice}), and we give conditions under which lattice cuts alone reproduce the cost-to-go function (Section~\ref{sec:periods}).

\subsection{Definition and Computation}\label{sec:lattice-def}

A \emph{lattice} $\Gamma\subseteq\Z^n$ is a set of integer vectors that is closed under addition and subtraction. We use lattices of full rank, for which the number of cosets $s+\Gamma$, called the \emph{index} $[\Z^n:\Gamma]$, is finite. In Example~\ref{ex:batch}, we have $\Gamma=m\Z$ and the cosets are the $m$ remainder classes. Let $\Gamma_t$ denote the lattice used for the states of stage $t=1,\ldots,T$, the choice of which is the subject of Section~\ref{sec:lattice}. We set $\Gamma_{T+1}=\Z^n$, so that $D_{T+1}$ is a single coset, because no cuts are generated for the states that follow the last stage. A terminal cost with a batch structure of its own can be modeled as an additional stage. We write $[s]_t=s+\Gamma_t$ for the coset of $s$, $\mathcal G_t$ for the set of cosets that meet $D_t$, and $D_g$ for the states of $D_t$ in coset $g$.

We fix a stage $t\ge 2$ and an outcome $w$, and we write $F=F_{tw}$. For a slope $\pi\in\R^n$ and a coset $g$, the largest intercept that keeps the affine function $\pi^\top s+\beta$ below $F$ on $D_g$ is
\begin{equation}\label{eq:coset-intercept}
 \beta_g(\pi)=\min_{z\in D_g}\bigl\{F(z)-\pi^\top z\bigr\}.
\end{equation}
A \emph{lattice cut} with slope $\pi$ is the function $h(s)=\pi^\top s+\beta_{[s]_t}(\pi)$. This cut is valid because $z=s$ is feasible in \eqref{eq:coset-intercept} for the coset of $s$. It is tight at a state $\hat s$ whenever $\hat s$ attains the minimum in \eqref{eq:coset-intercept} for its own coset, that is, whenever $F$ has a supporting affine function at $\hat s$ relative to the states of that coset. An affine cut is the special case $\Gamma_t=\Z^n$, in which there is a single coset.

Three observations make lattice cuts practical. First, the intercept \eqref{eq:coset-intercept} is the optimal value of one mixed-integer program (MIP), namely, the stage problem in which the incoming state is a variable $z$ that is restricted to $D_g$ by congruence constraints. Therefore, no enumeration of $D_g$ is involved. Second, a good slope is available from the linear programming (LP) relaxation of the stage problem at $\hat s$, and this slope can be improved by a subgradient method because $\pi\mapsto\pi^\top\hat s+\beta_{[\hat s]_t}(\pi)$ is concave and every minimizer $z_\pi$ of \eqref{eq:coset-intercept} gives the supergradient $\hat s-z_\pi$. Third, the intercepts of the other cosets need not be computed at once, because a single lower bound $b_0\le\min_g\beta_g(\pi)$ serves all cosets until the algorithm asks for more.

In the model of stage $t-1$, the cut is imposed on the next state $y$. The model contains one binary indicator $\delta_g$ for each coset $g\in\mathcal G_t$, with $\sum_g\delta_g=1$, and the indicators are linked to $y$ by the congruence constraints that define the cosets. With these indicators, the cut is the linear constraint
\begin{equation}\label{eq:lc-row}
 \theta_w\ \ge\ \pi^\top y+\sum_{g\in\mathcal G_t}\beta_g(\pi)\,\delta_g ,
\end{equation}
where $\theta_w$ represents the value of outcome $w$. Because the indicators belong to the stage and not to the cut, the number of binary variables in the model stays the same however many lattice cuts are added.

\subsection{Choice of the Lattice}\label{sec:lattice}

The lattice has to meet two requirements. It must be fine enough to separate the batch remainders of its stage, and the lattices of consecutive stages must fit together, so that the coset of the next state is determined by the coset of the current state and the action.

For the first requirement, suppose that stage $t$ contains a batch balance $Bz=b+Es-Wa$, where $z\in\Z^q$ counts batches, $B\in\Z^{q\times q}$ is nonsingular, $E\in\Z^{q\times n}$, and $b-Wa$ is integral for every action $a$. A change $v$ of the state leaves all batch remainders unchanged if and only if $Ev\in B\Z^q$. The Smith normal form $UBR=\operatorname{diag}(d_1,\ldots,d_q)$, where $U$ and $R$ are unimodular \citep{newman1972}, turns this condition into the congruences $(UEv)_i\equiv 0\pmod{d_i}$. Collecting the congruences of all batch balances of stage $t$, we obtain the \emph{stage lattice} $\Lambda_t$. In Example~\ref{ex:batch}, we have $B=m$, $E=1$, and $\Lambda_t=m\Z$. For a recipe that turns each batch into the product quantities given by the columns of a matrix $Y$, the stage lattice is $Y\Z^2$, whose cosets are not aligned with the coordinate axes. If a balance has more rows than batch variables, so that $B$ is not square, then the Smith normal form has rows without an invariant factor, and each such row requires an integer linear form $c^\top v$ of the state change to vanish instead of a congruence. On the bounded set $D_t$, the form $c^\top s$ takes values in an interval of some length $M-1$, and therefore, for differences of states in $D_t$, the equation $c^\top v=0$ is equivalent to the congruence $c^\top v\equiv 0\pmod M$. We use this congruence, so that the stage lattice keeps full rank. Among our benchmark families, this case occurs only in Bike-Sharing (Appendix~\ref{app:benchmarks}).

For the second requirement, suppose that the next state depends on the current state through $y=A_{tw}s+G_{tw}a+d_{tw}$, where $A_{tw}$ is an integer matrix.

\begin{definition}\label{def:compatible}
Lattices $\Gamma_1,\ldots,\Gamma_T$ of full rank are \emph{compatible} if $\Gamma_t\subseteq\Lambda_t$ for all $t$ and $A_{tw}\Gamma_t\subseteq\Gamma_{t+1}$ for all $t<T$ and $w\in\Omega_t$.
\end{definition}

Compatibility means that two states in the same coset at stage $t$ have the same batch remainders now and, under the same action, reach the same coset at stage $t+1$. The following proposition shows that, among all compatible families, one is coarsest and therefore has the fewest cosets.

\begin{proposition}\label{prop:compatible}
The family defined by
\begin{equation}\label{eq:compatible}
\Gamma_T=\Lambda_T,\qquad
\Gamma_t=\Lambda_t\cap\bigcap_{w\in\Omega_t} A_{tw}^{-1}\Gamma_{t+1}
\quad (t=T-1,\ldots,1)
\end{equation}
is compatible, and it contains every other compatible family stage by stage. In particular, it has the smallest index $[\Z^n:\Gamma_t]$ at every stage.
\end{proposition}

We compute the family \eqref{eq:compatible} with Algorithm~\ref{alg:lattice}. The algorithm requires only the matrices of the model and is run once, before any stage problem is solved. For each stage, it returns a unimodular matrix $V_t$ and moduli $q_{t,1},\ldots,q_{t,r_t}$ such that the coset of a state $s$ is identified by the label $\bigl((V_t^{-1}s)_i \bmod q_{t,i}\bigr)_{i=1}^{r_t}$. These labels provide the congruence constraints used in \eqref{eq:coset-intercept} and \eqref{eq:lc-row}. Although the intersections in \eqref{eq:compatible} could refine the lattice at every stage, so that the index would grow with the horizon, this does not happen for a stationary transition matrix and a common modulus, for which the recursion stops changing after at most $n$ steps (Lemma~\ref{lem:observation-phases}(iv)) and the index is bounded independently of $T$.

\begin{algorithm}[t]
\caption{(Compatible Lattices from the Stage Matrices)}\label{alg:lattice}
\begin{singlespace}\small
\textbf{Input:} batch balances $(B,E)$ of every stage; transition matrices $A_{tw}$.\par
\textbf{Output:} $V_t$ and $q_{t,1},\ldots,q_{t,r_t}$ for $t=1,\ldots,T$.
\begin{pseudocode}
\pcl{\kw{for} $t=1,\ldots,T$ \kw{do}}
\pcl[1]{\kw{for} every batch balance $(B,E)$ of stage $t$ \kw{do}}
\pcl[2]{\label{ln:lat-snf}Compute the Smith normal form $UBR=\operatorname{diag}(d_1,\ldots,d_q)$}
\pcl[2]{Record the congruence $(UE)_i\,s\equiv 0\pmod{d_i}$ for each $d_i>1$}
\pcl[1]{\kw{end}}
\pcl[1]{$\Lambda_t\leftarrow$ the lattice defined by the recorded congruences}
\pcl{\kw{end}}
\pcl{\label{ln:lat-last}$\Gamma_T\leftarrow\Lambda_T$}
\pcl{\label{ln:lat-back}\kw{for} $t=T-1,\ldots,1$ \kw{do}}
\pcl[1]{Start from the congruences of $\Lambda_t$}
\pcl[1]{\kw{for} every $w\in\Omega_t$ and every congruence $c^\top y\equiv 0\pmod q$ of $\Gamma_{t+1}$ \kw{do}}
\pcl[2]{\kw{if} $c^\top A_{tw}\,s\equiv 0\pmod q$ does not follow from the congruences already present \kw{then}}
\pcl[3]{Add the congruence $c^\top A_{tw}\,s\equiv 0\pmod q$}
\pcl[2]{\kw{end}}
\pcl[1]{\kw{end}}
\pcl[1]{$\Gamma_t\leftarrow$ the lattice defined by the resulting congruences}
\pcl{\label{ln:lat-backend}\kw{end}}
\pcl{\kw{for} $t=1,\ldots,T$ \kw{do}}
\pcl[1]{\label{ln:lat-label}Bring the congruences of $\Gamma_t$ to a common modulus and compute the Smith normal form of their coefficient matrix, which yields $V_t$ and the moduli $q_{t,i}>1$ of the coset label}
\pcl{\kw{end}}
\end{pseudocode}
\end{singlespace}
\end{algorithm}

\subsection{Exactness of Lattice Cuts}\label{sec:periods}

We now give a condition under which lattice cuts alone, without any further mechanism, reproduce the cost-to-go function. For a slope vector $\alpha\in\R^n$ and a coset $g'\in\mathcal G_{t+1}$, let
\[
 K^\alpha_{tw}(s,g')=\min\bigl\{c_{tw}(s,a,y)+\alpha^\top y:\ (a,y)\in X_{tw}(s),\ [y]_{t+1}=g'\bigr\}
\]
denote the cost of the cheapest way to move from state $s$ into coset $g'$ when the next state is valued at $\alpha^\top y$, with the value $+\infty$ if $g'$ cannot be reached. The following proposition characterizes the stages that pass a continuation value of the form of a linear function of the state plus a function of its coset on to the stage value, and it shows that, when all stages have this property, one lattice cut per stage and outcome is exact.

\begin{proposition}\label{prop:phasevalue}
Fix $t$, $w$, $\alpha$, and $\beta\in\R^n$. The following two statements are equivalent. (a)~For every function $\phi:\mathcal G_{t+1}\to\R$, the stage value under the continuation value $\alpha^\top y+\phi([y]_{t+1})$ has the form $\beta^\top s+\psi([s]_t)$ on $D_t$ for some function $\psi$ of the coset. (b)~For every $g'\in\mathcal G_{t+1}$, the difference $K^\alpha_{tw}(s,g')-\beta^\top s$ depends on $s$ only through the coset $[s]_t$. Consequently, if $V_{T+1}(y)=\alpha_{T+1}^\top y+\phi_{T+1}([y]_{T+1})$ and (b) holds at every $(t,w)$ with $\alpha=\alpha_{t+1}$ and some $\beta_{tw}$, where $\alpha_t=\sum_{w}p_{tw}\beta_{tw}$, then $V_t(s)=\alpha_t^\top s+\phi_t([s]_t)$ for all $t$, and the lattice cut of $Q_{tw}$ with slope $\beta_{tw}$ equals $Q_{tw}$ on all of $D_t$.
\end{proposition}

Condition (b) holds, for example, when the transitions are linear as above, the lattices are compatible, the set of feasible actions depends on the state only through its coset, and the stage cost is a linear function of the state plus a function of the coset and the action. In this case, we have $\beta_{tw}=\ell_{tw}+A_{tw}^\top\alpha_{t+1}$, where $\ell_{tw}$ is the linear cost coefficient of the state, the cost-to-go function is a linear trend plus a function of the coset, as in Figure~\ref{fig:problem}(b), and one lattice cut per stage and outcome is exact everywhere. Among the benchmark families of Section~\ref{sec:experiments}, Shift and Rotate are of this kind by construction because their state bounds never bind (Appendix~\ref{app:benchmarks}). In the other families, bounds and capacities make the set of feasible actions depend on the state itself, and the gap that lattice cuts then leave is measured by the bound of Section~\ref{sec:ideal-pool}.

We next show that, for a class of batch-balance models, the lattice of Algorithm~\ref{alg:lattice} cannot be replaced by a coarser one. To this end, we subtract the linear trend and consider the remaining periodic part $R_t$ of the cost-to-go function, which is defined on $\Z^n$ by $R_{T+1}=0$ and
\begin{equation}\label{eq:period-model}
 R_t(x)=\sum_{w\in\Omega_t} p_{tw}R_{tw}(x),\qquad
 R_{tw}(x)=\min_{a\in U_{tw}(C_{tw}x)}
 \bigl\{r_{tw}(C_{tw}x,a)+R_{t+1}(A_{tw}x+b_{tw}(a))\bigr\}.
\end{equation}
Here, $C_{tw}x$ collects the batch remainders of state $x$ (formally, $C_{tw}$ is a homomorphism from $\Z^n$ into a finite abelian group), $U_{tw}(g)$ is a finite nonempty set of actions available at remainder $g$, $r_{tw}\ge 0$ is the cost, and $b_{tw}(a)\in\Z^n$ is the shift of the state caused by action $a$. Let
\begin{equation}\label{eq:observer-kernel}
 K_{T+1}=\Z^n,\qquad
 K_t=\bigcap_{w\in\Omega_t}\bigl(\ker C_{tw}\cap A_{tw}^{-1}K_{t+1}\bigr)
 \quad (t=T,\ldots,1),
\end{equation}
which is the recursion \eqref{eq:compatible} with $\Lambda_t=\bigcap_{w}\ker C_{tw}$. We make the following assumption on this model.

\begin{assumption}\label{assump:period}
For all $t$ and $w$, the following conditions hold: (a)~at remainder zero there is an action $a_0\in U_{tw}(0)$ with $r_{tw}(0,a_0)=0$ and $b_{tw}(a_0)\in K_{t+1}$; (b)~there is $\epsilon_{tw}>0$ such that every pair $(g,a)$ with $r_{tw}(g,a)<\epsilon_{tw}$ has $g=0$ and $b_{tw}(a)\in K_{t+1}$.
\end{assumption}

Assumption~\ref{assump:period}(a) states that a state without mismatch can remain without mismatch at no cost, and Assumption~\ref{assump:period}(b) states that every mismatch, and every action that creates one, has a positive price. For the batch balance of Section~\ref{sec:sep}, Assumption~\ref{assump:period}(b) holds when the mismatch cost $\lambda$ and the cost of any nonzero adjustment are positive. Assumption~\ref{assump:period}(a) is a separate requirement on the model, namely, that it contains an action that costs nothing and keeps a matched state matched. The following theorem shows that, under this assumption, $K_t$ is the set of periods of $R_t$.

\begin{theorem}\label{thm:optimal-period}
Under Assumption~\ref{assump:period}, for every $t$, a vector $h\in\Z^n$ satisfies $R_t(x+h)=R_t(x)$ for all $x\in\Z^n$ if and only if $h\in K_t$. Moreover, $R_t(x)=0$ if and only if $x\in K_t$. Consequently, $R_t$ is a function of the coset of a lattice $L$ if and only if $L\subseteq K_t$, and $[\Z^n:K_t]$ is the smallest number of cosets with which a lattice can represent $R_t$.
\end{theorem}

Theorem~\ref{thm:optimal-period} is a statement about the periodic part $R_t$ on all of $\Z^n$. It identifies $K_t$ as the period lattice of $R_t$, and any coarser lattice merges cosets on which $R_t$ differs. We did not verify Assumption~\ref{assump:period} for our benchmark families of Section~\ref{sec:experiments}, and we do not use Theorem~\ref{thm:optimal-period} to explain a computational result (Appendix~\ref{app:results-class}). A numerical check of the conclusion of the theorem on four small instances is reported in Appendix~\ref{app:add-theory}. The amount that a cut on a coarser lattice $L$ loses on the bounded set $D_t$ depends on the states that $D_t$ contains. If $D_t$ contains the points $0$, $h$, and $kh$ for some $h\in L$ of order $k\ge 2$ modulo $K_t$, then the cut has to interpolate between states on which $R_t$ differs, as the affine cut of Proposition~\ref{prop:loss} does. Corollary~\ref{cor:price-period} in Appendix~\ref{ec:additional} bounds the resulting loss from below, and the bound remains positive when closing cuts with moderate penalties are added. We note that Theorem~\ref{thm:optimal-period} does not state that lattice cuts on $K_t$ are exact on $D_t$, which is the subject of Proposition~\ref{prop:phasevalue}. In Appendix~\ref{ec:box}, we determine the smallest index of a lattice for which lattice cuts are tight at every state for every function on an integer box, i.e., a guarantee that does not depend on the structure of the model.

\section{Closing Cuts and the Attainable Bound}\label{sec:cuts}

In this section, we introduce closing cuts, which remove the nonconvexity that remains inside a coset (Section~\ref{sec:closing-def}), and we characterize the best bound that lattice cuts and closing cuts can deliver for a given family of lattices (Section~\ref{sec:ideal-pool}).

\subsection{Closing Cuts}\label{sec:closing-def}

Although a lattice cut removes the nonconvexity between cosets, nonconvexity can remain inside a coset: at a state $\hat s$ near a bound, a capacity, or a fixed-charge threshold, the values of the other states of the coset may lie below every affine function through $F(\hat s)$, as in Figure~\ref{fig:problem}(c). We close this gap with a penalty for the distance from $\hat s$. For nonnegative penalty vectors $\rho^+,\rho^-\in\R^n_+$, we define
\begin{equation}\label{eq:aniso-penalty}
 \Phi(y;\hat s,\rho^+,\rho^-)=\sum_{i=1}^n\bigl[\rho_i^+(y_i-\hat s_i)_+ + \rho_i^-(\hat s_i-y_i)_+\bigr],
\end{equation}
where $(u)_+=\max\{u,0\}$. The penalty is zero at the anchor $\hat s$ and grows linearly with the distance, at a rate that may differ by coordinate and direction.

\begin{definition}\label{def:lattice-cuts}
Fix $t\ge 2$, $w\in\Omega_t$, a slope $\pi$, an anchor $\hat s\in D_t$, and penalties $\rho^+,\rho^-$. Let $b_g$ be a lower bound on $\min_{z\in D_g}\{F_{tw}(z)-\pi^\top z+\Phi(z;\hat s,\rho^+,\rho^-)\}$ for each coset $g$ of a set $\mathcal G\subseteq\mathcal G_t$, and let $b_0$ be a lower bound on the same expression minimized over all of $D_t$. Set $\bar b(g)=\max\{b_g,b_0\}$ for $g\in\mathcal G$ and $\bar b(g)=b_0$ otherwise.
\begin{enumerate}[label=(\roman*),leftmargin=2em,itemsep=1pt,topsep=3pt]
\item A \emph{lattice cut} has $\rho^+=\rho^-=0$ and reads $\theta_w\ge\pi^\top y+\bar b([y]_t)$.
\item A \emph{closing cut} reads $\theta_w\ge\pi^\top y+\bar b([y]_t)-\Phi(y;\hat s,\rho^+,\rho^-)$.
\item An \emph{isotropic closing cut} is a closing cut with $\rho_i^+=\rho_i^-=\rho$ for all $i$, in which case $\Phi(y;\hat s,\rho^+,\rho^-)=\rho\|y-\hat s\|_1$.
\end{enumerate}
\end{definition}

Definition~\ref{def:lattice-cuts} allows lower bounds in place of the exact intercepts \eqref{eq:coset-intercept}, so that a MIP that is stopped early still yields a valid cut, and it allows intercepts to be known for some cosets only. All three cuts are valid by the argument given for \eqref{eq:coset-intercept}. We next state how large the penalty must be for a closing cut to be tight at its anchor. Let $f(z)=F_{tw}(z)-\pi^\top z$, let $\hat g$ denote the coset of $\hat s$, and let $d(z)$ be either $\|z-\hat s\|_1$ or $\Phi(z;\hat s,\omega^+,\omega^-)$ with fixed positive weights $\omega^\pm$. For $\rho\ge 0$, let
\begin{equation}\label{eq:closing-cut}
 b_g^\rho=\min_{z\in D_g}\bigl\{f(z)+\rho\,d(z)\bigr\},\qquad
 h^\rho(y)=\pi^\top y+b_{[y]_t}^\rho-\rho\,d(y).
\end{equation}
The following theorem shows that the cut $h^\rho$ is valid for every penalty and gives the smallest penalty at which it is tight at the anchor.

\begin{theorem}\label{thm:closure}
For every $\rho\ge 0$, we have $h^\rho\le F_{tw}$ on $D_t$. The cut is tight at the anchor, that is, $h^\rho(\hat s)=F_{tw}(\hat s)$, if and only if $\rho\ge\rho_{\hat g}^*$, where
\begin{equation}\label{eq:least-penalty}
 \rho_{\hat g}^*=
 \max\!\left(\{0\}\,\cup
 \left\{\frac{f(\hat s)-f(z)}{d(z)}:\ z\in D_{\hat g},\ z\ne\hat s\right\}\right).
\end{equation}
Replacing $\Gamma_t$ by a finer lattice cannot increase $\rho_{\hat g}^*$.
\end{theorem}

The threshold \eqref{eq:least-penalty} is the steepest descent of $f$ from $\hat s$ to another state of the same coset, per unit of distance, and it explains why the lattice and the closing cut work well together. Without a lattice, the coset of $\hat s$ is all of $D_t$ and contains the neighbors of $\hat s$ with other remainders. The penalty must then be steep enough to bridge the periodic oscillation of the value function at distance one, and it must be steep at every state. With the lattice, the nearest competitors are a full period away and have the same remainder. Therefore, the threshold is small and, at many states, zero. We quantify this difference in Section~\ref{sec:exp-long}. The penalty is not chosen in advance. Starting from zero, our algorithm increases $\rho_i^+$ or $\rho_i^-$ only in the coordinates and directions in which a state of the coset violates the cut, and it stops after finitely many steps (Proposition~\ref{prop:dinkelbach}).

In the model of stage $t-1$, the term $-\Phi$ is concave in $y$ and requires binary variables that indicate on which side of the anchor each coordinate lies. We show in Section~\ref{sec:pooled-splits} how all closing cuts of a stage share these variables.

\subsection{The Bound a Lattice Can Deliver}\label{sec:ideal-pool}

We next characterize the best bound that can be obtained at the root for a given family of lattices $\Gamma=\{\Gamma_t\}$. To separate this bound from the progress of any particular run, suppose that every stage $t\ge 2$ had \emph{all} valid lattice cuts and all valid closing cuts with penalties $\rho_i^\pm\le\rho$ at its disposal, built recursively from stage $T$ backward, and let $W_\Gamma(\rho)$ denote the resulting bound at $s_1$. We show in Theorem~\ref{thm:root-price} that $W_\Gamma(\rho)$ is the optimal value of an LP.

The LP describes a relaxed decision maker who, on arriving at a state $s$ under outcome $w$, may act as if it were at other states $z$ of the same coset, provided that these states average to $s$. Let $S_t(s)=D_t\cap(s+\Gamma_t)$ for $t\ge 2$ and $S_1(s_1)=\{s_1\}$, and let $c^0_{tw}(z,y)$ denote the minimum cost of moving from $z$ to $y$ under outcome $w$. The variables are the probability $d_t(s)$ of arriving at $s$ and the flow $f_{twszy}\ge 0$ that arrives at $s$, acts as if at $z\in S_t(s)$, and moves to $y$. These variables form the polytope $\mathcal P_\Gamma$ given by
\begin{equation}\label{eq:occupation}
\begin{aligned}
 &d_1(s)=\mathbf 1\{s=s_1\},\qquad d_{t+1}(y)=\sum_{w,s,z}f_{twszy}\quad\text{for all } t,y,\\
 &\sum_{z,y} f_{twszy}=p_{tw}\,d_t(s),\qquad \sum_{z,y}(z-s)\,f_{twszy}=0\quad\text{for all } t,w,s.
\end{aligned}
\end{equation}
For $v=(d,f)\in\mathcal P_\Gamma$, let $J(v)=\sum f_{twszy}\,c^0_{tw}(z,y)+\sum_{y}d_{T+1}(y)V_{T+1}(y)$ denote the expected cost, and let $D(v)=\sum f_{twszy}\,\|z-s\|_1$ denote the expected displacement between the true state and the state acted upon.

\begin{theorem}\label{thm:root-price}
Let $\mathcal V_\Gamma$ be the set of vertices of $\mathcal P_\Gamma$. For every $\rho\ge 0$, we have
\begin{equation}\label{eq:root-price-value}
 W_\Gamma(\rho)=\min_{v\in\mathcal P_\Gamma}\bigl\{J(v)+\rho\,D(v)\bigr\}.
\end{equation}
The function $W_\Gamma$ is nondecreasing, concave, and piecewise linear, and $W_\Gamma(\rho)=V^\star$ if and only if
\begin{equation}\label{eq:root-price-threshold}
 \rho\ \ge\ \rho_\Gamma^\star=\max\!\left(\{0\}\cup
 \left\{\frac{V^\star-J(v)}{D(v)}:\ v\in\mathcal V_\Gamma,\ D(v)>0\right\}\right).
\end{equation}
Refining the lattices cannot decrease $W_\Gamma(\rho)$ and cannot increase $\rho_\Gamma^\star$.
\end{theorem}

Theorem~\ref{thm:root-price} has the following interpretation. A plan without displacement is an ordinary policy and costs at least $V^\star$. A plan with displacement can be cheaper because it mixes the costs of several states, and \eqref{eq:root-price-value} charges $\rho$ per unit of displacement for this. Lattice cuts alone correspond to $\rho=0$, in which case displacement is free but confined to a coset. The finer the lattice, the less room there is for displacement, and when the condition of Proposition~\ref{prop:phasevalue} holds, displacement inside a coset gains nothing and $W_\Gamma(0)=V^\star$. The threshold $\rho_\Gamma^\star$ is finite and can be bounded from the cost data of the stages (Corollary~\ref{cor:data-cap}).

Theorem~\ref{thm:root-price} also yields a decomposition of the gap of any run. Suppose that, after $k$ iterations, the algorithm has certified the bound $L_k$ and that the penalties of its closing cuts are at most $\bar\rho$. Then, we have $L_k\le W_\Gamma(\bar\rho)\le V^\star$ and
\begin{equation}\label{eq:error-decomp}
 V^\star - L_k \ =\ \underbrace{V^\star - W_\Gamma(\bar\rho)}_{\text{representation loss}}
 \ +\ \underbrace{W_\Gamma(\bar\rho) - L_k}_{\text{finite-pool loss}} .
\end{equation}
The representation loss depends only on the lattice and the penalty cap. It is not removed by any number of iterations, but only by a finer lattice or a larger cap. In contrast, the finite-pool loss is due to the cuts not yet generated and shrinks as cuts accumulate. In Appendix~\ref{ec:additional}, we express the finite-pool loss through the reduced costs of the LP \eqref{eq:root-price-value} and give a test for $L_k=W_\Gamma(\bar\rho)$ (Theorem~\ref{thm:finite-pool-root-gap}). In Appendices~\ref{ec:local-certs} and~\ref{ec:horizon-certs}, we give further certificates that are computed from the data collected during a run, for a single stage and for the whole horizon, respectively.

\section{Lattice-Cut Decomposition}\label{sec:algorithm}

In this section, we present the \LCD\ algorithm (Section~\ref{sec:cut-subproblems}), describe a pooled formulation that keeps the stage models small (Section~\ref{sec:pooled-splits}), and establish finite convergence (Section~\ref{sec:root}).

\subsection{The Algorithm}\label{sec:cut-subproblems}

\LCD\ is stated in Algorithm~\ref{alg:lcd}. The algorithm follows the forward and backward passes of stochastic dual dynamic programming \citep{pereira1991} and differs in what happens at a visited state.

\begin{algorithm}[t]
\caption{(Lattice-Cut Decomposition)}\label{alg:lcd}
\begin{singlespace}\small
\textbf{Input:} the problem data; lattices from Algorithm~\ref{alg:lattice}; tolerances; gate $\gamma\in(0,1]$.\par
\textbf{Output:} bounds $L_0\le L_1\le\cdots$ and cuts.
\begin{pseudocode}
\pcl{\label{ln:init}Initialize: build one model per stage with coset indicators, compute LP lower bounds on the cost-to-go functions and use them as initial cuts, and solve stage 1 for $L_0$}
\pcl{\kw{for} $k=1,2,\ldots$ \kw{do}}
\pcl[1]{\label{ln:forward}Forward pass: sample scenario paths from $s_1$, solve the stage problems along them with the current cuts, and record the visited states of stages $2,\ldots,T$}
\pcl[1]{\label{ln:backward}\kw{for} $t=T,\ldots,2$, every visited state $\hat s$ of stage $t$, and every $w\in\Omega_t$ \kw{do}}
\pcl[2]{\label{ln:solve}Solve the stage problem at $\hat s$}
\pcl[2]{\kw{if} the current cuts are not tight at $\hat s$ \kw{then}}
\pcl[3]{\label{ln:lattice}Compute a lattice cut at $\hat s$ (slope, intercept of the coset of $\hat s$, lower bound $b_0$) and add it to stage $t-1$, merging it with an existing cut of the same slope}
\pcl[3]{\label{ln:gate}\kw{if} a state of the coset of $\hat s$ lies below the lattice cut at $\hat s$, and $\hat s$ has been visited before or the lattice cut closed less than the fraction $\gamma$ of the gap \kw{then}}
\pcl[4]{\label{ln:closing}Raise the penalties $\rho^\pm$ in the violated directions until the closing cut is tight at $\hat s$, and add it to stage $t-1$}
\pcl[3]{\label{ln:gateend}\kw{end}}
\pcl[2]{\kw{end}}
\pcl[1]{\label{ln:backend}\kw{end}}
\pcl[1]{\label{ln:root}Solve stage 1}
\pcl[1]{\label{ln:escape}\kw{while} its solution lies in a coset of stage 2 without a computed intercept and the bound rises \kw{do}}
\pcl[2]{Compute that intercept and solve stage 1 again}
\pcl[1]{\label{ln:escapeend}\kw{end}}
\pcl[1]{\label{ln:bound}$L_k\leftarrow$ the value obtained}
\pcl[1]{\label{ln:stop}\kw{if} the time limit or the stopping criterion is met \kw{then} stop}
\pcl{\kw{end}}
\end{pseudocode}
\end{singlespace}
\end{algorithm}

\begin{figure}[t]
\centering
\resizebox{0.94\textwidth}{!}{\begin{tikzpicture}[font=\footnotesize, >={Stealth[length=2.1mm]}, line width=0.5pt,
  box/.style={draw, rounded corners=2pt, align=center, minimum height=1.0cm, inner sep=4pt},
  dec/.style={draw, diamond, aspect=2.1, align=center, inner sep=0.8pt},
  lab/.style={font=\scriptsize, inner sep=1.5pt},
  dot/.style={circle, fill, inner sep=1.1pt}]
\def\Xa{5.95}\def\Xb{10.35}\def\Xc{12.75}
\def\ya{-2.85}\def\yb{-4.85}\def\yc{-6.85}
\node[box, minimum width=4.4cm] (p1) at (0,0) {Stage matrices:\\batch balances, transitions};
\node[box, minimum width=4.4cm] (p2) at (5.55,0) {Compatible lattices\\$\Gamma_1,\ldots,\Gamma_T$ (Algorithm 1)};
\node[box, minimum width=4.4cm] (p3) at (11.1,0) {One model per stage:\\coset indicators, LP bounds};
\draw[->] (p1) -- (p2); \draw[->] (p2) -- (p3);
\node[lab, anchor=south west] at ($(p1.north west)+(0,0.03)$) {\emph{Preparation (once, before optimization)}};
\node[box, minimum width=3.5cm] (f) at (0,\ya) {Forward pass:\\sample paths, record\\visited states};
\node[box, minimum width=3.5cm] (r) at (0,\yb) {Solve stage 1:\\certified bound $L_k$};
\node[dec] (d3) at (0,\yc) {Stop?};
\node[box, minimum width=3.5cm] (o) at (0,-8.45) {Bound $L_k$ and policy};
\node[box, minimum width=4.0cm] (b1) at (\Xa,\ya) {Solve the stage problem\\at visited state $\hat s$,\\outcome $w$};
\node[dec] (d1) at (\Xb,\ya) {Cuts tight\\at $\hat s$?};
\node[box, minimum width=4.0cm] (b2) at (\Xa,\yb) {Lattice cut: one slope,\\intercept of the\\coset of $\hat s$};
\node[dec] (d2) at (\Xb,\yb) {Gap left inside\\the coset?};
\node[box, minimum width=4.0cm] (b3) at (\Xa,\yc) {Closing cut anchored\\at $\hat s$ (pooled\\split variables)};
\coordinate (busT) at (\Xc,\ya+1.12);
\draw[->] (p3.south) -- ++(0,-0.4) -| (f.north);
\draw[->] (f) -- (b1);
\draw[->] (b1) -- (d1);
\draw[->] (d1.south) -- ++(0,-0.28) -| node[lab, pos=0.1, above] {no} (b2.north);
\draw[->] (b2) -- (d2);
\draw[->] (d2.south) -- ++(0,-0.28) -| node[lab, pos=0.1, above] {yes} (b3.north);
\draw (d1.east) -- node[lab, above] {yes} (\Xc,\ya) node[dot] {};
\draw (d2.east) -- node[lab, above] {no} (\Xc,\yb) node[dot] {};
\draw[->] (b3.east) -- (\Xc,\yc) -- (busT) -| (b1.north);
\node[lab, anchor=north] at ($(\Xb,\ya)+(-0.2,1.09)$) {next visited state};
\node[draw, dashed, rounded corners=3pt, inner sep=7pt, fit=(b1)(d1)(b2)(d2)(b3)(busT)] (bw) {};
\node[lab, anchor=north east] at ($(bw.south east)+(0,-0.02)$) {\emph{Backward pass: stages $T,\ldots,2$, every visited state $\hat s$ and outcome $w$}};
\draw[->, dashed] (bw.west |- r) -- node[lab, above, align=center] {all states\\done} (r.east);
\draw[->] (r) -- (d3);
\draw[->] (d3) -- node[lab, right] {yes} (o);
\draw[->] (d3.west) -- (-2.5,\yc) |- node[lab, pos=0.22, left] {no} (f.west);
\end{tikzpicture}}
\caption{Flow Chart of Lattice-Cut Decomposition}\label{fig:flow}
\end{figure}
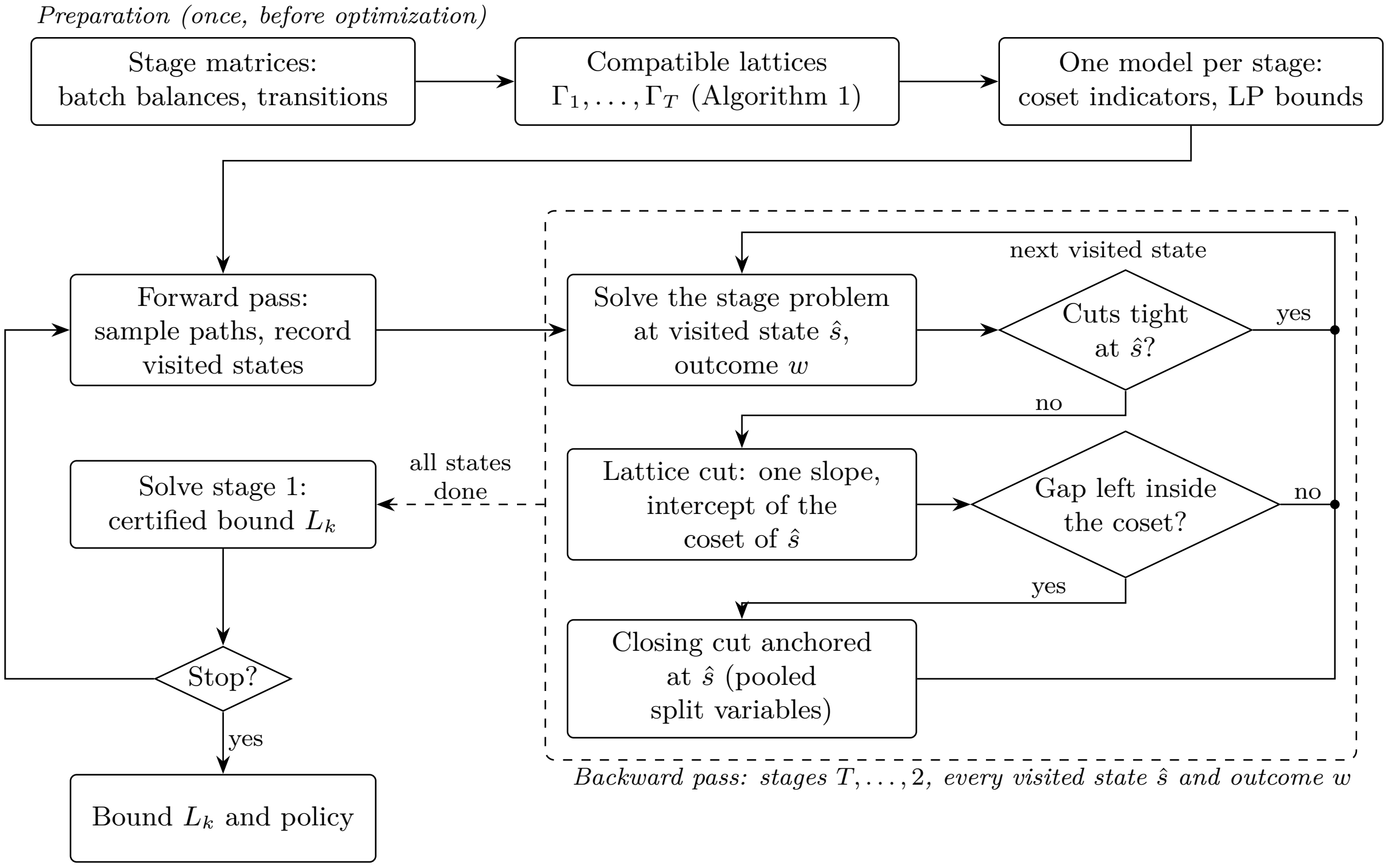

The forward pass in line~\ref{ln:forward} and the loop of the backward pass in line~\ref{ln:backward} of Algorithm~\ref{alg:lcd} are those of stochastic dual dynamic programming. The steps that are specific to \LCD\ are the coset indicators of line~\ref{ln:init}, the lattice cut of line~\ref{ln:lattice}, the gated closing cut of lines~\ref{ln:gate}--\ref{ln:gateend}, and the loop at the root in lines~\ref{ln:escape}--\ref{ln:escapeend}. Figure~\ref{fig:flow} shows the overall structure of \LCD. In this figure, the dashed box contains the steps of the backward pass, which are carried out for every visited state and outcome. We discuss the three components of our algorithm in turn.

\emph{Preparation (line~\ref{ln:init}).} We compute the lattices with Algorithm~\ref{alg:lattice}. Because the outcomes of a stage usually differ only in coefficients, \LCD\ builds one model per stage and switches between outcomes by updating coefficients, which also means that a cut is installed once and not once per outcome. Solving the LP relaxations of the stage models backward gives initial lower bounds on the cost-to-go functions, and a few passes with ordinary Benders cuts \citep{benders1962} on the relaxation provide a starting approximation.

\emph{Backward pass (lines~\ref{ln:backward}--\ref{ln:backend}).} At a visited state $\hat s$ and outcome $w$, \LCD\ solves the stage problem with the incoming state fixed at $\hat s$ (line~\ref{ln:solve}). If the existing cuts already attain this value, then nothing is added. Otherwise, \LCD\ generates a lattice cut as described in Section~\ref{sec:lattice-cuts} (line~\ref{ln:lattice}). Specifically, the slope is taken from the LP dual, the intercept of the coset of $\hat s$ is obtained from the program \eqref{eq:coset-intercept}, and the slope is improved by subgradient steps if the cut is not yet tight. In Appendix~\ref{sec:direction}, we characterize the gain in the value of a cut that is obtained when the slope may move in one additional direction. When a cut with the same slope exists, the new intercepts are merged into it, so that the number of rows does not grow. Merging is exact for identical slopes. In our implementation, we also merge slopes that agree to within $10^{-8}$ in every coordinate, and we reuse stored intercepts for slopes that agree to within $10^{-10}$. An intercept that is attached to a slope $\pi$ other than the slope $\pi'$ for which it was computed is first lowered by $\max\{0,\max_{y\in[\ell_t,u_t]}(\pi-\pi')^\top y\}$ and rounded down, so that every cut remains valid (Appendix~\ref{app:methods-ours}); Proposition~\ref{prop:merge-tolerance} in Appendix~\ref{ec:additional} bounds the effect of merging without this shift. A closing cut follows (line~\ref{ln:closing}) only if a state of the same coset is proven to lie below the lattice cut at $\hat s$ by more than a tolerance, so that no lattice cut with this slope can be tight. In dimension three and higher, \LCD\ postpones the closing cut on the first visit when the lattice cut has already closed a fraction $\gamma$ of the gap, and it generates the closing cut if the state is visited again. The value of $\gamma$ controls the tradeoff between the time of an iteration and the tightness of its cuts: a smaller value postpones more closing cuts, which makes an iteration faster but leaves it less complete.

\emph{Root (lines~\ref{ln:root}--\ref{ln:bound}).} Solving stage 1 gives the certified bound $L_k$. If the solution moves to a coset of stage 2 for which the active cuts hold only the common lower bound $b_0$, then \LCD\ computes the intercept of this coset and solves stage 1 again, until the bound stops rising (lines~\ref{ln:escape}--\ref{ln:escapeend}). This step is important for lattice cuts because a cut with a missing intercept on one coset allows the stage-1 solution to move to that coset.

\subsection{Pooled Formulation of the Stage Models}\label{sec:pooled-splits}

Cuts with nonconvex terms bring binary variables into the model of the preceding stage, and in a direct formulation the number of these variables grows with the number of cuts. \LCD\ avoids this growth for both of its cut families.

For lattice cuts, the coset indicators are shared, as noted after \eqref{eq:lc-row}. When a stage has at most 64 cosets, \LCD\ uses one indicator per coset. When a stage has more cosets, \LCD\ encodes the coset label of Algorithm~\ref{alg:lattice} in binary digits, which requires $\sum_j\lceil\log_2 q_{t,j}\rceil$ variables per stage, and it creates a continuous indicator only for the cosets that appear in a cut.

For closing cuts, the straightforward formulation writes $y_i-\hat s_i=v_i^+-v_i^-$ for each cut and each coordinate and uses one binary variable to force $v_i^+v_i^-=0$, that is,
\begin{equation}\label{eq:sign-binary}
 0\le v_i^+ \le (u_{t,i}-\hat s_i)\,z_i,\qquad
 0\le v_i^- \le (\hat s_i-\ell_{t,i})(1-z_i),\qquad z_i\in\{0,1\}.
\end{equation}
With $C$ closing cuts, this formulation adds up to $nC$ binary variables to every model that carries the cuts. However, the binary variable in \eqref{eq:sign-binary} only indicates whether $y_i$ exceeds $\hat s_i$, and this information can be shared by all cuts whose anchor has the same $i$th coordinate. Binary variables indexed by coordinate and value in this way appear in \citet{deng2024}, who show that rectified linear cuts on integer variables $x_i\in\{0,\ldots,B_i\}$ need at most $\sum_i(B_i-1)$ additional binary variables. We use this indexing for all closing cuts and all outcomes of a stage. For each coordinate $i$ and each threshold value $v$ that occurs as the $i$th coordinate of an anchor, with $\ell_{t,i}<v<u_{t,i}$, the stage model receives one triple $(p_{i,v},q_{i,v},z_{i,v})$ with the constraints
\begin{equation}\label{eq:pooled-split}
\begin{aligned}
 &y_i-p_{i,v}+q_{i,v}=v,\qquad
 0\le p_{i,v}\le(u_{t,i}-v)\,z_{i,v},\qquad
 0\le q_{i,v}\le(v-\ell_{t,i})(1-z_{i,v}),\\
 &z_{i,v}\in\{0,1\},\qquad
 z_{i,v_2}\le z_{i,v_1}\quad\text{for thresholds } v_1<v_2 .
\end{aligned}
\end{equation}
The last constraint orders the thresholds of a coordinate like a staircase: if $y_i$ exceeds $v_2$, then it also exceeds $v_1$. The following theorem establishes the validity of closing cuts in this formulation and bounds the number of binary variables that they require.

\begin{theorem}\label{thm:pooled-closing}
Fix $\hat s\in D_t$, a slope $\pi$, penalties $\rho^\pm\in\R^n_+$, and intercepts $\bar b$ as in Definition~\ref{def:lattice-cuts}.
\begin{enumerate}[label=(\roman*),leftmargin=2em,itemsep=1pt,topsep=3pt]
\item The closing cut $h(y)=\pi^\top y+\bar b([y]_t)-\Phi(y;\hat s,\rho^+,\rho^-)$ satisfies $h\le F_{tw}$ on $D_t$. Moreover, if $b_{\mathrm{LP}}+\pi_{\mathrm{LP}}^\top z\le F_{tw}(z)$ is a cut from the LP relaxation, $\underline F\le F_{tw}(\hat s)$ is a lower bound, and $\Delta=\underline F-b_{\mathrm{LP}}-\pi_{\mathrm{LP}}^\top\hat s>0$, then $\pi=\pi_{\mathrm{LP}}$, $\bar b\equiv b_{\mathrm{LP}}+\Delta$, and $\rho_i^\pm=\Delta$ give a valid closing cut with $h(\hat s)=\underline F$ without solving a further MIP.
\item If $\hat s_i=\ell_{t,i}$, then replacing $\pi_i$ by $\pi_i-\rho_i^+$ and adding $\rho_i^+\hat s_i$ to $\bar b$ leaves $h$ unchanged on $[\ell_t,u_t]$ and removes the penalty on coordinate $i$; if $\hat s_i=u_{t,i}$, then the same holds with $\pi_i+\rho_i^-$ and $-\rho_i^-\hat s_i$. Coordinates at a bound therefore need no binary variable.
\item For integer $y_i\in[\ell_{t,i},u_{t,i}]$, the constraints \eqref{eq:pooled-split} force $p_{i,v}=(y_i-v)_+$ and $q_{i,v}=(v-y_i)_+$, each triple is an ideal formulation, and $-\Phi(y;\hat s,\rho^+,\rho^-)=-\sum_i\bigl[\rho_i^+p_{i,\hat s_i}+\rho_i^-q_{i,\hat s_i}\bigr]$ for the coordinates not removed in (ii). The number of binary variables on coordinate $i$ is at most $\max\{0,u_{t,i}-\ell_{t,i}-1\}$, whatever the number of closing cuts and of outcomes.
\end{enumerate}
\end{theorem}

Theorem~\ref{thm:pooled-closing}(i) includes an inexpensive closing cut that lifts an LP cut to a known value. Theorem~\ref{thm:pooled-closing}(ii) and (iii) bound the size of the stage models, that is, the number of binary variables that closing cuts can add depends on the ranges of the state variables and not on the number of cuts. The count in Theorem~\ref{thm:pooled-closing}(iii) is the one given by \citet{deng2024} for rectified linear cuts. We regard the pooled formulation as an implementation technique for closing cuts and not as a contribution in its own right. In the 300-second runs of Section~\ref{sec:exp-equal} in which closing cuts are generated, the pooled formulation needs between 14 and 128 times fewer binary variables than \eqref{eq:sign-binary} would need for the same cuts (e.g., 301 instead of 4,334 on the 24-stage Case-Pack-Truck instance).

\subsection{Convergence}\label{sec:root}

The following theorem establishes that \LCD\ converges to the optimal value in finitely many iterations.

\begin{theorem}\label{thm:root}
Let Assumption~\ref{assump:msip} hold and let all subproblems be solved exactly, with zero tolerances. Suppose that cut generation is deterministic given the stage, outcome, state, current cuts of the next stage, and number of visits; that cuts are never removed; that ties in the forward pass are broken by a fixed rule; and that every scenario path is sampled with positive probability, independently across iterations. Then, with probability one, \LCD\ reaches $L_k=V^\star$ after finitely many iterations, and the policy defined by its cuts is optimal from $s_1$.
\end{theorem}

The proof of Theorem~\ref{thm:root} is by backward induction over the stages, in the spirit of the argument of \citet{zou2019} for binary states, and we outline it here. The approximation of the terminal cost never changes. Suppose that the approximation of stage $t+1$ no longer changes after some iteration. From then on, the stage value $F_{tw}$ is a fixed function on the finite set $D_t$. A state that the forward passes visit only finitely often no longer matters after its last visit. A state that is visited infinitely often receives, at the latest on its second visit after that iteration, a cut that is tight there, namely, the lattice cut if it is tight and otherwise a closing cut, whose penalties are raised in finitely many steps (Proposition~\ref{prop:dinkelbach}) until it is tight at its anchor (Theorem~\ref{thm:closure}). Later visits find the cuts tight and add nothing, and each of the finitely many missing intercepts is computed at most once. Therefore, the approximation of stage $t$ stops changing as well. Once all stages have stabilized, the forward pass follows one fixed policy. With probability one, every scenario path is sampled infinitely often, and therefore the cuts are tight at every state that this policy reaches. The bound $L_k$ then equals the expected cost of the policy, which is at least $V^\star$, and the validity of the cuts gives $L_k\le V^\star$.

Theorem~\ref{thm:root} also holds for the basic variant of our method that we use in part of our computational study and refer to as \LCDiso. This variant generates lattice cuts and isotropic closing cuts with the formulation \eqref{eq:sign-binary} at every state where a gap remains ($\gamma=1$), and it uses neither the shared stage models nor the merged cuts.

\section{Computational Study}\label{sec:experiments}

We report the results of our computational study, in which we evaluate lattice cuts and compare \LCD\ with existing methods. In Section~\ref{sec:exp-setup}, we describe the test instances, the methods, and the implementation. In Section~\ref{sec:exp-decomp}, we investigate whether the lattice removes the representation loss of affine cuts. In Section~\ref{sec:exp-long}, we investigate whether this yields tight certified bounds on long horizons. In Section~\ref{sec:exp-equal}, we compare \LCD\ with existing methods under equal time limits, and in Section~\ref{sec:exp-regimes}, we examine the cases in which the lattice does not help.

\subsection{Test Instances, Methods, and Implementation}\label{sec:exp-setup}

\subsubsection{Test Instances}\label{sec:exp-instances}

We conduct our study on the 15 problem families listed in Table~\ref{tab:families-lattices}. For each family, the table gives the batch structure that defines its lattice, the ranges of the number of state variables ($n$) and of the number of stages ($T$) over its instances, and the number of cosets per stage. The four families of Group~I are synthetic models built around a batch balance, with state dimensions and horizons that can be scaled freely. Group~II contains operational models with commitments over time windows, and Group~III contains batch inventory, workforce, production, rebalancing, and revenue management models. The number of cosets in the last column of Table~\ref{tab:families-lattices} is that of the lattice computed by Algorithm~\ref{alg:lattice}. This number ranges from 4 to more than $10^{100}$, and we find in Section~\ref{sec:exp-regimes} that it is one of the factors that decide where lattice cuts are useful.

\begin{table}[t]
\centering
\caption{Summary of the Characteristics of the Benchmark Families}\label{tab:families-lattices}
\begin{singlespace}\footnotesize
\setlength{\tabcolsep}{4pt}
\begin{tabular*}{\textwidth}{@{\extracolsep{\fill}}llrrl@{}}
\toprule
Family & Batch structure & $n$ & $T$ & Cosets per stage \\
\midrule
\multicolumn{5}{@{}l}{\emph{Group I: synthetic batch models}}\\
Shift & window sums modulo 6 & 6--6,144 & 3--2,048 & 36 \\
Rotate & crew differences modulo 6 & 6--6,144 & 3--2,048 & 36 \\
Odd-Cycle & parities along an odd cycle & 6--6,144 & 3--2,048 & 4 \\
Shear & prefix sums modulo 8, 4, and 2 & 6--6,144 & 3--2,048 & 64 \\
\multicolumn{5}{@{}l}{\emph{Group II: commitments over time windows}}\\
Hub-Consolidation & daily commitments modulo 6 & 130--2,050 & 52--260 & 36 \\
Sterile-Trays & two queues, modulo 6 and 8 & 256--2,048 & 84--336 & 2,304 \\
Pooled-Capacity & pooled supply modulo 6 & 2--80 & 5--25 & 1,296 \\
Fulfillment-Network & route commitments modulo 6 & 528--2,112 & 52--260 & $36^{D}$, $D$ routes \\
\multicolumn{5}{@{}l}{\emph{Group III: batch inventory, workforce, production, and revenue}}\\
Case-Pack & stock modulo pack size $p_i$ & 2--64 & 6--48 & $\prod_i p_i$ \\
Workforce-Cohort & staff modulo cohort size $c_i$ & 2--72 & 6--48 & $\prod_i c_i$ \\
Co-Product & recipe lattice $Y_g\Z^2$, $|\det Y_g|\in\{5,7\}$ & 2--48 & 6--48 & $\prod_g|\det Y_g|$ \\
Bike-Sharing & transfers in lots of $k$ bikes & 2--25 & 6--48 & $k^{S+1}$ or $Mk^{S}$, $S$ stations \\
Serial-Inventory & two echelons, modulo $q_{1i}$ and $q_{2i}$ & 4--120 & 5--48 & $\prod_i q_{1i}q_{2i}$ \\
Lead-Time-Pipeline & pipeline stock modulo $b_i$ & 2--144 & 6--48 & $18^3\prod_i b_i$ \\
Airline-RM & seat blocks of several sizes & 2--288 & 5--24 & up to $2.8\times10^{110}$ \\
\bottomrule
\end{tabular*}
\end{singlespace}
\figurenotes{Here, $n$ denotes the number of state variables and $T$ denotes the number of stages. Case-Pack instances with a fixed charge per truck are called Case-Pack-Truck. For Bike-Sharing, the declared lattice has $k^{S+1}$ cosets and the computed lattice has $Mk^S$ cosets, where $M$ is one plus the total capacity (Appendix~\ref{app:benchmarks}).}
\end{table}

\subsubsection{Methods}\label{sec:exp-methods}

We consider the following five methods.
\begin{enumerate}[leftmargin=*,itemsep=1pt,topsep=3pt]
\item \LCD: Lattice-Cut Decomposition as stated in Algorithm~\ref{alg:lcd}.
\item \LCDiso: The basic variant of our method described in Section~\ref{sec:root}.
\item Closing cuts without the lattice: \LCDiso\ with the trivial lattice $\Gamma_t=\Z^n$, i.e., affine cuts with isotropic closing cuts. We run this variant to isolate the effect of the lattice.
\item SDDiP: Stochastic dual dynamic integer programming with strengthened Benders and integer optimality cuts \citep{zou2019}.
\item SDDmiP: The normalized rectified linear cuts alternated with LP cuts of \citet{bansal2026}.
\end{enumerate}
We compare with SDDiP and SDDmiP because they represent the two lines of work described in Section~\ref{sec:intro}: SDDiP expands every integer state into binary digits, whereas SDDmiP keeps the original state and adds a nonconvex term to each cut. In the legends of Figures~\ref{fig:long} and~\ref{fig:anytime}, SDDmiP appears as SDDmiP~(Norm-A), and in the legend of Figure~\ref{fig:anytime}, SDDiP appears as SDDiP~(SB+I). In Section~\ref{sec:exp-regimes} and Appendix~\ref{app:additional}, we also run \LCD\ with the trivial lattice, in order to isolate the effect of the lattice in the full algorithm.

\subsubsection{Implementation Details}\label{sec:exp-implementation}

All methods are implemented in Python~3.14 with Gurobi~13.0. SDDiP and SDDmiP are our own implementations of the published algorithms, so that all methods share the same stage models and solver interface. Every run uses one thread of an eight-core desktop computer, and its time limit includes the construction of the models. We list the default parameters of \LCD\ and \LCDiso, including the gate $\gamma$, in Table~\ref{tab:ec-params} in Appendix~\ref{ec:benchmarks}. We give the benchmark formulations in Appendix~\ref{app:benchmarks}, describe the compared methods as implemented and their settings in Appendix~\ref{app:methods}, and report the results instance by instance in Appendix~\ref{app:results}. In the long runs of Section~\ref{sec:exp-long} and the 45-second runs of Section~\ref{sec:exp-regimes}, the lattices were declared from the batch structure of each model and checked against the stage matrices, whereas in the 300-second runs of Section~\ref{sec:exp-equal} and the catalog runs of Section~\ref{sec:exp-regimes}, they are computed by Algorithm~\ref{alg:lattice} within the time limit.

The runs with time limits of 20 to 300 seconds stop only at the time limit and receive no information about the optimal value. In the long runs of Section~\ref{sec:exp-long}, we stopped a run of \LCDiso, with or without the lattice, as soon as its bound reached the reference optimum to within a relative tolerance of $\max\{10^{-6},\,1.05\times10^{-7}(T-1)\}$, which covers the accumulated safety margins. This rule shortens these runs and cannot raise their bounds. \LCD\ and \LCDiso\ solve their MIPs to a relative gap of $10^{-9}$ and subtract a safety margin from every dual bound, so that the certified bounds remain valid under solver tolerances. Therefore, their bounds stop slightly below the optimum, at a relative distance of about $10^{-6}$. Reference optimal values $V^\star$ are obtained from exact dynamic programming over the state space or from the extensive form. Across all runs of our methods with a known $V^\star$, no certified bound exceeds $V^\star$. When \LCD\ merges or reuses intercepts across nearby slopes, each intercept is lowered as described in Appendix~\ref{app:methods-ours}, so that every cut remains valid. In addition, we evaluated all cuts that \LCD\ generates on five small instances at all states against exact dynamic programming and found no cut above the exact value (Appendix~\ref{app:add-validity}). Appendix~\ref{app:additional} reports this audit together with further experiments that we made after the study of this section, and we refer to them where they bear on a result.

\subsubsection{Performance Measures}\label{sec:exp-measures}

We compare the methods by the certified lower bound $L$ that they obtain at the root. When the optimal value $V^\star$ is known, we report the relative gap, which is calculated as $(V^\star-L)/|V^\star|$ and stated as a fraction in scientific notation. We also report the mean cost of the policy defined by the cuts of each method on sampled scenario paths. We emphasize that these sample means estimate the costs of the policies. They are not upper bounds on the optimal value, and a tight lower bound does not by itself certify the policy that the cuts define.

\subsection{Representation Loss With and Without the Lattice}\label{sec:exp-decomp}

\begin{figure}[t]
\centering
\includegraphics[width=0.94\textwidth]{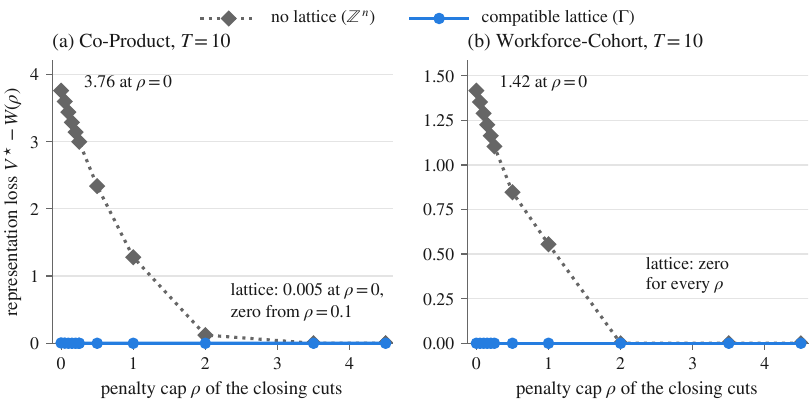}
\caption{\coloronline Representation Loss as a Function of the Penalty Cap on a Co-Product Instance (a) and a Workforce-Cohort Instance (b), With and Without the Lattice}\label{fig:root-gap-decomp}
\figurenotes{The values of $V^\star-W(\rho)$ are exact and are obtained from dynamic programming over the full state space. The loss at $\rho=0$ is the loss that lattice cuts (blue) or affine cuts (gray) leave when no closing cut is used.}
\end{figure}

We first examine whether the lattice removes the representation loss of affine cuts. Two 10-stage instances are small enough to evaluate the bound $W_\Gamma(\rho)$ of Theorem~\ref{thm:root-price} exactly by dynamic programming: a Co-Product instance with one recipe (169 states per stage, 5 cosets, $V^\star=168.26$) and a Workforce-Cohort instance with one cohort (38 states per stage, 6 cosets, $V^\star=468.67$). In Figure~\ref{fig:root-gap-decomp}, we plot the representation loss $V^\star-W(\rho)$ against the penalty cap for the lattice of Algorithm~\ref{alg:lattice} and for the trivial lattice.

On the Co-Product instance, affine cuts leave a loss of 3.76, or 2.2\% of the optimal value, and closing cuts remove this loss only when their penalties may reach 4.5. With the lattice, the loss at $\rho=0$ is 0.005, and it disappears at $\rho=0.1$. Therefore, the lattice removes 99.9\% of the loss of affine cuts and lowers the penalty needed for the rest by a factor of at least 35. On the Workforce-Cohort instance, the loss of affine cuts is 1.42, whereas with the lattice the loss is zero at $\rho=0$, that is, lattice cuts alone are exact on this instance.

Actual runs behave as predicted by the decomposition \eqref{eq:error-decomp}. On the Workforce-Cohort instance, the bound of affine cuts is within 0.01 of $W_{\Z^n}(0)$ after eight iterations and 1.43 below the optimum, and by Theorem~\ref{thm:root-price} it cannot pass $W_{\Z^n}(0)$. \LCDiso\ reaches the optimum to within the safety margin in 4.7 seconds and generates no closing cut. Closing cuts without the lattice also reach the optimum. However, they need 44 closing cuts and 23 seconds. On the Co-Product instance, lattice cuts alone approach $W_\Gamma(0)$ and, as expected from Theorem~\ref{thm:root-price}, cannot pass it: after 30 minutes, their bound is 0.007 below $W_\Gamma(0)$ and 0.012 below the optimum. \LCDiso\ passes $W_\Gamma(0)$ with four closing cuts, the largest penalty of which is 0.48, and ends at a relative gap of $5.5\times10^{-7}$.

\subsection{Long Horizons}\label{sec:exp-long}

\begin{figure}[t]
\centering
\includegraphics[width=\textwidth]{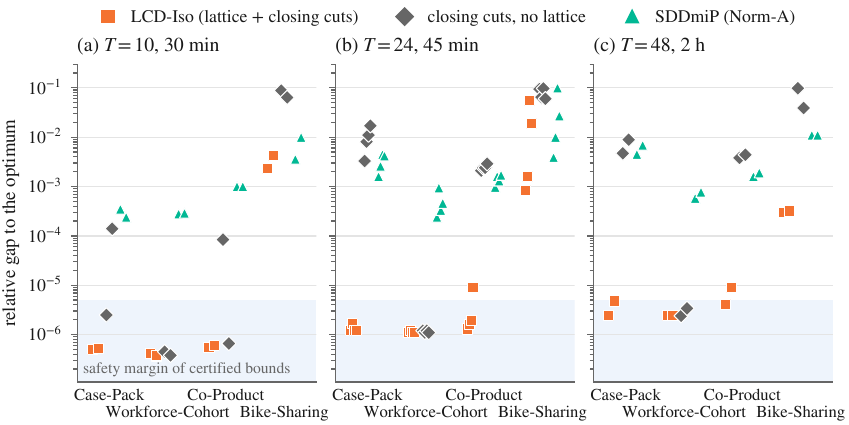}
\caption{\coloronline Relative Gap to the Optimum of \LCDiso, Closing Cuts Without the Lattice, and SDDmiP on 32 Long-Horizon Instances with $T=10$ (a), $T=24$ (b), and $T=48$ (c)}\label{fig:long}
\figurenotes{Each marker represents one instance and one method. The certified bounds of our methods carry a safety margin and cannot come closer to the optimum than the shaded band.}
\end{figure}

We next investigate whether lattice cuts yield tight certified bounds on long horizons. We consider 32 instances with one item, cohort, recipe, or station ($n=2$) from four families, with 10, 24, and 48 stages and between 4 and 36 cosets per stage. The time limits are 30 minutes, 45 minutes, and 2 hours, respectively, and the optimal values are known from dynamic programming. Figure~\ref{fig:long} shows the relative gap of the certified bound of \LCDiso, of closing cuts without the lattice, and of SDDmiP.

\LCDiso\ certifies a bound within $5\times10^{-6}$ of the optimum on 22 of the 32 instances, namely, on all 16 Case-Pack and Workforce-Cohort instances and on 6 of the 8 Co-Product instances, and the other two Co-Product instances are within $10^{-5}$. The horizon has little effect. In particular, at 48 stages, the gaps of \LCDiso\ on these instances are still below $10^{-5}$.

The comparison with closing cuts without the lattice shows what the lattice contributes, because the two runs differ in the lattice only. The bound of \LCDiso\ is higher than that of the same algorithm without the lattice on 21 of the 32 instances and tied with it on the other 11 (Table~\ref{tab:scope}). On the 12 Case-Pack and Co-Product instances with 24 and 48 stages, the gap without the lattice lies between $2.1\times10^{-3}$ and $1.7\times10^{-2}$, whereas the gap of \LCDiso\ lies between $1.2\times10^{-6}$ and $9.1\times10^{-6}$ (Table~\ref{tab:d-long-one}). Therefore, on these instances, the lattice lowers the gap by factors between about 300 and 14,000. On the eight Workforce-Cohort instances, both variants reach the optimum to within the safety margin, and the lattice replaces the closing cuts: \LCDiso\ generates none, and the same algorithm without the lattice generates up to 527 on an instance. Over the 32 runs, \LCDiso\ generates 86 closing cuts in total, all of them on Co-Product instances, and the 16 Case-Pack and Workforce-Cohort instances reach the optimum with lattice cuts alone. In contrast, the same algorithm without the lattice generates significantly more closing cuts, 50,509 in total, and comes within $5\times10^{-6}$ of the optimum on only 10 instances. This is the effect described after Theorem~\ref{thm:closure}: without the lattice, closing cuts have to bridge the periodic oscillation at every visited state.

On the eight Bike-Sharing instances, the gaps of \LCDiso\ remain between $3\times10^{-4}$ and $6\times10^{-2}$, although they are smaller than those of the same algorithm without the lattice, which lie between $4\times10^{-2}$ and $1\times10^{-1}$. Because the representation loss of lattice cuts is zero on these instances (Table~\ref{tab:d-long-one}), this gap is finite-pool loss: the iterations are slow, and the bound is still rising at the time limit. With three hours, the gap of the first 24-stage instance falls from $8\times10^{-4}$ to $1.5\times10^{-5}$.

Figure~\ref{fig:long} also shows SDDmiP. With the settings of Section~\ref{sec:exp-implementation}, its bound is lower than that of \LCDiso\ on all 32 instances, and its gap on the 24 Case-Pack, Workforce-Cohort, and Co-Product instances lies between $2\times10^{-4}$ and $7\times10^{-3}$. However, the distance between the two methods in the figure is not a measure of the strength of their cuts, because it depends on a setting of SDDmiP, which generates a rectified linear cut only when the stage value exceeds its approximation by more than 0.1 (Appendix~\ref{app:methods}). At 10 stages, the gaps of SDDmiP correspond to absolute differences of 0.10 to 0.25, which is the size of this threshold. In additional runs of 90 seconds with thresholds between $10^{-2}$ and $10^{-5}$ (Appendix~\ref{app:add-sddmip}), the gap of SDDmiP on the 24-stage Workforce-Cohort instance falls from $2.4\times10^{-4}$ to at most $8.1\times10^{-7}$, which is the level of our methods, whereas its gaps on the Case-Pack, Co-Product, and Bike-Sharing instances do not fall. Repeating the long runs on all 24 Case-Pack, Workforce-Cohort, and Co-Product instances and four Bike-Sharing instances with a threshold of $10^{-3}$, a bundle tolerance of $10^{-4}$, and a stage gap of $10^{-9}$ confirms this pattern (Table~\ref{tab:x-long-sddmip}): SDDmiP reaches a relative gap of at most $5\times10^{-6}$ on five of the 24 instances against 22 for \LCDiso, and \LCDiso\ outperforms SDDmiP on 13 instances and ties on 11 (16--11--1 across all 28). SDDmiP was run with a stage gap of $10^{-4}$, and repeating eight of its runs with a stage gap of $10^{-9}$ changed its gaps by a factor of at most 2.5 (Table~\ref{tab:d-sens}).

The better bounds do not translate into better policies on these instances as far as sampling can tell. On 200 sampled paths, the mean costs of the policies of \LCDiso\ and SDDmiP differ by less than the half-width of the 95\% confidence interval on 23 of the 24 instances with 24 and 48 stages. In additional runs with limits of 30 and 45 seconds, we therefore evaluate the policies of \LCD\ and SDDmiP exactly on eight small instances, by a recursion over all states that a policy can reach (Appendix~\ref{app:add-policy}). The policy of \LCD\ is optimal on five of these instances. The policy of SDDmiP, with the settings of Section~\ref{sec:exp-implementation}, is optimal on two of them, the Workforce-Cohort instances, and loses up to 0.08\% of the optimal value on the other three. On the Case-Pack-Truck instance, neither policy is optimal, and that of \LCD\ loses less (1.2\% against 2.6\%). On the Bike-Sharing and Serial-Inventory instances, however, on which the bound of \LCD\ is still far from the optimum after 30 seconds, the policy of \LCD\ is worse than that of SDDmiP (losses of 13\% and 15\% against 1.4\% and 7.6\%).

\subsection{Comparison Under Equal Time Limits}\label{sec:exp-equal}

We next compare \LCD\ with the other methods under equal time limits. In Table~\ref{tab:c1}, we present the bounds and policy costs of four methods after 300 seconds on 12 instances, which comprise fixed-charge and multi-item instances of Group~III and synthetic instances of Group~I with 12 and 96 state variables. For each instance, we report the number of state variables ($n$), the number of stages ($T$), the optimal value where it is known ($V^\star$), the certified lower bound of each method (Certified lower bound), and the mean cost of the policy of each method on sampled paths (Mean policy cost). Figure~\ref{fig:anytime} shows the bound as a function of time for six of these instances, and Table~\ref{tab:d-300} in Appendix~\ref{app:results-short} reports the iterations and cuts behind the bounds of all 12.

\begin{table}[t]
\centering
\caption{Certified Lower Bounds and Mean Policy Costs of \LCD, \LCDiso, SDDmiP, and SDDiP on 12 Instances (Time Limit of 300 Seconds)}\label{tab:c1}
\begin{singlespace}\footnotesize
\setlength{\tabcolsep}{1.7pt}
\begin{tabular*}{\textwidth}{@{\extracolsep{\fill}}lrrr rrrr rrrr@{}}
\toprule
 & & & & \multicolumn{4}{c}{Certified lower bound} & \multicolumn{4}{c}{Mean policy cost} \\
\cmidrule(lr){5-8}\cmidrule(l){9-12}
Instance & $n$ & $T$ & $V^\star$ & \LCD & \LCDiso & SDDmiP & SDDiP & \LCD & \LCDiso & SDDmiP & SDDiP \\
\midrule
Case-Pack-Truck & 2 & 10 & 887.7 & \tablebold{885.2} & 878.0 & 879.9 & 726.6 & \tablebold{891} & 893 & 893 & 930 \\
Case-Pack-Truck & 2 & 24 & 2,058.7 & \tablebold{1,967.6} & 1,943.4 & 1,919.1 & 1,000.3 & \tablebold{2,074} & 2,105 & 2,089 & 2,256 \\
Case-Pack-Truck (seed 2) & 2 & 24 & 1,900.3 & 1,786.6 & \tablebold{1,806.6} & 1,687.8 & 879.6 & \tablebold{1,952} & 1,954 & 1,966 & 2,157 \\
Case-Pack-Truck (4 items) & 8 & 6 & 1,014.8 & 948.8 & \tablebold{954.1} & 923.4 & 582.2 & 1,049 & \tablebold{1,033} & 1,062 & 1,277 \\
Serial-Inventory & 5 & 24 & 344.8 & \tablebold{274.3} & 257.6 & 272.1 & 186.0 & 416 & \tablebold{374} & 375 & 495 \\
Co-Product, 2 recipes & 4 & 24 & --- & 4,439.6 & 4,439.1 & \tablebold{4,441.1} & 4,432.2 & 4,467 & 4,467 & \tablebold{4,463} & 4,470 \\
Rotate & 12 & 24 & --- & \tablebold{714.6} & 683.3 & 663.1 & 679.1 & \tablebold{717} & 770 & 841 & 829 \\
Shear & 12 & 24 & --- & \tablebold{686.3} & 495.1 & 663.6 & --- & 1,230 & --- & --- & 1,770 \\
Shift & 96 & 64 & --- & 13,925.8 & \tablebold{14,132.6} & 13,676.1 & 13,498.2 & 14,523 & \tablebold{14,447} & 14,576 & 14,668 \\
Rotate & 96 & 64 & --- & 13,022.9 & \tablebold{13,144.6} & 12,871.2 & 13,130.8 & --- & 13,431 & 13,544 & 13,555 \\
Odd-Cycle & 96 & 64 & --- & \tablebold{12,833.5} & 12,436.2 & 11,150.4 & 8,287.8 & --- & --- & --- & --- \\
Shear & 96 & 64 & --- & \tablebold{12,785.5} & --- & 3,514.5 & --- & --- & --- & --- & 39,453 \\
 
\bottomrule
\end{tabular*}
\end{singlespace}
\figurenotes{All runs use one thread, and the time limit of 300 seconds includes model construction. \boldnote{Bold numbers indicate the highest bound and, where all four policies were evaluated, the lowest mean cost (before rounding). }The mean policy costs are computed on common sampled paths (100 paths for Case-Pack-Truck with $n=2$, and 32 paths otherwise). A dash indicates that no bound was obtained, that the optimal value is unknown, or that the policy evaluation did not finish within 90 seconds.}
\end{table}

\begin{figure}[t]
\centering
\includegraphics[width=\textwidth]{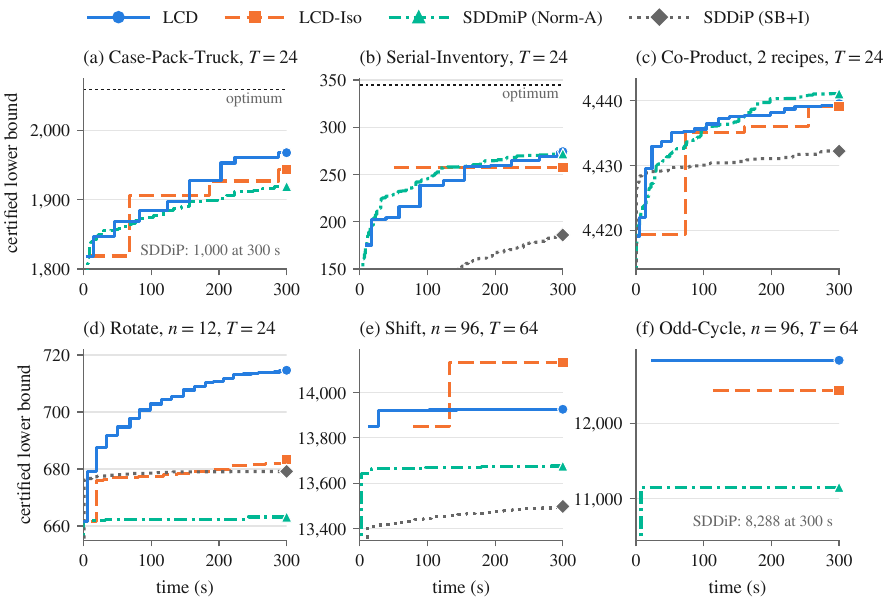}
\caption{\coloronline Certified Lower Bound over Time of \LCD, \LCDiso, SDDmiP, and SDDiP on Six Instances of Table~\ref{tab:c1}}\label{fig:anytime}
\figurenotes{The vertical axes are cut off from below. A method whose final bound lies below the axis is listed in the panel with its value at 300 seconds.}
\end{figure}

\LCD\ has the highest bound on 7 of the 12 instances, \LCDiso\ on 4, and SDDmiP on 1. \LCD\ outperforms SDDmiP on 11 instances and SDDiP on 11 (11--0--1 against SDDmiP with threshold $10^{-3}$ and 12--0--0 against SDDiP with Lagrangian cuts). We observe three patterns in Figure~\ref{fig:anytime}. On the instances with fixed charges (panels (a) and (b)), \LCD\ and SDDmiP rise at a similar rate at first, and \LCD\ pulls ahead after one to two minutes. Although the iterations of \LCD\ are significantly more expensive (8 against 48 in panel (a)), each one raises the bound more. On Rotate with 12 state variables (panel (d)), the condition of Proposition~\ref{prop:phasevalue} holds and lattice cuts are exact. Accordingly, \LCD\ reaches a bound of 714.6 with no closing cut, whereas the other methods stay between 663 and 683. The mean cost of the policy of \LCD\ on 32 sampled paths is 717.3. This Monte Carlo estimate is consistent with a small gap. On two of the four instances with 96 state variables, Shift (panel (e)) and Rotate, \LCDiso\ needs more than a minute for its first iteration and then overtakes \LCD. We believe this is because the mechanisms that make the iterations of \LCD\ fast, namely, intercepts computed only on demand and postponed closing cuts, leave each iteration less complete, which does not always pay off with many state variables. SDDiP iterates quickly. However, its bounds are weak whenever batches matter, as in panel (a), where it reaches 1,000 against 1,968 for \LCD.

The policy costs follow the bounds only loosely. On the eight instances on which all four policies were evaluated, \LCD\ has the lowest mean cost on four, \LCDiso\ on three, and SDDmiP on one. Where the bounds differ significantly, so do the policies: on Rotate with 12 state variables, the policy of \LCD\ costs 717 and the others cost between 770 and 841. We note that the runs of Table~\ref{tab:c1} are single runs. In additional runs of \LCD\ with three seeds on the 12 instances of Table~\ref{tab:c1} and on six instances at 180 seconds (Appendix~\ref{app:add-lattice}), the standard deviation of the bound over the seeds is below 0.6\% of its mean on Table~\ref{tab:c1} and up to 2.7\% at 180 seconds. Therefore, the small differences in Table~\ref{tab:c1} (e.g., between \LCD\ and SDDmiP on the Co-Product instance and between \LCD\ and \LCDiso\ on the Case-Pack-Truck instance with four items, which are 0.03\% and 0.6\% of the bound) are within the variation of a single run and do not rank the methods.

\subsection{Scope and Limits}\label{sec:exp-regimes}

Finally, we examine the cases in which the lattice does not help. Table~\ref{tab:scope} summarizes pairwise comparisons on all instance sets, including those on which our method loses. For each instance set, we report the number of instances or configurations (Number), the time limit, and the numbers of wins, ties, and losses on the certified bound of the method named in the panel against the same algorithm without the lattice and against each of the other methods. We identify three limits, and we then examine the effect of the lattice within \LCD.

\begin{table}[t]
\centering
\caption{Wins, Ties, and Losses of \LCDiso\ (Panel A) and \LCD\ (Panel B) on the Certified Bound}\label{tab:scope}
\begin{singlespace}\footnotesize
\setlength{\tabcolsep}{3pt}
\begin{tabular*}{\textwidth}{@{\extracolsep{\fill}}lrlcccc@{}}
\toprule
Instance set & Number & Time limit & vs.\ no lattice & vs.\ \LCDiso & vs.\ SDDmiP & vs.\ SDDiP \\
\midrule
\multicolumn{7}{@{}l}{\emph{Panel A: \LCDiso}}\\
One item, $T\in\{10,24,48\}$ & 32 & 30 min--2 h & 21--11--0 & & 32--0--0 & \\
Multiple items, $T\in\{10,24\}$ & 11 & 30--45 min & 5--4--2 & & 10--0--1 & \\
Truck fixed charge, $T\in\{10,24\}$ & 4 & 30--45 min & 4--0--0 & & 1--0--3 & \\
Full size, $T=6$, $n$ up to 144 & 15 & 15 min & 3--4--8 & & 8--1--6 & \\
\multicolumn{7}{@{}l}{\emph{Panel B: \LCD}}\\
Instances of Table~\ref{tab:c1} & 12 & 300 s & 2--1--9 & 8--0--4 & 11--0--1 & 11--0--1 \\
Harder instances of Panel A & 48 & 45 s & 16--8--24 & 39--3--6 & 33--3--12 & 47--0--1 \\
Catalog of 163 configurations & 163 & 20--25 s & 24--14--53 & 68--3--20 & 46--5--50 & 85--1--5 \\
\bottomrule
\end{tabular*}
\end{singlespace}
\figurenotes{Entries are wins--ties--losses of the method named in the panel. Two bounds tie if they differ by less than $10^{-4}$ in relative terms. A run without a bound loses against a run with one. In the catalog, neither method of a pair has a bound on 72, 72, 62, and 72 configurations, respectively, and these configurations are not counted. The column ``vs.\ no lattice'' compares the method of each panel with the same algorithm with the trivial lattice. On the eight Workforce-Cohort instances of the first row, the comparison with SDDmiP depends on its threshold (Section~\ref{sec:exp-long}; Table~\ref{tab:x-long-sddmip}).}
\end{table}

\emph{Nonconvexity that is not periodic in the chosen lattice.} In Case-Pack-Truck, the fixed charge per truck depends on the number of packs ordered, which the lattice of stock modulo pack size does not capture. On the 10-stage instance, the representation loss of lattice cuts is 1.6\% of the optimal value, and closing cuts must do the remaining work. With isotropic closing cuts and long time limits, \LCDiso\ outperforms the same algorithm without the lattice on all four truck instances. However, it outperforms SDDmiP on only one of them (Panel~A). The directional and pooled closing cuts of \LCD\ change this outcome at 300 seconds (Table~\ref{tab:c1}). A finer lattice also helps. In an exact evaluation of the 10-stage instance, the loss at $\rho=0$ falls from 1.63\% for the stock modulo the pack size to 0.96\% for the inventory position modulo 8 and to zero for the inventory position modulo 16. We note that this is a computational finding, because the fixed charge is outside the model of Theorem~\ref{thm:optimal-period}.

\emph{Too many cosets.} A lattice cut pays off when the cosets that it has priced are visited again. With one item, there are a handful of cosets. With $K$ items, the number of cosets is a product over the items and soon exceeds the number of states that a run ever visits (e.g., Airline-RM with 144 state variables has more than $10^{55}$ cosets). An intercept computed for one coset is then rarely used a second time, and \LCD\ behaves like a method with Lagrangian cuts and a more expensive subproblem. In instrumented 45-second runs of \LCD, the computed intercepts are used again between 69 and 12,558 times in a run on six instances with at most 64 cosets per stage and at most three times in a run on four instances with 1,296 cosets or more (Appendix~\ref{app:add-reuse}). On the 15 full-size six-stage instances, \LCDiso\ outperforms SDDmiP on 8, ties on 1, and loses on 6. On these instances, the lattice does not pay: \LCDiso\ outperforms the same algorithm without the lattice on 3 instances, ties on 4, and loses on 8. We note, however, that the number of cosets alone does not decide the outcome. Among the 48 instances of the 45-second runs, \LCD\ outperforms SDDmiP on 21 of the 27 instances with at most 64 cosets per stage, on 6 of the 7 with 65 to 4,096, on 3 of the 11 with 4,097 to $10^{8}$, and on all 3 with more (Table~\ref{tab:d-buckets}).

\emph{Short time limits in high dimension.} An iteration of \LCD\ solves MIPs with a free incoming state to a tight gap, which costs more than an iteration of SDDmiP. Within 45 seconds, \LCD\ outperforms SDDmiP on 33 of 48 instances, and on the 12 instances on which it loses, which have 2 to 24 state variables, \LCD\ completes between one and eight iterations. The policies do not follow the bounds on these 48 instances: on the 47 instances on which both policies were evaluated, the mean cost of the policy of \LCD\ is lower than that of SDDmiP on 16, tied on 4, and higher on 27 (Table~\ref{tab:d-buckets}). The catalog run pushes the short limit to the extreme, with 163 configurations that have up to 6,144 state variables and 2,048 stages under 20- or 25-second limits. On 62 of these configurations, none of the four methods obtains a bound because building the models already exhausts the limit. Where \LCD\ has a bound, it has completed at most one iteration in 75 of 91 runs. Therefore, this run measures start-up cost more than the strength of the cuts. \LCD\ outperforms SDDmiP on 46 configurations and loses on 50. \LCD\ wins on Co-Product (6 wins, 3 ties, no loss) and Bike-Sharing (9 wins) and loses on the multi-item configurations of Case-Pack (9 losses) and Pooled-Capacity (8 losses of 9).

\emph{The lattice within \LCD.} In Panel~B of Table~\ref{tab:scope}, \LCD\ outperforms the same algorithm without the lattice on 2 of 12 instances at 300 seconds, 16 of 48 at 45 seconds, and 24 of 91 in the catalog, and we made additional runs under both short and longer time limits (Appendix~\ref{app:add-lattice}). On 16 instances with a limit of 45 seconds, the lattice raises the bound of \LCD\ on 5, lowers it on 8, and leaves it tied on 3. The five are the one-item Case-Pack instance, the one-recipe Co-Product instance, one Serial-Inventory instance, Shift, and Rotate, on all of which lattice cuts alone are exact or nearly exact. On six instances with multiple items, recipes, or stations and a limit of 180 seconds, the mean bound over three seeds is higher without the lattice on four, by more than twice the standard deviation on three of them, although these four instances have only 27 to 48 cosets per stage. On ten instances with 1,296 cosets or more and a limit of 60 seconds, disabling the lattice raises the bound on seven. In the same runs, \LCD\ outperforms \LCDiso\ on 13 of the 16 instances with the lattice and on 13 without it. This suggests that, under short time limits, the advantage of \LCD\ over \LCDiso\ comes from the components in which the two algorithms differ, namely, the directional and pooled closing cuts and the implementation of Section~\ref{sec:algorithm}, and not from the lattice. With a limit of 1,800 seconds across 25 instances (Table~\ref{tab:x-crossover}), the comparison on the 11 multiple-item instances of Panel~A moves from 3--1--7 at 180 seconds to 2--6--3 at 1,800 seconds with roughly half as many closing cuts, and on the 48-stage Case-Pack and Co-Product instances \LCD\ reaches the optimum with 7 and 2,739 closing cuts against 5,453 and 9,498 without the lattice, whereas on the four truck instances it ends at 0--2--2.

Taken together, our results suggest using lattice cuts when lattice cuts alone are exact or nearly exact, that is, when the representation loss $V^\star-W_\Gamma(0)$ of Section~\ref{sec:ideal-pool} is small, and when the time limit allows the cosets to recur. This is the case for the one-item Case-Pack, Workforce-Cohort, and Co-Product instances and for Shift and Rotate. In this regime, the lattice lowers the gap of the same algorithm by up to four orders of magnitude or makes its closing cuts unnecessary. With multiple items, recipes, or stations, the results are mixed even when a stage has only a few dozen cosets, and with thousands of cosets or more, the lattice lowers the bound more often than it raises it within our time limits. For such models, we expect that coarser lattices or lattices on a subset of the state variables are a natural compromise, which we leave for future work.

\section{Conclusion}\label{sec:conclusion}

We introduce lattice cuts for multistage stochastic integer programs whose nonconvexity comes from batches. A lattice cut has one slope and one intercept per coset of a lattice that is computed from the constraint matrices of the stages, and it reproduces the periodic part of the cost-to-go function globally without adding a binary variable per cut. Closing cuts handle the nonconvexity that remains inside a coset. We show that a linear program characterizes the bound that a lattice and a penalty cap can deliver and separates the gap of a run into a representation loss and a finite-pool loss, and we prove that the resulting algorithm converges finitely with probability one.

In our computational study, the lattice removes nearly all of the representation loss of affine cuts on instances that can be checked exactly. On horizons of up to 48 stages, it brings the certified bound of the same algorithm closer to the optimum by two to four orders of magnitude on case-pack and co-production instances, with significantly fewer closing cuts than the algorithm needs without the lattice. The study also shows the limits of our approach, which requires that cosets recur, a condition that fails when many items multiply the number of cosets, and whose iterations are too expensive for very short time limits in high dimension. Under these conditions, our algorithm often attains a higher bound without the lattice than with it.

These limits suggest three directions for future work. One direction is to choose coarser lattices deliberately, trading representation loss for recurrence, with the linear program of Theorem~\ref{thm:root-price} as a guide. The other two are to detect fixed-charge multiples automatically and include them in the lattice, and to extend our bound characterization to risk-averse objectives.

\clearpage
\appendix
\pdfbookmark[0]{Appendices}{appendix-start}
\section*{Appendices}
\renewcommand{\theequation}{\Alph{section}.\arabic{equation}}
\renewcommand{\theHequation}{\Alph{section}.\arabic{equation}}
\renewcommand{\thetable}{\Alph{section}.\arabic{table}}
\renewcommand{\theHtable}{\Alph{section}.\arabic{table}}
\renewcommand{\thefigure}{\Alph{section}.\arabic{figure}}
\renewcommand{\theHfigure}{\Alph{section}.\arabic{figure}}
\renewcommand{\thetheorem}{\Alph{section}.\arabic{theorem}}
\renewcommand{\theHtheorem}{\Alph{section}.\arabic{theorem}}
\renewcommand{\theproposition}{\Alph{section}.\arabic{proposition}}
\renewcommand{\theHproposition}{\Alph{section}.\arabic{proposition}}
\renewcommand{\thelemma}{\Alph{section}.\arabic{lemma}}
\renewcommand{\theHlemma}{\Alph{section}.\arabic{lemma}}
\renewcommand{\thecorollary}{\Alph{section}.\arabic{corollary}}
\renewcommand{\theHcorollary}{\Alph{section}.\arabic{corollary}}
\renewcommand{\theremark}{\Alph{section}.\arabic{remark}}
\renewcommand{\theHremark}{\Alph{section}.\arabic{remark}}
\newcommand{\ecresetcounters}{\setcounter{equation}{0}\setcounter{table}{0}\setcounter{figure}{0}\setcounter{theorem}{0}%
  \setcounter{proposition}{0}\setcounter{lemma}{0}\setcounter{corollary}{0}\setcounter{remark}{0}}

\ecresetcounters
\section{Additional Results}\label{ec:additional}

In this appendix, we state the results to which the main paper refers without stating them in full. Their proofs are given in Appendix~\ref{ec:proofs}. Throughout the appendices, lattice cuts, closing cuts, and isotropic closing cuts are those of Definition~\ref{def:lattice-cuts}, and \textnormal{LCD-Iso} is the basic variant described in Section~\ref{sec:root}. We use $g_t(s)$ to denote the coset $[s]_t$ of a state $s$.

\subsection{Results for Section~\ref{sec:lattice-cuts}}

\paragraph{Stabilization of the backward recursion.}
We first verify that the backward pullbacks remain bounded in $T$ under a stationary matrix $A\in\mathbb Z^{n\times n}$ and remainders $Cs\bmod m$ ($C\in\mathbb Z^{p\times n}$, $m\ge 2$). To this end, we stack $C,CA,\ldots,CA^{h-1}$ into $\mathcal O_h$ over $R=\mathbb Z/m\mathbb Z$, and we set $K_h=\{s\in\mathbb Z^n:\mathcal O_hs\equiv 0\pmod m\}$, $G_h=\mathbb Z^n/K_h$, and $M_h=\operatorname{row}_R(\mathcal O_h\bmod m)$ ($K_0=\mathbb Z^n$, $M_0=\{0\}$). The following lemma provides a minimal coset label and shows that the recursion stabilizes after finitely many steps.

\begin{lemma}\label{lem:observation-phases}
The following statements hold.
\begin{enumerate}[label=(\roman*),leftmargin=2em,itemsep=1pt,topsep=3pt]
\item If an additive label $\phi$ on $\mathbb Z^n$ determines $\mathcal O_hs\bmod m$, then $\ker\phi\subseteq K_h$ and the quotient map $s\mapsto s+K_h$ factors uniquely through $\phi$. Therefore, no such label has fewer than $|G_h|$ values.
\item Given the Smith form $U\mathcal O_h V$ with invariant factors $e_1\mid\cdots\mid e_r$, let $q_i=m/\gcd(m,e_i)$. Then, a minimal residue label and the exact index are
\begin{equation}\label{eq:phaseorder}
 \lambda_h(s)=\bigl((V^{-1}s)_i\bmod q_i\bigr)_{i=1}^r,\qquad
 |G_h|=[\mathbb Z^n:K_h]=\prod_{i=1}^r q_i .
\end{equation}
\item The stabilization depth $h_*=\min\{h\ge 0:K_{h+1}=K_h\}$ is the first $h$ with $\operatorname{row}_R(CA^h\bmod m)\subseteq M_h$, which implies $K_j=K_{h_*}$ for all $j\ge h_*$ and $[K_h:K_{h+1}]\mid[K_{h-1}:K_h]$.
\item If $d_{\rm sc}$ is the least degree of a monic polynomial $P$ over $R$ with $CP(A)=0$ and $\nu$ is the minimum number of generators of $M_{h_*}$, then $h_*\le d_{\rm sc}\le \nu\le\min\{n,ph_*\}$ (with equality when $p=1$).
\end{enumerate}
\end{lemma}

\paragraph{Sharp form of Theorem~\ref{thm:optimal-period}.}
In the setting of \eqref{eq:period-model}--\eqref{eq:observer-kernel} and under Assumption~\ref{assump:period}, the periodic part $R_t$ is bounded below by a function $B_t$ that is positive outside $K_t$. We define $B_{T+1}\equiv 0$ and $B_t(x)=\sum_{w\in\Omega_t}p_{tw}B_{tw}(x)$, where
\begin{equation}\label{eq:observer-stop}
 B_{tw}(x)=\epsilon_{tw}\quad (C_{tw}x\ne 0),\qquad
 B_{tw}(x)=\min\{\epsilon_{tw},\,B_{t+1}(A_{tw}x)\}\quad (C_{tw}x=0).
\end{equation}
The following theorem establishes this lower bound, shows that it is sharp, and identifies $K_t$ as the period lattice of $R_t$.

\begin{theorem}\label{thm:optimal-period-full}
Under the neutral-action and positive-price conditions on \eqref{eq:period-model}, we have
\begin{equation}\label{eq:value-period}
 R_t(x)\ge B_t(x),\qquad
 \{x\in\mathbb Z^n:R_t(x)=0\}=K_t
 =\bigl\{h\in\mathbb Z^n:R_t(x+h)=R_t(x)\ \forall\,x\in\mathbb Z^n\bigr\},
\end{equation}
and a finite-action instance of \eqref{eq:period-model} attains $R_{tw}\equiv B_{tw}$ on $\mathbb Z^n$. Therefore, $[\mathbb Z^n:K_t]$ is the minimum lattice index whose cosets can represent $R_t$ on $\mathbb Z^n$.
\end{theorem}

\noindent Theorem~\ref{thm:optimal-period} of the main paper consists of the second and third statements of \eqref{eq:value-period}, and the proof in Appendix~\ref{ec:proofs} establishes Theorem~\ref{thm:optimal-period-full}. The following corollary bounds from below the loss of cuts on a lattice that is coarser than the period lattice.

\begin{corollary}\label{cor:price-period}
If a sublattice $L\subseteq\mathbb Z^n$ contains a vector $h\notin K_t$ of order $k\ge 2$ in $\mathbb Z^n/K_t$ with $0,h,kh\in D_t$ obeying \eqref{eq:period-model}, then any expected upper envelope $\underline R_t$ (net of linear terms) of valid coset or closing cuts on $L$ with directional or $\ell_1$ penalties $\le\rho$ satisfies
\begin{equation}\label{eq:period-r9}
 R_t(h)-\underline R_t(h)\ \ge\
 \bigl[B_t(h)-2(1-1/k)\rho\|h\|_1\bigr]_+ .
\end{equation}
\end{corollary}

\subsection{Results for Section~\ref{sec:cuts}}

\paragraph{Gain from refining the lattice.}
Let $\mathcal M(S,\hat s)$ denote the polytope of distributions $\mu$ on $S\subseteq D_t$ ($\hat s\in S$) with conditional barycenter $\sum_{z\in S}\mu_z z=\hat s$, and we define the barycentric envelope
\begin{equation}\label{eq:mixing}
 E_S(\hat s)=\min_{\mu\in\mathcal M(S,\hat s)}
                       \sum_{z\in S}\mu_z F(z).
\end{equation}
By LP duality, $E_S(\hat s)=\max_{\pi,b}\{\pi^\top\hat s+b:\pi^\top z+b\le F(z)\ \forall\,z\in S\}$ is the tightest value of a lattice cut at $\hat s$ over $S$, where $\mu$ represents geometric copy interpolation rather than the scenario probability $p_{tw}$. The following theorem gives the exact gain that is obtained by refining the lattice.

\begin{theorem}\label{thm:mixture}
For full-rank lattices $\Gamma'\subseteq\Gamma\subseteq\mathbb Z^n$, let $S=D_t\cap(\hat s+\Gamma)$ and $S'=D_t\cap(\hat s+\Gamma')$. Then, we have $E_{D_t}(\hat s)\le E_S(\hat s)\le E_{S'}(\hat s)\le F(\hat s)$, and for every affine minorant $u\le F$ on $S$ attaining $u(\hat s)=E_S(\hat s)$, we have $E_{S'}(\hat s)-E_S(\hat s)=\min_{\mu\in\mathcal M(S',\hat s)}\sum_{z\in S'}\mu_z(F(z)-u(z))$, which is strictly positive if and only if every optimizer for $E_S(\hat s)$ puts positive mass on $S\setminus S'$.
\end{theorem}

\paragraph{Least penalty when the slope is optimized as well.}
Let $f$ and $d$ be as in Theorem~\ref{thm:closure}. On a set $S\ni\hat s$, the contact residual
\begin{equation}\label{eq:residual}
 R_S(\rho)=\max_{z\in S}\bigl\{f(\hat s)-f(z)-\rho d(z)\bigr\}
\end{equation}
equals $F(\hat s)-h^\rho(\hat s)$ on $S=D_{\hat g}$ and reaches zero at $\rho_{\hat g}^*$. Optimizing $(\pi,\rho)$ jointly on $S=D_t\cap(\hat s+\Gamma_t)$ ($\mathcal M=\mathcal M(S,\hat s)$) defines $G_S(\rho)=\max_{\pi\in\mathbb R^n}\min_{z\in S}\{F(z)+\pi^\top(\hat s-z)+\rho d(z)\}$. The following theorem gives the least closing penalty when the slope is optimized as well.

\begin{theorem}\label{thm:local-price}
The function $G_S(\rho)$ is nondecreasing, concave, and piecewise linear on $\rho\ge 0$, and it attains $F(\hat s)$ if and only if $\rho\ge\rho_S^\circ$, where
\begin{equation}\label{eq:joint-penalty}
 \rho^\circ_S=
 \max_{\substack{\mu\in\mathcal M\\ \mathbb E_\mu d(Z)>0}}
 \frac{[F(\hat s)-\mathbb E_\mu F(Z)]_+}{\mathbb E_\mu d(Z)}
\end{equation}
($0$ if $\mathcal M=\{\delta_{\hat s}\}$). When $\rho_S^\circ>0$, the maximum in \eqref{eq:joint-penalty} is attained on at most $\dim(\operatorname{aff} S)+1$ points of $S\setminus\{\hat s\}$.
\end{theorem}

Because the function $t\mapsto\max\{\rho_i^+ t,-\rho_i^- t\}$ is sublinear, the penalty \eqref{eq:aniso-penalty} satisfies $\Phi(z;a,\rho^+,\rho^-)-\Phi(s;a,\rho^+,\rho^-)\le\Phi(z;s,\rho^+,\rho^-)\le\bar\rho\|z-s\|_1$ for any anchor $a\in D_t$ whenever $\max_i\max\{\rho_i^+,\rho_i^-\}\le\bar\rho$.

\paragraph{Single-stage form of the complete cut family.}
Suppose that every stage $t\ge 2$ admits \emph{all} valid lattice cuts, closing cuts, and isotropic closing cuts with penalties $\rho_i^\pm\le\rho$ against recursive continuations (and that stage~1 is solved at $s_1$). Then, the single-stage supremum at $s$ on $S=D_t\cap(s+\Gamma_t)$ over all anchors $a\in D_t$, slopes $\pi\in\mathbb R^n$, and $\rho_i^\pm\le\rho$ equals
\begin{equation}\label{eq:priced-envelope}
 \min_{\mu\in\mathcal M(S,s)}
 \sum_{z\in S}\mu_z\bigl(F(z)+\rho\|z-s\|_1\bigr),
\end{equation}
which is attained by centering an isotropic closing cut at $a=s$, whereas off-anchor closing cuts are bounded through the sublinearity of $\Phi$. The following corollary bounds the threshold $\rho_\Gamma^\star$ of Theorem~\ref{thm:root-price} explicitly in terms of the stage data.

\begin{corollary}\label{cor:data-cap}
Suppose that $\ell_{tw}\le c_{tw}\le u_{tw}$ and that $H_{T+1}$ bounds the oscillation of $V_{T+1}$ on $D_{T+1}$. Let $H_{tw}=u_{tw}-\ell_{tw}+H_{T+1}+\sum_{j=t+1}^T\sum_{b\in\Omega_j}p_{jb}(u_{jb}-\ell_{jb})$ and $\eta_t(s)=\min\{\mathbb E_\mu\|Z-s\|_1:\mu\in\mathcal M(S_t(s)\setminus\{s\},s)\}$ ($+\infty$ if empty, $H/(+\infty)=0$). Then, we have $\rho_\Gamma^\star\le\max_{t\ge 2,\,w\in\Omega_t,\,s\in D_t}H_{tw}/\eta_t(s)$ ($0$ if $T=1$).
\end{corollary}

\begin{remark}[Hierarchy of Lower Bounds and Error Decomposition]\label{rem:hierarchy}
For any finite cut pool that certifies $L_k$ with penalties $\le\bar\rho$ and valid continuation floors, backward monotonicity gives $W_{\mathbb Z^n}(0)\le W_\Gamma(0)\le W_\Gamma(\bar\rho)\le W_\Gamma(\rho_\Gamma^\star)=V^\star$ and $L_k\le W_\Gamma(\bar\rho)$. These inequalities split $V^\star-L_k=[V^\star-W_\Gamma(\bar\rho)]+[W_\Gamma(\bar\rho)-L_k]$ into the \emph{representation loss} and the \emph{finite-pool loss} (Theorem~\ref{thm:finite-pool-root-gap}).
\end{remark}

\subsection{Results for Section~\ref{sec:algorithm}}

\paragraph{Initial lower bounds.}
\textnormal{LCD} keeps one parametric model per stage. Solving these $T$ models as continuous LPs over $[\ell_\tau,u_\tau]$ yields $m_{\tau w}=\min\{c_{\tau w}(s,a,y):s\in[\ell_\tau,u_\tau],\,(a,y)\in X_{\tau w}^{\mathrm{LP}}(s)\}$ and $E_{T+1}=\min_{[\ell_{T+1},u_{T+1}]}V_{T+1}^{\mathrm{LP}}$, which define the Bellman-expectation continuation floors
\begin{equation}\label{eq:theta-lb}
 E_\tau=\sum_{w\in\Omega_\tau}p_{\tau w}m_{\tau w}+E_{\tau+1},\qquad
 \underline\theta_{\tau w}=m_{\tau w}+E_{\tau+1},\qquad
 \underline\theta_\tau=\min_{w\in\Omega_\tau}m_{\tau w}+E_{\tau+1}
 \quad (\tau=T,\ldots,1).
\end{equation}

\paragraph{Raising the penalties of a closing cut.}
We trigger pointwise closure only when a recorded coset copy $(z,U_z)$ ($z\in D_{\hat g}$, $U_z\ge F_{tw}(z)$) satisfies $L_{\hat s}-(\pi^\top\hat s+U_z-\pi^\top z)>\varepsilon_{\mathrm{close}}(1+|U_{\hat s}^{\min}|)$, which avoids false triggers from solver tolerances. On first visits in $n\ge 3$, a progress gate defers closure if the new lattice cut closes a fraction $\ge\gamma\in(0,1]$ of $L_{\hat s}-q_{tw}^-(\hat s)$, whereas closure is enforced on repeat visits. When closure is triggered, structural scales $w_i\ge 1$ that are read from the transition rows (e.g., $w_i=1$ for unit stock and $w_i=p_i$ for size-$p_i$ packs) guide the directional lifting under \eqref{eq:aniso-penalty}. Starting from $(\rho^{+,0},\rho^{-,0})$ seeded by violating copies, step $m$ minimizes $f(z)+\Phi(z;\hat s,\rho^{+,m},\rho^{-,m})$ on $D_{\hat g}$ ($f(z)=F_{tw}(z)-\pi^\top z$). If $\delta_m=f(\hat s)-f(z_m)-\Phi(z_m;\hat s,\rho^{+,m},\rho^{-,m})>0$, then, with $d_{w,m}=\sum_{j=1}^n w_j|z_{m,j}-\hat s_j|$, the penalties are updated to $\rho_i^{\pm,m+1}=\rho_i^{\pm,m}+w_i(\delta_m/d_{w,m})\mathbf 1\{\pm(z_{m,i}-\hat s_i)>0\}$. For isotropic closing cuts, that is, $\rho_i^\pm=\rho$, this update reduces to Dinkelbach's rule $\rho_{m+1}=[f(\hat s)-f(z_m)]/\|z_m-\hat s\|_1$ \citep{dinkelbach1967}. The following proposition shows that the lifting terminates after finitely many solves and that stopping it early still yields a valid cut.

\begin{proposition}\label{prop:dinkelbach}
Fix $F_{tw}$, $\pi\in\mathbb R^n$, and $\hat s\in D_{\hat g}$. Under isotropic closing cuts from $\rho_0\in[0,\rho_{\hat g}^*]$, exact coset solves generate $\rho_0<\rho_1<\cdots<\rho_K=\rho_{\hat g}^*$ and verify $h^{\rho_K}(\hat s)=F_{tw}(\hat s)$ in $K+1\le 1+|\{\|z-\hat s\|_1:z\in D_{\hat g}\setminus\{\hat s\}\}|\le 1+\sum_{i=1}^n(u_{t,i}-\ell_{t,i})$ solves. Under coordinate-scaled closing cuts ($w_i\ge 1$), exact solves reach $h(\hat s)=F_{tw}(\hat s)$ in $\le 1+|\{z\in D_{\hat g}:f(z)<f(\hat s)\}|\le|D_{\hat g}|$ solves while keeping $\rho_i^\pm=0$ on unviolated rays. Any early stop $(\rho^+,\rho^-)\in\mathbb R_+^{2n}$ with certified dual bounds $b_{\hat g},b_0$ yields a valid cut on $D_t$.
\end{proposition}

\paragraph{Stored intercepts and coset indicators.}
Initialized with $\underline\theta_t$, the child pools $\mathcal C_{tw}$ ($t\ge 2$) define the multicut approximation
\begin{equation}\label{eq:multicut}
 \underline V_t(y)=\sum_{w\in\Omega_t} p_{tw}\max_{h\in\mathcal C_{tw}} h(y).
\end{equation}
\textnormal{LCD} stores a three-level intercept hierarchy $(b_0,\{b^{(1)}_{r_1}:r_1\in\mathcal R_1\},\{b_g:g\in\mathcal G\})$ over global, first-invariant-factor ($(V_t^{-1}z)_1\equiv r_1\pmod{q_{t,1}}$), and coset fibers. It sets $\bar b(g)=\max(\{b_0\}\cup\{b^{(1)}_{g_1}:g_1\in\mathcal R_1\}\cup\{b_g:g\in\mathcal G\})$ and merges lattice cuts with a repeated slope in place through componentwise maxima. With the coset indicator $\mu_g(y)$, each parent cut reads
\begin{equation}\label{eq:parent-row}
 \theta_w \ \ge\ \pi^\top y + b_0 + \sum_{g\in\mathcal G:\,\bar b(g)>b_0}\bigl(\bar b(g)-b_0\bigr)\mu_g(y) - \Phi(y;\hat s,\rho^+,\rho^-) .
\end{equation}
Using the Smith remainders $(V_t^{-1}y)_j=q_{t,j}\xi_j+r_j$ ($0\le r_j<q_{t,j}$), stages with $[\mathbb Z^n:\Gamma_t]\le M_{\mathrm{oh}}$ share one-hot selectors $\delta_g\in\{0,1\}$ ($\sum_{g\in\mathcal G_t}\delta_g=1$, $\mu_g=\delta_g$). When $[\mathbb Z^n:\Gamma_t]>M_{\mathrm{oh}}$, \textnormal{LCD} expands $r_j$ into $B_t=\sum_{j=1}^{r_t}\lceil\log_2 q_{t,j}\rceil$ shared bits $\zeta_{j\ell}\in\{0,1\}$ and pools one $\mu_g\in[0,1]$ per active coset through $\mu_g+\operatorname{Ham}(\zeta,\operatorname{bits}(g))\ge 1$. The following proposition establishes the validity of the LP floors, of the merging of intercepts, and of the pooled encodings.

\begin{proposition}\label{prop:phase-hierarchy}
(i)~In \eqref{eq:theta-lb}, $Q_{\tau w}\ge\underline\theta_{\tau w}\ge\underline\theta_\tau$ and $V_\tau\ge E_\tau\ge\underline\theta_\tau$ on $D_\tau$, and $b_{\mathrm{LP}}(\pi)$ is a valid global floor. (ii)~Across $D_t$, $\pi^\top y+\bar b([y]_t)\le F_{tw}(y)\le Q_{tw}(y)$; in-place maxima and partition lifting $b_0\leftarrow\max\{b_0,\min_{g\in\mathcal G_t}\bar b(g)\}$ preserve validity and tighten one-hot LP relaxations without new rows. (iii)~For lattice cuts ($\Phi\equiv 0$), shared one-hot selectors attain the simplex hull with zero per-cut variables, and pooled $\mu_g\in[0,1]$ under $\zeta\in\{0,1\}^{B_t}$ evaluates \eqref{eq:parent-row} exactly on $D_t$.
\end{proposition}

\paragraph{Effect of the merging tolerance.}
Proposition~\ref{prop:phase-hierarchy}(ii) covers the merging of intercepts that belong to the same slope. The implementation also attaches an intercept to a slope that differs slightly from the slope for which the intercept was computed (Section~\ref{sec:cut-subproblems}). The following proposition bounds the effect of such differences on the cuts and on the bound $L_k$. For $t\ge 2$, let $\sigma_t=\max_{y\in D_t}\|y\|_1$.

\begin{proposition}\label{prop:merge-tolerance}
Let $\varepsilon\ge 0$. Suppose that the pools are generated as in Algorithm~\ref{alg:lcd}, except that every intercept of a cut with slope $\pi$ is a lower bound on the program that defines it with a slope $\pi'$ in place of $\pi$, where $\|\pi-\pi'\|_\infty\le\varepsilon$ and $\pi'$ may differ from intercept to intercept. Then, every cut $h\in\mathcal C_{tw}$ with $t\ge 2$ satisfies
\[
 h(y)\ \le\ Q_{tw}(y)+\varepsilon\sum_{j=t}^{T}\sigma_j\qquad\text{for all } y\in D_t,
\]
and the bound of every iteration satisfies $L_k\le V^\star+\kappa$ with $\kappa=\varepsilon\sum_{t=2}^{T}\sigma_t$.
\end{proposition}

For $\varepsilon=0$, the proposition reduces to the validity of the cuts. The constant $\kappa$ does not depend on the number of cuts or of iterations. It grows linearly with $\varepsilon$ and with the horizon, because an error in a cut of stage $t+1$ enters the intercepts of the cuts of stage $t$.

\paragraph{Gain at the root.}
For a root outcome $w\in\Omega_1$ with feasible successors $Y_w\subseteq D_2$, prior objective $J_w(y)$, minimum $j_w=\min_{Y_w}J_w$, and stage-2 update $\Delta(y)\ge 0$, the certified root-bound gain satisfies
\begin{equation}\label{eq:rootgain}
 L^{+} - L \ =\ \sum_{w\in\Omega_1} p_{1w}\min_{y\in Y_w}
        \bigl\{J_w(y)-j_w+\Delta(y)\bigr\},
\end{equation}
which is positive under $w$ if and only if $\Delta(y)>0$ across $\arg\min_{Y_w}J_w$ (Corollary~\ref{cor:root-gain-coverage}).

\paragraph{Finite-pool loss.}
We suppress the iteration index $k$. After stage~1 is solved at $s_1$, let $q_{tw}(s)=\max_{h\in\mathcal C_{tw}}h(s)$ ($t\ge 2$), $q_{1w}(s_1)=F_{1w}(s_1)$, $u_t(s)=\sum_{w\in\Omega_t}p_{tw}q_{tw}(s)$, $u_{T+1}=V_{T+1}$, and $L=u_1(s_1)$. With the active slopes $\pi_{tws}$ ($\pi_{1ws_1}=0$), we define on each arc $e=(t,w,s,z,y)$ of \eqref{eq:occupation} the reduced cost $g_e=g_{twszy}=c^0_{tw}(z,y)+u_{t+1}(y)+\bar\rho\|z-s\|_1-q_{tw}(s)-\pi_{tws}^\top(z-s)$. The following theorem expresses the finite-pool loss through these reduced costs and provides a test for the equality $L=W_\Gamma(\bar\rho)$.

\begin{theorem}\label{thm:finite-pool-root-gap}
The potentials $(u,q)$ and the active slopes $\pi$ are dual feasible for the occupancy LP of Theorem~\ref{thm:root-price} at cap $\bar\rho$, with dual objective $L$, reduced costs $g_e\ge 0$, and $J(v)+\bar\rho D(v)-L=\sum_e f_e g_e$ for all $v=(d,f)\in\mathcal P_\Gamma$, so that $W_\Gamma(\bar\rho)-L=\min_{v\in\mathcal P_\Gamma}\sum_e f_e g_e$. Let $A_{T+1}=D_{T+1}$ and, backward for $t=T,\ldots,2$, let $B_{tws}=\{z\in S_t(s):\exists\,y\in Y_{tw}(z)\cap A_{t+1},\ g_{twszy}=0\}$ and $A_t=\{s\in D_t:s\in\operatorname{conv}B_{tws}\ \forall\,w\in\Omega_t\}$. Then, the equality $L=W_\Gamma(\bar\rho)$ holds if and only if every $w\in\Omega_1$ admits an optimal successor $y\in Y_{1w}(s_1)\cap A_2$ ($g_{1ws_1s_1y}=0$).
\end{theorem}

\begin{remark}[Computable Finite-Data Certificates]\label{rem:ec-certificates}
The feasible transitions that are retained during a run provide LP and convex-piecewise diagnostics without MIP solves at unvisited states (Appendices~\ref{ec:local-certs} and~\ref{ec:horizon-certs}). These diagnostics are ray cut coverage at $\le\min\{K+1,2R+2\}$ kink neighbors (Proposition~\ref{prop:directional-coverage}), single-stage envelope LPs (Proposition~\ref{prop:finite-pool-certificate}), stage-2 regional root-gain tests (Corollary~\ref{cor:root-gain-coverage}), recorded-transition upper intervals $U_{\mathcal D}\ge W_\Gamma(\bar\rho)$ (Proposition~\ref{prop:record-root-certificate}), and quotient Markov replay on cosets $V^\star-L_k\le J(\mu^k)-L_k$ (Proposition~\ref{prop:phase-policy-replay}).
\end{remark}

\ecresetcounters
\section{Proofs of Statements in the Main Paper}\label{ec:proofs}

In this appendix, we prove the statements of the main paper and of Appendix~\ref{ec:additional}.

\subsection{Proofs for Sections~\ref{sec:model} and~\ref{sec:lattice-cuts}}

\begin{proof}[Proof of Proposition~\ref{prop:compatible}]
We first show that the family \eqref{eq:compatible} is compatible. By construction, we have $\Gamma_t\subseteq\Lambda_t$ and $A_{tw}\Gamma_t\subseteq\Gamma_{t+1}$ for all $t<T$ and $w\in\Omega_t$, so that $\{\Gamma_t\}_{t=1}^T$ is admissible and compatible. Each $\Gamma_t$ has finite index in $\mathbb Z^n$ by backward induction, because $\Gamma_T=\Lambda_T$ has finite index, each pullback $A_{tw}^{-1}\Gamma_{t+1}$ is the kernel of the induced homomorphism $v\mapsto A_{tw}v+\Gamma_{t+1}$ into the finite quotient group $\mathbb Z^n/\Gamma_{t+1}$, and a finite intersection of finite-index sublattices has finite index. Next, we show that this family contains every other compatible family. If $\{K_t\}_{t=1}^T$ is any compatible family satisfying $K_t\subseteq\Lambda_t$, then $K_T\subseteq\Lambda_T=\Gamma_T$. Assuming $K_{t+1}\subseteq\Gamma_{t+1}$, compatibility implies $A_{tw}K_t\subseteq K_{t+1}\subseteq\Gamma_{t+1}$ for every $w\in\Omega_t$, and therefore $K_t\subseteq\Lambda_t\cap\bigcap_{w\in\Omega_t} A_{tw}^{-1}\Gamma_{t+1}=\Gamma_t$, which completes the induction.
\end{proof}

\begin{proof}[Proof of Lemma~\ref{lem:observation-phases}]
An additive homomorphism $\phi$ on $\mathbb Z^n$ recovers $\mathcal O_h s\bmod m$ if and only if $\ker\phi\subseteq K_h$, because this inclusion makes $\phi(s)\mapsto\mathcal O_h s\bmod m$ well defined, whereas recoverability forces every $s\in\ker\phi$ to share the zero state's observations. The assignment $\phi(s)\mapsto s+K_h$ from the image of $\phi$ onto $G_h$ is then well defined, surjective, and unique, so that the image of $\phi$ has at least $|G_h|$ elements and $G_h$ itself achieves the minimum label cardinality. With $x=V^{-1}s$ in unimodular Smith coordinates \citep{newman1972}, the condition $\mathcal O_h s\equiv 0\pmod m$ becomes $e_i x_i\equiv 0\pmod m$ for $i=1,\ldots,r$, which is equivalent to $q_i\mid x_i$. The matrix $V^{-1}$ maps $\mathbb Z^n$ bijectively onto $\mathbb Z^n$. Therefore, the coordinates $x_i$ vary independently, which establishes \eqref{eq:phaseorder}.

Up to unimodular coordinate change, the same Smith decomposition identifies $M_h$ with the direct sum of cyclic ideals $e_i R\subseteq R$, of cardinality $\prod_{i=1}^r q_i = |G_h|$ even for composite $m$. Because $K_{h+1}\subseteq K_h$ and $M_h\subseteq M_{h+1}$ with $|M_h|=|G_h|$, the equality $K_{h+1}=K_h$ holds if and only if $M_{h+1}=M_h$, that is,
\[
 CA^h\equiv \sum_{j=0}^{h-1}B_jCA^j \pmod m,\qquad B_j\in R^{p\times p}.
\]
We have $M_{j+1}=\operatorname{row}_R(C)+M_j A$, and therefore equality at $h_*$ implies $M_{j}=M_{h_*}$ and $K_j=K_{h_*}$ for all $j\ge h_*$. By the Cayley--Hamilton theorem over the commutative ring $R$, $A^n$ is an $R$-linear combination of $I,A,\ldots,A^{n-1}$, and therefore $h_*\le n$. The natural projection $G_h\to G_{h-1}$ has uniform fiber cardinality $g_h=[K_{h-1}:K_h]=|M_h/M_{h-1}|$. Right multiplication by $A$ sends $M_{h-1}$ into $M_h$ and maps the rows of $\mathcal O_h$ onto the newly adjoined rows of $\mathcal O_{h+1}$ modulo $M_h$, which induces a surjective homomorphism $M_h/M_{h-1}\twoheadrightarrow M_{h+1}/M_h$. Therefore, $g_{h+1}$ divides $g_h$.

For the chain of degree bounds, any monic scalar polynomial $P(x)=x^d-\sum_{j=0}^{d-1}b_j x^j\in R[x]$ with $CP(A)=0$ yields a special block recurrence with $B_j=b_j I_p$, which gives $h_*\le d_{\rm sc}$ (with equality precisely when $CA^{h_*}$ lies in the scalar $R$-span of $C,CA,\ldots,CA^{h_*-1}$). Next, we stack a minimal generating set of $\nu$ rows for $M_{h_*}$ into $H\in R^{\nu\times n}$. The invariance $M_{h_*}A\subseteq M_{h_*}$ yields $HA=TH$ for some $T\in R^{\nu\times\nu}$. Applying the monic characteristic polynomial $\chi_T\in R[x]$ of degree $\nu$ gives $H\chi_T(A)=\chi_T(T)H=0$. The rows of $C\bmod m$ belong to $\operatorname{row}_R(H)$, and it follows that $C\chi_T(A)=0$ and $d_{\rm sc}\le\nu$. Finally, $\mathcal O_{h_*}$ has $ph_*$ rows, and its Smith form produces at most $n$ nontrivial cyclic factors ($q_i>1$), which generate $M_{h_*}$, and therefore $\nu\le\min\{n,ph_*\}$. To verify that $\nu$ equals the exact count $r_+$ of factors with $q_i>1$, we pick a prime $\ell$ dividing the smallest nontrivial order $q_{r_+}$. Because $q_{i+1}\mid q_i$, the prime $\ell$ divides every $q_i$ for $1\le i\le r_+$, so that $\dim_{\mathbb F_\ell}(M_{h_*}/\ell M_{h_*})=r_+$, which prevents generation by fewer than $r_+$ elements. When $C\equiv 0\pmod m$, the empty label and $P(x)=1$ give $h_*=d_{\rm sc}=\nu=0$.
\end{proof}

\begin{proof}[Proof of Proposition~\ref{prop:phasevalue}]
We first show that statement (b) implies statement (a). Partitioning the feasible pairs $(a,y)\in X_{tw}(s)$ by their successor coset $g'=g_{t+1}(y)$ expresses the stage minimum minus $\beta^\top s$ as $\min_{g'\in\mathcal G_{t+1}}\{K^\alpha_{tw}(s,g')-\beta^\top s+\phi(g')\}$, which reduces to $\psi_\phi(g_t(s))$ whenever $K^\alpha_{tw}(s,g')-\beta^\top s=\bar K_{tw}(g_t(s),g')$. Conversely, suppose that the stage minimum minus $\beta^\top s$ depends only on $g_t(s)$ for every finite $\phi$. We fix a target coset $g'\in\mathcal G_{t+1}$, set $\phi(g')=0$, and set $\phi(\tilde g)=M$ for all $\tilde g\ne g'$. Letting $M\to+\infty$ drives the adjusted minimum to $K^\alpha_{tw}(s,g')-\beta^\top s$ (including $+\infty$ if $g'$ is unreachable from $s$). Agreement at two states in the same coset for all finite $M$ therefore forces identical values of $K^\alpha_{tw}(\cdot,g')-\beta^\top(\cdot)$ for every $g'$. Taking expectations over $w\in\Omega_t$ and inducting backward from $T+1$ yields the horizon-wide representation.
\end{proof}

\begin{proof}[Proof of Theorem~\ref{thm:optimal-period}]
We first show that every $h\in K_t$ is a period of $R_t$ and that $R_t$ vanishes on $K_t$. Let $h\in K_t$. For any $x\in\mathbb Z^n$ and $w\in\Omega_t$, we have $C_{tw}(x+h)=C_{tw}x$, so that $x$ and $x+h$ share the same feasible action set, whereas their successors under any $a$ differ by $A_{tw}h\in K_{t+1}$. Backward induction from $R_{T+1}\equiv 0$ proves $R_{tw}(x+h)=R_{tw}(x)$ and $R_t(x+h)=R_t(x)$. At $x\in K_t$, selecting the zero-cost action $a_0$ with $b_{tw}(a_0)\in K_{t+1}$ keeps $A_{tw}x+b_{tw}(a_0)\in K_{t+1}$, so that backward induction and nonnegativity give $R_t(x)=0$ on $K_t$.

Next, we show that $R_t\ge B_t$. Assume inductively that $R_{t+1}\ge B_{t+1}$. At stage $t$ and outcome $w$, any action with $r_{tw}(g,a)\ge\epsilon_{tw}$ already incurs at least $\epsilon_{tw}\ge B_{tw}(x)$ because $R_{t+1}\ge0$. Any action with $r_{tw}(g,a)<\epsilon_{tw}$ is available only when $g=C_{tw}x=0$ and satisfies $b_{tw}(a)\in K_{t+1}$. By the $K_{t+1}$-periodicity of $R_{t+1}$, its total cost is at least $R_{t+1}(A_{tw}x+b_{tw}(a))=R_{t+1}(A_{tw}x)\ge B_{t+1}(A_{tw}x)$. Minimizing over $a$ yields $R_{tw}(x)\ge B_{tw}(x)$ and therefore $R_t(x)\ge B_t(x)$. Because $p_{tw}>0$ and $\epsilon_{tw}>0$, \eqref{eq:observer-stop} gives $B_t(x)=0$ if and only if $C_{tw}x=0$ and $B_{t+1}(A_{tw}x)=0$ for all $w\in\Omega_t$, which by induction and \eqref{eq:observer-kernel} characterizes $\{B_t=0\}=K_t$. Together with $R_t\ge B_t\ge0$ and $R_t|_{K_t}=0$, this proves $\{R_t=0\}=K_t$. Finally, any period $h$ of $R_t$ satisfies $R_t(h)=R_t(0)=0$, and therefore $h\in K_t$, which completes the proof of \eqref{eq:value-period}.

To prove sharpness, we construct a single model in class \eqref{eq:period-model} by offering at $g=0$ a free action with $b=0$, and at every $g$ paid actions of cost $\epsilon_{tw}$ whose offsets $b$ range over $\{0,\ldots,m-1\}^n$. Choosing $b$ so that $A_{tw}x+b\equiv 0\pmod m$ lands in $m\mathbb Z^n\subseteq K_{t+1}$, where the future cost is zero. Therefore, the best paid option always costs $\epsilon_{tw}$, whereas the free option at $g=0$ costs $R_{t+1}(A_{tw}x)$, which matches \eqref{eq:observer-stop} identically.
\end{proof}

\begin{proof}[Proof of Corollary~\ref{cor:price-period}]
We define a two-point random variable $Z$ that takes the value $kh$ with probability $1/k$ and the value $0$ with probability $1-1/k$, so that $\mathbb E[Z]=h$. Because $h\in L$, all three points $0,h,kh$ lie in the same $L$-coset, whereas $0\in K_t$ and $kh\in K_t$ imply $R_{tw}(0)=R_{tw}(kh)=0$ for every $w\in\Omega_t$ by Theorem~\ref{thm:optimal-period}. On this $L$-coset, every valid cut in outcome $w$'s pool has the form $\ell(x)-r\|x-a\|_1$ with $\ell$ affine and $0\le r\le\rho$. By affinity of $\ell$ and the reverse triangle inequality $\|Z-a\|_1-\|h-a\|_1\ge -\|Z-h\|_1$, the cut's value at $h$ exceeds its expectation under $Z$ by at most $r\,\mathbb E\|Z-h\|_1=2r(1-1/k)\|h\|_1\le 2(1-1/k)\rho\|h\|_1$. Because validity at $0$ and $kh$ forces the cut's expectation under $Z$ to be nonpositive, taking the supremum over cuts in outcome $w$ and averaging over $w\in\Omega_t$ yields $\underline R_t(h)\le 2(1-1/k)\rho\|h\|_1$. Combining this with $R_t(h)\ge B_t(h)$ and $\underline R_t(h)\le R_t(h)$ gives \eqref{eq:period-r9}.
\end{proof}

\subsection{Proofs for Section~\ref{sec:cuts}}

\begin{proof}[Proof of Theorem~\ref{thm:mixture}]
Because $\hat s\in S'\subseteq S\subseteq D$, the feasible distribution polytopes satisfy $\delta_{\hat s}\in\mathcal M(S',\hat s)\subseteq\mathcal M(S,\hat s)\subseteq\mathcal M(D,\hat s)$, where $\delta_{\hat s}$ is the Dirac mass at $\hat s$, which establishes the chain of inequalities. Both finite linear programs attain primal and dual optima. For any $\mu\in\mathcal M(S',\hat s)$, the affinity of $u$ gives $\sum_{z\in S'}\mu_z u(z)=u(\hat s)=E_S(\hat s)$, and subtracting this constant from $\sum_{z\in S'}\mu_z F(z)$ yields the exact gain identity. Because $F(z)-u(z)\ge0$ on $S$ and $\mathcal M(S',\hat s)=\{\mu\in\mathcal M(S,\hat s):\operatorname{supp}(\mu)\subseteq S'\}$, the minimum over $\mathcal M(S',\hat s)$ vanishes if and only if some minimizer of $\sum_{z\in S}\mu_z F(z)$ over $\mathcal M(S,\hat s)$ is supported entirely inside $S'$.
\end{proof}

\begin{proof}[Proof of Theorem~\ref{thm:closure}]
We first show validity. Evaluating the definition of $b_{[s]_t}^\rho$ in \eqref{eq:closing-cut} at the copy $z=s$ gives $b_{[s]_t}^\rho\le f(s)+\rho d(s)=F(s)-\pi^\top s+\rho d(s)$, which rearranges to $h^\rho(s)\le F(s)$ on $D_t$. At $s=\hat s$, we have $d(\hat s)=0$, so that $h^\rho(\hat s)=\pi^\top\hat s+b_{\hat g}^\rho$ equals $F(\hat s)$ if and only if $z=\hat s$ minimizes $f(z)+\rho d(z)$ over $z\in D_{\hat g}$, that is, $f(\hat s)\le f(z)+\rho d(z)$ for every $z\in D_{\hat g}\setminus\{\hat s\}$. Solving these finitely many inequalities yields \eqref{eq:least-penalty}, which is finite because $D_{\hat g}$ is finite. Finally, refining the anchor coset shrinks $D_{\hat g}$ and deletes competitors from the maximum in \eqref{eq:least-penalty}.
\end{proof}

\begin{proof}[Proof of Theorem~\ref{thm:local-price}]
By LP duality applied to the copy-interpolation problem \eqref{eq:mixing} with perturbed objective $F(z)+\rho d(z)$, we have $G_S(\rho)=\min_{\mu\in\mathcal M}\mathbb E_\mu[F(Z)+\rho d(Z)]$, which is the pointwise minimum of finitely many affine functions of $\rho$ with nonnegative slopes $\mathbb E_\mu d(Z)\ge0$. Because $\delta_{\hat s}\in\mathcal M$ achieves the objective value $F(\hat s)$, the equality $G_S(\rho)=F(\hat s)$ holds if and only if $\mathbb E_\mu F(Z)+\rho\,\mathbb E_\mu d(Z)\ge F(\hat s)$ for all $\mu\in\mathcal M$ with $\mathbb E_\mu d(Z)>0$, which yields \eqref{eq:joint-penalty}. For the support bound, we observe that removing any probability mass on $\hat s$ and renormalizing leaves both $[F(\hat s)-\mathbb E_\mu F(Z)]_+/\mathbb E_\mu d(Z)$ and the zero-mean condition $\mathbb E_\mu(Z-\hat s)=0$ invariant. Setting $v_z=z-\hat s$ and homogenizing via Charnes--Cooper variables $\eta_z=\mu_z/\mathbb E_\mu d(Z)$ for $z\in S\setminus\{\hat s\}$ rewrites $\rho^\circ_S>0$ as the maximum of $\sum_{z\ne\hat s}\eta_z(F(\hat s)-F(z))$ over the bounded polytope $\{\eta\ge0:\sum_{z\ne\hat s}\eta_z v_z=0,\ \sum_{z\ne\hat s}\eta_z d(z)=1\}$ (bounded because $d(z)\ge1$ for $z\ne\hat s$). At an extreme-point optimizer, the active columns $(v_z,d(z))$ are linearly independent in $(\operatorname{aff} S-\hat s)\times\mathbb R$, bounding the support cardinality by $\dim(\operatorname{aff} S)+1$.
\end{proof}

\begin{proof}[Proof of Theorem~\ref{thm:root-price}]
We first prove \eqref{eq:root-price-value}. Applying the single-stage priced envelope \eqref{eq:priced-envelope} inductively from $V_{T+1}$ down to $t=2$, followed by physical optimization at $s_1$, identifies $W_\Gamma(\rho)$ with the optimal expected cost of a dynamic control problem that, upon observing $(s,w)$ at $t\ge2$, selects a mean-$s$ copy distribution on $S_t(s)$, pays $\rho\|z-s\|_1$, and selects $y\in Y_{tw}(z)$ at cost $c^0_{tw}(z,y)$. Every Markov randomized policy for this relaxed process induces state-arrival and transition masses $(d,f)$ satisfying \eqref{eq:occupation}. Conversely, normalizing any $(d,f)\in\mathcal P_\Gamma$ at pairs $(t,w,s)$ with $d_t(s)>0$ recovers admissible conditional copy and successor distributions achieving objective $J(v)+\rho D(v)$. Because $D_t$ and $\Omega_t$ are finite and probability masses sum to $1$ at each stage, $\mathcal P_\Gamma$ is a nonempty compact polytope, and therefore the linear functional $J+\rho D$ attains its minimum at a vertex $v\in\mathcal V_\Gamma$, which proves \eqref{eq:root-price-value}. Concavity, monotonicity, and piecewise linearity follow immediately.

Next, we establish the threshold \eqref{eq:root-price-threshold}. Flows $v\in\mathcal P_\Gamma$ with $D(v)=0$ satisfy $f_{twszy}=0$ whenever $z\ne s$, and they correspond bijectively to randomized physical MSIP policies, whose minimum expected cost over $\mathcal P_\Gamma\cap\{D(v)=0\}$ equals $V^\star$ by Bellman optimality. Because $\{v\in\mathcal P_\Gamma:D(v)=0\}$ is a nonempty face of $\mathcal P_\Gamma$, we have $W_\Gamma(\rho)\le V^\star$ for all $\rho\ge0$, with equality $W_\Gamma(\rho)=V^\star$ if and only if $J(v)+\rho D(v)\ge V^\star$ holds at every vertex $v\in\mathcal V_\Gamma$ with $D(v)>0$, which yields \eqref{eq:root-price-threshold}. Finally, replacing $\Gamma_t$ by $\Gamma'_t\subseteq\Gamma_t$ restricts $\mathcal P_{\Gamma'}$ to the face of $\mathcal P_\Gamma$ with $f_{twszy}=0$ for $z\notin s+\Gamma'_t$, which still contains the entire zero-displacement face $\{D(v)=0\}$. Minimizing over a smaller polytope containing the physical optimum weakly raises $W_\Gamma(\rho)$ and shrinks the vertex set in \eqref{eq:root-price-threshold}.
\end{proof}

\noindent\emph{Extension to anisotropic directional gauges \eqref{eq:aniso-penalty}.} For any componentwise nonnegative vectors $\rho^+,\rho^-\in\mathbb R_+^n$ and triples $z,s,c\in\mathbb R^n$, applying $(a+b)_+\le a_++b_+$ coordinatewise to $z_i-c_i=(z_i-s_i)+(s_i-c_i)$ and $c_i-z_i=(c_i-s_i)+(s_i-z_i)$ yields $(z_i-c_i)_+-(s_i-c_i)_+\le(z_i-s_i)_+$ and $(c_i-z_i)_+-(c_i-s_i)_+\le(s_i-z_i)_+$, and therefore
\[
 \Phi(z;c,\rho^+,\rho^-)-\Phi(s;c,\rho^+,\rho^-)\ \le\ \Phi(z;s,\rho^+,\rho^-).
\]
When $(\rho^+,\rho^-)=\rho(w^+,w^-)$ for fixed positive directional weights $w_i^+,w_i^->0$ \citep[cf.][]{ahmed2022}, the gauge $d_{w^+,w^-}(z)=\Phi(z;\hat s,w^+,w^-)$ satisfies $d_{w^+,w^-}(\hat s)=0$ and $d_{w^+,w^-}(z)\ge\min_i\min\{w_i^+,w_i^-\}>0$ on $\mathbb Z^n\setminus\{\hat s\}$ (keeping the Charnes--Cooper polytope in Theorem~\ref{thm:local-price} bounded, while Theorem~\ref{thm:closure} at fixed $\pi$ requires positivity only on $\{z\in D_{\hat g}:f(z)<f(\hat s)\}$). Consequently, Theorems~\ref{thm:closure}--\ref{thm:root-price} hold verbatim with $d(z)$ and $\|z-s\|_1$ replaced by $\Phi(z;\hat s,w^+,w^-)$ and $\Phi(z;s,w^+,w^-)$. Specifically, centering at $c=s$ with $(\rho^+,\rho^-)=\rho(w^+,w^-)$ attains the directional single-stage envelope $\min_{\mu\in\mathcal M(S,s)}\sum_{z\in S}\mu_z\{F(z)+\rho\,\Phi(z;s,w^+,w^-)\}$, whereas sublinearity ensures that every valid closing cut centered at an arbitrary $c\in D_t$ with $0\le\rho_i^\pm\le\rho w_i^\pm$ is dominated by that envelope (likewise extending Corollary~\ref{cor:price-period} to directional cuts with $0\le\rho_i^\pm\le\rho$ via $\Phi(Z;a,\rho^+,\rho^-)-\Phi(h;a,\rho^+,\rho^-)\le\rho\|Z-h\|_1$).

\begin{proof}[Proof of Corollary~\ref{cor:data-cap}]
Stagewise exogeneity of outcomes and the immediate cost bounds $\ell_{jb}\le c_{jb}\le u_{jb}$ imply $\max_{s\in D_t}Q_{tw}(s)-\min_{s\in D_t}Q_{tw}(s)\le H_{tw}$ for every $(t,w)$. Applying the single-stage slope-optimized threshold \eqref{eq:joint-penalty} to the true outcome value $F=Q_{tw}$ at any $s\in D_t$ ($t\ge2$) after stripping mass at $s$ bounds the numerator $[Q_{tw}(s)-\mathbb E_\mu Q_{tw}(Z)]_+$ above by $H_{tw}$ and the denominator $\mathbb E_\mu\|Z-s\|_1$ below by $\eta_t(s)$ (note that $\eta_t(s)\ge1$ whenever finite, because every nonzero integer displacement has $\ell_1$ norm at least $1$). Consequently, for $\rho$ at or above the displayed maximum, the complete cut family reproduces $Q_{tw}(s)$ exactly at every nonroot state and outcome by backward induction from $V_{T+1}$, forcing $W_\Gamma(\rho)=V^\star$.
\end{proof}

\subsection{Proofs for Section~\ref{sec:algorithm}}

\begin{proof}[Proof of Proposition~\ref{prop:dinkelbach}]
We first establish finite shell-monotonicity for isotropic Dinkelbach search \citep{dinkelbach1967}. At any iterate $\rho_m$ where the maximum residual $\max_{z\in D_{\hat g}}\{a_z-\rho_m d_z\}$ on $D_{\hat g}$ is strictly positive, a maximizing copy $z=z_m$ satisfies $a_z-\rho_m d_z>0$. We necessarily have $z\ne\hat s$ because $a_{\hat s}=d_{\hat s}=0$, and therefore $d_z=\|z-\hat s\|_1\ge 1$. Setting $\rho_{m+1}=a_z/d_z$ yields $\rho_{m+1}>\rho_m$, resets the residual of $z$ to $a_z-\rho_{m+1}d_z=0$, and satisfies $\rho_{m+1}\le\rho_{\hat g}^*$ by definition of $\rho_{\hat g}^*$ in \eqref{eq:least-penalty}.

Suppose that a copy $z'\in D_{\hat g}$ causes another strictly positive update at $\rho_{m+1}$. Optimality of $z$ at $\rho_m$ gives $a_{z'}-\rho_m d_{z'}\le a_z-\rho_m d_z$, whereas positivity at $\rho_{m+1}$ gives $a_{z'}-\rho_{m+1}d_{z'}>0=a_z-\rho_{m+1}d_z$. Subtracting the first relation from the second yields
\[
 (\rho_{m+1}-\rho_m)(d_z-d_{z'}) > 0,
\]
and therefore $d_{z'}<d_z$ because $\rho_{m+1}>\rho_m$. Consequently, every positive penalty update consumes a strictly smaller positive $\ell_1$ distance in $D_{\hat g}$, even when multiple copies tie in an intermediate coset subproblem. Because $D_{\hat g}$ realizes at most $|\{d_z:z\in D_{\hat g},\,d_z>0\}|$ distinct positive $\ell_1$ distances, at most that many updates can occur before $\max_{z\in D_{\hat g}}\{a_z-\rho d_z\}=0$ (which implies $\rho\ge\rho_{\hat g}^*$ and therefore $\rho=\rho_{\hat g}^*$), and one final solve verifies contact, which gives $K+1\le 1+|\{d_z:z\in D_{\hat g},\,d_z>0\}|$. When $D_t\subseteq\prod_{i=1}^n\{\ell_i,\ldots,\ell_i+M_i\}$, positive integer $\ell_1$ distances lie in $\{1,\ldots,\sum_{i=1}^n M_i\}$.

Next, we consider coordinate-scaled sparse directional lifting under \eqref{eq:aniso-penalty} with positive matrix scales $w_i\ge 1$ ($i=1,\ldots,n$). Let $\mathcal V_{\hat g}(\pi)=\{z\in D_{\hat g}:f(z)<f(\hat s)\}\subseteq D_{\hat g}\setminus\{\hat s\}$ collect the strictly violating copies in $D_{\hat g}$. At step $m\ge 0$ with nonnegative vectors $(\rho^{+(m)},\rho^{-(m)})\in\mathbb R_+^n\times\mathbb R_+^n$, let $z_m\in\arg\min_{z\in D_{\hat g}}\{f(z)+\Phi(z;\hat s,\rho^{+(m)},\rho^{-(m)})\}$. If $f(z_m)+\Phi(z_m;\hat s,\rho^{+(m)},\rho^{-(m)})\ge f(\hat s)$, then $z=\hat s$ (where $\Phi(\hat s;\hat s,\cdot,\cdot)=0$) already attains the coset minimum $f(\hat s)$ and the search terminates with exact contact at $\hat s$. Otherwise, the deficit $\Delta_m=f(\hat s)-f(z_m)-\Phi(z_m;\hat s,\rho^{+(m)},\rho^{-(m)})>0$ forces $z_m\in\mathcal V_{\hat g}(\pi)$, so that $\delta_m=z_m-\hat s\in\mathbb Z^n\setminus\{0\}$ has positive weighted norm $d_w(\delta_m)=\sum_{i=1}^n w_i|\delta_{m,i}|\ge 1$. Updating
\[
 \rho_i^{+(m+1)}=\rho_i^{+(m)}+\eta_m w_i\mathbf 1_{\{\delta_{m,i}>0\}},\qquad
 \rho_i^{-(m+1)}=\rho_i^{-(m)}+\eta_m w_i\mathbf 1_{\{\delta_{m,i}<0\}}
\]
with step size $\eta_m\ge\Delta_m/d_w(\delta_m)>0$ ensures $(\rho^{+(m+1)},\rho^{-(m+1)})\ge(\rho^{+(m)},\rho^{-(m)})$ componentwise and increases $\Phi(z_m;\hat s,\cdot,\cdot)$ by $\eta_m\sum_{i:\delta_{m,i}>0}w_i\delta_{m,i}+\eta_m\sum_{i:\delta_{m,i}<0}w_i(-\delta_{m,i})=\eta_m d_w(\delta_m)\ge\Delta_m$, which yields $f(z_m)+\Phi(z_m;\hat s,\rho^{+(m+1)},\rho^{-(m+1)})\ge f(\hat s)$. Because $\Phi(z;\hat s,\rho^+,\rho^-)$ is componentwise nondecreasing in $(\rho^+,\rho^-)\ge 0$, every copy $z_r$ ($0\le r\le m$) eliminated at or before step $m$ satisfies $f(z_r)+\Phi(z_r;\hat s,\rho^{+(m')},\rho^{-(m')})\ge f(\hat s)$ for all subsequent steps $m'>m$ and can never again produce a positive deficit. Therefore, at most $|\mathcal V_{\hat g}(\pi)|$ updates occur before exact contact is verified on solve $K+1\le 1+|\mathcal V_{\hat g}(\pi)|\le|D_{\hat g}|$, whereas $\rho_i^+$ (resp.\ $\rho_i^-$) remains zero on every coordinate $i$ with $z_i\le\hat s_i$ (resp.\ $z_i\ge\hat s_i$) across encountered violators. Finally, at any iterate $(\rho^+,\rho^-)\ge 0$ (including isotropic $\rho_i^+=\rho_i^-=\rho$), substituting $z=s$ into the certified lower bounds $b_{\hat g}\le\min_{z\in D_{\hat g}}\{f(z)+\Phi(z;\hat s,\rho^+,\rho^-)\}$ and $b_0\le\min_{z\in D_t}\{f(z)+\Phi(z;\hat s,\rho^+,\rho^-)\}$ yields $b_0+\max\{0,b_{\hat g}-b_0\}\mathbf 1_{\{[s]_t=\hat g\}}-\Phi(s;\hat s,\rho^+,\rho^-)\le F_{tw}(s)-\pi^\top s$, which establishes global validity on $D_t$.
\end{proof}

\begin{proof}[Proof of Proposition~\ref{prop:phase-hierarchy}]
\emph{(i) Bellman-expectation LP floors and LP-certified global floor skips.}
Let $X_{\tau w}^{\mathrm{LP}}([\ell_\tau,u_\tau])$ denote the continuous LP relaxation of stage $(\tau,w)$ over the incoming state box $[\ell_\tau,u_\tau]\supseteq D_\tau$, and let $m_{\tau w}\le\min_{(s,a,y)\in X_{\tau w}^{\mathrm{LP}}([\ell_\tau,u_\tau])}c_{\tau w}(s,a,y)$ and $E_{T+1}\le\min_{y\in D_{T+1}}V_{T+1}(y)$ be certified lower bounds. We prove $Q_{\tau w}(s)\ge\underline\theta_{\tau w}\ge\underline\theta_\tau$ and $V_\tau(s)\ge E_\tau\ge\underline\theta_\tau$ in \eqref{eq:theta-lb} by backward induction on $\tau=T+1,T,\ldots,1$. At $\tau=T+1$, $V_{T+1}(y)\ge E_{T+1}$ holds on $D_{T+1}$ by definition. Suppose that $V_{\tau+1}(y)\ge E_{\tau+1}$ holds for every $y\in D_{\tau+1}$. For any $s\in D_\tau$ and $w\in\Omega_\tau$, relatively complete recourse and attained finite minima provide $(a^\star,y^\star)\in X_{\tau w}(s)$ with $y^\star\in D_{\tau+1}$ and $Q_{\tau w}(s)=c_{\tau w}(s,a^\star,y^\star)+V_{\tau+1}(y^\star)$. Because $(s,a^\star,y^\star)\in X_{\tau w}^{\mathrm{LP}}([\ell_\tau,u_\tau])$, we have $c_{\tau w}(s,a^\star,y^\star)\ge m_{\tau w}$, and the induction hypothesis gives $V_{\tau+1}(y^\star)\ge E_{\tau+1}$. Summing yields
\[
 Q_{\tau w}(s)\ \ge\ m_{\tau w}+E_{\tau+1}\ =\ \underline\theta_{\tau w}
 \ \ge\ \min_{\tilde w\in\Omega_\tau}m_{\tau\tilde w}+E_{\tau+1}\ =\ \underline\theta_\tau.
\]
Taking expectations over $w\in\Omega_\tau$ ($\sum_{w\in\Omega_\tau}p_{\tau w}=1$, $p_{\tau w}>0$) gives $V_\tau(s)=\sum_{w\in\Omega_\tau}p_{\tau w}Q_{\tau w}(s)\ge\sum_{w\in\Omega_\tau}p_{\tau w}m_{\tau w}+E_{\tau+1}=E_\tau\ge\underline\theta_\tau$, which completes the induction. Expanding $E_{\tau+1}$ shows that
\[
 \underline\theta_\tau - \Bigl(\sum_{j=\tau}^T\min_{w\in\Omega_j}m_{jw}+E_{T+1}\Bigr)
 \ =\ \sum_{j=\tau+1}^T\Bigl(\sum_{w\in\Omega_j}p_{jw}m_{jw}-\min_{w\in\Omega_j}m_{jw}\Bigr)\ \ge\ 0.
\]
Similarly, at any backward step $(t,w)$ with current continuation relaxation $F_{tw}^{\mathrm{LP}}$ over $[\ell_t,u_t]\supseteq D_t$, integrality relaxation implies $F_{tw}^{\mathrm{LP}}(z)\le F_{tw}(z)$ for every $z\in D_t$, so that
\[
 b_{\mathrm{LP}}(\pi)=\min_{z\in[\ell_t,u_t]}\bigl\{F_{tw}^{\mathrm{LP}}(z)-\pi^\top z\bigr\}
 \ \le\ \min_{z\in D_t}\bigl\{F_{tw}(z)-\pi^\top z\bigr\}
 \ \le\ \beta_{\hat g}(\pi).
\]
Setting $b_0=b_{\mathrm{LP}}(\pi)$ (minus LP optimality tolerance) therefore yields a certified global intercept on $D_t$ with zero MIP solves, and when $b_{\hat g}(\pi)-b_{\mathrm{LP}}(\pi)\le\varepsilon$, skipping the free-incoming MIP loses at most $\varepsilon$ on unvisited cosets.

\emph{(ii) Pointwise Bellman validity, in-place merging monotonicity, and LP dominance.}
We fix $(t,w)$, a slope $\pi\in\mathbb R^n$, and a valid continuation $\underline V_{t+1}\le V_{t+1}$ so that $F_{tw}\le Q_{tw}$ on $D_t$. We write $\phi_{t,1}(z)=(V_t^{-1}z)_1\bmod q_{t,1}\in\{0,\ldots,q_{t,1}-1\}$ for the first Smith residue coordinate. Let $b_0\le\min_{z\in D_t}\{F_{tw}(z)-\pi^\top z\}$, let $b^{(1)}_{r_1}\le\min_{z\in D_t:\phi_{t,1}(z)=r_1}\{F_{tw}(z)-\pi^\top z\}$ for $r_1\in\mathcal R_1\subseteq\{0,\ldots,q_{t,1}-1\}$, and let $b_g\le\min_{z\in D_g}\{F_{tw}(z)-\pi^\top z\}$ for $g\in\mathcal G\subseteq\mathcal G_t$, where $+\infty$ is assigned only when the corresponding integer feasible set is proven empty. Consider any physical state $y\in D_t$ with label $\hat g=[y]_t=(\hat g_1,\ldots,\hat g_{r_t})$. Because $y$ belongs simultaneously to $D_t$, $\{z\in D_t:\phi_{t,1}(z)=\hat g_1\}$, and $D_{\hat g}$, all three sets are nonempty, which ensures $\bar b(\hat g)=\max\bigl(\{b_0\}\cup\{b^{(1)}_{\hat g_1}:\hat g_1\in\mathcal R_1\}\cup\{b_{\hat g}:\hat g\in\mathcal G\}\bigr)<+\infty$. Substituting $z=y$ into each defining minimization gives $\bar b(\hat g)\le F_{tw}(y)-\pi^\top y$, and therefore $\pi^\top y+\bar b([y]_t)\le F_{tw}(y)\le Q_{tw}(y)$.

Suppose that a second valid lattice cut at the same slope $\pi$ supplies intercepts $(b_0',\{b^{\prime(1)}_{r_1}\},\{b_g'\})$ certified against $\tilde F_{tw}$ with $F_{tw}\le\tilde F_{tw}\le Q_{tw}$. Because both sets of intercepts minorize $Q_{tw}(z)-\pi^\top z$ on the identical domain subsets (and minorize $\tilde F_{tw}(z)-\pi^\top z$ because $F_{tw}\le\tilde F_{tw}$), their componentwise maxima $\tilde b_0=\max\{b_0,b_0'\}$, $\tilde b^{(1)}_{r_1}=\max\{b^{(1)}_{r_1},b^{\prime(1)}_{r_1}\}$, and $\tilde b_g=\max\{b_g,b_g'\}$ preserve pointwise lower-bound validity while monotonically raising $\bar b(g)$. Whenever certified bounds cover all factor-1 residues $\{0,\ldots,q_{t,1}-1\}$ or all cosets partitioning $D_t$, the disjoint union $D_t=\bigsqcup_{r_1=0}^{q_{t,1}-1}\{z\in D_t:\phi_{t,1}(z)=r_1\}=\bigsqcup_{g\in\mathcal G_t}D_g$ implies $\min_{z\in D_t}\{\tilde F_{tw}(z)-\pi^\top z\}\ge\max\{\min_{r_1}\tilde b^{(1)}_{r_1},\,\min_{g\in\mathcal G_t}\bar b(g)\}$, which validates partition lifting of $\tilde b_0$. In the parent LP relaxation over one-hot weights $\delta\in\Delta(\mathcal G_t)$, Jensen's inequality applied to the convex maximum gives
\[
 \sum_{g\in\mathcal G_t}\Bigl(\max_{1\le k\le K}\bar b^{(k)}(g)\Bigr)\delta_g
 \ \ge\ \max_{1\le k\le K}\sum_{g\in\mathcal G_t}\bar b^{(k)}(g)\delta_g,
\]
which shows that merging $K$ same-slope cuts into a single row strictly strengthens fractional parent LP relaxations whenever distinct cuts elevate distinct cosets.

\emph{(iii) Exactness of one-hot and stage-pooled binary-digit parent formulations.}
When $[\mathbb Z^n:\Gamma_t]\le M_{\mathrm{oh}}$, fixing $\delta_g=0$ whenever $\bar b(g)=+\infty$ excludes no state in $D_t$ (because $D_g=\varnothing$), and for every $y\in D_t$ the unique feasible binary selector $\delta_g=\mathbf 1\{g=[y]_t\}$ evaluates \eqref{eq:parent-row} at $\Phi\equiv 0$ to $\pi^\top y+\bar b([y]_t)\le Q_{tw}(y)$. Because the simplex $\{\delta\ge 0:\sum_g\delta_g=1,\ \delta_g=0\ (\bar b(g)=\infty)\}$ has exclusively binary extreme points, the encoding achieves the convex-hull LP relaxation over one-hot weights with zero per-cut variables.

When $[\mathbb Z^n:\Gamma_t]>M_{\mathrm{oh}}$, let $\zeta=(\zeta_{j\ell})\in\{0,1\}^{B_t}$ ($B_t=\sum_{j=1}^{r_t}\lceil\log_2 q_{t,j}\rceil$) encode the Smith residues $r_j(y)=\sum_{\ell=0}^{\lceil\log_2 q_{t,j}\rceil-1}2^\ell\zeta_{j\ell}\in\{0,\ldots,q_{t,j}-1\}$. Injectivity of binary expansion maps each $g=(g_1,\ldots,g_{r_t})$ to a unique bit string $\operatorname{bits}(g)\in\{0,1\}^{B_t}$ (and $g_1$ to $\operatorname{bits}_1(g_1)$ on $\zeta_1$). If $\bar b(g)=+\infty$, the no-good inequality $\operatorname{Ham}(\zeta,\operatorname{bits}(g))\ge 1$ cuts off only $\zeta=\operatorname{bits}(g)$, which corresponds to no $y\in D_t$. For finite intercepts, each stage-pooled variable $\mu_g\in[0,1]$ (and optional factor-1 $\mu^{(1)}_{g_1}\in[0,1]$) obeys $\mu_g+\operatorname{Ham}(\zeta,\operatorname{bits}(g))\ge 1$ (resp.\ $\mu^{(1)}_{g_1}+\operatorname{Ham}(\zeta_1,\operatorname{bits}_1(g_1))\ge 1$), and each cut installs the telescoped row
\[
 \theta_w \ \ge\ \pi^\top y + b_0 + \sum_{g_1\in\mathcal R_1}G_1(g_1)\mu^{(1)}_{g_1}
   + \sum_{g\in\mathcal G}G_{\mathrm{fine}}(g)\mu_g,
\]
with nonnegative gains $G_1(g_1)=[b^{(1)}_{g_1}-b_0]_+\ge 0$ and $G_{\mathrm{fine}}(g)=[\bar b(g)-(b_0+G_1(g_1))]_+\ge 0$ (reducing directly to \eqref{eq:parent-row} with $G_1\equiv 0$ and $G_{\mathrm{fine}}(g)=[\bar b(g)-b_0]_+$ when factor-1 indicators are folded into $\mathcal G$). At any integer state $y\in D_t$ with $[y]_t=\hat g=(\hat g_1,\ldots,\hat g_{r_t})$, we have $\zeta=\operatorname{bits}(\hat g)$, so that $\operatorname{Ham}(\zeta,\operatorname{bits}(\hat g))=0$ forces $\mu_{\hat g}=1$ (and $\mu^{(1)}_{\hat g_1}=1$), whereas for every $g\ne\hat g$ (and $g_1\ne\hat g_1$) the Hamming distance is at least $1$, which makes $1-\operatorname{Ham}\le 0$, so that the lower bound $0\le\mu_g\le 1$ is tight at $0$. We stress that the single assignment $\mu_g=\mathbf 1\{g=\hat g\}$ and $\mu^{(1)}_{g_1}=\mathbf 1\{g_1=\hat g_1\}$ depends solely on $y$ and is simultaneously feasible across all pooled Hamming constraints. Because all coefficients $G_1(g_1),G_{\mathrm{fine}}(g)\ge 0$ and outcome weights $p_{tw}>0$ are nonnegative in minimization, $\mu_g=\mathbf 1\{g=\hat g\}$ and $\mu^{(1)}_{g_1}=\mathbf 1\{g_1=\hat g_1\}$ minimize every cut's right-hand side simultaneously, evaluating to $\pi^\top y+b_0+G_1(\hat g_1)+G_{\mathrm{fine}}(\hat g)=\pi^\top y+\bar b([y]_t)\le F_{tw}(y)\le Q_{tw}(y)$.
\end{proof}

\begin{proof}[Proof of Proposition~\ref{prop:merge-tolerance}]
For $t\ge 2$, let $e_t$ denote the supremum of $h(y)-Q_{tw}(y)$ over all $w\in\Omega_t$, all cuts $h$ that belong to $\mathcal C_{tw}$ at any time of the run, and all $y\in D_t$, or zero if this supremum is negative, and let $e_{T+1}=0$. We show that $e_t\le e_{t+1}+\varepsilon\sigma_t$ for $t=T,\ldots,2$, which gives the first claim. Fix a cut $h$ of $\mathcal C_{tw}$ with slope $\pi$, anchor $\hat s$, and penalties $\rho^\pm$, and fix a state $y\in D_t$. The value $h(y)$ equals $\pi^\top y+b-\Phi(y;\hat s,\rho^+,\rho^-)$, where $b$ is the intercept of the coset of $y$, the common lower bound $b_0$, or the minimum of the intercepts of all cosets. In each case, by the hypothesis, there is a slope $\pi'$ with $\|\pi-\pi'\|_\infty\le\varepsilon$ such that $b\le\tilde F(y)-\pi'^\top y+\Phi(y;\hat s,\rho^+,\rho^-)$, where $\tilde F(y)$ is the optimal value of the stage problem at $y$ and $w$ with the cuts of stage $t+1$ that were present when $b$ was computed. Hence,
\[
 h(y)\ \le\ \tilde F(y)+(\pi-\pi')^\top y\ \le\ \tilde F(y)+\varepsilon\|y\|_1\ \le\ \tilde F(y)+\varepsilon\sigma_t .
\]
By the definition of $e_{t+1}$, the approximation of $V_{t+1}$ in this stage problem exceeds $V_{t+1}$ by at most $e_{t+1}$ at every state of $D_{t+1}$, and for $t=T$ it equals $V_{T+1}$. Therefore, $\tilde F(y)\le Q_{tw}(y)+e_{t+1}$, and we obtain $h(y)\le Q_{tw}(y)+e_{t+1}+\varepsilon\sigma_t$. An intercept that is the maximum of several intercepts satisfies the same inequality, because each of them does. For the second claim, $L_k$ is a lower bound on the optimal value of the stage-1 problem with the cuts of stage 2, in which the approximation of $V_2$ exceeds $V_2$ by at most $e_2\le\kappa$ at every state. Therefore, $L_k\le V^\star+\kappa$.
\end{proof}

\begin{proof}[Proof of Theorem~\ref{thm:pooled-closing}]
\emph{(i) Global validity, finite anchor exactness, and zero-solve analytic LP closure.}
We fix $(t,w)$, an anchor $\hat s\in D_t$ with $\hat g=[\hat s]_t$, a slope $\pi\in\mathbb R^n$, and directional penalties $\rho^+,\rho^-\in\mathbb R_+^n$ in \eqref{eq:aniso-penalty}. Let $b_0\le\min_{z\in D_t}\{F_{tw}(z)-\pi^\top z+\Phi(z;\hat s,\rho^+,\rho^-)\}$ and $b_g\le\min_{z\in D_g}\{F_{tw}(z)-\pi^\top z+\Phi(z;\hat s,\rho^+,\rho^-)\}$ ($g\in\mathcal G\ni\hat g$) be certified lower bounds, and set $\bar b(g)=\max\{b_g,b_0\}$ for $g\in\mathcal G$ and $\bar b(g)=b_0$ otherwise. We note that $b_0$ may be certified through a MIP dual bound, through the penalized LP relaxation over $[\ell_t,u_t]\supseteq D_t$, or analytically from any valid global LP cut $b_{\mathrm{LP}}+\pi_{\mathrm{LP}}^\top z\le F_{tw}(z)$ via
\[
\begin{aligned}
 \min_{z\in D_t}\bigl\{F_{tw}(z)-\pi^\top z+\Phi(z;\hat s,\rho^+,\rho^-)\bigr\}
 \ge b_{\mathrm{LP}}+\sum_{i=1}^n\min_{z_i\in\{\ell_{t,i},\hat s_i,u_{t,i}\}}
 \bigl\{(\pi_{\mathrm{LP},i}-\pi_i)z_i\\[-3pt]
 {}+\rho_i^+(z_i-\hat s_i)_++\rho_i^-(\hat s_i-z_i)_+\bigr\},
\end{aligned}
\]
because each coordinate summand is piecewise affine on $[\ell_{t,i},\hat s_i]$ and $[\hat s_i,u_{t,i}]$ and therefore attains its minimum on $\{\ell_{t,i},\hat s_i,u_{t,i}\}$. For every $y\in D_t$, substituting $z=y\in D_{[y]_t}\subseteq D_t$ into the definitions of $b_0$ and $b_{[y]_t}$ gives $\bar b([y]_t)\le F_{tw}(y)-\pi^\top y+\Phi(y;\hat s,\rho^+,\rho^-)$, which rearranges to $h(y)=\pi^\top y+\bar b([y]_t)-\Phi(y;\hat s,\rho^+,\rho^-)\le F_{tw}(y)\le Q_{tw}(y)$. At $y=\hat s$, we have $\Phi(\hat s;\hat s,\rho^+,\rho^-)=0$. By Proposition~\ref{prop:dinkelbach}, sparse coordinate-scaled directional lifting achieves $\min_{z\in D_{\hat g}}\{F_{tw}(z)-\pi^\top z+\Phi(z;\hat s,\rho^+,\rho^-)\}=F_{tw}(\hat s)-\pi^\top\hat s$ (and therefore $h(\hat s)=F_{tw}(\hat s)$) in at most $1+|\{z\in D_{\hat g}:F_{tw}(z)-\pi^\top z<F_{tw}(\hat s)-\pi^\top\hat s\}|\le|D_{\hat g}|$ separations. Finally, given $\underline F_{tw}(\hat s)\le F_{tw}(\hat s)$ and a global LP cut $b_{\mathrm{LP}}+\pi_{\mathrm{LP}}^\top z\le F_{tw}(z)$ with $\Delta=[\underline F_{tw}(\hat s)-b_{\mathrm{LP}}-\pi_{\mathrm{LP}}^\top\hat s]_+>0$, setting $\pi=\pi_{\mathrm{LP}}$, $\bar b(g)\equiv b_{\mathrm{LP}}+\Delta$, and $\rho^+=\rho^-=\Delta\mathbf 1$ requires zero MIP solves and satisfies $h(\hat s)=\pi_{\mathrm{LP}}^\top\hat s+b_{\mathrm{LP}}+\Delta=\underline F_{tw}(\hat s)\le F_{tw}(\hat s)$ at $\hat s$, whereas for every $y\in D_t\setminus\{\hat s\}\subseteq\mathbb Z^n$, integrality implies $\|y-\hat s\|_1\ge 1$ and therefore
\[
 h(y)=\pi_{\mathrm{LP}}^\top y+b_{\mathrm{LP}}+\Delta\bigl(1-\|y-\hat s\|_1\bigr)
 \ \le\ \pi_{\mathrm{LP}}^\top y+b_{\mathrm{LP}}\ \le\ F_{tw}(y)\ \le\ Q_{tw}(y).
\]

\emph{(ii) Exact zero-binary slope and intercept absorption on box boundaries.}
We write $[\ell_t,u_t]=\prod_{i=1}^n[\ell_{t,i},u_{t,i}]\supseteq D_t$ and partition the boundary anchor coordinates into $I_\ell=\{i:\hat s_i\le\ell_{t,i}\}$ and $I_u=\{i\notin I_\ell:\hat s_i\ge u_{t,i}\}$. For every $y\in[\ell_t,u_t]$, the following holds:
\begin{itemize}[leftmargin=1.8em,itemsep=1pt,topsep=2pt]
\item if $i\in I_\ell$, then $y_i\ge\ell_{t,i}\ge\hat s_i$, so $(y_i-\hat s_i)_+=y_i-\hat s_i$, $(\hat s_i-y_i)_+=0$, and $\pi_i y_i-\rho_i^+(y_i-\hat s_i)_+-\rho_i^-(\hat s_i-y_i)_+=(\pi_i-\rho_i^+)y_i+\rho_i^+\hat s_i$;
\item if $i\in I_u$, then $y_i\le u_{t,i}\le\hat s_i$, so $(y_i-\hat s_i)_+=0$, $(\hat s_i-y_i)_+=\hat s_i-y_i$, and $\pi_i y_i-\rho_i^+(y_i-\hat s_i)_+-\rho_i^-(\hat s_i-y_i)_+=(\pi_i+\rho_i^-)y_i-\rho_i^-\hat s_i$.
\end{itemize}
Replacing $\pi_i$ by $\tilde\pi_i=\pi_i-\rho_i^+\mathbf 1_{\{i\in I_\ell\}}+\rho_i^-\mathbf 1_{\{i\in I_u\}}$, shifting every finite intercept $\bar b(g)$ to $\tilde b(g)=\bar b(g)+\sum_{i\in I_\ell}\rho_i^+\hat s_i-\sum_{i\in I_u}\rho_i^-\hat s_i$, and setting $\tilde\rho_i^+=\tilde\rho_i^-=0$ on $I_\ell\cup I_u$ leaves $h(y)$ identically unchanged at every point of $\prod_{i=1}^n[\ell_{t,i},u_{t,i}]$ (both continuous and integer) while adding zero split variables, zero binary indicators, and zero auxiliary rows for coordinates in $I_\ell\cup I_u$.

\emph{(iii) Ideal convex-hull MIP representation of stage-pooled staircase splits \eqref{eq:pooled-split}.}
We fix an integer coordinate $i\in\{1,\ldots,n\}$ with finite integer box bounds $\ell_{t,i}<u_{t,i}$ and first consider a single interior threshold $v\in\{\ell_{t,i}+1,\ldots,u_{t,i}-1\}$ in \eqref{eq:pooled-split}:
\[
 y_i - p_{i,v} + q_{i,v} = v,\qquad
 0\le p_{i,v}\le(u_{t,i}-v)z_{i,v},\qquad
 0\le q_{i,v}\le(v-\ell_{t,i})(1-z_{i,v}),\qquad
 z_{i,v}\in[0,1].
\]
For any point $(y_i,p,q,z)$ satisfying these linear constraints, we define the four weights
\[
 \alpha_1=\frac{p}{u_{t,i}-v},\qquad
 \alpha_2=z-\frac{p}{u_{t,i}-v},\qquad
 \alpha_3=\frac{q}{v-\ell_{t,i}},\qquad
 \alpha_4=(1-z)-\frac{q}{v-\ell_{t,i}}.
\]
By the variable bounds on $(p,q,z)$, we have $\alpha_r\ge 0$ ($r=1,2,3,4$) and $\sum_{r=1}^4\alpha_r=z+(1-z)=1$. Direct substitution verifies the convex decomposition
\[
 (y_i,p,q,z)
 =\alpha_1(u_{t,i},\,u_{t,i}-v,\,0,\,1)
 +\alpha_2(v,\,0,\,0,\,1)
 +\alpha_3(\ell_{t,i},\,0,\,v-\ell_{t,i},\,0)
 +\alpha_4(v,\,0,\,0,\,0),
\]
whose four extreme points all have binary $z\in\{0,1\}$ and project onto the three vertices $(\ell_{t,i},0,v-\ell_{t,i})$, $(v,0,0)$, and $(u_{t,i},u_{t,i}-v,0)$ of $\operatorname{conv}\{(y_i,(y_i-v)_+,(v-y_i)_+):y_i\in\{\ell_{t,i},\ldots,u_{t,i}\}\}$ (because $(y_i-v)_+$ and $(v-y_i)_+$ are affine on $[\ell_{t,i},v]$ and $[v,u_{t,i}]$). Therefore, each individual threshold split is locally ideal (sharp disjunctive).

Now let $\mathcal V_{t,i}=\{v_1<v_2<\cdots<v_m\}\subseteq\{\ell_{t,i}+1,\ldots,u_{t,i}-1\}$ collect all active interior integer thresholds pooled on coordinate $i$ across cuts and child outcomes at stage $t-1$, linked by the bidiagonal staircase inequalities $z_{i,v_{r+1}}\le z_{i,v_r}$ ($r=1,\ldots,m-1$). Consider any integer $y_i\in\{\ell_{t,i},\ldots,u_{t,i}\}$. Setting $z_{i,v_r}=\mathbf 1_{\{y_i>v_r\}}\in\{0,1\}$, $p_{i,v_r}=(y_i-v_r)_+$, and $q_{i,v_r}=(v_r-y_i)_+$ satisfies \eqref{eq:pooled-split} for each $r$ and obeys $z_{i,v_{r+1}}\le z_{i,v_r}$ because $y_i>v_{r+1}\implies y_i>v_r$. Conversely, for any binary assignment $z_{i,v_r}\in\{0,1\}$ feasible with $y_i$ in \eqref{eq:pooled-split}, the following holds:
\begin{itemize}[leftmargin=1.8em,itemsep=1pt,topsep=2pt]
\item if $z_{i,v_r}=1$, then $0\le q_{i,v_r}\le 0$ forces $q_{i,v_r}=0$, so $p_{i,v_r}=y_i-v_r\ge 0$, requiring $y_i\ge v_r$ and yielding $p_{i,v_r}=(y_i-v_r)_+$, $q_{i,v_r}=(v_r-y_i)_+$;
\item if $z_{i,v_r}=0$, then $0\le p_{i,v_r}\le 0$ forces $p_{i,v_r}=0$, so $q_{i,v_r}=v_r-y_i\ge 0$, requiring $y_i\le v_r$ and again yielding $p_{i,v_r}=(y_i-v_r)_+$, $q_{i,v_r}=(v_r-y_i)_+$ (including at $y_i=v_r$, where both $z_{i,v_r}\in\{0,1\}$ force $p_{i,v_r}=q_{i,v_r}=0$).
\end{itemize}
Therefore, $(p_{i,v_r},q_{i,v_r})=((y_i-v_r)_+,(v_r-y_i)_+)$ is uniquely and simultaneously exact across all $r=1,\ldots,m$ at every integer state $y\in D_t$, which allows every closing cut in $\bigcup_{w\in\Omega_t}\mathcal C_{tw}$ to reuse $\tilde\rho_i^+ p_{i,\hat s_i}+\tilde\rho_i^- q_{i,\hat s_i}$ without restricting $D_t$ or relaxing $\Phi(y;\hat s,\tilde\rho^+,\tilde\rho^-)$ at integer points. On fractional LP relaxations, $z_{i,v_{r+1}}\le z_{i,v_r}$ eliminates contradictory cross-cut indicator configurations ($z_{i,v_2}>z_{i,v_1}$ for $v_1<v_2$) and propagates bound-fixing across prefix/suffix threshold chains during branching. Finally, on the single parametric stage-$(t-1)$ model, each active interior pair $(i,v)$ creates exactly one binary $z_{i,v}$, which caps the total number of stage-$(t-1)$ split binaries on coordinate $i$ at $|\mathcal V_{t,i}|\le\max\{0,u_{t,i}-\ell_{t,i}-1\}$, independent of both $|\mathcal C_{tw}|$ and $|\Omega_{t-1}|$.
\end{proof}

\begin{proof}[Proof of Theorem~\ref{thm:root}]
We prove almost-sure finite stabilization and exactness simultaneously for immediate unpooled isotropic closure (\textnormal{LCD-Iso}, $\gamma=1$, $\rho_i^\pm=\rho$) and progress/repeat-visit gated directional pooled closure (\textnormal{LCD}, Algorithm~\ref{alg:lcd}, $\gamma\in(0,1]$ with closure enforced on repeat visits). First, parametric outcome switching on the persistent stage-$t$ solver modifies only the outcome-dependent right-hand sides, bounds, objective, and matrix coefficients, and therefore solves the same MIP $\min_{(a,y)\in X_{tw}(s)}\{c_{tw}(s,a,y)+\underline V_{t+1}^{(k)}(y)\}$ as dedicated per-outcome models. Proposition~\ref{prop:phase-hierarchy} and Theorem~\ref{thm:pooled-closing} ensure that every installed or in-place merged cut satisfies $h(y)\le F_{tw}^{(k)}(y)\le Q_{tw}(y)$ across $D_t$ and that $q_{tw}^{(k)}(s)=\max_{h\in\mathcal C_{tw}^{(k)}}h(s)$ and $\underline V_t^{(k)}(s)=\sum_{w\in\Omega_t}p_{tw}q_{tw}^{(k)}(s)$ are monotonically nondecreasing in $k$ on $D_t$.

We next show by backward induction on $t=T+1,T,\ldots,2$ that, along every sample path, the stage continuation $\underline V_t^{(k)}:D_t\to\mathbb R$ stabilizes permanently at some finite iteration $K_t<\infty$. At $t=T+1$, $\underline V_{T+1}^{(k)}\equiv V_{T+1}$ is fixed at $K_{T+1}=0$. Suppose that $\underline V_{t+1}^{(k)}\equiv\underline V_{t+1}^\star$ on $D_{t+1}$ for all $k\ge K_{t+1}$. Then, $F_{tw}^{(k)}\equiv F_{tw}^\star$ is fixed on $D_t$ for all $w\in\Omega_t$ and $k\ge K_{t+1}$. We partition the finite set $D_t=D_t^{\mathrm{fin}}\sqcup D_t^{\mathrm{inf}}$ into the states that are visited finitely many times and the states that are visited infinitely many times by forward passes over $k\ge 1$. The set $D_t^{\mathrm{fin}}$ is finite, and therefore there exists a finite iteration $K'_t\ge K_{t+1}$ after which forward passes visit exclusively states in $D_t^{\mathrm{inf}}$.

We fix any $(\hat s,w)\in D_t^{\mathrm{inf}}\times\Omega_t$ and consider the nonincreasing residual $\delta_k(\hat s)=F_{tw}^\star(\hat s)-q_{tw}^{(k)}(\hat s)\ge 0$ for $k\ge K'_t$. Under immediate exact closure ($\gamma=1$), the first visit to $\hat s$ after $K'_t$ installs a cut attaining $F_{tw}^\star(\hat s)$, which gives $\delta_k(\hat s)=0$. Under \textnormal{LCD} ($\gamma\in(0,1]$ with closure enforced on repeat visits), if $\delta_{k-1}(\hat s)>0$ on a visit $k\ge K'_t$ to $\hat s$, either that visit is already a repeat visit (triggering exact directional closure immediately via Theorem~\ref{thm:pooled-closing}(i) and setting $\delta_k(\hat s)=0$), or it is the first visit to $\hat s$, in which case $\hat s\in D_t^{\mathrm{inf}}$ is necessarily visited a second time at some $k'>k$, where the repeat-visit rule enforces exact directional closure whenever $\delta_{k'-1}(\hat s)>0$, which drives $\delta_{k'}(\hat s)=0$ in at most $|D_{\hat g}|$ directional steps. (Even if closure were gated purely by fractional progress $<\gamma\,\delta_{k-1}(\hat s)$ without a repeat-visit flag, whenever $E_{D_{\hat g}}(\hat s)<F_{tw}^\star(\hat s)$ every zero-penalty lattice cut satisfies $h(\hat s)\le E_{D_{\hat g}}(\hat s)$, so that $\delta(\hat s)\ge F_{tw}^\star(\hat s)-E_{D_{\hat g}}(\hat s)>0$ can contract by $(1-\gamma)$ at most finitely many visits before closure fires.) Because $D_t^{\mathrm{inf}}\times\Omega_t$ is finite, there exists a finite iteration $K''_t\ge K'_t$ such that $q_{tw}^{(k)}(\hat s)=F_{tw}^\star(\hat s)$ holds for all $(\hat s,w)\in D_t^{\mathrm{inf}}\times\Omega_t$ and all $k\ge K''_t$. For all $k>K''_t$, every backward visit at stage $t$ lands in $D_t^{\mathrm{inf}}$ with zero separation residual and is skipped by the separation gate ($\varepsilon_{\mathrm{sep}}=0$ in exact arithmetic, or reproduces already attained cut values), whereas stage-2 escape tightening and lazy coset evaluations can certify each of the finitely many pairs $(h,g)\in\mathcal C_{tw}^{(K''_t)}\times\mathcal G_t$ against fixed $F_{tw}^\star$ at most once. Therefore, $\underline V_t^{(k)}\equiv\underline V_t^\star$ freezes permanently on $D_t$ at a finite iteration $K_t\ge K''_t$.

Setting $K_\star=\max_{2\le t\le T}K_t<\infty$, all continuation models $\underline V_t^\star$ are frozen for $k\ge K_\star$, so that the fixed forward tie-breaking rule induces a single stationary Markov policy $\pi^\star$ starting from $s_1$. Scenario paths are i.i.d.\ from a law assigning positive probability to every positive-probability path over the finite tree. Therefore, the Borel--Cantelli lemma ensures that with probability one every positive-probability path is sampled infinitely often. Consequently, with probability one, every state $\hat s$ reached with positive probability under $\pi^\star$ at stage $t\ge 2$ belongs to $D_t^{\mathrm{inf}}$ and therefore satisfies $q_{tw}^\star(\hat s)=F_{tw}^\star(\hat s)$ for all $w\in\Omega_t$, which gives $\underline V_t^\star(\hat s)=\sum_{w\in\Omega_t}p_{tw}F_{tw}^\star(\hat s)$. Telescoping the greedy minimizing transitions of $\pi^\star$ backward from $\underline V_{T+1}^\star=V_{T+1}$ down to $t=1$ yields $L_{K_\star}=J(\pi^\star)$. Global cut validity gives $L_{K_\star}\le V_1(s_1)$, and feasibility of $\pi^\star$ gives $V_1(s_1)\le J(\pi^\star)=L_{K_\star}$, which proves $L_{K_\star}=V_1(s_1)=J(\pi^\star)$ almost surely.
\end{proof}

\begin{proof}[Proof of Theorem~\ref{thm:finite-pool-root-gap}]
We first establish dual feasibility. In the occupancy system \eqref{eq:occupation}, the incoming state masses $d_t(s)$ may be treated as unrestricted variables because nonnegativity of arc flows $f_e\ge 0$ and flow conservation automatically enforce $d_t(s)\ge 0$. Dualizing \eqref{eq:occupation} with multipliers $\nu_t(s)$ for incoming state mass conservation, $\eta_{tw}(s)$ for scenario mass conservation $\sum_{z,y}f_{twszy}=p_{tw}d_t(s)$, and $\beta_{tws}\in\mathbb R^n$ for the conditional barycenter equation $\sum_{z,y}(z-s)f_{twszy}=0$ maximizes $\nu_1(s_1)$ subject to
\[
\begin{aligned}
 \nu_t(s)&=\sum_{w\in\Omega_t} p_{tw}\eta_{tw}(s),\qquad \nu_{T+1}(y)=V_{T+1}(y),\\
 \eta_{tw}(s)+\beta_{tws}^{\top}(z-s)
 &\le c^0_{tw}(z,y)+\bar\rho\|z-s\|_1+\nu_{t+1}(y)
 \qquad\text{for all admissible arcs }(t,w,s,z,y).
\end{aligned}
\]
We substitute $(\nu,\eta,\beta)=(u,q,\pi)$, which satisfies all state-potential equalities with dual objective $u_1(s_1)=L$. On root arcs $(1,w,s_1,s_1,y)$, exact stage-1 optimality $q_{1w}(s_1)=F_{1w}(s_1)\le c^0_{1w}(s_1,y)+u_2(y)$ gives $g_e\ge 0$. On any nonroot arc $(t,w,s,z,y)$ with $z\in S_t(s)$, let $h\in\mathcal C_{tw}$ be the chosen cut active at $s$, written on the shared coset $[z]_t=[s]_t$ in the general directional form $h(z)=\pi^\top z+b-\Phi(z;c,\rho^+,\rho^-)$ with $\pi=\pi_{tws}$, coset intercept $b=\bar b([s]_t)$, center $c$, and directional penalties $0\le\rho_i^+,\rho_i^-\le\bar\rho$ (recovering isotropic closing cuts when $\rho_i^+=\rho_i^-=r\in[0,\bar\rho]$ and $\Phi(z;c,\rho^+,\rho^-)=r\|z-c\|_1$, and lattice cuts when $\rho_i^\pm=0$). From $q_{tw}(s)=h(s)=\pi^\top s+b-\Phi(s;c,\rho^+,\rho^-)$, we obtain $q_{tw}(s)+\pi^\top(z-s)=h(z)+\Phi(z;c,\rho^+,\rho^-)-\Phi(s;c,\rho^+,\rho^-)$. Adding and subtracting $F_{tw}(z)$ decomposes $g_{twszy}$ into three terms:
\[
\begin{aligned}
 g_{twszy}={}&\bigl[c^0_{tw}(z,y)+u_{t+1}(y)-F_{tw}(z)\bigr]
                  +\bigl[F_{tw}(z)-h(z)\bigr]\\
             &+\Bigl[\bar\rho\|z-s\|_1-\bigl(\Phi(z;c,\rho^+,\rho^-)-\Phi(s;c,\rho^+,\rho^-)\bigr)\Bigr].
\end{aligned}
\]
The first bracket is nonnegative by definition of the frozen stage minimum $F_{tw}(z)$, and the second is nonnegative by validity of $h$ for $F_{tw}$. For the third bracket, directional sublinearity $\Phi(z;c,\rho^+,\rho^-)-\Phi(s;c,\rho^+,\rho^-)\le\Phi(z;s,\rho^+,\rho^-)$ combined with $0\le\rho_i^+,\rho_i^-\le\bar\rho$ yields
\[
 \bar\rho\|z-s\|_1-\bigl(\Phi(z;c,\rho^+,\rho^-)-\Phi(s;c,\rho^+,\rho^-)\bigr)
 \ \ge\ \sum_{i=1}^n\Bigl[(\bar\rho-\rho_i^+)(z_i-s_i)_++(\bar\rho-\rho_i^-)(s_i-z_i)_+\Bigr]
 \ \ge\ 0.
\]
Therefore, $g_e\ge 0$ on every arc, which establishes dual feasibility.

Next, we derive the identity for the gap. For any primal occupancy flow $v=(d,f)\in\mathcal P_\Gamma$, summing $f_e g_e$ across all arcs eliminates the slope terms $\pi_{tws}^\top(z-s)$ via the conditional barycenter constraints $\sum_{z,y}(z-s)f_{twszy}=0$. Scenario mass balance gives $\sum_e f_e q_{tw}(s)=\sum_{t=1}^T\sum_s d_t(s)u_t(s)$, whereas interstage flow conservation gives $\sum_e f_e u_{t+1}(y)=\sum_{t=1}^T\sum_y d_{t+1}(y)u_{t+1}(y)$. Telescoping their difference across $t=1,\ldots,T$ leaves $\sum_y d_{T+1}(y)V_{T+1}(y)-u_1(s_1)=J(v)+\bar\rho D(v)-L$ after adding the immediate cost and penalty sums. Minimizing over $v\in\mathcal P_\Gamma$ yields $W_\Gamma(\bar\rho)-L$ by Theorem~\ref{thm:root-price}.

Finally, we prove the exactness test. Because $\mathcal P_\Gamma$ is a nonempty compact polytope and $g_e\ge 0$, the minimum $\min_{v\in\mathcal P_\Gamma}\sum_e f_e g_e$ equals zero if and only if $\mathcal P_\Gamma$ admits a flow supported exclusively on zero-slack arcs ($g_e=0$). At any stage-$t$ state $s$ with $d_t(s)>0$ in a zero-slack flow, normalizing $f_{twszy}/(p_{tw}d_t(s))$ for each $w\in\Omega_t$ ($p_{tw}>0$) yields a probability distribution on $S_t(s)$ with barycenter $s$ whose supported arcs have $g_{twszy}=0$ and reach successors $y$ with $d_{t+1}(y)>0$. Backward induction from $A_{T+1}=D_{T+1}$ places every supported copy in $B_{tws}$ and therefore $s\in\operatorname{conv}B_{tws}$ for all $w$, which gives $s\in A_t$. At $t=1$, $d_1(s_1)=1$ and $S_1(s_1)=\{s_1\}$ require a zero-slack root successor in $A_2$ for each $w\in\Omega_1$. Conversely, if each $w\in\Omega_1$ has $y_{1w}\in Y_{1w}(s_1)\cap A_2$ with $g_{1ws_1s_1y_{1w}}=0$, sending mass $p_{1w}$ to $y_{1w}$ at stage~1 and inductively splitting mass $p_{tw}d_t(s)$ at each reached $s\in A_t$ ($t\ge 2$) via a barycentric convex combination in $\operatorname{conv}B_{tws}$ routed to witnessing successors in $A_{t+1}$ along $g_{twszy}=0$ arcs constructs a feasible flow in $\mathcal P_\Gamma$ with $\sum_e f_e g_e=0$.
\end{proof}

\ecresetcounters
\section{Universal Box Support and Multiplier Subspace Refinement}\label{ec:box}

\subsection{Universal Affine Support on Integer Boxes}\label{app:box}

In binary-state SDDiP \citep{zou2019}, every state in $\{0,1\}^n$ is an extreme point of $[0,1]^n$, which guarantees that any finite stage value function admits an exact globally affine lower support at every state. On a general bounded integer box, a sublattice $\Gamma\subseteq\mathbb Z^n$ can reproduce this extreme-point geometry without state binarization by partitioning the box into coset fibers \citep[cf.][]{er2025}. In this appendix, we determine the exact minimum lattice index $[\mathbb Z^n:\Gamma]$ that is required when a single lattice must provide exact zero-penalty lattice cuts ($\rho=0$) for \emph{every} real-valued function on the box, which establishes the zero-penalty limit noted after Theorem~\ref{thm:mixture} and Corollary~\ref{cor:data-cap}.

Let $D\subset\mathbb Z^n$ be finite. We say that a full-rank sublattice $\Gamma\subseteq\mathbb Z^n$ admits \emph{universal pointwise support} on $D$ if, for every function $F:D\to\mathbb R$ and every anchor $s\in D$, there exists an affine function on the fiber $D\cap(s+\Gamma)$ that minorizes $F$ throughout $D\cap(s+\Gamma)$ and attains $F(s)$ at $s$. Geometrically, universal pointwise support holds if and only if every $s\in D$ is a vertex of $\operatorname{conv}(D\cap(s+\Gamma))$. Every vertex of a finite polytope is exposed, and positively scaling an exposing linear functional therefore yields an affine minorant passing through $(s,F(s))$ below all other points of the fiber. Conversely, if $s$ is a nontrivial convex combination of other points in $D\cap(s+\Gamma)$, then the indicator function equal to one at $s$ and zero elsewhere admits no affine support attaining its value at $s$. The following theorem gives the smallest index of a lattice that admits universal pointwise support on an integer box.

\begin{theorem}\label{thm:box}
Let $D=\prod_{i=1}^n\{l_i,\ldots,l_i+M_i\}$ with integers $l_i\in\mathbb Z$ and widths $M_i\in\mathbb Z_+$, and let $k_i=\lfloor M_i/2\rfloor+1$ for $i=1,\ldots,n$. Then, we have
\begin{equation}\label{eq:box-index}
 \min_{\Gamma\ \mathrm{admitting\ universal\ pointwise\ support\ on}\ D}
 [\mathbb Z^n:\Gamma] \ =\ \prod_{i=1}^n k_i .
\end{equation}
\end{theorem}
\begin{proof}
We may take $l_i=0$ for all $i$, because translating $D$ by $-l$ leaves lattice fibers and affine supports invariant. We first show that the index in \eqref{eq:box-index} is attained. Consider the rectangular lattice $\Gamma_0=\prod_{i=1}^n k_i\mathbb Z$, which has index $[\mathbb Z^n:\Gamma_0]=\prod_{i=1}^n k_i$. Because $2k_i=2\lfloor M_i/2\rfloor+2>M_i$, each coordinate $z_i$ within any fixed fiber $D\cap(s+\Gamma_0)$ takes at most two distinct values in $\{0,\ldots,M_i\}$. We fix an anchor $s\in D$ and, for each coordinate $i$, define $\nu_i=1/k_i$ if $s_i$ is the larger of two distinct fiber levels in coordinate $i$, $\nu_i=-1/k_i$ if $s_i$ is the smaller of two levels, and $\nu_i=0$ if coordinate $i$ takes a single value on the fiber. Then, for every copy $z\in D\cap(s+\Gamma_0)$, each differing coordinate contributes $(z_i-s_i)\nu_i=-1$, which gives $\nu^\top(z-s)=-|\{i:z_i\ne s_i\}|$. Setting $\Delta=\max\{0,\,F(s)-\min_{z\in D\cap(s+\Gamma_0)}F(z)\}$ (or replacing the minimum by any known finite lower bound on the fiber), the affine function $u(z)=F(s)+\Delta\,\nu^\top(z-s)$ satisfies $u(s)=F(s)$ and, for every $z\in D\cap(s+\Gamma_0)\setminus\{s\}$, $u(z)\le F(s)-\Delta\le F(z)$. Therefore, $\Gamma_0$ admits universal pointwise support on $D$.

Conversely, suppose that $[\mathbb Z^n:\Gamma]<\prod_{i=1}^n k_i$. By the pigeonhole principle, two distinct integer points in $\prod_{i=1}^n\{0,\ldots,k_i-1\}$ map to the same coset in $\mathbb Z^n/\Gamma$. Their nonzero difference $h\in\Gamma\setminus\{0\}$ satisfies $|h_i|\le k_i-1=\lfloor M_i/2\rfloor$, and therefore $2|h_i|\le M_i$ for every $i=1,\ldots,n$. We define $b\in\mathbb Z^n$ coordinatewise by $b_i=\max\{0,-2h_i\}$. Then, we have $0\le b_i,b_i+h_i,b_i+2h_i\le M_i$ for all $i$, so that the three collinear points $b$, $b+h$, and $b+2h$ all lie in $D\cap(b+\Gamma)$. We pick an index $j$ with $h_j\ne 0$. The single-stage modular balance
\[
 2|h_j|v + e^+ - e^- = z_j - b_j,\qquad v\in\mathbb Z,\quad e^+,e^-\in\mathbb Z_+,
\]
with objective $e^++e^-$ evaluates to $0$ at both $z=b$ and $z=b+2h$ and to $|h_j|\ge 1$ at the midpoint $z=b+h=\frac{1}{2}b+\frac{1}{2}(b+2h)$. Convex interpolation of the two endpoints forces every affine minorant on $D\cap(b+\Gamma)$ to be nonpositive at $b+h$, which rules out pointwise support there.
\end{proof}

For binary ($M_i=1$) or fixed ($M_i=0$) coordinates, we have $k_i=1$, which recovers $[\mathbb Z^n:\mathbb Z^n]=1$. On multilevel integer boxes, \eqref{eq:box-index} counts the ambient cosets needed for a worst-case, function-independent guarantee. In contrast, when stage congruence equations govern the periodic structure of $F$ (Section~\ref{sec:lattice}), far coarser compatible lattices suffice, and any residual nonperiodic gap is closed pointwise by closing cuts with direction-dependent or isotropic penalties without refining $\Gamma$ to index $\prod_i k_i$. We note that interstage compatibility \eqref{eq:compatible} imposes an independent coupling across stages beyond single-stage box geometry.

\subsection{Marginal Value of a New Multiplier Direction}\label{sec:direction}

In moderate-to-high state dimension, Lattice-Cut Decomposition (\textnormal{LCD}, Algorithm~\ref{alg:lcd}) refines the LP trial slope through Polyak copy-displacement steps and a low-dimensional active subspace spanned by recent cut slopes and copy subgradients \citep{chen2022} before deciding whether directional closure is needed (Section~\ref{sec:cut-subproblems}). On a fixed finite set $S\subset\mathbb Z^n$ of states in one coset that contains the visited anchor $s\in S$, we characterize the exact directional return to enlarging the current multiplier subspace by one orthogonal unit vector, which separates restricted-subspace slack from intrinsic nonconvexity within the coset.

Let $U\subseteq\mathbb R^n$ be a linear subspace of multipliers (after translating the origin to the current reference slope if affine), let $P_U$ denote the orthogonal projection onto $U$, and let $\Delta(S)$ denote the probability simplex on $S$. For any candidate direction $v\perp U$ with $\|v\|_2\le 1$ and step budget $\tau\ge 0$, we define the restricted dual cut height
\[
 H_v(\tau)=\max_{\pi\in U,\,|\lambda|\le\tau}
 \min_{z\in S}\bigl\{F(z)+(\pi+\lambda v)^\top(s-z)\bigr\}.
\]
For $\mu\in\Delta(S)$, let $d_\mu=\mathbb E_\mu Z-s$ denote its copy barycenter discrepancy, and let $\mathcal M_U=\{\mu\in\Delta(S):P_U d_\mu=0\}$ denote the set of distributions that are balanced along $U$. Let $D^*=\{d_\mu:\mu\in\arg\min_{\mathcal M_U}\mathbb E_\mu F(Z)\}$ collect the barycenter discrepancies of all optimal distributions at $\tau=0$. Because $\mathcal M_U$ is a nonempty compact polytope (containing unit mass at $s$) and $\mu\mapsto d_\mu$ is affine, the set $D^*$ is a nonempty compact convex subset of $U^\perp$. The following theorem gives the largest first-order gain in the cut height that one additional multiplier direction can provide.

\begin{theorem}\label{thm:direction}
The maximum right-hand derivative of the attainable cut height at $\tau=0$ is
\begin{equation}\label{eq:best-direction}
 \max_{v\perp U,\ \|v\|_2\le 1} H_v'(0+)
 \ =\ \operatorname{dist}_2(0,D^*).
\end{equation}
If $\operatorname{dist}_2(0,D^*)>0$, then the unique Euclidean projection $d^*$ of the origin onto $D^*$ yields an attaining unit direction $v^*=d^*/\|d^*\|_2$, which strictly increases $H_{v^*}(\tau)>H_{v^*}(0)$ for every $\tau>0$. If $\operatorname{dist}_2(0,D^*)=0$, then the subspace $U$ already attains the unrestricted full-slope coset envelope $E_S(s)$ at $s$.
\end{theorem}
\begin{proof}
We dualize the finite inner minimization over $z\in S$ with probability weights $\mu\in\Delta(S)$ and exchange $\max_{\pi\in U,|\lambda|\le\tau}$ with $\min_{\mu\in\Delta(S)}$ by LP duality. The free multipliers $\pi\in U$ enforce $P_U d_\mu=0$, whereas $|\lambda|\le\tau$ contributes $\max_{|\lambda|\le\tau}\lambda v^\top(-d_\mu)=\tau|v^\top d_\mu|$. Therefore, we have
\[
 H_v(\tau)=\min_{\mu\in\mathcal M_U}
 \bigl\{\mathbb E_\mu F(Z)+\tau |v^\top d_\mu|\bigr\}.
\]
We partition $\mathcal M_U$ into the two subpolytopes $\{\mu\in\mathcal M_U:v^\top d_\mu\ge 0\}$ and $\{\mu\in\mathcal M_U:v^\top d_\mu\le 0\}$, on each of which the objective is affine in $\mu$ for fixed $\tau$. Minimizing over their finitely many vertices expresses $H_v(\tau)$ as a concave piecewise-linear function of $\tau\ge 0$ whose right-hand derivative at $\tau=0$ equals the minimum slope over $\tau=0$ minimizers:
\[
 H_v'(0+)=\min_{d\in D^*}|v^\top d|.
\]
For every $v\perp U$ with $\|v\|_2\le 1$, Cauchy--Schwarz gives $\min_{d\in D^*}|v^\top d|\le |v^\top d^*|\le \|d^*\|_2=\operatorname{dist}_2(0,D^*)$. When $d^*\ne 0$, the variational inequality for Euclidean projection onto the convex set $D^*$ gives $(d-d^*)^\top d^*\ge 0$, and therefore $d^\top d^*\ge \|d^*\|_2^2>0$, for all $d\in D^*$. Setting $v^*=d^*/\|d^*\|_2\in U^\perp$ yields $(v^*)^\top d\ge \|d^*\|_2$ across $D^*$, which achieves equality in \eqref{eq:best-direction}. Because $|(v^*)^\top d_\mu|>0$ on the entire $\tau=0$ optimal face of $\mathcal M_U$, compactness implies $H_{v^*}(\tau)>H_{v^*}(0)$ for every $\tau>0$. Finally, $\operatorname{dist}_2(0,D^*)=0$ holds if and only if $0\in D^*$, that is, some $\tau=0$ minimizer $\mu\in\mathcal M_U$ already satisfies $d_\mu=0$, so that $\min_{\mathcal M_U}\mathbb E_\mu F(Z)=\min_{\mathcal M(S,s)}\mathbb E_\mu F(Z)=E_S(s)$.
\end{proof}

Equation~\eqref{eq:best-direction} distinguishes two sources of gap at $s$ and motivates the Polyak-plus-subspace dual step in \textnormal{LCD}. When $\operatorname{dist}_2(0,D^*)>0$ (for instance, at $U=\{0\}$ when prior cuts share a zero LP slope), a step along $v^*$ gives a strict first-order gain within the fixed coset. The direction $d^*$ is a convex combination of the displacements $z-s$ of the copies that support an optimal distribution. A step along the displacement of a single copy need not give a first-order gain. For example, for $U=\{0\}$, $s=(0,0)$, $S=\{(0,0),(1,1),(1,-1)\}$, $F(s)=1$, and $F=0$ at the other two points, we have $D^*=\operatorname{conv}\{(1,1),(1,-1)\}$ and $d^*=(1,0)$, whereas $H_v'(0+)=0$ for $v=(1,1)/\sqrt2$. For this reason, the dual step of \textnormal{LCD} does not stop at steps along single displacements and optimizes over the span of several of them (Table~\ref{tab:ec-params}). In contrast, when $0\in D^*$, any remaining gap $F(s)-E_S(s)>0$ survives optimization over all affine slopes in $\mathbb R^n$ and requires either lattice refinement (Theorem~\ref{thm:mixture}) or positive directional closure (Theorem~\ref{thm:local-price}).

\ecresetcounters
\section{Single-Stage Certificates}\label{ec:local-certs}

In this appendix, we establish the three local finite-data certificates summarized in Remark~\ref{rem:ec-certificates}: directional cut accuracy along a lattice segment (Proposition~\ref{prop:directional-coverage}), LP certification of the complete single-stage bound from feasible state evaluations (Proposition~\ref{prop:finite-pool-certificate}), and regional stage-2 coverage conditions for a stage-1 bound gain (Corollary~\ref{cor:root-gain-coverage}).

\subsection{Directional Coverage Along a Lattice Segment}

Whenever a feasible action found at $\hat s$ remains feasible along a lattice direction $d\in\Gamma_t\setminus\{0\}$, it provides a constructive upper bound on $F_{tw}$ along that discrete ray. Because downstream closing cuts introduce piecewise-linear kinks at their centers, the coordinate crossings identify a small finite test set that certifies cut accuracy over the entire integer segment. The following proposition makes this precise.

\begin{proposition}\label{prop:directional-coverage}
Fix $t<T$, an outcome $w\in\Omega_t$, and the multicut continuation $\underline V_{t+1}$ in \eqref{eq:multicut}, and write $F=F_{tw}$ as in \eqref{eq:current}. Let $h\le F$ be a valid lattice or closing cut with slope $\pi$, penalty $\rho\ge 0$, and center $\hat s$. Consider the discrete lattice segment $z_k=\hat s+kd$ for $d\in\Gamma_t\setminus\{0\}$ and integers $0\le k\le K$ ($K\ge 1$). Suppose that a common action (or affine action sequence) is feasible at every $z_k$, incurs immediate cost $c_0+k\gamma$, and produces successor $y_k=\hat y+kv$ with $v\in\Gamma_{t+1}$. Let
\[
 U_k=c_0+k\gamma+\underline V_{t+1}(y_k),\qquad C(k)=U_k-h(z_k).
\]
Then, we have $0\le F(z_k)-h(z_k)\le C(k)$ for all $k=0,\ldots,K$, and the certificate slack $C(k)-[F(z_k)-h(z_k)]=U_k-F(z_k)$ equals the suboptimality of the reused action at $z_k$.

For each downstream stage-$(t+1)$ outcome $q\in\Omega_{t+1}$ and installed closing cut $j\in\mathcal C_{t+1,q}$ with positive penalty $\rho_{qj}>0$ and center $c_{qj}$, collect the distinct interior ray parameters where a moving successor coordinate crosses the corresponding center coordinate:
\[
 B=(0,K)\cap
 \left\{\frac{c_{qj,i}-\hat y_i}{v_i}:\ \rho_{qj}>0,\ v_i\ne 0\right\},
 \qquad R=|B|.
\]
For every tolerance $\varepsilon\ge 0$, the integer indices $\{k\in\{0,\ldots,K\}:C(k)\le\varepsilon\}$ form a union of at most $R+1$ intervals of consecutive integers. Moreover, with $P=\{0,K\}\cup\{\lfloor b\rfloor,\lceil b\rceil:b\in B\}$, we have
\[
 \max_{k=0,\ldots,K} C(k)=\max_{k\in P}C(k),\qquad
 |P|\le\min\{K+1,\,2R+2\}.
\]
\end{proposition}
\begin{proof}
Primal feasibility of the reused action at $z_k$ gives $F(z_k)\le U_k$, and cut validity gives $h(z_k)\le F(z_k)$, which yields $0\le F(z_k)-h(z_k)\le C(k)$ and $C(k)-[F(z_k)-h(z_k)]=U_k-F(z_k)$.

We have $d\in\Gamma_t$ and $v\in\Gamma_{t+1}$, and therefore every $z_k$ shares the coset $[\hat s]_t$ and every $y_k$ shares the successor coset $[\hat y]_{t+1}$. Fixing those coset intercepts in $h$ and in every installed cut in $\mathcal C_{t+1,q}$, we extend the algebraic formula for $C(k)$ to real $x\in[0,K]$ through $z(x)=\hat s+xd$ and $y(x)=\hat y+xv$. Because $x\ge 0$ and $h$ is centered at $\hat s$, we have $\|z(x)-\hat s\|_1=x\|d\|_1$ (and likewise $\Phi(z(x);\hat s,\rho^+,\rho^-)=x\,\Phi(\hat s+d;\hat s,\rho^+,\rho^-)$ under \eqref{eq:aniso-penalty}), so that $h(z(x))$ is affine on $[0,K]$. We partition $[0,K]$ into $R+1$ closed subintervals through $\{0,K\}\cup B$. On the interior of each subinterval, no coordinate $y_i(x)=\hat y_i+xv_i$ changes sign relative to $c_{qj,i}$, and therefore, by continuity, every downstream penalty term $-\rho_{qj}\|y(x)-c_{qj}\|_1$ (or $-\Phi(y(x);c_{qj},\rho_{qj}^+,\rho_{qj}^-)$, where $B$ need only collect coordinates with $\rho_{qj,i}^++\rho_{qj,i}^->0$) is affine on the closed subinterval. Therefore, each outcome envelope $\max_{j\in\mathcal C_{t+1,q}}h_{qj}(y(x))$ is convex on each subinterval, as is $C(x)$ after taking nonnegative expectations, adding $c_0+x\gamma$, and subtracting affine $h(z(x))$.

On each subinterval, $\{x:C(x)\le\varepsilon\}$ is a convex real interval whose integer points are consecutive, which yields at most $R+1$ discrete intervals overall. On any closed subinterval $[a,b]$ containing integers, its extreme integers $\lceil a\rceil,\lfloor b\rfloor$ lie in $P$, and convexity bounds $C(k)\le\max\{C(\lceil a\rceil),C(\lfloor b\rfloor)\}$, which proves $\max_{0\le k\le K}C(k)=\max_{k\in P}C(k)$ and $|P|\le\min\{K+1,2R+2\}$.
\end{proof}

When all downstream cuts are zero-penalty lattice cuts ($\rho_{qj}=0$), we have $B=\varnothing$ and $R=0$ even if $h$ has $\rho>0$. In this case, $P=\{0,K\}$, so that $C(0)\le\varepsilon$ and a single evaluation $C(K)\le\varepsilon$ certify the entire segment.

\subsection{Certifying the Complete Single-Stage Bound from Feasible Actions}

At a fixed nonroot stage $t$, outcome $w$, and continuation model, we write $F=F_{tw}$. Feasible actions that are recorded across visits supply upper bounds $U_z\ge F(z)$ on $Z\subseteq S=D_t\cap(s+\Gamma_t)$, which yield a small LP upper bound on the complete-family single-stage operator \eqref{eq:priced-envelope}. The following proposition states this bound and characterizes when it coincides with the value of the current cut pool.

\begin{proposition}\label{prop:finite-pool-certificate}
Let $S=D_t\cap(s+\Gamma_t)$, and denote the complete-family single-stage value \eqref{eq:priced-envelope} at penalty cap $\bar\rho\ge 0$ by $H_{\bar\rho}F(s)=\min_{\mu\in\mathcal M(S,s)}\sum_{z\in S}\mu_z(F(z)+\bar\rho\|z-s\|_1)$. Let $H_{\mathrm{pool}}(s)$ be the maximum at $s$ of a nonempty finite pool of cuts valid for $F$ on $S$, each having penalty in $[0,\bar\rho]$. For any subset $Z\subseteq S$ with $s\in\operatorname{conv}Z$ and known finite upper bounds $U_z\ge F(z)$ ($z\in Z$), define
\[
 M_Z(s)=\min_{\lambda\in\mathcal M(Z,s)}
          \sum_{z\in Z}\lambda_z\bigl(U_z+\bar\rho\|z-s\|_1\bigr).
\]
Then, we have $H_{\mathrm{pool}}(s)\le H_{\bar\rho}F(s)\le M_Z(s)$. Enlarging $Z$ or tightening $U_z$ weakly decreases $M_Z(s)$, and if $P_s$ denotes the minimal face of $\operatorname{conv}Z$ containing $s$, then an optimizer uses at most $\dim P_s+1$ copies.

Next, consider any cut $h_*$ that is active in the pool at $s$, written on $S$ as $h_*(z)=\pi_*^\top z+b_*-r_*\|z-c_*\|_1$ with $r_*\in[0,\bar\rho]$. Let $E_*=\{z\in Z:U_z=h_*(z)\}$ and $L_*(s)=\{z\in\mathbb R^n:z_i=s_i\text{ whenever }s_i\ne c_{*,i}\}$. Then, the following conditions are necessary and sufficient for $M_Z(s)=H_{\mathrm{pool}}(s)$:
\[
\begin{array}{ll}
 \bar\rho=0\,:             & s\in\operatorname{conv}E_*,\\[2pt]
 0\le r_*<\bar\rho\,:      & s\in E_*,\\[2pt]
 r_*=\bar\rho>0\,:         & s\in\operatorname{conv}\bigl(E_*\cap L_*(s)\bigr).
\end{array}
\]
\end{proposition}
\begin{proof}
Consider the affine function $\ell_s(z)=h_*(s)+\pi_*^\top(z-s)$. For every $z\in S$ and $u\ge F(z)$, we have
\[
\begin{aligned}
 u+\bar\rho\|z-s\|_1-\ell_s(z)&=\bigl[u-h_*(z)\bigr]+D_*(z;s),\\
 D_*(z;s)&=\bar\rho\|z-s\|_1-r_*\bigl(\|z-c_*\|_1-\|s-c_*\|_1\bigr)
          \ge(\bar\rho-r_*)\|z-s\|_1\ge 0,
\end{aligned}
\]
by cut validity and the reverse triangle inequality. Because $\ell_s$ has mean $h_*(s)=H_{\mathrm{pool}}(s)$ under every distribution in $\mathcal M(S,s)$, taking $u=F(z)$ on $S$ gives $H_{\mathrm{pool}}(s)\le H_{\bar\rho}F(s)$, and restricting to $\mathcal M(Z,s)$ with $u=U_z$ gives $H_{\bar\rho}F(s)\le M_Z(s)$ along with
\[
 M_Z(s)-H_{\mathrm{pool}}(s)=
 \min_{\lambda\in\mathcal M(Z,s)}
 \sum_{z\in Z}\lambda_z\bigl[U_z-h_*(z)+D_*(z;s)\bigr].
\]
Enlarging $Z$ or lowering $U_z$ weakly decreases this minimum, and any extreme-point optimizer in the minimal face $P_s$ has affinely independent support of size at most $\dim P_s+1$.

The gap vanishes if and only if some $\lambda\in\mathcal M(Z,s)$ is supported on $\{z\in E_*:D_*(z;s)=0\}$. At $\bar\rho=0$, $D_*\equiv 0$ gives $s\in\operatorname{conv}E_*$. At $r_*<\bar\rho$, $D_*(z;s)\ge(\bar\rho-r_*)\|z-s\|_1$ forces $z=s$, that is, $s\in E_*$. At $r_*=\bar\rho>0$, $D_*(z;s)=0$ requires coordinatewise equality $(z_i-s_i)(s_i-c_{*,i})\ge 0$. Whenever $s_i\ne c_{*,i}$, the zero mean forces $z_i=s_i$ across the support, which places all support points in $E_*\cap L_*(s)$, whereas every $z\in L_*(s)$ satisfies $\|z-c_*\|_1-\|s-c_*\|_1=\|z-s\|_1$ and $D_*(z;s)=0$.
\end{proof}

\subsection{Regional Cut Coverage Required for a Root-Bound Gain}

We consider an update of the stage-2 pools $\{\mathcal C_{2q}\}_{q\in\Omega_2}$ with downstream continuations fixed, and we recall \eqref{eq:rootgain}. Let $Q_q^-,Q_q^+$ denote the old and updated envelopes for $F_q=F_{2q}$. Given new cuts $h$ with certified error $F_q(y)-h(y)\le B_{qh}(y)$ on $C_{qh}\subseteq D_2$ and certified old shortfalls $\ell_q(y)\le F_q(y)-Q_q^-(y)$, we define
\[
 d_q(y)=\max\Bigl(\{0\}\cup
       \bigl\{\ell_q(y)-B_{qh}(y):h\text{ is new for }q,\ y\in C_{qh}\bigr\}\Bigr),
 \qquad d(y)=\sum_{q\in\Omega_2} p_{2q}d_q(y),
\]
so that $0\le d(y)\le\Delta(y)=\sum_{q\in\Omega_2}p_{2q}(Q_q^+(y)-Q_q^-(y))$. The following corollary states the coverage that these certified increases must have for the root bound to rise.

\begin{corollary}\label{cor:root-gain-coverage}
In \eqref{eq:rootgain}, let $a_w(y)=J_w(y)-j_w\ge 0$ for $w\in\Omega_1$ and $y\in Y_w$. For any certified increase $0\le d\le\Delta$, the root bound rises by at least $\sum_{w\in\Omega_1} p_{1w}G_w(d)$, where $G_w(d)=\min_{y\in Y_w}\{a_w(y)+d(y)\}$. For $\tau_w\ge 0$, we have
\[
 G_w(d)\ge\tau_w\quad\Longleftrightarrow\quad
 d(y)\ge\tau_w-a_w(y)\quad\text{for every }y\in Y_w\text{ with }a_w(y)<\tau_w .
\]
With $\kappa_w(b)=\min\{a_w(y):y\in Y_w,\ d(y)\le b\}$ ($+\infty$ if empty), we have $G_w(d)=\min_{b\in d(Y_w)}\{b+\kappa_w(b)\}$ and, in particular, $G_w(b\mathbf 1_C)=\min\{b,\,\min_{y\in Y_w\setminus C}a_w(y)\}$ for $b>0$ and $C\subseteq D_2$.
\end{corollary}
\begin{proof}
Monotonicity of \eqref{eq:rootgain} gives $L^+-L\ge\sum_{w\in\Omega_1}p_{1w}G_w(d)$. The inequality $a_w(y)+d(y)\ge\tau_w$ holds automatically whenever $a_w(y)\ge\tau_w$, because $d\ge 0$. For any $b\in d(Y_w)$ and $y\in Y_w$ with $d(y)\le b$, $G_w(d)\le a_w(y)+b$ implies $G_w(d)\le b+\kappa_w(b)$, whereas taking $b^*=d(y^*)$ at $y^*\in\arg\min_{Y_w}(a_w+d)$ gives the reverse inequality. For $d=b\mathbf 1_C$, the levels $\{0,b\}$ have $\{d\le 0\}=Y_w\setminus C$ and $\kappa_w(b)=\min_{Y_w}a_w=0$.
\end{proof}

Restricting the old scenario-$w$ root MIP to $\{y\in Y_w:\widetilde d(y)\le b\}$ at finitely many levels $b<\tau_w$ and checking $\underline m_w(b)\ge\overline j_w+\tau_w-b$ (for solver lower bound $\underline m_w(b)$ and upper bound $\overline j_w\ge j_w$) certifies $G_w(\widetilde d)\ge\tau_w$, and rounding $d$ down to $\widetilde d$ with $0\le d-\widetilde d\le\varepsilon_w$ loses at most $\sum_{w\in\Omega_1}p_{1w}\varepsilon_w$.

\ecresetcounters
\section{Horizon-Wide Certificates}\label{ec:horizon-certs}

In this appendix, we extend finite-data certification across the horizon through recorded-action occupancy flows (Proposition~\ref{prop:record-root-certificate}) and quotient Markov replay of frozen greedy policies on cosets (Proposition~\ref{prop:phase-policy-replay}).

\subsection{Root Intervals from Recorded Feasible Actions}

For each $(t,w)$, let $\mathcal D_{tw}$ be a finite set of feasible transition records $i=(z_i,a_i,y_i,c_i)$ with $(a_i,y_i)\in X_{tw}(z_i)$ and $c_i=c_{tw}(z_i,a_i,y_i)\ge c^0_{tw}(z_i,y_i)$, and let $I_{tw}(s)=\{i\in\mathcal D_{tw}:z_i\in S_t(s)\}$ ($S_1(s_1)=\{s_1\}$). Let $\mathcal P_{\mathcal D}$ be the polytope of record flows $(d,x)$ satisfying \eqref{eq:occupation} with outgoing masses $x_{twsi}\ge 0$ on $i\in I_{tw}(s)$, and define
\[
 \mathsf C_{\mathcal D}(d,x)=\sum_{t,w,s}\sum_{i\in I_{tw}(s)} x_{twsi}
       \bigl(c_i+\bar\rho\|z_i-s\|_1\bigr)
       +\sum_{y\in D_{T+1}} d_{T+1}(y)V_{T+1}(y),
\]
and $U_{\mathcal D}=\min_{(d,x)\in\mathcal P_{\mathcal D}}\mathsf C_{\mathcal D}(d,x)$ (with $U_{\mathcal D}=+\infty$ if $\mathcal P_{\mathcal D}=\varnothing$). The following proposition bounds the finite-pool loss by an interval that is computed from the recorded actions.

\begin{proposition}\label{prop:record-root-certificate}
With superscript $k$ denoting iteration $k$ in Theorem~\ref{thm:finite-pool-root-gap}, define $\widehat g^k_{twsi}=c_i+u^k_{t+1}(y_i)+\bar\rho\|z_i-s\|_1-q^k_{tw}(s)-(\pi^k_{tws})^\top(z_i-s)$ for $i\in I_{tw}(s)$. Whenever $\mathcal P_{\mathcal D}\ne\varnothing$, we have $\widehat g^k_{twsi}\ge 0$ and
\[
 0 \ \le\ W_\Gamma(\bar\rho)-L_k
 \ \le\ U_{\mathcal D}-L_k
 \ =\ \min_{(d,x)\in\mathcal P_{\mathcal D}}\sum_{t,w,s}\sum_{i\in I_{tw}(s)}x_{twsi}\widehat g^k_{twsi},
\]
where $U_{\mathcal D}$ and its minimizing flows in $\mathcal P_{\mathcal D}$ are independent of $k$. Let $\widehat A_{T+1}=D_{T+1}$ and $\widehat A_t=\{s:s\in\operatorname{conv}\{z_i:i\in I_{tw}(s)\text{ qualifying},\ y_i\in\widehat A_{t+1}\}\ \forall\,w\in\Omega_t\}$ for $t=T,\ldots,2$. Then, the existence for each $w\in\Omega_1$ of a qualifying $i\in I_{1w}(s_1)$ with $y_i\in\widehat A_2$ characterizes $U_{\mathcal D}<+\infty$ when all records qualify, and characterizes $U_{\mathcal D}=L_k$ when only records with $\widehat g^k_{twsi}=0$ qualify.
\end{proposition}
\begin{proof}
Aggregating $f_{twszy}=\sum_{i\in I_{tw}(s):(z_i,y_i)=(z,y)}x_{twsi}$ maps $(d,x)\in\mathcal P_{\mathcal D}$ into $v=(d,f)\in\mathcal P_\Gamma$ with $c_i\ge c^0_{tw}(z_i,y_i)$, which gives $W_\Gamma(\bar\rho)\le\mathsf C_{\mathcal D}(d,x)$, $\widehat g^k_{twsi}=g^k_e+(c_i-c^0_{tw}(z_i,y_i))\ge 0$, and $\sum x_{twsi}\widehat g^k_{twsi}=\mathsf C_{\mathcal D}(d,x)-L_k$ by Theorem~\ref{thm:finite-pool-root-gap}. Minimizing over $\mathcal P_{\mathcal D}$ gives the interval. The bound $U_{\mathcal D}$ and its minimizers are independent of $k$, because $k$ only subtracts the constant $L_k$. Normalizing positive outgoing masses yields barycentric distributions on qualifying copies, and routing barycentric weights forward from $s_1$ constructs a flow in $\mathcal P_{\mathcal D}$, which proves both hull tests because $\widehat g^k_{twsi}\ge 0$.
\end{proof}

For nested records $\mathcal D_k\subseteq\mathcal D_{k+1}$ and $L_{k+1}\ge L_k$, the interval widths $E_k=U_{\mathcal D_k}-L_k$ satisfy $E_k-E_{k+1}=(L_{k+1}-L_k)+(U_{\mathcal D_k}-U_{\mathcal D_{k+1}})$, which separates the cut gain from the primal record enrichment.

\subsection{Replay of a Frozen Policy on Cosets}

We fix an iteration $k$ and a physical Markov policy $\mu^k$ that greedily minimizes $c_{tw}(s,a,y)+u^k_{t+1}(y)$ over $X_{tw}(s)$, with successor-coset law $P_{tw}^{\mu^k}(g'\mid s)$ and optimal successor-coset set $\mathcal H_{tw}^k(s)\subseteq\mathcal G_{t+1}$. The following proposition characterizes when this policy can be replayed on cosets and bounds its cost.

\begin{proposition}\label{prop:phase-policy-replay}
The policy $\mu^k$ admits an exact coset kernel $\bar P_{tw}(\cdot\mid g)$ across all representatives $s$ of $g$ if and only if $P_{tw}^{\mu^k}(\cdot\mid s)$ is constant on each coset $g$. Such an optimal rule exists if and only if
\[
 I_{tw}^k(g)=\bigcap_{s\in D_t:\,g_t(s)=g}\mathcal H_{tw}^k(s)\ \ne\ \varnothing
\]
for every occupied $g\in\mathcal G_t$ and $w\in\Omega_t$, with attainable common kernels equal to distributions supported on $I_{tw}^k(g)$. Under $\bar P_{tw}$, the coset distribution $\omega_t(g)=\Pr_{\mu^k}\{g_t(s_t)=g\}$ satisfies $\omega_1(g)=\mathbf 1_{\{g=g_1(s_1)\}}$ and $\omega_{t+1}(g')=\sum_{g,w}\omega_t(g)p_{tw}\bar P_{tw}(g'\mid g)$. Moreover, if $r^k_{tw}(s)=F^k_{tw}(s)-q^k_{tw}(s)\le\bar r^k_{twg}$ on positive-probability states of coset $g$ ($t\ge 2$, with $q^k_{1w}(s_1)=F^k_{1w}(s_1)$), then
\[
 0 \ \le\ V^\star - L_k \ \le\ J(\mu^k) - L_k
 \ =\ \mathbb E_{\mu^k}\sum_{t=2}^T r^k_{t,w_t}(s_t)
 \ \le\ \sum_{t=2}^T\sum_{w\in\Omega_t}\sum_{g\in\mathcal G_t} p_{tw}\omega_t(g)\bar r^k_{twg}.
\]
\end{proposition}
\begin{proof}
Invariance of $P_{tw}^{\mu^k}(\cdot\mid s)=\bar P_{tw}(\cdot\mid g)$ within each coset $g$ is necessary and sufficient for every history ending in $g$ to share kernel $\bar P_{tw}(\cdot\mid g)$, which gives the forward recursion by stagewise independence. Common optimal rules must be supported in $\bigcap_{g_t(s)=g}\mathcal H_{tw}^k(s)=I_{tw}^k(g)$, and mixing representative minimizers realizes any distribution on $I_{tw}^k(g)$. Greedy optimality gives $c_{t,w_t}+u^k_{t+1}(s_{t+1})=q^k_{t,w_t}(s_t)+r^k_{t,w_t}(s_t)$ with $\mathbb E[q^k_{t,w_t}(s_t)\mid s_t]=u^k_t(s_t)$, $r^k_{1,w_1}(s_1)=0$, and $u^k_{T+1}=V_{T+1}$. Taking expectations telescopes intermediate values to $J(\mu^k)-L_k=\mathbb E_{\mu^k}\sum_{t=2}^T r^k_{t,w_t}(s_t)$, and bounding $r^k_{tw}\le\bar r^k_{twg}$ yields the coset-marginal bound.
\end{proof}

\ecresetcounters
\section{Computational Details}\label{ec:benchmarks}

In this appendix, we list the default parameters of \textnormal{LCD} and \textnormal{LCD-Iso} in Table~\ref{tab:ec-params}, and we break the catalog run of Section~\ref{sec:exp-regimes} down by family in Table~\ref{tab:ec-catalog}. The catalog has 16 size configurations for each family of Group~I and 9 for each of the other families, 163 in total, and each configuration was evaluated at one instance seed. The four methods had 20 seconds each on 18 configurations and 25 seconds each on the remaining 145. We count a bound of zero returned by a run that completed no iteration as no bound.

\begin{table}[!htbp]
\centering
\caption{Default Algorithmic, Subproblem, and Numerical Parameters in \textnormal{LCD} (Algorithm~\ref{alg:lcd})}\label{tab:ec-params}
\vspace{3pt}
\begin{singlespace}\scriptsize
\setlength{\tabcolsep}{4.0pt}
\renewcommand{\arraystretch}{1.14}
\begin{tabular}{@{}>{\raggedright\arraybackslash}p{1.32in}>{\raggedright\arraybackslash}p{1.22in}>{\raggedright\arraybackslash}p{3.62in}@{}}
\toprule
Parameter / Mechanism & Default Value & Operational Role in \textnormal{LCD} (and \textnormal{LCD-Iso} Comparison) \\
\midrule
One-hot threshold $M_{\mathrm{oh}}$ & $64$ ($256$ in \textnormal{LCD-Iso}) & Shared one-hot $\delta_g\in\{0,1\}$ if $[\mathbb Z^n:\Gamma_t]\le M_{\mathrm{oh}}$; binary residue digits $\zeta_{j\ell}$ with stage-pooled $\mu_g\in[0,1]$ otherwise (Proposition~\ref{prop:phase-hierarchy}(iii)). \\
Separation gate $\varepsilon_{\mathrm{sep}}$ & $10^{-7}$ (rel.) & Skips backward separation at $(t,w,\hat s)$ when $\max_{h\in\mathcal C_{tw}}h(\hat s)\ge U_{\hat s}-\varepsilon_{\mathrm{sep}}(1+|U_{\hat s}|)$. \\
Dual-refinement tol.\ $\varepsilon_{\mathrm{tight}}$ & $10^{-6}$ (rel.) & Triggers $2$ Polyak copy-displacement steps plus $\le 5$ active-subspace level-bundle iterations \citep{chen2022} (span $6$; $\le 1$ once $\ge 8$ closing cuts exist). \\
LP floor skip tol.\ $\varepsilon_{\mathrm{floor}}$ & $10^{-5}$ (rel.) & Certifies $b_0=b_{\mathrm{LP}}(\pi)$ with zero full-box MIPs when $b_{\hat g}(\pi)-b_{\mathrm{LP}}(\pi)\le\varepsilon_{\mathrm{floor}}(1+|b_{\hat g}|)$; in-place merge at $\|\pi-\pi'\|_\infty\le 10^{-8}$ with box-shift subtraction. \\
Closure miss tol.\ $\varepsilon_{\mathrm{close}}$ & $10^{-6}$ rel.\ / $10^{-7}$ abs. & Triggers directional lifting \citep{dinkelbach1967} ($\le 8$ solves in \textnormal{LCD}; $\le 12$ isotropic solves in \textnormal{LCD-Iso}) on proven coset shortfall. \\
First-visit progress gate $\gamma$ & $0.5$ ($\gamma=1$ in \textnormal{LCD-Iso}) & Defers first-visit closure ($n\ge 3$) when the lattice cut closes a fraction $\ge\gamma$ of shortfall ($0.75$ for $3\le n\le 10$, $0.25$ for $n>10$); unconditional on repeat visits; $0$-MIP LP closure if $<1.25\,\mathrm{s}$ left. \\
Warm-up \& escape caps & $\le 30\%$ budget & Corner-probed LP warm-up ($\le 6$ paths, $\le 7.5\,\mathrm{s}$; $1$-pass exit if $\pi^{\mathrm{LP}}\equiv 0$); up to $36$ stage-2 escape-coset rounds ($\le 6.0\,\mathrm{s}$). \\
MIP \& feasibility tols. & $10^{-9}$ & Gurobi \texttt{MIPGap}, \texttt{IntFeasTol}, \texttt{FeasibilityTol}, and \texttt{OptimalityTol} on integer stage and coset subproblems. \\
Dual-bound safety margin & $10^{-7}$ rel.\ / $10^{-6}$ abs. & Subtracted from solver \texttt{ObjBound} before installing any intercept; post-hoc validity gate $L_k\le V^\star+10^{-5}$, tie tol.\ $10^{-4}\max(1,|a|,|b|)$. \\
Thread \& timer discipline & $1$ CPU thread & \texttt{Params.Threads=1} audited per solve; wall-clock timer includes instance load, Smith form, and model build; no information about $V^\star$ in the runs with limits of 20 to 300 seconds. \\
\bottomrule
\end{tabular}
\end{singlespace}
\end{table}

\begin{table}[!htbp]
\centering
\caption{Comparison of the Bounds of \textnormal{LCD} and SDDmiP by Family in the Catalog Runs (Time Limits of 20 and 25 Seconds)}\label{tab:ec-catalog}
\begin{singlespace}\footnotesize
\begin{tabular}{@{}lrrrrr@{}}
\toprule
Family & Configurations & \textnormal{LCD} higher & Tie & \textnormal{LCD} lower & Neither has a bound \\
\midrule
Shift & 16 & 3 & 0 & 3 & 10 \\
Rotate & 16 & 1 & 2 & 3 & 10 \\
Odd-Cycle & 16 & 3 & 0 & 0 & 13 \\
Shear & 16 & 3 & 0 & 0 & 13 \\
Hub-Consolidation & 9 & 4 & 0 & 0 & 5 \\
Sterile-Trays & 9 & 2 & 0 & 2 & 5 \\
Pooled-Capacity & 9 & 1 & 0 & 8 & 0 \\
Fulfillment-Network & 9 & 1 & 0 & 2 & 6 \\
Case-Pack & 9 & 0 & 0 & 9 & 0 \\
Workforce-Cohort & 9 & 3 & 0 & 6 & 0 \\
Co-Product & 9 & 6 & 3 & 0 & 0 \\
Bike-Sharing & 9 & 9 & 0 & 0 & 0 \\
Serial-Inventory & 9 & 4 & 0 & 5 & 0 \\
Lead-Time-Pipeline & 9 & 2 & 0 & 7 & 0 \\
Airline-RM & 9 & 4 & 0 & 5 & 0 \\
\midrule
All families & 163 & 46 & 5 & 50 & 62 \\
 
\bottomrule
\end{tabular}
\end{singlespace}
\end{table}

\ecresetcounters
\section{Benchmark Formulations}\label{app:benchmarks}
In this appendix, we give the formulations of the 15 benchmark families of Table~\ref{tab:families-lattices} and, for each family, the label that identifies the cosets of its stage lattices. All models minimize the expected total cost, have bounded integer states, and have zero terminal cost.

Section~\ref{sec:experiments} uses two collections of instances. The \emph{catalog} has 163 size configurations: 16 for each of the four families of Group~I and 9 for each of the other 11 families. Each configuration is run with one instance seed. The \emph{long-run instances} belong to Group~III. They comprise 32 instances with one item, cohort, recipe, or station and 10, 24, or 48 stages (Section~\ref{sec:exp-long}), and 30 instances with multiple items, with the truck fixed charge, or of full size with 6 stages (Table~\ref{tab:scope}, Panel~A). The 45-second runs use 48 instances from the same generators, and the 300-second runs use six instances of Group~III and six of Group~I. Appendix~\ref{app:results} lists every instance with its size.

\subsection{Group I. Synthetic Batch Models}\label{axbench:modular}

Shift, Rotate, Odd-Cycle, and Shear have the same project selection problem and differ in how the state evolves and how batches are charged. We index coordinates from $0$ to $n-1$. The incoming state $x$ and outgoing state $y$ are nonnegative integer vectors. After observing the stage outcome, the planner accepts projects $a_i\in\{0,1\}$, at a cost $c_i$ for each project declined. The root outcome is known. At each later stage one of five equiprobable joint outcomes $w=0,\ldots,4$ is drawn independently of earlier stages. It specifies rejection costs $c_i\in\{1+w,\ldots,9+w\}$, resource weights, and shocks $\xi_i\in\{0,1\}$, and the root shocks are zero. Once the instance seed is fixed, each support point is a fixed vector, which is observed before $a$ is chosen.

Two overlapping crew partitions constrain the accepted projects. Let $b=\min\{8,\max\{2,\lfloor n/2\rfloor\}\}$ and $r=60+2w$. The first partition divides the indices into consecutive blocks of at most $b$, and the second shifts each index of these blocks by $\lfloor b/2\rfloor$ modulo $n$. For each block $G$ in partition $j$, project weights $h_{ji}\in\{1,2,3\}$ give
\[
 \sum_{i\in G}h_{ji}a_i\le
 \max\left\{\max_{i\in G}h_{ji},\left\lfloor r\sum_{i\in G}h_{ji}/100\right\rfloor\right\},
 \qquad
 \sum_i a_i\le\max\{1,\lfloor n(r-8)/100\rfloor\}.
\]
The last inequality applies once to the full crew. Let $x^0$ denote the initial state, and let $v=x-x^0$. The common transition and one-stage cost are
\begin{equation}\label{axbench:projecttransition}
 y=x^0+Av+a+\xi+d,\qquad
 \sum_i c_i(1-a_i)+\phi(x)+k(d)+h_t^\top x .
\end{equation}
Here $d$ is a family-specific adjustment, $\phi$ is the incoming-state batch charge, and $k(d)$ denotes the cost of the controls producing $d$. Each model minimizes the expected undiscounted sum of its stage costs with zero terminal cost. We write $\delta_m(u)=\min_{z\in\mathbb Z}|u-mz|$.

\paragraph{Shift.}
Here $x^0=0$, $(Av)_{2j}=v_{2j+1}$, and $(Av)_{2j+1}=0$. Coordinate zero has upper bound 10, other even coordinates have upper bound 4, and odd coordinates have upper bound 2. The adjustment $d=q e_0$ costs $q$, with $q\in\{0,\ldots,6\}$. The batch cost is $\phi(x)=5\delta_6(\sum_j x_{2j})$. At stage $t=1,\ldots,T$, $(h_t)_i$ is one if $i+\mathrm{seed}+t-1$ is divisible by seven and zero otherwise.

\paragraph{Rotate.}
Here $x^0=0$, $(Av)_i=v_{(i+1)\bmod n}$, and $d=h_t=0$. Each incoming coordinate has upper bound $2(t-1)$ and each outgoing coordinate has upper bound $2t$. Billing compares the first two positions of each triple:
\[
 \phi(x)=5\delta_6\left(\sum_{j=0}^{n/3-1}(x_{3j}-x_{3j+1})\right).
\]

\paragraph{Odd-Cycle.}
Here $x^0=4\mathbf 1$ and the state box is $[0,9]^n$. With $p=(2j+2)\bmod n$, the next pair receives $(Av)_{2j}=-v_{p+1}$ and $(Av)_{2j+1}=v_p-v_{p+1}$. Let $B$ denote the unsigned vertex--edge incidence matrix of a triangle on $0,1,2$ with the path $2,3,\ldots,n-1$ attached. The billing equation is
\[
 Bz+(e^+-e^-+u)e_0=x,\qquad
 z\in\mathbb Z^n,\quad e^+,e^-,u\in\{0,1\}.
\]
It charges $5(e^++e^-)+2u$. The image of $B$ is the set of integer vectors with even coordinate sum, so that minimizing this charge gives $\phi(x)=2\delta_2(\sum_i x_i)$. The repair $u$ affects only the bill. For state balancing, $d_i=2(b_i^+-b_i^-)$, where $b_i^\pm\in\{0,\ldots,8\}$ cost $(1+i\bmod 2)(b_i^++b_i^-)$. Here $h_t=0$.

\paragraph{Shear.}
Here $x^0=8\mathbf 1$, the box is $[0,23]^n$, and $(Av)_i=v_i+2v_{i+1}$ with $v_n=0$. The billing matrix has columns $e_j-e_{j+1}$ for $j<n-1$ and $8e_{n-1}$. The equation $Bz+(e^+-e^-)e_{n-1}=x$ charges $5(e^++e^-)$ with $0\le e^\pm\le4$ integral, which gives $\phi(x)=5\delta_8(\sum_i x_i)$. The balancing controls and their costs are as for Odd-Cycle, with $d_i=8(b_i^+-b_i^-)$ and $h_t=0$. These expressions for $\phi$ eliminate only the current billing variables, whereas future project choices and box constraints remain in the stage problems.

\paragraph{The stage lattices.}
Let $h=T-t+1$ be the number of billing stages remaining. If the current bill observes $Cx\bmod m$, then two states under the same actions and shocks have identical billing residues over those stages if and only if their difference belongs to
\begin{equation}\label{axbench:obskernel}
 \Gamma(h)=\{v\in\mathbb Z^n: CA^jv\equiv0\pmod m,\ 0\le j<h\}.
\end{equation}
The empty observation gives $\Gamma(0)=\mathbb Z^n$. The following coordinates describe the quotient explicitly. In each update, we give the image of $Ax$, and adding the label of $x^0-Ax^0+a+\xi+d$ from \eqref{axbench:projecttransition} gives the label of the outgoing state.
\begin{itemize}
 \item \emph{Shift.} Let $E=\sum_j x_{2j}$ and $O=\sum_j x_{2j+1}$. The label is $(E,O)\bmod6$, and the update is $(O,0)$. With one stage remaining, only $E\bmod6$ is retained. The indices for one stage and for two or more stages are $6$ and $36$, respectively.
 \item \emph{Rotate.} Let $S_j=\sum_{i\equiv j\,(\mathrm{mod}\,3)}x_i$, $P=S_0-S_1$, and $Q=S_1-S_2$. The label is $(P,Q)\bmod6$, with update $(Q,-P-Q)$. Only $P$ is needed at $h=1$. The indices are $6$ and $36$.
 \item \emph{Odd-Cycle.} With $E,O$ as for Shift, the label is $(P,Q)=(E+O,E)\bmod2$, with update $(Q,P+Q)$. Only $P$ is needed at $h=1$. The indices are $2$ and $4$.
 \item \emph{Shear.} Let $S_j=\sum_{i=j}^{n-1}x_i$. The ordered coordinates are $(S_0\bmod8,S_1\bmod4,S_2\bmod2)$, with update $(S_0+2S_1,S_1+2S_2,S_2)$ in the same moduli. We keep the first $\min(h,3)$ coordinates. The indices for $h=1,2$ and $h\ge3$ are $8,32,64$.
\end{itemize}

For each family, the catalog crosses $n\in\{96,384,1536,6144\}$ with $T\in\{64,256,1024,2048\}$. The smaller instances with $(n,T)=(12,24)$ in Section~\ref{sec:exp-equal} use the same stage generator. For a fixed family, dimension, and seed, changing $T$ retains the common prefix of the outcome supports.

\subsection{Group II. Commitments over Time Windows}\label{axbench:windows}

All four families minimize expected cumulative cost, with bounded nonnegative integer states and actions, binary setup decisions where indicated, and zero terminal value. The first stage has a known outcome. At each later stage, Hub-Consolidation, Fulfillment-Network, and Sterile-Trays have five equiprobable joint outcomes. The branch sets the common load level (65, 85, 100, 120, or 145 percent), and the seed fixes item arrivals at that level using a seven-day seasonal pattern. Positive arrivals contain one to five units. Pooled-Capacity has three equiprobable joint outcomes at each stage after the first.

\paragraph{Hub-Consolidation and Fulfillment-Network.}
Each item $i$ belongs to a destination $g(i)$. At a destination $g$, the state records inventories $I_i$ and aggregate freight volumes $(Q_{g0},Q_{g1})$ committed for today and tomorrow. After observing demand $D_i$, the planner produces $p_i$, pays for a setup $z_i\in\{0,1\}$, outsources $o_i$ units, and promises $b_{ij}$ units for dispatch in $j=0,1,2$ days. All promised units leave inventory now. Let $v_i$ denote the unit volume, let $d_i$ denote the maximum permitted delay, and let $B_{gj}=\sum_{i:g(i)=g}v_i b_{ij}$. The constraints are
\begin{align*}
 o_i+\sum_j b_{ij}&=D_i, &
 I'_i&=I_i+p_i-\sum_j b_{ij}, &
 p_i&\le\bar p_i z_i,\\
 Q_{g0}+B_{g0}&\le6 f_g, &
 Q'_{g0}&=Q_{g1}+B_{g1}, &
 Q'_{g1}&=B_{g2}.
\end{align*}
Here $f_g$ is the number of six-unit freight bundles purchased. Bookings vanish for $j>\min\{d_i,T-t\}$. Each destination has four line capacities and one labor capacity. For each resource $r$, the load is $\sum_i(\alpha_{ri}p_i+\beta_{ri}z_i)$ over the items using it. Nonnegative inventory and queue bounds enforce storage and booking limits. Initial commitments are zero. The stage objective charges production, setups, outgoing inventory, outsourcing, booking delays, and freight:
\[
 \sum_i\left(c^p_i p_i+c^z_i z_i+h_i I'_i+c^o_i o_i+
                 c^d_i\sum_j j b_{ij}\right)+\sum_g F_g f_g .
\]
Hub-Consolidation has one destination and $F_0=24$. Fulfillment-Network has $D$ destinations, each stocking the same 64-product catalog, with local inventory, demand seasonality, capacities, and freight fees. Two additional upstream resources couple production and setups across all destinations. The catalog sizes are $K=128,512,2048$ items for Hub and $D=8,16,32$ destinations ($K=64D$ item--destination pairs) for Network. Each width is crossed with $T=52,104,260$. The state has $K+2D$ coordinates, with $D=1$ for Hub.

Freight billing observes today's committed volume $Q_{g0}\bmod6$ on each route. Tomorrow's commitment becomes today's queue at the next stage. Therefore, for $h\ge2$, the lattice is the kernel of the $D$ pairs $(Q_{g0},Q_{g1})\bmod6$. Inventory coordinates are unrestricted in this kernel. Under the bookings above, each pair updates to $(Q_{g1}+B_{g1},B_{g2})\bmod6$. At $h=1$, only $Q_{g0}\bmod6$ is retained, and at $h=0$, no label remains. The respective indices are $36^D$, $6^D$, and $1$.

\paragraph{Sterile-Trays.}
For tray type $i$, $u_i$ jobs are due now, and $v_i$ jobs are due in the next shift. The branch reveals intakes $(A_i,B_i)$ that become available next shift, due then or one shift later. The planner outsources $o_i\le\min\{u_i,\bar o_i\}$ of today's jobs and processes $a_i\le v_i$ of tomorrow's jobs early. Local processing and the queue transitions are
\[
 l_i=u_i-o_i+a_i,\qquad u'_i=v_i-a_i+A_i,\qquad v'_i=B_i .
\]
Items use one of two sterilizer groups $g(i)$ with batch sizes $(m_1,m_2)=(6,8)$ and one of four skills $s(i)$. Integer cycle counts $f_g$, machine overtime $e$, and skill overtime $e_s$ satisfy
\[
 \sum_{i:g(i)=g}l_i\le m_g f_g,\qquad
 2f_1+3f_2\le H+e,\qquad
 \sum_{i:s(i)=s}\tau_i l_i\le H_s+e_s .
\]
The cost is $\sum_i(c_i^o o_i+c_i^a a_i)+30f_1+36f_2+10e+9\sum_s e_s$. Initial queues are drawn in $\{0,1,2\}$, and intake is suppressed whenever its due date would exceed $T$. The terminal queues are empty. The catalog crosses the sizes $K=128,512,1024$ with the horizons $T=84,168,336$, and the state dimension is $2K$.

For group $g$, let $U_g=\sum_{i:g(i)=g}u_i$ and $V_g=\sum_{i:g(i)=g}v_i$. With fixed early-processing and outsourcing decisions, a difference in the current local load is a difference in $U_g$. The two-stage billing label consists of $(U_g,V_g)\bmod m_g$ for each of the two groups. The queue balances update it to $(V_g+\sum_{i:g(i)=g}(A_i-a_i),\sum_{i:g(i)=g} B_i)\bmod m_g$. The lattice is the kernel of these labels. Its index is $(6\cdot8)^2=2{,}304$ for $h\ge2$, falls to $48$ at $h=1$ when only the two $U_g$ residues are retained, and equals $1$ at the terminal state.

\paragraph{Pooled-Capacity.}
For each product $i$, the state is inventory $I_i$ and installed capacity $K_i$. Incoming capacity produces $y_iK_i$ units automatically. New capacity $k_i$ becomes productive next period. After observing demand $D_i$ and resource budgets, the planner chooses shipments $s_i$, lost sales $\ell_i$, discards $v_i$, installations $k_i$, and binary setups $z_i$, with
\[
 I'_i=I_i+y_iK_i-s_i-v_i,\qquad K'_i=K_i+k_i,\qquad
 s_i+\ell_i=D_i,\qquad k_i\le2z_i.
\]
Products belong to two shipping pools. For pool $g$, an integer load count $f_g$ and binary opening $o_g$ satisfy
\[
 0\le6 f_g-\sum_{i:g(i)=g}s_i\le5,\qquad f_g\le\bar f_g o_g.
\]
Capital limits $\sum_i(\alpha_i k_i+\beta_i z_i)$, installation crews limit $\sum_i z_i$, and handling capacity limits $f_1+2f_2+2\sum_i z_i$. Inventories and capacities are bounded, with $K'_i\le6$. Costs comprise holding and carrying costs on $(I',K')$, lost sales, discards, installations, setups, and pool-specific load and opening fees. The catalog crosses 12, 24, and 40 products with 8, 14, and 25 stages and uses heterogeneous yields $y_i\in\{1,3\}$ in both pools. The state dimension is twice the number of products.

Let $U_g=\sum_{i:g(i)=g}(I_i+y_iK_i)$ and $J_g=\sum_{i:g(i)=g}y_iK_i$. Eliminating shipments from the six-unit billing row gives the incoming supply observation $U_g\bmod6$. Its next observation under zero action is $(U_g+J_g)\bmod6$. Therefore, for $h\ge2$, the lattice is the kernel of $(U_g,J_g)\bmod6$ for both pools. Let $R_g=\sum_{i:g(i)=g}(s_i+v_i)$ and $N_g=\sum_{i:g(i)=g}y_i k_i$. The label updates to $(U_g+J_g-R_g+N_g,J_g+N_g)\bmod6$. The quotient has index $6^4=1{,}296$. With $h=1$, only the two $U_g$ remain (index $36$), and the terminal label has index $1$.

\subsection{Group III. Batch Inventory, Workforce, Production, and Revenue}
\label{axops:groupiii}

Each model below minimizes expected total cost over a finite horizon, with zero terminal value. A prime denotes the outgoing state. The current outcome is observed before the local decisions. States and actions are nonnegative integers with the capacities and upper bounds of the instance, and setup indicators are binary. In the catalog and in the instances of the long runs, the first outcome is known, and each later stage has three equiprobable joint outcomes, which are sampled independently across stages. A demand vector is one outcome, even if several products have different demands. The full-size Airline-RM network keeps the supports and probabilities of its source instance.

For a queue \(P_i=(P_{i1},\ldots,P_{iL})\), we write \(\mathcal S(P_i,u_i)=(P_{i2},\ldots,P_{iL},u_i)\), where the first entry is due now. Queue entries count packs, classes, or pallets as specified below, whereas in Lead-Time-Pipeline they count physical units. Each label map \(\phi_t\) below acts on the incoming integer state, and two states lie in the same coset when their difference belongs to \(\Gamma_t=\ker\phi_t\). The terminal map is zero. Outside Lead-Time-Pipeline, the maps are the same in all stages.

\paragraph{Case-Pack and Case-Pack-Truck.}
For product \(i\), the state is shelf stock \(I_i\) and \(L\) periods of in-transit packs \(P_i\). Each pack contains \(p_i\) units. After observing demand \(d_i\), the planner orders \(x_i\) packs, serves \(s_i\) units, loses \(\ell_i\) sales, and discards \(e_i\) units. With \(z\) trucks and \(Q\) packs per truck, we have
\begin{align*}
 I_i' &= I_i+p_iP_{i1}-s_i-e_i, &
 P_i'&=\mathcal S(P_i,x_i), &
 s_i+\ell_i&=d_i,\\
 0\le I_i'&\le\bar I_i, &
 0\le x_i&\le\bar x_i, &
 \sum_i x_i&\le Qz .
\end{align*}
The stage cost is \(Fz+\sum_i(c_ip_ix_i+h_i I_i'+\lambda_i\ell_i+r_ie_i)\). Case-Pack-Truck charges \(F=25\) per truck with \(Q=8\) packs per truck, and Case-Pack sets \(F=0\). The catalog configurations of this family include the truck charge and use \(L=1\). The label is \(\phi(I,P)=(I_i\bmod p_i)_i\), of index \(\prod_i p_i\). For the 10-stage Case-Pack-Truck instance with one item and \(p=4\), Section~\ref{sec:exp-regimes} also evaluates the finer labels \(I+4P\bmod 8\) and \(I+4P\bmod 16\) of the inventory position.

\paragraph{Workforce-Cohort.}
For skill \(i\), \(S_i\) is staff on hand and \(P_i\) counts training classes by time to graduation. Each class has \(a_i\) workers. The planner starts \(x_i\) classes, lays off \(e_i\) workers, and splits the observed workload \(d_i\) between staff service \(s_i\) and outsourcing \(o_i\). The balances are
\[
 S_i'=S_i+a_iP_{i1}-e_i,\qquad P_i'=\mathcal S(P_i,x_i),\qquad
 s_i\le S_i',\qquad s_i+o_i=d_i .
\]
Service uses the end-of-period workforce, and therefore it does not deplete \(S_i'\). The stage cost is \(\sum_i(w_iS_i'+c_ix_i+\lambda_i o_i+r_i e_i)\). The model imposes \(\sum_i x_i\le H\), where \(H\) is the shared trainer limit, alongside \(x_i\le\bar x_i\) and \(S_i'\le\bar S_i\). The variant without a trainer limit drops the first of these constraints. The catalog uses \(L=2\) and the shared limit. The label is \(\phi(S,P)=(S_i\bmod a_i)_i\), of index \(\prod_i a_i\).

\paragraph{Co-Product.}
Each group \(g\) has two inventories \(I_g\in\mathbb Z_+^2\) and two recipes. Column \(r\) of \(Y_g\in\mathbb Z_+^{2\times2}\) is the two-product output of one batch of recipe \(r\). Production \(x_g\) is available for service in the same stage, which gives
\[
 I_g'=I_g+Y_gx_g-s_g-e_g,\qquad s_g+\ell_g=d_g,\qquad
 0\le x_g\le\bar x_g,\qquad 0\le I_g'\le\bar I_g,\qquad
 \sum_{g,r}x_{gr}\le H .
\]
Here, \(e_g\) is disposal and \(\ell_g\) is lost demand. The cost is \(\sum_g(c_g^\top x_g+h_g^\top I_g'+\lambda_g^\top\ell_g+r_g^\top e_g)\). The shared plant limit \(H\) counts batches. The recipe lattice \(Y_g\mathbb Z^2\) has index \(5\) or \(7\). The label is the tuple of cosets \((I_g+Y_g\mathbb Z^2)_g\). Equivalently, the difference \(\Delta I\) of two states belongs to the lattice when \(Y_g^{-1}\Delta I_g\in\mathbb Z^2\) for every group \(g\).

\paragraph{Bike-Sharing.}
The state gives stocks \(B_i\) at stations \(i=1,\ldots,S\) and stock \(D\) at the depot. An outcome specifies rental requests \(d_i\) and returns \(u_i\) together. The planner delivers \(a_i\) racks and picks up \(b_i\) racks of \(k\) bikes each. Stocks immediately after the transfers satisfy
\[
 N_i=B_i+k(a_i-b_i)\in[0,\bar B_i],\qquad
 N_D=D-k\sum_i(a_i-b_i)\ge0,\qquad \sum_i(a_i+b_i)\le H .
\]
The planner then serves \(s_i\) rentals from \(N_i\), docks \(v_i\) returns, and sends the overflow \(o_i\) to the depot. With lost rentals \(\ell_i\) and depot spill \(e\), the service and end-of-period balances are
\[
 s_i+\ell_i=d_i,\quad s_i\le N_i,\quad v_i+o_i=u_i,\quad
 B_i'=N_i-s_i+v_i\le\bar B_i,\quad D'=N_D+\sum_i o_i-e\le\bar D .
\]
The cost is \(c_m\sum_i(a_i+b_i)+c_\ell\sum_i\ell_i+c_o\sum_i o_i+c_e e\).

Bike-Sharing is the one family whose batch balance does not define a lattice of full rank. The transfer columns \(k(\mathbf e_i-\mathbf e_D)\) have zero sum, and their integer span is the kernel of
\[
 \phi(B,D)=\biggl(D+\sum_i B_i,\ (B_i\bmod k)_{i=1}^S\biggr)
 \ \in\ \mathbb Z\oplus(\mathbb Z/k\mathbb Z)^S .
\]
This span has rank \(S\) in \(\mathbb Z^{S+1}\), because a transfer preserves the residues \(B_i\bmod k\) and the total number of bikes, and the total is an integer, not a residue. Rentals, returns, and spill enter after the transfer and change the total. We obtain lattices of full rank in two ways. Where the lattice is computed by Algorithm~\ref{alg:lattice} (the catalog runs), the total is taken modulo \(M=1+\bar D+\sum_i\bar B_i\). On the state box, the total lies in \(\{0,\ldots,M-1\}\), so that two states of the box have the same total if and only if their totals agree modulo \(M\). The lattice \(\{v:\ v_i\equiv 0\ (\mathrm{mod}\ k)\text{ for all stations } i,\ \sum_j v_j\equiv 0\ (\mathrm{mod}\ M)\}\) therefore has full rank and index \(Mk^S\), and on the box, its cosets are the classes of \(\ker\phi\). The index is large, for example \(273\cdot 5^{6}\) with six stations, and only the cosets that appear in a cut receive an indicator (Section~\ref{sec:pooled-splits}). Where the lattice is declared (the long runs and the 45-second runs), we use the coarser label \(((B_i\bmod k)_i,\,D\bmod k)\), whose lattice has full rank and index \(k^{S+1}\). It contains \(\ker\phi\), and therefore it merges classes that the transfers distinguish. Cuts on it remain valid, as they are for any lattice.

\paragraph{Serial-Inventory.}
Product \(i\) has store stock \(S_i\), warehouse stock \(W_i\), a store-bound queue \(A_i\) of packs of \(q_{1i}\) units, and a supplier queue \(P_i\) of pallets of \(q_{2i}\) units, where \(q_{1i}\) divides \(q_{2i}\). The planner ships \(x_i\) packs from the warehouse and orders \(y_i\) pallets. The two-echelon balances are
\begin{align*}
 S_i'&=S_i+q_{1i}A_{i1}-s_i-e_i^S, &
 A_i'&=\mathcal S(A_i,x_i),\\
 W_i'&=W_i+q_{2i}P_{i1}-q_{1i}x_i-e_i^W, &
 P_i'&=\mathcal S(P_i,y_i),\\
 s_i+\ell_i&=d_i, & \sum_i y_i&\le Qz .
\end{align*}
Here, \(z\) counts supplier trucks, and \(e_i^S,e_i^W\) are disposals. The cost is \(Fz+\sum_i(c_iq_{2i}y_i+c_i^x x_i+h_i^S S_i'+h_i^W W_i' +\lambda_i\ell_i+r_i(e_i^S+e_i^W))\). Instances labeled ``truck charge'' in Appendix~\ref{app:results} have \(F>0\), and the other instances of the long runs set \(F=0\). The catalog includes the truck charge and uses store and supplier lead times \(1\) and \(2\). Let \(E_i=S_i+q_{1i}\sum_k A_{ik}+W_i\). The label is \(\phi(S,A,W,P)=(S_i\bmod q_{1i},E_i\bmod q_{2i})_i\), of index \(\prod_i q_{1i}q_{2i}\). Supplier pallets add multiples of \(q_{2i}\) to the full echelon position and drop out of its residue.

\paragraph{Lead-Time-Pipeline.}
Product \(i\) has inventory \(I_i\) and a queue \(P_i\) of physical units with lead time \(L_i\). The planner produces \(x_i\) lots of \(b_i\) units, with setup \(z_i\in\{0,1\}\). For each machine \(j\), \(a_{ji}\) is usage per unit, \(s_{ji}\) is setup time, and \(K_j\) is capacity. With service \(u_i\), lost demand \(\ell_i\), and disposal \(e_i\), the constraints are
\begin{align*}
 I_i'&=I_i+P_{i1}-u_i-e_i, &
 P_i'&=\mathcal S(P_i,b_ix_i), & u_i+\ell_i&=d_i,\\
 x_i&\le\bar x_i z_i, &
 \sum_i(a_{ji}b_ix_i+s_{ji}z_i)&\le K_j &&(\text{each }j).
\end{align*}
Two integer bills charge for outgoing inventory. For bundle size \(m=6\) and volume weights \(v_i\), the bill counts are
\[
 H=\left\lceil\frac{\sum_i I_i'}m\right\rceil,\qquad
 V=\left\lceil\frac{\sum_i v_i I_i'}m\right\rceil .
\]
The stage cost is \(c_H H+c_V V+ \sum_i(c_i b_i x_i+f_i z_i+h_iI_i'+\lambda_i\ell_i+r_i e_i)\). The catalog uses \(L_i=3\) and weights \(1\) and \(3\).

For the lattice, let \(h\) denote the number of stages remaining, including the current one. We define the inventory position \(Y_i=I_i+\sum_{k=1}^{L_i}P_{ik}\) and the unforced stock after \(j\) queue shifts \(X_{ij}=I_i+\sum_{k=1}^{\min(j,L_i)}P_{ik}\). For \(h\ge1\), the label collects
\[
 (Y_i\bmod b_i)_i,\qquad
 \left(
   \sum_i X_{ij}\bmod 6,\ \sum_i v_i X_{ij}\bmod 6
 \right)_{1\le j\le\min(h,\max_i L_i)} .
\]
Each pair is the pullback of the two billing rows through \(j\) transitions, and the sequence stabilizes after three shifts.

\paragraph{Airline-RM.}
For fare class \(j\), the state records cumulative bookings \(B_j\) and cancellations \(C_j\), initially zero. After observing requests \(d_j\), the planner accepts \(b_j\) bookings and records \(c_j\) cancellations. With fare \(f_j\), cancellation numerator \(\alpha_j\), leg incidence \(A_{ij}\), and leg capacity \(K_i\), the model minimizes the stage cost \(\sum_j f_j(c_j-b_j)\) subject to
\begin{align*}
 B_j' &= B_j+b_j, & C_j'&=C_j+c_j,\\
 -10 &\le 20C_j'-\alpha_j B_j'\le10, &
 \sum_j A_{ij}(B_j'-C_j')&\le K_i \qquad(\text{each leg }i).
\end{align*}
The closed strip retains both possible integers at half-integer ties. Net bookings occupy seats, and cancellations refund the full fare. The network instances impose \(b_j\le\min(d_j,100)\), \(c_j\le100\), and instance-specific cumulative bounds. With \(m_j=20/\gcd(20,\alpha_j)\), the label is \(\phi(B,C)=(B_j\bmod m_j)_j\), of index \(\prod_j m_j\).

The six-stage network instance of the full-size comparison (Table~\ref{tab:scope}, Panel~A) has 72 fare classes, 12 legs, and no requests in the first stage. The catalog keeps the incidence, fare, and cancellation profiles of a network of this kind with 36, 72, or 144 fare classes. It subsamples or duplicates classes, perturbs the requests within their booking windows, and rescales the capacity of each leg to its expected net requests.

\ecresetcounters
\section{Compared Methods and Experimental Settings}\label{app:methods}

In this appendix, we describe the methods of Section~\ref{sec:experiments} as they were run. All of them are written in Python and call Gurobi through its Python interface. They read the same stage models, with the same constraints, bounds, and integrality requirements, and differ only in how they approximate the cost-to-go functions. Because SDDiP and SDDmiP are our own implementations of published algorithms, we state what we implemented, which settings we used, and what we did not verify, so that the reader can judge the comparison.

\subsection{Common Elements}\label{app:methods-common}

We fix a stage $t$, an outcome $w$, and the current approximation $\underline V_{t+1}$ of the next stage. With $z$ denoting a copy of the incoming state, we write
\[
 f_w(z,a,y)=c_{tw}(z,a,y)+\underline V_{t+1}(y),\qquad
 F_w(s)=\min_{(a,y)\in X_{tw}(s)} f_w(s,a,y).
\]
Every method alternates a forward pass and a backward pass. The forward pass samples scenario paths and solves the stage problems along them as MIPs. \textnormal{LCD} and \textnormal{LCD-Iso} sample three paths per iteration (\textnormal{LCD} one in its first iteration), and SDDmiP and SDDiP sample one. The backward pass visits the states of these paths from the last stage to the second and generates cuts for every outcome of the stage. After each iteration, a method reports as its lower bound the dual bound that Gurobi returns for the first-stage problem, not the value of the incumbent. The reported bound is therefore valid for any stage gap.

\paragraph{Reference values.}
An optimal value $V^\star$ is computed either by exact dynamic programming over the state space, that is, by evaluating the recursion \eqref{eq:bellman} at every state of every stage, or by solving the extensive form, which has one copy of the stage decisions for every node of the scenario tree, with Gurobi. Instances for which neither is possible have no reference value, and the tables then report the bounds themselves.

\paragraph{Policies.}
After training, the cuts of a method are frozen. Its policy solves, in every stage, the stage problem with these cuts and implements the integer solution. The policies of all methods are evaluated on the same sampled scenario paths. A mean over sampled paths is an estimate of the cost of a policy and is not an upper bound on $V^\star$.

\subsection{\textnormal{LCD}, \textnormal{LCD-Iso}, and the Control Without a Lattice}\label{app:methods-ours}

Table~\ref{tab:methods-ours} lists the three variants of our method that appear in Section~\ref{sec:experiments}. They share the cut definitions of Definition~\ref{def:lattice-cuts} and differ in the classes over which an intercept is computed, in the penalty of a closing cut, and in the formulation of the stage models. Table~\ref{tab:ec-params} in Appendix~\ref{ec:benchmarks} gives their parameters. Figure~\ref{fig:flow} in Section~\ref{sec:cut-subproblems} shows the overall structure of \textnormal{LCD}.

\begin{table}[H]
\centering
\caption{Variants of Our Method in the Computational Study}\label{tab:methods-ours}
\begin{singlespace}\footnotesize
\begin{tabular}{@{}p{0.17\linewidth}p{0.17\linewidth}p{0.27\linewidth}p{0.30\linewidth}@{}}
\toprule
Variant & Intercepts & Closing cuts & Stage models \\
\midrule
\textnormal{LCD} (Algorithm~\ref{alg:lcd}) & one per coset of $\Gamma_t$, computed on demand & directional penalties $\rho^\pm$, shared binary variables \eqref{eq:pooled-split}, postponed on a first visit ($\gamma<1$) & one model per stage, shared by all outcomes; cuts of equal slope are merged \\
\addlinespace[2pt]
\textnormal{LCD-Iso} & one per coset of $\Gamma_t$ & one penalty $\rho$ for all coordinates, formulation \eqref{eq:sign-binary}, at every state with a gap ($\gamma=1$) & one model per stage and outcome; no merging \\
\addlinespace[2pt]
\textnormal{LCD-Iso} without the lattice & one for all states ($\Gamma_t=\Z^n$) & as \textnormal{LCD-Iso} & as \textnormal{LCD-Iso} \\
\bottomrule
\end{tabular}
\end{singlespace}
\end{table}

Both variants solve their MIPs to a relative gap of $10^{-9}$. They use the dual bound of every such program, reduced by a safety margin of $10^{-7}$ in relative and $10^{-6}$ in absolute terms, as the intercept of a cut, so that a cut remains valid when a program is stopped at its gap.

\paragraph{Merging of cuts in \textnormal{LCD}.}
When a new lattice cut has the slope of an existing one, \textnormal{LCD} keeps one row and takes the larger intercept coset by coset. This is exact for identical slopes (Proposition~\ref{prop:phase-hierarchy}). The implementation uses three tolerances. It merges two cuts whenever their slopes agree to within $10^{-8}$ in every coordinate, it reuses a stored intercept or lower bound $b_0$ for a slope that agrees with the slope of the stored value to within $10^{-10}$, and it reuses the value of the LP relaxation for a slope that agrees with the LP slope to within $10^{-12}$. In each of these steps, an intercept that is attached to a slope $\pi$ other than the slope $\pi'$ for which it was computed is first lowered by $\max\{0,\sum_{i=1}^n\max\{(\pi_i-\pi_i')\ell_{t,i},(\pi_i-\pi_i')u_{t,i}\}\}$, which bounds $(\pi-\pi')^\top y$ on the box $[\ell_t,u_t]$, and rounded down. With this rule, every cut remains valid across repeated merges, and the constant $\kappa$ of Proposition~\ref{prop:merge-tolerance} is not needed. Across all 568 runs of \textnormal{LCD} and its control without the lattice, 55,068 of 501,629 merges lowered an intercept, by at most $4.4\times10^{-7}$, and no reported bound exceeded a known $V^\star$ (Appendix~\ref{app:add-validity}).

\paragraph{States visited by \textnormal{LCD}.}
In instances with at most 12 state variables, \textnormal{LCD} also visits states that do not lie on the sampled paths. When a lattice cut is not tight at its state, the outgoing states in the solutions of the programs \eqref{eq:coset-intercept} are recorded. In an iteration that follows one without an improvement of the bound, up to two further paths start at recorded states and are continued by sampling, and their states are added to the backward pass. These states only add cuts and do not affect validity.

\subsection{SDDmiP}\label{app:methods-sddmip}

SDDmiP is the variant of \citet{bansal2026} that alternates cuts from the LP relaxation with normalized rectified linear cuts \citep{bansal2026b}. Our implementation keeps the original integer state, with one epigraph variable $\theta_w$ per outcome of the next stage. At a visited state $\hat s$ and an outcome $w$, with $\hat\theta_w$ the current value of the approximation at $\hat s$, it proceeds as follows.

\emph{LP cut.} The LP relaxation of the stage problem is solved with the incoming state fixed at $\hat s$. If its value $v^{\mathrm{LP}}$ exceeds $\hat\theta_w$, the cut $\theta_w\ge v^{\mathrm{LP}}+\pi^\top(y-\hat s)$ is added, where $\pi$ is the dual vector of the fixing constraints, and no further cut is generated for this outcome in this iteration.

\emph{Rectified linear cut.} Otherwise, the stage problem is solved as a MIP at $\hat s$, with value $\hat f$. If its dual bound exceeds $\hat\theta_w$ by more than $0.1$, a rectified linear cut is generated. Let $d^+(z)=(z-\hat s)_+$ and $d^-(z)=(\hat s-z)_+$, taken coordinate by coordinate, and for multipliers $\lambda=(\lambda^+,\lambda^-)$ and $\lambda_0\ge 0$ let
\[
 H_w(\lambda,\lambda_0)=\min_{z,\ (a,y)\in X_{tw}(z)}
 \bigl\{\lambda_0\,f_w(z,a,y)+(\lambda^+)^\top d^+(z)+(\lambda^-)^\top d^-(z)\bigr\},
\]
where the copy $z$ of the incoming state is continuous within the bounds of the state and all other variables keep their integrality. Every lower bound $\underline H$ on $H_w(\lambda,\lambda_0)$ gives the valid cut
\[
 \lambda_0\,\theta_w+(\lambda^+)^\top d^+(y)+(\lambda^-)^\top d^-(y)\ \ge\ \underline H .
\]
The multipliers are chosen to maximize $H_w(\lambda,\lambda_0)-\lambda_0\hat\theta_w$ subject to the normalization
\[
 h^\top\lambda+\bigl(\hat f-\hat\theta_w+0.1\bigr)\lambda_0\le 1,\qquad 0\le\lambda_0\le 1,
\]
and to finite bounds on $\lambda$. The vector $h$ is the CoreMid choice of \citet{bansal2026b}: for a state variable with range $[0,u_i]$, $(h_i^+,h_i^-)$ equals $(u_i/2,0)$ if $\hat s_i=0$, $(0,u_i/2)$ if $\hat s_i=u_i$, and $(r_i,r_i)$ with $r_i=\hat s_i(1-\hat s_i/u_i)/2$ otherwise. The maximization uses a level bundle method with level parameter $0.5$, tolerance $0.01$, and at most 200 iterations (300 for Airline-RM), each of which solves one MIP.

\emph{Cuts in the stage models.} In the model of the preceding stage, $d^+(y)$ and $d^-(y)$ are represented with one binary variable per state variable, as in \eqref{eq:sign-binary}. Cuts with the same anchor share these variables. All rectified linear cuts are kept, whereas an LP cut stays in the model as long as it is the highest LP cut at one of the visited states.

\emph{Settings.} The stage gap is $10^{-4}$. The original stopping rule, which ends a run when the bound has improved by less than 1\% over ten iterations, is disabled, so that a run uses its whole time limit.

\subsection{SDDiP}\label{app:methods-sddip}

SDDiP is stochastic dual dynamic integer programming \citep{zou2019} with strengthened Benders cuts and integer optimality cuts. Each state variable $s_i\in\{0,\ldots,u_i\}$ is replaced by binary variables $b_{ij}$ with $s_i=\sum_j 2^j b_{ij}$, and the bound $s_i\le u_i$ is kept. Let $b$ be the vector of all these variables and $\hat b$ its value at a visited state. Outgoing states are binary in every stage problem. In the relaxations that generate cuts, the copy $v$ of the incoming binary state is continuous in $[0,1]$.

\emph{Strengthened Benders cut.} The LP relaxation of the stage problem with $v$ fixed at $\hat b$ gives a dual vector $q_w$ of the fixing constraints. The intercept is then recomputed with the integrality requirements of the stage in place. Specifically, $A_w(q_w)$ is the minimum of the stage objective minus $q_w^\top v$ over the stage constraints with $v\in[0,1]$ free. The cut is
\[
 \theta\ \ge\ \sum_{w\in\Omega_t}p_{tw}\bigl(\underline A_w(q_w)+q_w^\top b\bigr),
\]
where $\underline A_w$ is the dual bound of the MIP that defines $A_w$.

\emph{Integer optimality cut.} Let $\ell_t$ denote a lower bound on the cost-to-go function of stage $t$ that holds at all states, and let $h$ denote the larger of $\ell_t$ and the probability-weighted sum of the dual bounds of the stage problems at $\hat b$. The cut is
\[
 \theta\ \ge\ h-(h-\ell_t)\Bigl[\sum_{j:\hat b_j=1}(1-b_j)+\sum_{j:\hat b_j=0}b_j\Bigr].
\]
It equals $h$ at $\hat b$ and is at most $\ell_t$ at every other binary state.

One cut of each kind is added per visited state and iteration, with a single epigraph variable per stage. Lagrangian cuts are not used. The stage gap is $10^{-4}$.

\subsection{Settings by Experiment}\label{app:methods-settings}

Table~\ref{tab:methods-design} lists the experiments of Section~\ref{sec:experiments} with their settings. Every run uses one thread, and its time limit includes reading the instance, computing the lattices where they are computed, and building the models.

\begin{table}[!htbp]
\centering
\caption{Experiments and Their Settings}\label{tab:methods-design}
\begin{singlespace}\footnotesize
\setlength{\tabcolsep}{3pt}
\begin{tabular}{@{}p{0.20\linewidth}rp{0.17\linewidth}p{0.13\linewidth}p{0.10\linewidth}p{0.15\linewidth}p{0.12\linewidth}@{}}
\toprule
Experiment & Inst. & Methods & Time limit & Lattices & Stopping & Policy paths \\
\midrule
One item, long runs (Section~\ref{sec:exp-long}) & 32 & \textnormal{LCD-Iso} with and without the lattice, SDDmiP & 30 min, 45 min, 2 h for $T=10,24,48$ & declared & ours at $V^\star$ or the limit; SDDmiP at the limit & 200 \\
\addlinespace[2pt]
Other long runs (Table~\ref{tab:scope}, Panel~A) & 30 & \textnormal{LCD-Iso}, SDDmiP & 30--45 min; 15 min for $T=6$ & declared & as above & not reported \\
\addlinespace[2pt]
300-second runs (Section~\ref{sec:exp-equal}) & 12 & \textnormal{LCD}, \textnormal{LCD-Iso}, SDDmiP, SDDiP & 300 s & computed & time limit only & 100 or 32 \\
\addlinespace[2pt]
45-second runs (Table~\ref{tab:scope}, Panel~B) & 48 & as above & 45 s & declared & time limit only & 32 \\
\addlinespace[2pt]
Catalog (Table~\ref{tab:scope}, Panel~B) & 163 & as above & 20 s or 25 s & computed & time limit only & 16 \\
\bottomrule
\end{tabular}
\end{singlespace}
\figurenotes{``Declared'' lattices are specified from the batch structure of the model (Appendix~\ref{app:benchmarks}), whereas ``computed'' lattices are obtained by Algorithm~\ref{alg:lattice} within the time limit. The stage gap is $10^{-9}$ for \textnormal{LCD} and \textnormal{LCD-Iso} and $10^{-4}$ for SDDmiP and SDDiP in all experiments. Policy evaluation has a separate time limit of 90, 20, and 12 seconds in the last three experiments, and a policy whose evaluation does not finish is reported as missing.}
\end{table}

\paragraph{Comparing bounds and policies.}
We count two bounds $L$ and $L'$ as tied if $|L-L'|\le 10^{-4}\max\{1,|L|,|L'|\}$. A run that ends without a bound loses against a run with a bound. A run that completes no iteration and returns the trivial bound zero is counted as a run without a bound. For two policies with costs $C(\omega^r)$ and $C'(\omega^r)$ on the common paths $\omega^1,\ldots,\omega^N$, we report the mean of the differences $d_r=C(\omega^r)-C'(\omega^r)$ with the interval $\bar d\pm t_{N-1,0.975}\,s_d/\sqrt N$, where $s_d$ is the sample standard deviation of the $d_r$.

\subsection{Relation to the Published Implementations}\label{app:methods-caveats}

The comparison methods are reimplementations. SDDmiP follows \citet{bansal2026,bansal2026b}. We consulted the authors' Julia package for algorithmic details, and no Julia code is executed. SDDiP follows Algorithm~2 of \citet{zou2019}. On three small test cases with known optimal values, namely, two two-stage server location cases and a three-stage lot-sizing case, our implementations of SDDiP and SDDmiP attain the optimal value to within $2\times10^{-8}$. We have not reproduced the computational results of these papers with our implementations, and we have not compared the speed of our implementations with that of the authors' codes. Therefore, the two methods should be read as our implementations of the published algorithms and not as the authors' codes. The following three settings bear on the comparison.

First, the comparison methods solve their stage problems to a gap of $10^{-4}$, whereas our methods solve them to a gap of $10^{-9}$. Table~\ref{tab:d-sens} in Appendix~\ref{app:results} repeats eight long runs with the gaps exchanged. The gaps of SDDmiP change by a factor of at most 2.5, and those of \textnormal{LCD-Iso} at a stage gap of $10^{-4}$ remain below those of SDDmiP at $10^{-9}$ on all eight instances.

Second, SDDmiP generates a rectified linear cut only if the stage value exceeds the current approximation by more than $0.1$ in absolute terms. This threshold limits the accuracy that SDDmiP can reach in Figure~\ref{fig:long}. The gaps of SDDmiP in Section~\ref{sec:exp-long} correspond to absolute differences of 0.10 to 0.25 at 10 stages, 0.25 to 4.4 at 24 stages, and 1.2 to 21 at 48 stages. At 10 stages, they are of the size of the threshold, so that Figure~\ref{fig:long} does not show what SDDmiP attains with a smaller one. In Appendix~\ref{app:add-sddmip}, we vary the threshold in runs of 45 to 90 seconds and repeat 28 long runs at threshold $10^{-3}$ (Table~\ref{tab:x-long-sddmip}). A smaller threshold lowers the gaps of SDDmiP on the Workforce-Cohort instances, to the level of our methods at 24 stages, and does not close its gaps on the 24- and 48-stage Case-Pack and 48-stage Co-Product instances.

Third, the stopping rule of SDDmiP is disabled, so that it uses its whole time limit. In the long runs, it continues after our methods have stopped at the reference value.

\ecresetcounters
\section{Detailed Results}\label{app:results}

In this appendix, we report the results behind Section~\ref{sec:experiments} instance by instance. The relative gap of a certified bound $L$ is $(V^\star-L)/|V^\star|$. Where no optimal value is known, the tables give the bounds themselves. A dash marks a run without a bound or a quantity that was not computed.

\subsection{Exact Representation Loss}\label{app:results-exact}

Table~\ref{tab:d-exact} contains the values plotted in Figure~\ref{fig:root-gap-decomp}. They are computed by dynamic programming over the full state space and involve no sampling and no cut generation.

\begin{table}[H]
\centering
\caption{Representation Loss $V^\star-W(\rho)$ With and Without the Lattice on Two 10-Stage Instances}\label{tab:d-exact}
\begin{singlespace}\footnotesize
\begin{tabular}{@{}rrrrr@{}}
\toprule
 & \multicolumn{2}{c}{Co-Product ($V^\star=168.2649$)} & \multicolumn{2}{c}{Workforce-Cohort ($V^\star=468.6667$)} \\
\cmidrule(lr){2-3}\cmidrule(l){4-5}
Penalty cap $\rho$ & No lattice & Lattice & No lattice & Lattice \\
\midrule
0 & 3.7588 & 0.0050 & 1.4167 & 0.0000 \\
0.05 & 3.5968 & 0.0023 & 1.3529 & 0.0000 \\
0.1 & 3.4389 & 0.0000 & 1.2891 & 0.0000 \\
0.15 & 3.2863 & 0.0000 & 1.2252 & 0.0000 \\
0.2 & 3.1413 & 0.0000 & 1.1635 & 0.0000 \\
0.25 & 2.9981 & 0.0000 & 1.1037 & 0.0000 \\
0.5 & 2.3378 & 0.0000 & 0.8469 & 0.0000 \\
1 & 1.2812 & 0.0000 & 0.5556 & 0.0000 \\
2 & 0.1224 & 0.0000 & 0.0000 & 0.0000 \\
3.5 & 0.0005 & 0.0000 & 0.0000 & 0.0000 \\
4.5 & 0.0000 & 0.0000 & 0.0000 & 0.0000 \\
 
\bottomrule
\end{tabular}
\end{singlespace}
\end{table}

\subsection{Long Runs}\label{app:results-long}

Table~\ref{tab:d-long-one} lists the 32 one-item instances of Section~\ref{sec:exp-long} and Figure~\ref{fig:long}, and Table~\ref{tab:d-long-other} lists the 30 further instances that Panel~A of Table~\ref{tab:scope} summarizes. Both tables include the control without the lattice, that is, \textnormal{LCD-Iso} with $\Gamma_t=\Z^n$. In these runs, \textnormal{LCD-Iso} and the control stop when their bound reaches the reference value to within the tolerance stated in Section~\ref{sec:exp-setup}, whereas SDDmiP runs to the time limit (Appendix~\ref{app:methods-settings}). In Table~\ref{tab:d-long-one}, gaps of \textnormal{LCD-Iso} of the order of $10^{-6}$ or below are the accumulated safety margins.

\begin{table}[H]
\centering
\caption{Relative Gaps and Effort in the Long Runs on 32 One-Item Instances}\label{tab:d-long-one}
\begin{singlespace}\scriptsize
\setlength{\tabcolsep}{2.5pt}
\begin{tabular*}{\textwidth}{@{\extracolsep{\fill}}lrrrrccccrrr@{}}
\toprule
 & & & & & Loss & \multicolumn{3}{c}{Relative gap} & \multicolumn{2}{c}{\textnormal{LCD-Iso}} & No lattice \\
\cmidrule(lr){7-9}\cmidrule(lr){10-11}\cmidrule(l){12-12}
Family & $T$ & Seed & Cosets & $V^\star$ & at $\rho=0$ & \textnormal{LCD-Iso} & No lattice & SDDmiP & Iter. & Closing & Closing \\
\midrule
Case-Pack & 10 & 0 & 4 & 702.70 & 0 & $5.0\times10^{-7}$ & $2.5\times10^{-6}$ & $3.5\times10^{-4}$ & 4 & 0 & 710 \\
Case-Pack & 10 & 1 & 6 & 609.82 & 0 & $5.2\times10^{-7}$ & $1.4\times10^{-4}$ & $2.4\times10^{-4}$ & 3 & 0 & 665 \\
Workforce-Cohort & 10 & 0 & 6 & 468.67 & 0 & $4.2\times10^{-7}$ & $4.5\times10^{-7}$ & $2.8\times10^{-4}$ & 3 & 0 & 35 \\
Workforce-Cohort & 10 & 1 & 6 & 347.89 & 0 & $3.7\times10^{-7}$ & $3.8\times10^{-7}$ & $2.9\times10^{-4}$ & 1 & 0 & 0 \\
Co-Product & 10 & 0 & 5 & 168.26 & $2.9\times10^{-5}$ & $5.5\times10^{-7}$ & $8.4\times10^{-5}$ & $1.0\times10^{-3}$ & 27 & 4 & 816 \\
Co-Product & 10 & 1 & 5 & 120.28 & 0 & $6.0\times10^{-7}$ & $6.6\times10^{-7}$ & $1.0\times10^{-3}$ & 6 & 2 & 632 \\
Bike-Sharing & 10 & 0 & 36 & 50.46 & 0 & $2.3\times10^{-3}$ & $8.8\times10^{-2}$ & $3.6\times10^{-3}$ & 19 & 0 & 838 \\
Bike-Sharing & 10 & 1 & 25 & 20.50 & 0 & $4.2\times10^{-3}$ & $6.4\times10^{-2}$ & $9.9\times10^{-3}$ & 29 & 0 & 928 \\
\addlinespace[2pt]
Case-Pack & 24 & 0 & 4 & 1,613.98 & 0 & $1.2\times10^{-6}$ & $3.3\times10^{-3}$ & $1.6\times10^{-3}$ & 17 & 0 & 1,602 \\
Case-Pack & 24 & 1 & 6 & 1,485.18 & 0 & $1.7\times10^{-6}$ & $8.1\times10^{-3}$ & $2.6\times10^{-3}$ & 8 & 0 & 1,606 \\
Case-Pack & 24 & 2 & 8 & 978.86 & 0 & $1.2\times10^{-6}$ & $1.1\times10^{-2}$ & $4.5\times10^{-3}$ & 14 & 0 & 1,542 \\
Case-Pack & 24 & 3 & 12 & 862.35 & 0 & $1.2\times10^{-6}$ & $1.7\times10^{-2}$ & $4.2\times10^{-3}$ & 10 & 0 & 1,384 \\
Workforce-Cohort & 24 & 0 & 6 & 1,053.33 & 0 & $1.1\times10^{-6}$ & $1.2\times10^{-6}$ & $2.4\times10^{-4}$ & 4 & 0 & 100 \\
Workforce-Cohort & 24 & 1 & 6 & 795.89 & 0 & $1.2\times10^{-6}$ & $1.1\times10^{-6}$ & $9.4\times10^{-4}$ & 2 & 0 & 53 \\
Workforce-Cohort & 24 & 2 & 4 & 801.07 & 0 & $1.1\times10^{-6}$ & $1.2\times10^{-6}$ & $3.3\times10^{-4}$ & 4 & 0 & 43 \\
Workforce-Cohort & 24 & 3 & 6 & 1,433.69 & 0 & $1.1\times10^{-6}$ & $1.1\times10^{-6}$ & $4.6\times10^{-4}$ & 4 & 0 & 270 \\
Co-Product & 24 & 0 & 5 & 412.22 & $2.6\times10^{-5}$ & $1.3\times10^{-6}$ & $2.1\times10^{-3}$ & $9.7\times10^{-4}$ & 55 & 5 & 1,791 \\
Co-Product & 24 & 1 & 5 & 306.79 & 0 & $1.6\times10^{-6}$ & $2.4\times10^{-3}$ & $1.6\times10^{-3}$ & 14 & 3 & 1,505 \\
Co-Product & 24 & 2 & 5 & 298.14 & $1.5\times10^{-5}$ & $1.9\times10^{-6}$ & $2.4\times10^{-3}$ & $1.3\times10^{-3}$ & 29 & 15 & 1,522 \\
Co-Product & 24 & 3 & 5 & 305.43 & $3.2\times10^{-5}$ & $9.1\times10^{-6}$ & $2.9\times10^{-3}$ & $1.7\times10^{-3}$ & 73 & 20 & 1,540 \\
Bike-Sharing & 24 & 0 & 36 & 219.77 & 0 & $8.4\times10^{-4}$ & $9.4\times10^{-2}$ & $3.9\times10^{-3}$ & 9 & 0 & 1,472 \\
Bike-Sharing & 24 & 1 & 25 & 95.16 & 0 & $1.6\times10^{-3}$ & $6.5\times10^{-2}$ & $1.0\times10^{-2}$ & 23 & 0 & 1,800 \\
Bike-Sharing & 24 & 2 & 16 & 4.31 & 0 & $5.5\times10^{-2}$ & $9.7\times10^{-2}$ & $1.0\times10^{-1}$ & 6 & 0 & 2,542 \\
Bike-Sharing & 24 & 3 & 16 & 28.32 & 0 & $1.9\times10^{-2}$ & $6.0\times10^{-2}$ & $2.7\times10^{-2}$ & 8 & 0 & 1,529 \\
\addlinespace[2pt]
Case-Pack & 48 & 0 & 4 & 3,247.65 & 0 & $2.4\times10^{-6}$ & $4.7\times10^{-3}$ & $4.5\times10^{-3}$ & 15 & 0 & 4,222 \\
Case-Pack & 48 & 1 & 6 & 3,022.50 & 0 & $4.8\times10^{-6}$ & $8.9\times10^{-3}$ & $6.9\times10^{-3}$ & 26 & 0 & 4,038 \\
Workforce-Cohort & 48 & 0 & 6 & 2,113.11 & 0 & $2.4\times10^{-6}$ & $2.4\times10^{-6}$ & $5.8\times10^{-4}$ & 9 & 0 & 527 \\
Workforce-Cohort & 48 & 1 & 6 & 1,602.56 & 0 & $2.4\times10^{-6}$ & $3.4\times10^{-6}$ & $7.7\times10^{-4}$ & 4 & 0 & 164 \\
Co-Product & 48 & 0 & 5 & 828.74 & $1.7\times10^{-5}$ & $4.1\times10^{-6}$ & $3.8\times10^{-3}$ & $1.6\times10^{-3}$ & 43 & 23 & 4,630 \\
Co-Product & 48 & 1 & 5 & 639.82 & $3.5\times10^{-5}$ & $9.0\times10^{-6}$ & $4.4\times10^{-3}$ & $1.9\times10^{-3}$ & 81 & 14 & 3,855 \\
Bike-Sharing & 48 & 0 & 36 & 412.10 & 0 & $3.0\times10^{-4}$ & $9.8\times10^{-2}$ & $1.1\times10^{-2}$ & 16 & 0 & 3,899 \\
Bike-Sharing & 48 & 1 & 25 & 267.25 & 0 & $3.2\times10^{-4}$ & $3.9\times10^{-2}$ & $1.1\times10^{-2}$ & 26 & 0 & 4,249 \\
 
\bottomrule
\end{tabular*}
\end{singlespace}
\figurenotes{The time limits are 30 minutes for $T=10$, 45 minutes for $T=24$, and 2 hours for $T=48$. ``Loss at $\rho=0$'' is the relative representation loss $(V^\star-W_\Gamma(0))/V^\star$ of lattice cuts, computed by dynamic programming. ``No lattice'' is \textnormal{LCD-Iso} with $\Gamma_t=\Z^n$. ``Iter.'' is the number of iterations, and ``Closing'' is the number of closing cuts generated.}
\end{table}

\begin{table}[H]
\centering
\caption{Results of \textnormal{LCD-Iso}, the Control Without the Lattice, and SDDmiP in the Long Runs on the Other Instance Sets}\label{tab:d-long-other}
\begin{singlespace}\scriptsize
\setlength{\tabcolsep}{3.2pt}
\begin{tabular}{@{}lrrrrrrrr@{}}
\toprule
Instance & $n$ & $T$ & Cosets & $V^\star$ & Loss at $\rho=0$ & \textnormal{LCD-Iso} & No lattice & SDDmiP \\
\midrule
\multicolumn{9}{@{}l}{\emph{Multiple items}}\\
Case-Pack, lead time 3 & 4 & 10 & 8 & --- & --- & 710.03 & 710.03 & 709.87 \\
Case-Pack, lead time 3 & 4 & 24 & 8 & --- & --- & 1,481.14 & 1,465.75 & 1,464.53 \\
Case-Pack, lead time 3 & 4 & 24 & 8 & 557.38 & 0 & $2.0\times10^{-6}$ & $4.7\times10^{-3}$ & $1.3\times10^{-3}$ \\
Co-Product & 4 & 10 & 35 & --- & --- & 1,840.16 & 1,840.16 & 1,839.95 \\
Co-Product & 4 & 24 & 35 & --- & --- & 4,442.89 & 4,443.44 & 4,441.51 \\
Co-Product & 4 & 24 & 35 & --- & --- & 3,422.00 & 3,422.01 & 3,420.48 \\
Bike-Sharing & 3 & 10 & 64 & --- & --- & 197.85 & 190.87 & 189.58 \\
Bike-Sharing & 3 & 24 & 64 & --- & --- & 489.87 & 483.37 & 482.39 \\
Serial-Inventory & 5 & 10 & 27 & 90.82 & 0 & $5.9\times10^{-7}$ & $6.5\times10^{-7}$ & $7.8\times10^{-4}$ \\
Serial-Inventory & 5 & 24 & 27 & 270.02 & 0 & $1.2\times10^{-2}$ & $2.9\times10^{-2}$ & $3.2\times10^{-2}$ \\
Serial-Inventory & 5 & 24 & 27 & 344.76 & 0 & $8.5\times10^{-2}$ & $6.1\times10^{-2}$ & $5.7\times10^{-2}$ \\
\multicolumn{9}{@{}l}{\emph{Truck fixed charge}}\\
Case-Pack-Truck & 2 & 10 & 4 & 887.74 & $1.6\times10^{-2}$ & $1.7\times10^{-3}$ & $6.4\times10^{-3}$ & $2.9\times10^{-4}$ \\
Case-Pack-Truck & 2 & 10 & 6 & 773.56 & $1.4\times10^{-2}$ & $1.5\times10^{-3}$ & $6.2\times10^{-3}$ & $2.2\times10^{-4}$ \\
Case-Pack-Truck & 2 & 24 & 4 & 2,058.71 & $2.4\times10^{-2}$ & $1.8\times10^{-2}$ & $2.3\times10^{-2}$ & $1.7\times10^{-2}$ \\
Case-Pack-Truck & 2 & 24 & 6 & 1,900.27 & $1.9\times10^{-2}$ & $1.2\times10^{-2}$ & $2.6\times10^{-2}$ & $1.3\times10^{-2}$ \\
\multicolumn{9}{@{}l}{\emph{Full size, $T=6$}}\\
Case-Pack-Truck & 8 & 6 & 3,072 & 1,014.83 & --- & $4.8\times10^{-2}$ & $5.1\times10^{-2}$ & $8.4\times10^{-2}$ \\
Case-Pack & 8 & 6 & 3,072 & 893.62 & --- & $4.3\times10^{-2}$ & $3.5\times10^{-2}$ & $6.3\times10^{-2}$ \\
Case-Pack & 16 & 6 & 7,962,624 & 1,217.98 & --- & $1.0\times10^{-1}$ & $9.5\times10^{-2}$ & $1.1\times10^{-1}$ \\
Workforce-Cohort & 6 & 6 & 600 & 3,879.89 & --- & $2.8\times10^{-7}$ & $2.5\times10^{-7}$ & $1.8\times10^{-4}$ \\
Workforce-Cohort, no trainer limit & 6 & 6 & 600 & 3,826.67 & --- & $1.9\times10^{-7}$ & $1.9\times10^{-7}$ & $3.9\times10^{-6}$ \\
Workforce-Cohort & 12 & 6 & 960,000 & 7,112.59 & --- & $8.7\times10^{-4}$ & $2.3\times10^{-7}$ & $2.0\times10^{-3}$ \\
Workforce-Cohort, no trainer limit & 12 & 6 & 960,000 & 7,027.81 & --- & $2.0\times10^{-7}$ & $2.0\times10^{-7}$ & $7.4\times10^{-4}$ \\
Co-Product & 6 & 6 & 125 & 1,337.78 & --- & $2.9\times10^{-7}$ & $5.9\times10^{-6}$ & $3.3\times10^{-4}$ \\
Co-Product & 12 & 6 & 21,875 & 2,564.28 & --- & $6.6\times10^{-3}$ & $6.7\times10^{-3}$ & $4.4\times10^{-3}$ \\
Bike-Sharing & 4 & 6 & 256 & 236.00 & --- & $2.3\times10^{-2}$ & $2.3\times10^{-2}$ & $1.7\times10^{-2}$ \\
Bike-Sharing & 7 & 6 & 279,936 & 406.57 & --- & $6.3\times10^{-2}$ & $5.8\times10^{-2}$ & $5.7\times10^{-2}$ \\
Serial-Inventory, truck charge & 12 & 6 & 5,832 & 364.98 & --- & $1.3\times10^{-1}$ & $6.8\times10^{-2}$ & $1.2\times10^{-1}$ \\
Serial-Inventory & 12 & 6 & 5,832 & 305.88 & --- & $1.1\times10^{-1}$ & $5.6\times10^{-2}$ & $8.2\times10^{-2}$ \\
Serial-Inventory & 24 & 6 & $7.2\times10^{7}$ & 766.56 & --- & $1.0\times10^{-1}$ & $8.6\times10^{-2}$ & $7.4\times10^{-2}$ \\
Airline-RM & 144 & 6 & $1.7\times10^{55}$ & -207,524.84 & --- & $2.2\times10^{-2}$ & $2.2\times10^{-2}$ & $3.0\times10^{-2}$ \\
 
\bottomrule
\end{tabular}
\end{singlespace}
\figurenotes{The last three columns report relative gaps where $V^\star$ is known and certified bounds otherwise. ``No lattice'' is \textnormal{LCD-Iso} with $\Gamma_t=\Z^n$. ``Loss at $\rho=0$'' is as in Table~\ref{tab:d-long-one} and is available where $V^\star$ comes from dynamic programming. The time limits are 30 minutes for $T=10$, 45 minutes for $T=24$, and 15 minutes for the full-size instances with $T=6$, whose reference values come from the extensive form.}
\end{table}

Counted with the tie rule of Appendix~\ref{app:methods-settings}, the bound of \textnormal{LCD-Iso} in Table~\ref{tab:d-long-other} is higher than that of the control without the lattice on 5 of the 11 instances with multiple items, tied on 4, and lower on 2. It is higher on all four instances with a truck fixed charge. On the 15 full-size instances, it is higher on 3, tied on 4, and lower on 8. The losses include the Workforce-Cohort instance with 12 state variables and 960,000 cosets per stage, on which the control reaches the optimum to within the safety margin and \textnormal{LCD-Iso} ends at a gap of $8.7\times10^{-4}$.

Table~\ref{tab:d-sens} addresses the different stage gaps of the two methods. It repeats the eight 24-stage runs with seeds 0 and 1 with the gaps exchanged, that is, with \textnormal{LCD-Iso} at $10^{-4}$ and SDDmiP at $10^{-9}$. We observe that the gaps of SDDmiP change by a factor between 0.8 and 2.5. With a stage gap of $10^{-4}$, \textnormal{LCD-Iso} no longer reaches the optimum to within the safety margin, and its gaps on the six instances other than Bike-Sharing lie between $1\times10^{-6}$ and $2.3\times10^{-5}$. However, they remain below those of SDDmiP at either stage gap on all eight instances. The threshold of $0.1$ below which SDDmiP generates no rectified linear cut (Appendix~\ref{app:methods-caveats}) is the same in all of these runs.

\begin{table}[H]
\centering
\caption{Relative Gaps With the Stage Gaps of the Two Methods Exchanged (24 Stages, 45 Minutes)}\label{tab:d-sens}
\begin{singlespace}\footnotesize
\begin{tabular}{@{}lrrccccr@{}}
\toprule
 & & & \multicolumn{2}{c}{\textnormal{LCD-Iso}} & \multicolumn{3}{c}{SDDmiP} \\
\cmidrule(lr){4-5}\cmidrule(l){6-8}
Family & Seed & $V^\star$ & gap $10^{-9}$ & gap $10^{-4}$ & gap $10^{-4}$ & gap $10^{-9}$ & Iter.\ at $10^{-9}$ \\
\midrule
Case-Pack & 0 & 1,613.98 & $1.2\times10^{-6}$ & $1.2\times10^{-5}$ & $1.6\times10^{-3}$ & $1.2\times10^{-3}$ & 168 \\
Case-Pack & 1 & 1,485.18 & $1.7\times10^{-6}$ & $3.1\times10^{-6}$ & $2.6\times10^{-3}$ & $2.6\times10^{-3}$ & 152 \\
Workforce-Cohort & 0 & 1,053.33 & $1.1\times10^{-6}$ & $1.1\times10^{-6}$ & $2.4\times10^{-4}$ & $1.6\times10^{-4}$ & 5,771 \\
Workforce-Cohort & 1 & 795.89 & $1.2\times10^{-6}$ & $2.3\times10^{-5}$ & $9.4\times10^{-4}$ & $3.8\times10^{-4}$ & 4,673 \\
Co-Product & 0 & 412.22 & $1.3\times10^{-6}$ & $2.9\times10^{-6}$ & $9.7\times10^{-4}$ & $6.2\times10^{-4}$ & 2,022 \\
Co-Product & 1 & 306.79 & $1.6\times10^{-6}$ & $2.3\times10^{-5}$ & $1.6\times10^{-3}$ & $1.9\times10^{-3}$ & 1,670 \\
Bike-Sharing & 0 & 219.77 & $8.4\times10^{-4}$ & $6.0\times10^{-4}$ & $3.9\times10^{-3}$ & $3.3\times10^{-3}$ & 291 \\
Bike-Sharing & 1 & 95.16 & $1.6\times10^{-3}$ & $2.5\times10^{-3}$ & $1.0\times10^{-2}$ & $1.0\times10^{-2}$ & 245 \\
 
\bottomrule
\end{tabular}
\end{singlespace}
\figurenote{The columns ``gap $10^{-9}$'' for \textnormal{LCD-Iso} and ``gap $10^{-4}$'' for SDDmiP are the settings of Table~\ref{tab:d-long-one}.}
\end{table}

\clearpage
\subsection{Runs With Time Limits of 300 and 45 Seconds}\label{app:results-short}

Table~\ref{tab:c1} in the main text gives the bounds and policy costs of the 300-second runs, and Table~\ref{tab:d-300} adds the effort behind them. Table~\ref{tab:d-45} gives the bounds and policy costs of the 48 instances with a limit of 45 seconds, which Panel~B of Table~\ref{tab:scope} summarizes. These 48 instances are the eight one-item instances with 48 stages, the 30 instances of Table~\ref{tab:d-long-other}, five two-item instances with 6 stages, and five full-size instances with 12 or 24 stages for which no optimal value is known.

\begin{table}[H]
\centering
\caption{Iterations and Cuts in the 300-Second Runs of Table~\ref{tab:c1}}\label{tab:d-300}
\begin{singlespace}\scriptsize
\setlength{\tabcolsep}{3pt}
\begin{tabular*}{\textwidth}{@{\extracolsep{\fill}}lrr rrr rrr rrr rr@{}}
\toprule
 & & & \multicolumn{3}{c}{\textnormal{LCD}} & \multicolumn{3}{c}{\textnormal{LCD-Iso}} & \multicolumn{3}{c}{SDDmiP} & \multicolumn{2}{c}{SDDiP} \\
\cmidrule(lr){4-6}\cmidrule(lr){7-9}\cmidrule(lr){10-12}\cmidrule(l){13-14}
Instance & $n$ & $T$ & Iter. & Cuts & Closing & Iter. & Cuts & Closing & Iter. & Cuts & ReLU & Iter. & Cuts \\
\midrule
Case-Pack-Truck & 2 & 10 & 16 & 813 & 482 & 7 & 1,898 & 191 & 51 & 597 & 512 & 126 & 2,268 \\
Case-Pack-Truck & 2 & 24 & 8 & 1,281 & 739 & 4 & 5,874 & 303 & 48 & 880 & 681 & 83 & 3,842 \\
Case-Pack-Truck & 2 & 24 & 7 & 1,236 & 635 & 4 & 5,265 & 195 & 55 & 669 & 508 & 92 & 4,232 \\
Case-Pack-Truck & 8 & 6 & 9 & 530 & 147 & 5 & 1,608 & 9 & 30 & 225 & 44 & 78 & 782 \\
Serial-Inventory & 5 & 24 & 10 & 2,017 & 164 & 1 & 6,369 & 0 & 77 & 1,107 & 427 & 115 & 5,320 \\
Co-Product & 4 & 24 & 19 & 1,865 & 1,520 & 3 & 2,266 & 19 & 89 & 954 & 833 & 121 & 5,590 \\
Rotate & 12 & 24 & 22 & 115 & 0 & 15 & 11,460 & 0 & 5 & 195 & 79 & 145 & 6,694 \\
Shear & 12 & 24 & 1 & 2,030 & 0 & 0 & 5,380 & 20 & 120 & 5,555 & 14 & --- & 10 \\
Shift & 96 & 64 & 9 & 331 & 0 & 3 & 20,795 & 0 & 5 & 675 & 40 & 23 & 2,898 \\
Rotate & 96 & 64 & 5 & 315 & 0 & 4 & 22,135 & 0 & 2 & 331 & 16 & 15 & 1,940 \\
Odd-Cycle & 96 & 64 & 1 & 3,086 & 0 & 0 & 6,203 & 28 & 3 & 1,145 & 7 & 4 & 510 \\
Shear & 96 & 64 & 0 & 737 & 0 & --- & 6,145 & 0 & 18 & 2,092 & 0 & --- & 2 \\
 
\bottomrule
\end{tabular*}
\end{singlespace}
\figurenotes{``Cuts'' is the total number of cuts generated, ``Closing'' the number of closing cuts among them, and ``ReLU'' the number of rectified linear cuts of SDDmiP. Rows are in the order of Table~\ref{tab:c1}.}
\end{table}

\begin{table}[H]
\centering
\caption{Certified Lower Bounds and Mean Policy Costs of \textnormal{LCD}, \textnormal{LCD-Iso}, SDDmiP, and SDDiP on 48 Instances (Time Limit of 45 Seconds)}\label{tab:d-45}
\begin{singlespace}\scriptsize
\setlength{\tabcolsep}{3pt}
\resizebox{\textwidth}{!}{%
\begin{tabular}{@{}lrrrr rrrr rrrr@{}}
\toprule
 & & & & & \multicolumn{4}{c}{Certified lower bound} & \multicolumn{4}{c}{Mean policy cost} \\
\cmidrule(lr){6-9}\cmidrule(l){10-13}
Instance & $n$ & $T$ & Seed & $V^\star$ & \textnormal{LCD} & \textnormal{LCD-Iso} & SDDmiP & SDDiP & \textnormal{LCD} & \textnormal{LCD-Iso} & SDDmiP & SDDiP \\
\midrule
Case-Pack-Truck & 2 & 10 & 0 & 887.7 & \tablebold{863.4} & 851.3 & 839.7 & 519.1 & 894 & 918 & 918 & 967 \\
Case-Pack-Truck & 2 & 10 & 1 & 773.6 & \tablebold{749.2} & 738.5 & 712.3 & 459.7 & 786 & 779 & 778 & 842 \\
Case-Pack-Truck & 2 & 24 & 0 & 2,058.7 & \tablebold{1,866.9} & 1,818.8 & 1,855.4 & 419.4 & 2,174 & 2,148 & 2,130 & 2,449 \\
Case-Pack-Truck & 2 & 24 & 1 & 1,900.3 & \tablebold{1,648.9} & 1,585.7 & 1,619.4 & 425.4 & 2,014 & 1,984 & 2,002 & 2,229 \\
Case-Pack & 2 & 48 & 0 & 3,247.7 & \tablebold{3,247.0} & 2,932.9 & 3,174.8 & 299.3 & 3,258 & 4,118 & 3,259 & 3,947 \\
Case-Pack & 2 & 48 & 1 & 3,022.5 & 2,376.0 & 2,432.5 & \tablebold{2,793.4} & 273.3 & 3,177 & 3,520 & 3,045 & 3,679 \\
Case-Pack-Truck & 4 & 6 & 0 & 446.4 & \tablebold{418.0} & 368.3 & 396.9 & 334.0 & 455 & 457 & 449 & 482 \\
Case-Pack & 4 & 6 & 0 & 353.5 & \tablebold{343.1} & 329.5 & 336.3 & 269.9 & 361 & 410 & 357 & 383 \\
Case-Pack, lead time 3 & 4 & 10 & 0 & --- & \tablebold{709.5} & 706.4 & 702.5 & 546.0 & 709 & 709 & 707 & 751 \\
Case-Pack, lead time 3 & 4 & 24 & 0 & --- & \tablebold{1,428.5} & 1,423.8 & 1,410.4 & 556.1 & 1,501 & 1,473 & 1,516 & 1,746 \\
Case-Pack, lead time 3 & 4 & 24 & 1 & 557.4 & \tablebold{511.4} & 476.8 & 444.0 & 370.0 & 588 & 589 & 589 & 628 \\
Case-Pack-Truck & 8 & 6 & 0 & 1,014.8 & 897.9 & \tablebold{901.4} & 895.9 & 493.8 & 1,053 & 1,039 & 1,036 & 1,159 \\
Case-Pack & 8 & 6 & 0 & 893.6 & \tablebold{850.5} & 834.5 & 825.0 & 431.9 & 906 & 907 & 910 & 1,051 \\
Case-Pack & 16 & 6 & 0 & 1,218.0 & 817.2 & 768.3 & \tablebold{1,057.4} & 384.3 & 1,667 & 1,438 & 1,420 & 1,709 \\
Case-Pack-Truck & 16 & 12 & 0 & --- & 972.8 & 1,067.0 & \tablebold{2,032.5} & 440.5 & --- & --- & 3,215 & --- \\
\addlinespace[2pt]
Workforce-Cohort & 2 & 48 & 0 & 2,113.1 & \tablebold{2,113.1} & 2,084.4 & 2,076.3 & 1,607.4 & 2,126 & 2,125 & 2,179 & 2,561 \\
Workforce-Cohort & 2 & 48 & 1 & 1,602.6 & \tablebold{1,602.6} & \tablebold{1,602.6} & 1,601.3 & 1,274.0 & 1,595 & 1,595 & 1,595 & 1,846 \\
Workforce-Cohort & 4 & 6 & 0 & 3,685.3 & \tablebold{3,685.3} & 3,680.2 & 3,684.7 & 2,875.5 & 3,683 & 3,687 & 3,683 & 3,799 \\
Workforce-Cohort & 6 & 6 & 0 & 3,879.9 & \tablebold{3,879.9} & \tablebold{3,879.9} & 3,879.2 & 2,414.1 & 3,883 & 3,883 & 3,883 & 4,138 \\
Workforce-Cohort, no trainer limit & 6 & 6 & 0 & 3,826.7 & \tablebold{3,826.7} & \tablebold{3,826.7} & \tablebold{3,826.7} & 3,063.1 & 3,828 & 3,826 & 3,826 & 4,057 \\
Workforce-Cohort & 12 & 6 & 0 & 7,112.6 & 7,063.4 & \tablebold{7,077.7} & 7,063.6 & 3,716.2 & 7,161 & 7,118 & 7,125 & 7,948 \\
Workforce-Cohort, no trainer limit & 12 & 6 & 0 & 7,027.8 & \tablebold{7,023.2} & 7,018.4 & 7,013.4 & 4,096.6 & 7,106 & 7,019 & 7,047 & 7,894 \\
Workforce-Cohort & 18 & 12 & 0 & --- & 12,468.0 & 12,358.7 & \tablebold{12,696.3} & 3,363.8 & 13,086 & 13,262 & 13,057 & 14,431 \\
\addlinespace[2pt]
Co-Product & 2 & 48 & 0 & 828.7 & 805.8 & \tablebold{826.6} & 782.6 & 801.6 & 834 & 831 & 833 & 843 \\
Co-Product & 2 & 48 & 1 & 639.8 & 608.1 & 578.6 & 587.8 & \tablebold{609.2} & 644 & 644 & 643 & 657 \\
Co-Product & 4 & 6 & 0 & 905.4 & 902.6 & 899.7 & \tablebold{905.1} & 901.9 & 892 & 892 & 887 & 889 \\
Co-Product & 4 & 10 & 0 & --- & 1,836.7 & 1,829.1 & \tablebold{1,838.1} & 1,835.4 & 1,830 & 1,825 & 1,825 & 1,831 \\
Co-Product & 4 & 24 & 0 & --- & \tablebold{4,430.5} & 4,419.3 & \tablebold{4,430.2} & 4,429.1 & 4,474 & 4,470 & 4,469 & 4,483 \\
Co-Product & 4 & 24 & 1 & --- & \tablebold{3,416.3} & 3,407.7 & 3,415.4 & 3,412.3 & 3,427 & 3,424 & 3,424 & 3,435 \\
Co-Product & 6 & 6 & 0 & 1,337.8 & \tablebold{1,337.8} & 1,337.2 & 1,337.3 & 1,335.7 & 1,335 & 1,337 & 1,335 & 1,339 \\
Co-Product & 12 & 6 & 0 & 2,564.3 & \tablebold{2,544.0} & 2,542.0 & 2,538.1 & 2,470.7 & 2,570 & 2,573 & 2,580 & 2,663 \\
Co-Product & 12 & 12 & 0 & --- & 5,643.3 & \tablebold{5,648.1} & \tablebold{5,648.5} & 5,575.1 & 5,671 & 5,672 & 5,697 & 5,830 \\
Co-Product & 24 & 24 & 0 & --- & \tablebold{23,618.6} & 23,607.2 & 23,558.4 & 23,544.8 & 23,784 & 23,796 & 23,809 & 24,102 \\
\addlinespace[2pt]
Bike-Sharing & 2 & 48 & 0 & 412.1 & \tablebold{343.0} & 337.3 & 342.4 & 42.0 & 443 & 435 & 435 & 483 \\
Bike-Sharing & 2 & 48 & 1 & 267.3 & \tablebold{253.2} & 252.9 & 252.5 & 6.5 & 277 & 279 & 280 & 302 \\
Bike-Sharing & 3 & 6 & 0 & 66.5 & \tablebold{62.2} & 59.7 & 59.7 & 24.3 & 76 & 70 & 70 & 79 \\
Bike-Sharing & 3 & 10 & 0 & --- & \tablebold{190.7} & 189.9 & 187.6 & 20.6 & 197 & 193 & 194 & 263 \\
Bike-Sharing & 3 & 24 & 0 & --- & \tablebold{482.7} & 482.2 & 474.8 & 14.4 & 481 & 479 & 482 & --- \\
Bike-Sharing & 4 & 6 & 0 & 236.0 & \tablebold{230.3} & 228.7 & 228.0 & 56.3 & 240 & 244 & 237 & 325 \\
Bike-Sharing & 7 & 6 & 0 & 406.6 & \tablebold{380.8} & 380.6 & 372.4 & 59.3 & 407 & 405 & 409 & 455 \\
Bike-Sharing & 13 & 24 & 0 & --- & \tablebold{2,035.4} & 2,032.1 & 1,796.9 & 70.7 & 2,486 & 2,470 & 2,491 & --- \\
\addlinespace[2pt]
Serial-Inventory & 5 & 10 & 0 & 90.8 & 80.7 & 61.5 & \tablebold{81.8} & 76.8 & 112 & 121 & 128 & 133 \\
Serial-Inventory & 5 & 24 & 0 & 270.0 & \tablebold{179.3} & 159.9 & 150.5 & 76.5 & 332 & 332 & 384 & 408 \\
Serial-Inventory & 5 & 24 & 1 & 344.8 & 219.1 & 186.5 & \tablebold{225.1} & 87.4 & 429 & 450 & 409 & 603 \\
Serial-Inventory, truck charge & 12 & 6 & 0 & 365.0 & 296.8 & 284.5 & \tablebold{307.6} & 145.1 & 465 & 441 & 385 & 527 \\
Serial-Inventory & 12 & 6 & 0 & 305.9 & 260.7 & 240.7 & \tablebold{268.1} & 124.2 & 351 & 362 & 328 & 390 \\
Serial-Inventory & 24 & 6 & 0 & 766.6 & 468.2 & 361.1 & \tablebold{692.2} & 221.5 & 1,074 & 972 & 823 & 1,553 \\
\addlinespace[2pt]
Airline-RM & 144 & 6 & 0 & -207,524.8 & \tablebold{-211,816.0} & -212,770.3 & -213,931.1 & -254,210.0 & -206,853 & -206,048 & -205,838 & -189,244 \\
 
\bottomrule
\end{tabular}}
\end{singlespace}
\figurenotes{\boldnote{Bold marks the highest bound and the bounds tied with it. }Policy costs are means over 32 common sampled paths, and a dash means that the evaluation did not finish within 20 seconds.}
\end{table}

Table~\ref{tab:d-buckets} breaks the comparisons of Table~\ref{tab:d-45} down by the number of cosets per stage. \textnormal{LCD} outperforms SDDmiP on 21 of the 27 instances with at most 64 cosets and on 3 of the 11 instances with 4,097 to $10^{8}$ cosets. However, it outperforms SDDmiP on all three instances with more than $10^{8}$ cosets, so that the number of cosets alone does not decide the comparison. The last column compares the mean policy costs of \textnormal{LCD} and SDDmiP on the 47 instances on which both policies were evaluated: the policy of \textnormal{LCD} has the lower mean cost on 16 instances and the higher on 27.

\begin{table}[H]
\centering
\caption{Wins, Ties, and Losses of \textnormal{LCD} on the 48 Instances of Table~\ref{tab:d-45} by the Number of Cosets per Stage}\label{tab:d-buckets}
\begin{singlespace}\footnotesize
\begin{tabular}{@{}lrccccc@{}}
\toprule
 & & \multicolumn{4}{c}{Certified bound} & Mean policy cost \\
\cmidrule(lr){3-6}\cmidrule(l){7-7}
Cosets per stage & Instances & vs.\ no lattice & vs.\ \textnormal{LCD-Iso} & vs.\ SDDmiP & vs.\ SDDiP & vs.\ SDDmiP \\
\midrule
at most 64 & 27 & 9--3--15 & 24--1--2 & 21--1--5 & 26--0--1 & 9--1--17 \\
65 to 4,096 & 7 & 1--5--1 & 4--2--1 & 6--1--0 & 7--0--0 & 1--3--3 \\
4,097 to $10^{8}$ & 11 & 4--0--7 & 8--0--3 & 3--1--7 & 11--0--0 & 3--0--7 \\
more than $10^{8}$ & 3 & 2--0--1 & 3--0--0 & 3--0--0 & 3--0--0 & 3--0--0 \\
\midrule
All & 48 & 16--8--24 & 39--3--6 & 33--3--12 & 47--0--1 & 16--4--27 \\
 
\bottomrule
\end{tabular}
\end{singlespace}
\figurenotes{Entries are wins--ties--losses of \textnormal{LCD}, with the tie rule of Appendix~\ref{app:methods-settings}. In the last column, the lower mean cost on the 32 common sampled paths wins, and one instance of the third row is missing because the evaluation of a policy did not finish.}
\end{table}

\clearpage
\subsection{Catalog}\label{app:results-catalog}

The catalog has 163 configurations. On 62 of them, none of the four methods obtains a bound within the limit of 20 or 25 seconds. Tables~\ref{tab:d-catalog-1} and~\ref{tab:d-catalog-2} list the other 101 configurations with the bounds of the four methods. Many of these bounds stem from the initial approximation or from a single iteration (Section~\ref{sec:exp-regimes}). Table~\ref{tab:ec-catalog} in Appendix~\ref{ec:benchmarks} summarizes the comparison of \textnormal{LCD} with SDDmiP by family.

\begin{table}[H]
\centering
\caption{Certified Bounds on the Catalog Configurations With at Least One Bound (Part 1)}\label{tab:d-catalog-1}
\begin{singlespace}\scriptsize
\begin{tabular}{@{}lrrrrrrr@{}}
\toprule
Family & $n$ & $T$ & Limit (s) & \textnormal{LCD} & \textnormal{LCD-Iso} & SDDmiP & SDDiP \\
\midrule
Shift & 96 & 64 & 25 & \tablebold{13,848.8} & \tablebold{13,848.8} & 13,663.9 & 13,416.8 \\
Shift & 96 & 256 & 25 & \tablebold{55,602.0} & --- & 54,835.5 & --- \\
Shift & 96 & 1,024 & 25 & --- & --- & \tablebold{127.0} & --- \\
Shift & 384 & 64 & 25 & \tablebold{54,470.3} & --- & 53,691.7 & --- \\
Shift & 384 & 256 & 25 & --- & --- & \tablebold{586.0} & --- \\
Shift & 1,536 & 64 & 25 & --- & --- & \tablebold{2,186.0} & --- \\
\addlinespace[2pt]
Rotate & 96 & 64 & 25 & 13,006.4 & 12,871.2 & 12,871.2 & \tablebold{13,128.3} \\
Rotate & 96 & 256 & 25 & \tablebold{51,637.8} & --- & \tablebold{51,633.9} & --- \\
Rotate & 96 & 1,024 & 25 & --- & --- & \tablebold{132.0} & --- \\
Rotate & 384 & 64 & 25 & \tablebold{50,507.1} & --- & \tablebold{50,503.2} & --- \\
Rotate & 384 & 256 & 25 & --- & --- & \tablebold{579.0} & --- \\
Rotate & 1,536 & 64 & 25 & --- & --- & \tablebold{2,258.0} & --- \\
\addlinespace[2pt]
Odd-Cycle & 96 & 64 & 25 & \tablebold{12,829.8} & 12,436.2 & 11,150.2 & --- \\
Odd-Cycle & 96 & 256 & 25 & \tablebold{51,517.8} & --- & 136.0 & --- \\
Odd-Cycle & 384 & 64 & 25 & \tablebold{50,127.8} & --- & 512.0 & --- \\
\addlinespace[2pt]
Shear & 96 & 64 & 25 & \tablebold{12,785.5} & --- & 777.3 & --- \\
Shear & 96 & 256 & 25 & \tablebold{51,508.4} & --- & 129.0 & --- \\
Shear & 384 & 64 & 25 & \tablebold{50,026.3} & --- & 538.0 & --- \\
\addlinespace[2pt]
Hub-Consolidation & 130 & 52 & 25 & \tablebold{73,635.4} & --- & 64,763.6 & --- \\
Hub-Consolidation & 130 & 104 & 25 & \tablebold{144,261.8} & --- & 126,719.9 & --- \\
Hub-Consolidation & 130 & 260 & 25 & \tablebold{104,694.9} & --- & 968.0 & --- \\
Hub-Consolidation & 514 & 52 & 25 & \tablebold{82,287.7} & --- & 4,510.0 & --- \\
\addlinespace[2pt]
Sterile-Trays & 256 & 84 & 25 & 28,845.1 & \tablebold{40,193.0} & 3,283.1 & --- \\
Sterile-Trays & 256 & 168 & 25 & \tablebold{27,774.2} & --- & 988.0 & --- \\
Sterile-Trays & 256 & 336 & 25 & --- & --- & \tablebold{988.0} & --- \\
Sterile-Trays & 1,024 & 84 & 25 & --- & --- & \tablebold{3,586.0} & --- \\
\addlinespace[2pt]
Pooled-Capacity & 24 & 8 & 20 & 7,548.5 & 7,734.9 & \tablebold{8,638.9} & 1,591.3 \\
Pooled-Capacity & 24 & 14 & 25 & 10,661.0 & 11,371.7 & \tablebold{12,753.8} & 1,691.7 \\
Pooled-Capacity & 24 & 25 & 25 & 17,853.6 & 12,413.5 & \tablebold{22,593.7} & 1,575.2 \\
Pooled-Capacity & 48 & 8 & 20 & 26,311.8 & 25,595.6 & \tablebold{27,712.9} & 4,968.6 \\
Pooled-Capacity & 48 & 14 & 25 & 29,821.8 & 27,111.9 & \tablebold{39,083.9} & 4,618.9 \\
Pooled-Capacity & 48 & 25 & 25 & 45,219.9 & 35,425.7 & \tablebold{47,605.6} & 4,466.0 \\
Pooled-Capacity & 80 & 8 & 25 & 31,432.3 & --- & \tablebold{33,731.3} & 6,686.8 \\
Pooled-Capacity & 80 & 14 & 25 & 45,112.4 & 33,518.1 & \tablebold{50,358.4} & 5,581.6 \\
Pooled-Capacity & 80 & 25 & 25 & \tablebold{71,906.4} & 52,831.2 & 56,339.2 & 3,386.8 \\
\addlinespace[2pt]
Fulfillment-Network & 528 & 52 & 25 & \tablebold{79,522.1} & --- & 4,417.0 & --- \\
Fulfillment-Network & 528 & 104 & 25 & --- & --- & \tablebold{4,417.0} & --- \\
Fulfillment-Network & 1,056 & 52 & 25 & --- & --- & \tablebold{7,515.0} & --- \\
\addlinespace[2pt]
Case-Pack & 16 & 12 & 20 & 715.1 & 638.7 & \tablebold{2,884.6} & 497.0 \\
Case-Pack & 16 & 24 & 25 & 633.7 & 638.2 & \tablebold{3,879.3} & 491.0 \\
Case-Pack & 16 & 48 & 25 & 566.6 & 555.5 & \tablebold{2,356.9} & 491.0 \\
Case-Pack & 32 & 12 & 20 & 1,081.0 & 1,106.5 & \tablebold{2,300.3} & 927.0 \\
Case-Pack & 32 & 24 & 25 & 983.6 & 986.0 & \tablebold{1,832.5} & 923.0 \\
Case-Pack & 32 & 48 & 25 & 982.3 & 984.0 & \tablebold{1,252.5} & 917.9 \\
Case-Pack & 64 & 12 & 25 & 2,250.1 & 2,256.7 & \tablebold{2,307.2} & 2,116.8 \\
Case-Pack & 64 & 24 & 25 & 2,256.7 & 2,252.7 & \tablebold{2,278.2} & 2,099.8 \\
Case-Pack & 64 & 48 & 25 & 2,251.7 & 2,221.0 & \tablebold{2,253.8} & 2,077.7 \\
\addlinespace[2pt]
Workforce-Cohort & 18 & 12 & 20 & 13,350.4 & 13,246.0 & \tablebold{13,786.1} & 2,835.0 \\
Workforce-Cohort & 18 & 24 & 25 & 25,307.8 & 22,735.1 & \tablebold{26,507.4} & 2,786.3 \\
Workforce-Cohort & 18 & 48 & 25 & 48,069.4 & 44,551.4 & \tablebold{49,661.5} & 2,708.0 \\
Workforce-Cohort & 36 & 12 & 20 & 27,143.3 & 25,805.6 & \tablebold{28,811.0} & 6,082.7 \\
 
\bottomrule
\end{tabular}
\end{singlespace}
\end{table}

\begin{table}[H]
\centering
\caption{Certified Bounds on the Catalog Configurations With at Least One Bound (Part 2)}\label{tab:d-catalog-2}
\begin{singlespace}\scriptsize
\begin{tabular}{@{}lrrrrrrr@{}}
\toprule
Family & $n$ & $T$ & Limit (s) & \textnormal{LCD} & \textnormal{LCD-Iso} & SDDmiP & SDDiP \\
\midrule
Workforce-Cohort & 36 & 24 & 25 & 49,692.4 & 47,638.6 & \tablebold{54,868.8} & 6,040.2 \\
Workforce-Cohort & 36 & 48 & 25 & \tablebold{96,886.8} & 92,345.6 & 79,861.6 & 5,876.9 \\
Workforce-Cohort & 72 & 12 & 20 & 54,006.1 & 52,739.8 & \tablebold{57,306.5} & 12,944.1 \\
Workforce-Cohort & 72 & 24 & 25 & \tablebold{103,674.0} & 101,365.1 & 96,602.9 & 12,537.9 \\
Workforce-Cohort & 72 & 48 & 25 & \tablebold{203,351.0} & 197,929.4 & 102,738.3 & 11,439.7 \\
\addlinespace[2pt]
Co-Product & 12 & 12 & 20 & \tablebold{4,974.1} & 4,965.3 & 4,949.4 & 4,958.4 \\
Co-Product & 12 & 24 & 25 & 10,421.5 & \tablebold{10,430.1} & 10,420.8 & 10,427.1 \\
Co-Product & 12 & 48 & 25 & 21,306.9 & 21,301.1 & 21,307.0 & \tablebold{21,323.6} \\
Co-Product & 24 & 12 & 20 & \tablebold{12,370.5} & 12,366.2 & 12,351.1 & 12,361.8 \\
Co-Product & 24 & 24 & 25 & \tablebold{26,672.7} & 26,667.3 & 26,661.2 & 26,640.0 \\
Co-Product & 24 & 48 & 25 & \tablebold{54,790.2} & \tablebold{54,793.1} & \tablebold{54,790.9} & 54,775.6 \\
Co-Product & 48 & 12 & 20 & 24,905.5 & \tablebold{24,915.2} & 24,891.5 & 24,903.1 \\
Co-Product & 48 & 24 & 25 & 52,478.4 & 52,495.5 & 52,466.8 & \tablebold{52,508.3} \\
Co-Product & 48 & 48 & 25 & 108,065.6 & 108,069.0 & 108,051.2 & \tablebold{108,111.9} \\
\addlinespace[2pt]
Bike-Sharing & 7 & 12 & 20 & \tablebold{368.2} & 367.5 & 313.4 & 54.6 \\
Bike-Sharing & 7 & 24 & 25 & 707.3 & \tablebold{714.1} & 587.4 & 50.0 \\
Bike-Sharing & 7 & 48 & 25 & 1,432.3 & \tablebold{1,448.4} & 1,040.4 & 40.0 \\
Bike-Sharing & 13 & 12 & 20 & \tablebold{703.4} & 687.4 & 634.8 & 59.3 \\
Bike-Sharing & 13 & 24 & 25 & \tablebold{1,508.8} & 1,499.9 & 1,303.8 & 45.0 \\
Bike-Sharing & 13 & 48 & 25 & \tablebold{3,043.1} & 2,990.0 & 2,525.5 & 41.3 \\
Bike-Sharing & 25 & 12 & 20 & \tablebold{1,510.8} & 1,466.3 & 1,278.6 & 110.0 \\
Bike-Sharing & 25 & 24 & 25 & 3,186.1 & \tablebold{3,274.0} & 2,699.0 & 96.3 \\
Bike-Sharing & 25 & 48 & 25 & \tablebold{6,228.4} & 6,160.9 & 4,716.3 & 78.3 \\
\addlinespace[2pt]
Serial-Inventory & 30 & 12 & 20 & 808.8 & 531.3 & \tablebold{1,345.5} & 345.5 \\
Serial-Inventory & 30 & 24 & 25 & 802.8 & 797.3 & \tablebold{1,293.0} & 95.0 \\
Serial-Inventory & 30 & 48 & 25 & 492.2 & 489.1 & \tablebold{905.4} & 327.0 \\
Serial-Inventory & 60 & 12 & 20 & 1,746.7 & 1,640.5 & \tablebold{2,116.6} & 398.4 \\
Serial-Inventory & 60 & 24 & 25 & 1,795.8 & 1,720.2 & \tablebold{1,844.7} & 190.0 \\
Serial-Inventory & 60 & 48 & 25 & \tablebold{1,850.3} & 1,847.9 & 1,798.8 & --- \\
Serial-Inventory & 120 & 12 & 20 & \tablebold{1,820.2} & 1,749.5 & 1,817.7 & 1,116.1 \\
Serial-Inventory & 120 & 24 & 25 & \tablebold{2,411.3} & 2,402.3 & 2,274.6 & 289.0 \\
Serial-Inventory & 120 & 48 & 25 & 1,686.5 & \tablebold{1,688.9} & 1,657.9 & 868.3 \\
\addlinespace[2pt]
Lead-Time-Pipeline & 48 & 12 & 20 & \tablebold{1,792.9} & 1,256.2 & 1,551.7 & 845.6 \\
Lead-Time-Pipeline & 48 & 24 & 25 & 1,353.8 & 985.1 & \tablebold{1,476.4} & 612.0 \\
Lead-Time-Pipeline & 48 & 48 & 25 & \tablebold{1,227.4} & 1,128.0 & 1,213.5 & --- \\
Lead-Time-Pipeline & 96 & 12 & 20 & 1,689.0 & 1,752.0 & \tablebold{2,739.7} & 1,172.9 \\
Lead-Time-Pipeline & 96 & 24 & 25 & 1,797.7 & 2,081.7 & \tablebold{2,887.5} & 1,161.0 \\
Lead-Time-Pipeline & 96 & 48 & 25 & 1,690.3 & 2,292.5 & \tablebold{2,366.9} & --- \\
Lead-Time-Pipeline & 144 & 12 & 25 & 4,370.4 & 2,170.0 & \tablebold{4,388.1} & 1,569.9 \\
Lead-Time-Pipeline & 144 & 24 & 25 & 1,978.0 & 2,136.8 & \tablebold{4,381.7} & 1,539.9 \\
Lead-Time-Pipeline & 144 & 48 & 25 & 2,208.0 & 2,634.5 & \tablebold{3,730.0} & --- \\
\addlinespace[2pt]
Airline-RM & 72 & 6 & 25 & \tablebold{-166,373.0} & -167,189.0 & -166,654.2 & -200,280.0 \\
Airline-RM & 72 & 12 & 25 & -292,589.7 & -296,546.0 & \tablebold{-290,757.7} & -358,300.0 \\
Airline-RM & 72 & 24 & 25 & \tablebold{-561,762.3} & -579,827.8 & -566,084.8 & -670,580.0 \\
Airline-RM & 144 & 6 & 25 & \tablebold{-335,255.2} & -336,598.5 & -336,600.7 & -396,120.0 \\
Airline-RM & 144 & 12 & 25 & -585,480.6 & -594,003.5 & \tablebold{-579,807.2} & -687,800.0 \\
Airline-RM & 144 & 24 & 25 & -1,085,836.5 & -1,188,681.3 & \tablebold{-1,068,745.6} & -1,254,000.0 \\
Airline-RM & 288 & 6 & 25 & \tablebold{-699,004.1} & -706,161.7 & -705,434.7 & -806,565.0 \\
Airline-RM & 288 & 12 & 25 & -1,206,237.4 & -1,251,639.0 & \tablebold{-1,197,529.0} & -1,396,133.0 \\
Airline-RM & 288 & 24 & 25 & -2,310,494.7 & -2,470,445.5 & \tablebold{-2,243,966.7} & -2,594,055.0 \\
 
\bottomrule
\end{tabular}
\end{singlespace}
\figurenotes{\boldnote{Bold marks the highest bound and the bounds tied with it. }A dash means that the method obtained no bound within the limit. Optimal values are not known for these configurations.}
\end{table}

\clearpage
\subsection{Which Results Apply to Which Family}\label{app:results-class}

Table~\ref{tab:d-class} states, for each family, what we know about the exactness of lattice cuts. The second column refers to the sufficient condition given after Proposition~\ref{prop:phasevalue}, which we checked on the formulations of Appendix~\ref{app:benchmarks}. The third column reports what was computed, where $W_\Gamma(0)$ is the bound of Theorem~\ref{thm:root-price} at $\rho=0$, evaluated by dynamic programming (Tables~\ref{tab:d-long-one} and~\ref{tab:d-long-other}). We did not verify Assumption~\ref{assump:period} for any benchmark family. Theorem~\ref{thm:optimal-period} enters the paper as a statement about the model class \eqref{eq:period-model} and is not used to explain a computational result. In Appendix~\ref{app:add-theory}, we check its conclusion numerically on four small instances.

\begin{table}[H]
\centering
\caption{Structural Statements and Computed Evidence on the Exactness of Lattice Cuts by Family}\label{tab:d-class}
\begin{singlespace}\footnotesize
\begin{tabular}{@{}p{0.22\linewidth}p{0.30\linewidth}p{0.42\linewidth}@{}}
\toprule
Family & Condition of Proposition~\ref{prop:phasevalue} & Computed evidence \\
\midrule
Shift, Rotate & holds by construction: the state bounds never bind & no closing cut in any run of either variant \\
\addlinespace[2pt]
Odd-Cycle, Shear & not established: balancing actions are limited by the state bounds & no closing cut in the runs of \textnormal{LCD}; optimal values unknown \\
\addlinespace[2pt]
Group II (four families) & not established & no closing cut in the runs of \textnormal{LCD}; optimal values unknown \\
\addlinespace[2pt]
Workforce-Cohort & not established & $W_\Gamma(0)=V^\star$ on all eight one-cohort instances, which reach $V^\star$ without a closing cut \\
\addlinespace[2pt]
Case-Pack & not established & $W_\Gamma(0)=V^\star$ on all eight one-item instances, which reach $V^\star$ without a closing cut \\
\addlinespace[2pt]
Co-Product & does not hold & $W_\Gamma(0)<V^\star$ on six of the eight one-recipe instances, with a relative loss of at most $3.5\times10^{-5}$; closing cuts on all eight \\
\addlinespace[2pt]
Case-Pack-Truck & does not hold: the fixed charge is not a function of the coset & $W_\Gamma(0)$ is 1.4\% to 2.4\% below $V^\star$ on the four one-item instances \\
\addlinespace[2pt]
Bike-Sharing & not established & $W_\Gamma(0)=V^\star$ on all eight one-station instances; the long runs do not reach $V^\star$ within the limit, so their gap is finite-pool loss \\
\addlinespace[2pt]
Serial-Inventory & not established & $W_\Gamma(0)=V^\star$ on the three one-item instances \\
\addlinespace[2pt]
Lead-Time-Pipeline, Airline-RM & not established & not evaluated exactly \\
\bottomrule
\end{tabular}
\end{singlespace}
\end{table}

\clearpage
\ecresetcounters
\section{Additional Experiments}\label{app:additional}

In this appendix, we report experiments that we made after the study of Section~\ref{sec:experiments} in order to examine specific questions: whether the cuts of \textnormal{LCD} are valid at every state (Appendix~\ref{app:add-validity}), how SDDmiP depends on its settings (Appendix~\ref{app:add-sddmip}), what the lattice contributes within \textnormal{LCD} (Appendix~\ref{app:add-lattice}), how often computed intercepts are reused (Appendix~\ref{app:add-reuse}), how good the policies are when they are evaluated exactly (Appendix~\ref{app:add-policy}), and what happens without batch structure (Appendix~\ref{app:add-nobatch}). In Appendix~\ref{app:add-theory}, we report numerical checks of two statements of Sections~\ref{sec:model} and~\ref{sec:lattice-cuts}.

These experiments differ from those of Section~\ref{sec:experiments} in three respects. First, their time limits are short, between 12 and 180 seconds, except in Tables~\ref{tab:x-long-sddmip} and~\ref{tab:x-crossover}. Second, where Table~\ref{tab:x-design} says so, \textnormal{LCD} was run with the corrected merging rule of Appendix~\ref{app:methods-ours}, whereas Tables~\ref{tab:x-settings}, \ref{tab:x-highcoset}, and~\ref{tab:x-nobatch} come from earlier runs without the slope shift (Proposition~\ref{prop:merge-tolerance}). Third, the runs of different instances were executed four at a time (six in Tables~\ref{tab:x-long-sddmip} and~\ref{tab:x-crossover}) on the eight-core computer, with all methods of an instance one after the other in one process and every run on one thread. Because of the first and third points, the bounds in a table of this appendix are comparable with each other, and we do not compare them with the tables of Appendix~\ref{app:results}. All runs are single runs with seed 0 unless stated otherwise, the instances are named as in Appendix~\ref{app:results}, and two bounds are counted as tied by the rule of Appendix~\ref{app:methods-settings}.

\begin{table}[H]
\centering
\caption{Additional Experiments and Their Settings}\label{tab:x-design}
\begin{singlespace}\footnotesize
\setlength{\tabcolsep}{3pt}
\begin{tabular}{@{}p{0.30\linewidth}rp{0.31\linewidth}p{0.12\linewidth}p{0.15\linewidth}@{}}
\toprule
Experiment & Inst. & Methods & Time limit & Merging in \textnormal{LCD} \\
\midrule
Audit of the cuts (Table~\ref{tab:x-audit}) & 5 & \textnormal{LCD} & 20 s & three rules \\
\addlinespace[2pt]
Merging rules under a time limit (Table~\ref{tab:x-merge}) & 24 & \textnormal{LCD} & 12--45 s & three rules \\
\addlinespace[2pt]
Threshold of SDDmiP (Tables~\ref{tab:x-sweep} and~\ref{tab:x-targets}) & 10 & \textnormal{LCD}, \textnormal{LCD-Iso}, SDDmiP with four thresholds & 90 s & corrected \\
\addlinespace[2pt]
Long-horizon calibrated SDDmiP (Table~\ref{tab:x-long-sddmip}) & 28 & SDDmiP with threshold $10^{-3}$ & 30--60 min & --- \\
\addlinespace[2pt]
Costs multiplied by 100 (Table~\ref{tab:x-scale}) & 6 & \textnormal{LCD}, SDDmiP with two thresholds & 60 s & corrected \\
\addlinespace[2pt]
Stage gap and thresholds of SDDmiP (Table~\ref{tab:x-settings}) & 10 & \textnormal{LCD} and SDDmiP with several settings & 45 s & uncorrected \\
\addlinespace[2pt]
Lattice on and off (Table~\ref{tab:x-ablation}) & 16 & \textnormal{LCD} and \textnormal{LCD-Iso}, each with and without the lattice & 45 s & corrected \\
\addlinespace[2pt]
Long-budget lattice crossover (Table~\ref{tab:x-crossover}) & 25 & \textnormal{LCD} with and without the lattice & 1,800 s & corrected \\
\addlinespace[2pt]
Three seeds (Table~\ref{tab:x-seeds}) & 6 & \textnormal{LCD} with and without the lattice, SDDmiP, SDDiP & 180 s & corrected \\
\addlinespace[2pt]
Many cosets (Table~\ref{tab:x-highcoset}) & 10 & \textnormal{LCD} with and without the lattice, SDDmiP, SDDiP & 60 s & uncorrected \\
\addlinespace[2pt]
Reuse of intercepts (Table~\ref{tab:x-reuse}) & 10 & \textnormal{LCD} & 45 s & corrected \\
\addlinespace[2pt]
Exact policy evaluation (Table~\ref{tab:x-policy}) & 8 & \textnormal{LCD} with and without the lattice, SDDmiP with two thresholds & 30 s or 45 s & corrected \\
\addlinespace[2pt]
Without batch structure (Table~\ref{tab:x-nobatch}) & 12 & \textnormal{LCD}, \textnormal{LCD-Iso} without the lattice, SDDmiP, SDDiP & 30 s & uncorrected \\
\bottomrule
\end{tabular}
\end{singlespace}
\figurenotes{``Corrected'' refers to the box-shift merging rule of Appendix~\ref{app:methods-ours}, and ``uncorrected'' refers to the tolerances analyzed in Proposition~\ref{prop:merge-tolerance}. Unless a table states otherwise, SDDmiP and SDDiP use the settings of Appendix~\ref{app:methods}.}
\end{table}

\subsection{Validity of the Cuts}\label{app:add-validity}

We examine the validity of the cuts of \textnormal{LCD} in two ways. First, we audit every cut at every state. We consider five instances that are small enough for dynamic programming over the full state space: the 10-stage Co-Product, Workforce-Cohort, and Case-Pack-Truck instances with seed 0, a 6-stage Case-Pack-Truck instance with pack size one, and a 4-stage lot-sizing instance with two products (Appendix~\ref{app:add-nobatch}). On each of them, we run \textnormal{LCD} for 20 seconds, and we then evaluate every lattice cut and closing cut that it has generated at every state of its stage and compare the value $h(y)$ with the exact value $Q_{tw}(y)$. Table~\ref{tab:x-audit} summarizes the audit for three merging rules: the uncorrected tolerances of Proposition~\ref{prop:merge-tolerance}, the corrected rule of Appendix~\ref{app:methods-ours}, and merging and reuse for identical slopes only.

\begin{table}[H]
\centering
\caption{Audit of All Cuts of \textnormal{LCD} at All States Against Exact Dynamic Programming on Five Instances}\label{tab:x-audit}
\begin{singlespace}\footnotesize
\setlength{\tabcolsep}{4.5pt}
\begin{tabular}{@{}lrrrrc@{}}
\toprule
Merging rule & Instances & Cuts & Evaluations & Above $Q_{tw}$ & Largest $h(y)-Q_{tw}(y)$ \\
\midrule
Tolerances as in the reported runs & 5 & 697 & 356,284 & 0 & $-1.0\times10^{-6}$ \\
Intercepts corrected for the slope difference & 5 & 1,005 & 459,292 & 0 & $-1.0\times10^{-6}$ \\
Identical slopes only & 5 & 1,014 & 465,568 & 0 & $-1.0\times10^{-6}$ \\
 
\bottomrule
\end{tabular}
\end{singlespace}
\figurenotes{``Evaluations'' is the number of pairs of a cut and a state of its stage, and ``Above $Q_{tw}$'' is the number of these pairs with $h(y)>Q_{tw}(y)$. The runs with the first rule and the runs with the other two rules are separate runs, which is why the numbers of cuts differ.}
\end{table}

We find no cut above the exact value at any state under any of the three rules. The largest value of $h(y)-Q_{tw}(y)$ is $-1.0\times10^{-6}$, which is the absolute part of the safety margin. As a control, we also audited a variant that takes the value of the incumbent instead of the dual bound of every MIP at a stage gap of $10^{-4}$. For this variant, 529 of 381,035 evaluations lie above the exact value, by up to $1.7\times10^{-2}$, so that the audit does detect invalid cuts.

Second, we compare the three merging rules under time limits. Table~\ref{tab:x-merge} reports the bounds of 24 runs with limits between 12 and 45 seconds, together with the number of merges in the run with the corrected rule, the number of merges whose correction is not zero, and the largest correction of an intercept.

\begin{table}[H]
\centering
\caption{Bounds of \textnormal{LCD} With Three Merging Rules on 24 Runs}\label{tab:x-merge}
\begin{singlespace}\scriptsize
\setlength{\tabcolsep}{3pt}
\begin{tabular*}{\textwidth}{@{\extracolsep{\fill}}lrrrr rrr rrc@{}}
\toprule
 & & & & & \multicolumn{3}{c}{Certified lower bound} & \multicolumn{3}{c}{Corrected rule} \\
\cmidrule(lr){6-8}\cmidrule(l){9-11}
Instance & $n$ & $T$ & Seed & Limit (s) & Uncorrected & Corrected & Identical slopes & Merges & Nonzero & Largest \\
\midrule
Case-Pack-Truck & 2 & 10 & 0 & 30 & 858.65 & 858.65 & 858.51 & 739 & 8 & $1.6\times10^{-8}$ \\
Case-Pack-Truck & 2 & 10 & 1 & 30 & 746.16 & 745.49 & 747.29 & 737 & 22 & $5.8\times10^{-8}$ \\
Case-Pack-Truck & 2 & 24 & 0 & 45 & 1,877.01 & 1,876.98 & 1,863.55 & 2,257 & 13 & $3.6\times10^{-8}$ \\
Case-Pack-Truck & 2 & 24 & 1 & 45 & 1,650.76 & 1,651.56 & 1,654.43 & 2,354 & 141 & $1.4\times10^{-13}$ \\
Shift & 12 & 24 & 6101 & 12 & 752.40 & 752.40 & 752.40 & 2,864 & 0 & 0 \\
Shift & 12 & 24 & 6102 & 12 & 736.20 & 736.40 & 736.80 & 2,946 & 0 & 0 \\
Shear & 12 & 24 & 6101 & 12 & 683.54 & 683.54 & 683.54 & 492 & 0 & 0 \\
Shear & 12 & 24 & 6102 & 12 & 675.23 & 675.23 & 675.23 & 500 & 0 & 0 \\
Odd-Cycle & 12 & 24 & 6101 & 12 & 677.25 & 677.25 & 677.25 & 147 & 4 & $1.2\times10^{-16}$ \\
Odd-Cycle & 12 & 24 & 6102 & 12 & 677.03 & 677.03 & 677.03 & 112 & 3 & $8.5\times10^{-16}$ \\
Rotate & 12 & 24 & 6101 & 12 & 680.80 & 680.80 & 680.80 & 253 & 0 & 0 \\
Rotate & 12 & 24 & 6102 & 12 & 707.60 & 707.60 & 707.60 & 233 & 0 & 0 \\
Shift & 96 & 64 & 6101 & 15 & 13,848.84 & 13,848.84 & 13,848.84 & 833 & 0 & 0 \\
Shear & 96 & 64 & 6101 & 15 & 12,785.51 & 12,785.51 & 12,785.51 & 523 & 7 & $8.6\times10^{-8}$ \\
Odd-Cycle & 96 & 64 & 6101 & 15 & 12,828.05 & 12,828.05 & 12,828.05 & 3 & 1 & $1.6\times10^{-17}$ \\
Rotate & 96 & 64 & 6101 & 15 & 12,871.21 & 12,871.21 & 12,871.21 & 270 & 0 & 0 \\
Shift & 96 & 64 & 6101 & 45 & 13,918.96 & 13,920.36 & 13,918.96 & 2,844 & 2 & $1.8\times10^{-15}$ \\
Shear & 96 & 64 & 6101 & 45 & 12,785.51 & 12,785.51 & 12,785.51 & 2,038 & 0 & 0 \\
Odd-Cycle & 96 & 64 & 6101 & 45 & 12,829.40 & 12,829.40 & 12,829.40 & 8 & 2 & $4.2\times10^{-16}$ \\
Rotate & 96 & 64 & 6101 & 45 & 13,015.00 & 13,014.91 & 13,013.45 & 620 & 0 & 0 \\
Co-Product & 4 & 24 & 0 & 30 & 4,430.06 & 4,430.06 & 4,429.89 & 846 & 78 & $4.5\times10^{-14}$ \\
Serial-Inventory & 5 & 24 & 1 & 30 & 204.54 & 204.54 & 204.93 & 642 & 102 & $1.2\times10^{-13}$ \\
Case-Pack-Truck & 8 & 6 & 0 & 25 & 898.11 & 909.92 & 910.70 & 331 & 82 & $2.0\times10^{-13}$ \\
Hub-Consolidation & 130 & 52 & 6101 & 25 & 73,969.51 & 73,969.51 & 73,969.51 & 25 & 0 & 0 \\
 
\bottomrule
\end{tabular*}
\end{singlespace}
\figurenotes{``Uncorrected'' uses the tolerances of Proposition~\ref{prop:merge-tolerance}. ``Nonzero'' is the number of merges in which the intercept was lowered, and ``Largest'' is the largest amount by which an intercept was lowered.}
\end{table}

The corrected and the uncorrected rule give the same bound to all digits in 15 of the 24 runs. In the other nine runs, the bounds differ by up to 1.3\% and in either direction: counted with the tie tolerance, the corrected rule is higher in four runs and lower in one. However, the corrections themselves are very small. Of the 22,617 merges in these runs, 465 require a correction, and the largest correction of an intercept is $8.6\times10^{-8}$. The differences between the bounds therefore do not measure the size of the corrections. They arise because a run under a time limit is sensitive to small perturbations: a change of $10^{-8}$ in one intercept can change the solution of a later stage problem and, with it, the states that the following iterations visit. The differences are of the size of the variation over seeds in Table~\ref{tab:x-seeds}.

\subsection{Settings of SDDmiP}\label{app:add-sddmip}

SDDmiP generates a rectified linear cut only if the stage value exceeds its approximation by more than a threshold of $0.1$ (Appendix~\ref{app:methods-caveats}). We examine the effect of this threshold in four experiments.

\emph{Threshold.} Table~\ref{tab:x-sweep} reports the relative gap after 90 seconds on ten one-item instances with a known optimal value, for \textnormal{LCD}, \textnormal{LCD-Iso}, and SDDmiP with four thresholds. With the threshold $10^{-1}$, SDDmiP has the settings of Section~\ref{sec:experiments}. With the thresholds $10^{-2}$, $10^{-3}$, and $10^{-5}$, we set the tolerance of its bundle method to $10^{-3}$, $10^{-4}$, and $10^{-4}$, respectively, and its stage gap to $10^{-9}$. Table~\ref{tab:x-targets} counts the instances on which each method reaches a given relative gap within the 90 seconds.

\begin{table}[H]
\centering
\caption{Relative Gap After 90 Seconds of \textnormal{LCD}, \textnormal{LCD-Iso}, and SDDmiP With Four Thresholds}\label{tab:x-sweep}
\begin{singlespace}\footnotesize
\setlength{\tabcolsep}{3.5pt}
\begin{tabular*}{\textwidth}{@{\extracolsep{\fill}}lrr cc cccc@{}}
\toprule
 & & & & & \multicolumn{4}{c}{SDDmiP with threshold} \\
\cmidrule(l){6-9}
Instance & $T$ & $V^\star$ & \textnormal{LCD} & \textnormal{LCD-Iso} & $10^{-1}$ & $10^{-2}$ & $10^{-3}$ & $10^{-5}$ \\
\midrule
Case-Pack & 10 & 702.70 & $5.1\times10^{-7}$ & $5.0\times10^{-7}$ & $4.5\times10^{-3}$ & $5.5\times10^{-3}$ & $6.2\times10^{-3}$ & $7.3\times10^{-3}$ \\
Case-Pack & 24 & 1,613.98 & $1.2\times10^{-6}$ & $1.0\times10^{-5}$ & $1.6\times10^{-2}$ & $1.7\times10^{-2}$ & $1.7\times10^{-2}$ & $1.7\times10^{-2}$ \\
Workforce-Cohort & 10 & 468.67 & $4.0\times10^{-7}$ & $4.2\times10^{-7}$ & $2.8\times10^{-4}$ & $1.3\times10^{-4}$ & $5.2\times10^{-5}$ & $5.1\times10^{-5}$ \\
Workforce-Cohort & 24 & 1,053.33 & $1.1\times10^{-6}$ & $1.1\times10^{-6}$ & $2.4\times10^{-4}$ & $8.1\times10^{-7}$ & $2.3\times10^{-7}$ & $1.7\times10^{-9}$ \\
Co-Product & 10 & 168.26 & $5.5\times10^{-7}$ & $7.9\times10^{-5}$ & $2.5\times10^{-3}$ & $3.3\times10^{-3}$ & $1.6\times10^{-2}$ & $1.1\times10^{-2}$ \\
Co-Product & 24 & 412.22 & $7.5\times10^{-5}$ & $1.5\times10^{-5}$ & $3.1\times10^{-2}$ & $3.7\times10^{-2}$ & $4.0\times10^{-2}$ & $3.9\times10^{-2}$ \\
Bike-Sharing & 10 & 50.46 & $1.2\times10^{-1}$ & $2.8\times10^{-2}$ & $6.4\times10^{-2}$ & $6.5\times10^{-2}$ & $9.6\times10^{-2}$ & $1.1\times10^{-1}$ \\
Bike-Sharing & 24 & 219.77 & $1.4\times10^{-1}$ & $1.9\times10^{-1}$ & $1.2\times10^{-1}$ & $1.4\times10^{-1}$ & $1.3\times10^{-1}$ & $1.5\times10^{-1}$ \\
Case-Pack-Truck & 10 & 887.74 & $1.5\times10^{-2}$ & $3.1\times10^{-2}$ & $3.8\times10^{-2}$ & $4.6\times10^{-2}$ & $4.5\times10^{-2}$ & $4.6\times10^{-2}$ \\
Case-Pack-Truck & 24 & 2,058.71 & $8.0\times10^{-2}$ & $7.4\times10^{-2}$ & $9.4\times10^{-2}$ & $9.8\times10^{-2}$ & $1.0\times10^{-1}$ & $1.0\times10^{-1}$ \\
 
\bottomrule
\end{tabular*}
\end{singlespace}
\figurenote{All instances have one item, cohort, recipe, or station ($n=2$) and seed 0.}
\end{table}

\begin{table}[H]
\centering
\caption{Number of the Ten Instances of Table~\ref{tab:x-sweep} on Which a Relative Gap Is Reached Within 90 Seconds}\label{tab:x-targets}
\begin{singlespace}\footnotesize
\begin{tabular}{@{}lrrrr@{}}
\toprule
 & \multicolumn{4}{c}{Relative gap of at most} \\
\cmidrule(l){2-5}
Method & $10^{-2}$ & $10^{-3}$ & $10^{-4}$ & $10^{-5}$ \\
\midrule
\textnormal{LCD} & 6 & 6 & 6 & 5 \\
\textnormal{LCD-Iso} & 6 & 6 & 6 & 3 \\
SDDmiP, threshold $10^{-1}$ & 4 & 2 & 0 & 0 \\
SDDmiP, threshold $10^{-2}$ & 4 & 2 & 1 & 1 \\
SDDmiP, threshold $10^{-3}$ & 3 & 2 & 2 & 1 \\
SDDmiP, threshold $10^{-5}$ & 3 & 2 & 2 & 1 \\
SDDmiP, best threshold per instance & 4 & 2 & 2 & 1 \\
 
\bottomrule
\end{tabular}
\end{singlespace}
\figurenote{The last row takes, for each instance, the threshold with which SDDmiP has the smallest gap after 90 seconds.}
\end{table}

On the 24-stage Workforce-Cohort instance, the gap of SDDmiP falls from $2.4\times10^{-4}$ with the threshold $10^{-1}$ to $8.1\times10^{-7}$, $2.3\times10^{-7}$, and $1.7\times10^{-9}$ with the three smaller thresholds. These values are below the gap of \textnormal{LCD}, $1.1\times10^{-6}$, whose bound stops at its safety margin. On the 10-stage Workforce-Cohort instance, the gap of SDDmiP falls from $2.8\times10^{-4}$ to $5.1\times10^{-5}$. On the other eight instances, the gap of SDDmiP is larger with each of the smaller thresholds than with $10^{-1}$. With a smaller threshold, SDDmiP generates rectified linear cuts at more states and completes fewer iterations on these instances (e.g., 107, 55, 38, and 30 iterations with the four thresholds on the 10-stage Co-Product instance). Counted on the bounds, \textnormal{LCD} outperforms SDDmiP with the best threshold of each instance on six instances, ties on the two Workforce-Cohort instances, and loses on the two Bike-Sharing instances.

\emph{Costs multiplied by 100.} As a second test, we multiply all costs of six instances by 100, which makes the threshold one hundred times smaller in relation to the objective, and we run \textnormal{LCD} and SDDmiP for 60 seconds at both scales (Table~\ref{tab:x-scale}). On the Workforce-Cohort instance, the gap of SDDmiP with the threshold $10^{-1}$ falls from $2.8\times10^{-4}$ to $1.6\times10^{-5}$. On the other five instances, it does not fall.

\begin{table}[H]
\centering
\caption{Relative Gap After 60 Seconds With the Original Costs ($1\times$) and With All Costs Multiplied by 100 ($100\times$)}\label{tab:x-scale}
\begin{singlespace}\footnotesize
\setlength{\tabcolsep}{3.5pt}
\begin{tabular*}{\textwidth}{@{\extracolsep{\fill}}lrr cc cc cc@{}}
\toprule
 & & & \multicolumn{2}{c}{\textnormal{LCD}} & \multicolumn{2}{c}{SDDmiP, threshold $10^{-1}$} & \multicolumn{2}{c}{SDDmiP, threshold $10^{-3}$} \\
\cmidrule(lr){4-5}\cmidrule(lr){6-7}\cmidrule(l){8-9}
Instance & $T$ & $V^\star$ & $1\times$ & $100\times$ & $1\times$ & $100\times$ & $1\times$ & $100\times$ \\
\midrule
Case-Pack & 10 & 702.70 & $5.1\times10^{-7}$ & $5.0\times10^{-7}$ & $7.4\times10^{-3}$ & $8.4\times10^{-3}$ & $7.4\times10^{-3}$ & $8.3\times10^{-3}$ \\
Workforce-Cohort & 10 & 468.67 & $4.0\times10^{-7}$ & $3.8\times10^{-7}$ & $2.8\times10^{-4}$ & $1.6\times10^{-5}$ & $5.2\times10^{-5}$ & $2.6\times10^{-4}$ \\
Co-Product & 10 & 168.26 & $2.1\times10^{-5}$ & $5.0\times10^{-7}$ & $3.7\times10^{-3}$ & $2.2\times10^{-2}$ & $2.3\times10^{-2}$ & $2.2\times10^{-2}$ \\
Bike-Sharing & 10 & 50.46 & $1.4\times10^{-1}$ & $1.5\times10^{-1}$ & $7.5\times10^{-2}$ & $1.0\times10^{-1}$ & $1.3\times10^{-1}$ & $1.5\times10^{-1}$ \\
Case-Pack-Truck & 10 & 887.74 & $1.7\times10^{-2}$ & $2.0\times10^{-2}$ & $4.7\times10^{-2}$ & $5.1\times10^{-2}$ & $5.7\times10^{-2}$ & $5.7\times10^{-2}$ \\
Co-Product & 24 & 412.22 & $5.4\times10^{-4}$ & $1.1\times10^{-3}$ & $3.8\times10^{-2}$ & $4.0\times10^{-2}$ & $4.2\times10^{-2}$ & $4.6\times10^{-2}$ \\
 
\bottomrule
\end{tabular*}
\end{singlespace}
\figurenote{$V^\star$ is the optimal value with the original costs. The thresholds keep their absolute values at both scales.}
\end{table}

The two experiments agree. The threshold of $0.1$ limits the accuracy of SDDmiP on the Workforce-Cohort instances, on which a smaller threshold lowers its gap within 90 seconds, by a factor of five at 10 stages and to the level of our methods at 24 stages, and it does not limit SDDmiP on the other instances within 60 or 90 seconds.

\emph{Long-horizon runs with threshold $10^{-3}$.} To examine whether a smaller threshold improves SDDmiP over long horizons, we ran SDDmiP with a threshold of $10^{-3}$, a bundle tolerance of $10^{-4}$, and a stage gap of $10^{-9}$ across all 24 Case-Pack, Workforce-Cohort, and Co-Product instances of Table~\ref{tab:d-long-one} and four Bike-Sharing instances, under time limits of 30 minutes ($T=10$), 45 minutes ($T=24$), and 60 minutes ($T=48$), stopping early if the relative gap reached $5\times10^{-7}$ (Table~\ref{tab:x-long-sddmip}). Across the 24 Case-Pack, Workforce-Cohort, and Co-Product instances, calibrated SDDmiP reaches a relative gap of at most $5\times10^{-6}$ on five instances (compared with 22 for \textnormal{LCD-Iso}), and \textnormal{LCD-Iso} certifies a higher bound on 13 instances and ties on 11 (16 higher, 11 tied, and 1 lower across all 28 instances). On the Workforce-Cohort instances, lowering the threshold reduces the gaps of SDDmiP to between $1.7\times10^{-7}$ and $2.0\times10^{-4}$. On the 24- and 48-stage Case-Pack instances and the 48-stage Co-Product instances, however, SDDmiP generates between 1,974 and 6,452 rectified linear cuts, completes fewer iterations, and ends with relative gaps between $3.5\times10^{-3}$ and $2.4\times10^{-2}$ (higher than at threshold $0.1$ at 48 stages).

\begin{table}[H]
\centering
\caption{Relative Gaps in Long-Horizon Runs of Calibrated SDDmiP (Threshold $10^{-3}$, Stage Gap $10^{-9}$)}\label{tab:x-long-sddmip}
\begin{singlespace}\scriptsize
\setlength{\tabcolsep}{3.2pt}
\begin{tabular*}{\textwidth}{@{\extracolsep{\fill}}lrrr ccc rr@{}}
\toprule
 & & & & \multicolumn{3}{c}{Relative gap} & \multicolumn{2}{c}{SDDmiP ($10^{-3}$)} \\
\cmidrule(lr){5-7}\cmidrule(l){8-9}
Family & $T$ & Seed & $V^\star$ & \textnormal{LCD-Iso} & SDDmiP ($0.1$) & SDDmiP ($10^{-3}$) & Iter. & ReLU \\
\midrule
Case-Pack & 10 & 0 & 702.70 & $5.0\times10^{-7}$ & $3.5\times10^{-4}$ & $1.1\times10^{-4}$ & 2,425 & 703 \\
Case-Pack & 10 & 1 & 609.82 & $5.2\times10^{-7}$ & $2.4\times10^{-4}$ & $1.1\times10^{-6}$ & 2,419 & 722 \\
Case-Pack & 24 & 0 & 1,613.98 & $1.2\times10^{-6}$ & $1.6\times10^{-3}$ & $3.5\times10^{-3}$ & 110 & 1,974 \\
Case-Pack & 24 & 1 & 1,485.18 & $1.7\times10^{-6}$ & $2.6\times10^{-3}$ & $8.2\times10^{-3}$ & 103 & 1,989 \\
Case-Pack & 24 & 2 & 978.86 & $1.2\times10^{-6}$ & $4.5\times10^{-3}$ & $9.6\times10^{-3}$ & 103 & 2,400 \\
Case-Pack & 24 & 3 & 862.35 & $1.2\times10^{-6}$ & $4.2\times10^{-3}$ & $2.4\times10^{-2}$ & 106 & 2,118 \\
Case-Pack & 48 & 0 & 3,247.65 & $2.4\times10^{-6}$ & $4.5\times10^{-3}$ & $1.3\times10^{-2}$ & 111 & 2,410 \\
Case-Pack & 48 & 1 & 3,022.50 & $4.8\times10^{-6}$ & $6.9\times10^{-3}$ & $2.0\times10^{-2}$ & 115 & 2,793 \\
\addlinespace[2pt]
Workforce-Cohort & 10 & 0 & 468.67 & $4.2\times10^{-7}$ & $2.8\times10^{-4}$ & $5.2\times10^{-5}$ & 9,757 & 45 \\
Workforce-Cohort & 10 & 1 & 347.89 & $3.7\times10^{-7}$ & $2.9\times10^{-4}$ & $1.7\times10^{-4}$ & 4,538 & 21 \\
Workforce-Cohort & 24 & 0 & 1,053.33 & $1.1\times10^{-6}$ & $2.4\times10^{-4}$ & $2.3\times10^{-7}$ & 41 & 260 \\
Workforce-Cohort & 24 & 1 & 795.89 & $1.2\times10^{-6}$ & $9.4\times10^{-4}$ & $2.0\times10^{-4}$ & 4,625 & 158 \\
Workforce-Cohort & 24 & 2 & 801.07 & $1.1\times10^{-6}$ & $3.3\times10^{-4}$ & $5.7\times10^{-7}$ & 3,698 & 135 \\
Workforce-Cohort & 24 & 3 & 1,433.69 & $1.1\times10^{-6}$ & $4.6\times10^{-4}$ & $1.7\times10^{-7}$ & 51 & 465 \\
Workforce-Cohort & 48 & 0 & 2,113.11 & $2.4\times10^{-6}$ & $5.8\times10^{-4}$ & $5.8\times10^{-5}$ & 4,753 & 1,878 \\
Workforce-Cohort & 48 & 1 & 1,602.56 & $2.4\times10^{-6}$ & $7.7\times10^{-4}$ & $9.2\times10^{-5}$ & 4,334 & 579 \\
\addlinespace[2pt]
Co-Product & 10 & 0 & 168.26 & $5.5\times10^{-7}$ & $1.0\times10^{-3}$ & $1.6\times10^{-6}$ & 3,390 & 490 \\
Co-Product & 10 & 1 & 120.28 & $6.0\times10^{-7}$ & $1.0\times10^{-3}$ & $8.8\times10^{-5}$ & 3,201 & 514 \\
Co-Product & 24 & 0 & 412.22 & $1.3\times10^{-6}$ & $9.7\times10^{-4}$ & $2.6\times10^{-5}$ & 990 & 2,183 \\
Co-Product & 24 & 1 & 306.79 & $1.6\times10^{-6}$ & $1.6\times10^{-3}$ & $7.8\times10^{-6}$ & 474 & 3,042 \\
Co-Product & 24 & 2 & 298.14 & $1.9\times10^{-6}$ & $1.3\times10^{-3}$ & $1.5\times10^{-4}$ & 647 & 2,603 \\
Co-Product & 24 & 3 & 305.43 & $9.1\times10^{-6}$ & $1.7\times10^{-3}$ & $6.4\times10^{-4}$ & 331 & 3,280 \\
Co-Product & 48 & 0 & 828.74 & $4.1\times10^{-6}$ & $1.6\times10^{-3}$ & $7.5\times10^{-3}$ & 182 & 5,761 \\
Co-Product & 48 & 1 & 639.82 & $9.0\times10^{-6}$ & $1.9\times10^{-3}$ & $1.1\times10^{-2}$ & 153 & 6,452 \\
\addlinespace[2pt]
Bike-Sharing & 10 & 0 & 50.46 & $2.3\times10^{-3}$ & $3.6\times10^{-3}$ & $9.9\times10^{-4}$ & 227 & 1,122 \\
Bike-Sharing & 10 & 1 & 20.50 & $4.2\times10^{-3}$ & $9.9\times10^{-3}$ & $4.9\times10^{-3}$ & 184 & 1,201 \\
Bike-Sharing & 24 & 0 & 219.77 & $8.4\times10^{-4}$ & $3.9\times10^{-3}$ & $2.7\times10^{-2}$ & 115 & 1,995 \\
Bike-Sharing & 24 & 1 & 95.16 & $1.6\times10^{-3}$ & $1.0\times10^{-2}$ & $3.3\times10^{-2}$ & 105 & 1,955 \\
 
\bottomrule
\end{tabular*}
\end{singlespace}
\figurenotes{Columns ``\textnormal{LCD-Iso}'' and ``SDDmiP ($0.1$)'' are from Table~\ref{tab:d-long-one}. ``SDDmiP ($10^{-3}$)'' uses a cut-improvement threshold of $10^{-3}$, a bundle tolerance of $10^{-4}$, and a stage gap of $10^{-9}$, with time limits of 30 minutes ($T=10$), 45 minutes ($T=24$), and 60 minutes ($T=48$).}
\end{table}

\emph{Stage gap and thresholds at 45 seconds.} Table~\ref{tab:x-settings} reports the bounds after 45 seconds of \textnormal{LCD} with stage gaps of $10^{-9}$ and $10^{-4}$ and of SDDmiP with four settings: stage gaps of $10^{-4}$ (as in Section~\ref{sec:experiments}), $10^{-6}$, and $10^{-9}$ with the threshold $0.1$, and a stage gap of $10^{-9}$ with a threshold of $10^{-6}$ and a bundle tolerance of $10^{-4}$.

\begin{table}[H]
\centering
\caption{Certified Lower Bounds After 45 Seconds of \textnormal{LCD} and SDDmiP With Different Stage Gaps and Thresholds}\label{tab:x-settings}
\begin{singlespace}\scriptsize
\setlength{\tabcolsep}{3pt}
\begin{tabular*}{\textwidth}{@{\extracolsep{\fill}}lrrr rr rrrr@{}}
\toprule
 & & & & \multicolumn{2}{c}{\textnormal{LCD}, stage gap} & \multicolumn{3}{c}{SDDmiP, threshold $0.1$, stage gap} & SDDmiP, \\
\cmidrule(lr){5-6}\cmidrule(lr){7-9}
Instance & $n$ & $T$ & $V^\star$ & $10^{-9}$ & $10^{-4}$ & $10^{-4}$ & $10^{-6}$ & $10^{-9}$ & threshold $10^{-6}$ \\
\midrule
Co-Product & 2 & 10 & 168.26 & 168.26 & 168.24 & 163.23 & 163.03 & 163.63 & 161.67 \\
Co-Product & 2 & 24 & 412.22 & 406.59 & 406.92 & 392.48 & 391.12 & 391.12 & 390.58 \\
Workforce-Cohort & 2 & 24 & 1,053.33 & 1,053.33 & 1,053.31 & 1,052.84 & 1,052.06 & 1,051.17 & 1,046.28 \\
Bike-Sharing & 2 & 24 & 219.77 & 185.40 & 185.58 & 185.52 & 184.56 & 182.45 & 181.10 \\
Case-Pack-Truck & 2 & 10 & 887.74 & 851.23 & 846.17 & 827.14 & 822.08 & 827.75 & 823.25 \\
Case-Pack-Truck & 2 & 24 & 2,058.71 & 1,850.58 & 1,848.84 & 1,842.46 & 1,838.66 & 1,840.33 & 1,832.89 \\
Serial-Inventory & 5 & 24 & 344.76 & 200.25 & 177.25 & 196.47 & 200.67 & 208.43 & 180.20 \\
Case-Pack-Truck & 8 & 6 & 1,014.83 & 900.00 & 911.18 & 894.50 & 894.51 & 894.51 & 894.51 \\
Rotate & 12 & 24 & --- & 685.88 & 685.88 & 661.98 & 662.04 & 661.95 & 662.07 \\
Shear & 12 & 24 & --- & 683.88 & 683.88 & 513.04 & 512.87 & 512.87 & 517.35 \\
 
\bottomrule
\end{tabular*}
\end{singlespace}
\figurenote{The last column has a stage gap of $10^{-9}$ and a bundle tolerance of $10^{-4}$. The Serial-Inventory instance has seed 1, and the other instances have seed 0 (6101 for Rotate and Shear).}
\end{table}

\textnormal{LCD} outperforms SDDmiP on nine of the ten instances at each of the three stage gaps of SDDmiP, on all ten with the threshold $10^{-6}$, and on eight when the best of the four settings of SDDmiP is taken for each instance. The stage gap has little effect on either method: with a stage gap of $10^{-4}$, the bound of \textnormal{LCD} is higher on three instances, lower on three, and tied on four. Lowering the threshold of SDDmiP to $10^{-6}$ lowers its bound on seven instances and raises it on two, which agrees with Table~\ref{tab:x-sweep}.

\subsection{The Lattice Within \textnormal{LCD}}\label{app:add-lattice}

The long runs of Section~\ref{sec:exp-long} isolate the effect of the lattice for \textnormal{LCD-Iso}. Here, we isolate it for \textnormal{LCD}, by running \textnormal{LCD} with the trivial lattice $\Gamma_t=\Z^n$ and all other components unchanged.

\emph{Lattice on and off at 45 seconds.} Table~\ref{tab:x-ablation} reports the bounds after 45 seconds of \textnormal{LCD-Iso} and \textnormal{LCD}, each without and with the lattice, on 16 instances.

\begin{table}[H]
\centering
\caption{Certified Lower Bounds After 45 Seconds of \textnormal{LCD-Iso} and \textnormal{LCD} Without and With the Lattice}\label{tab:x-ablation}
\begin{singlespace}\footnotesize
\setlength{\tabcolsep}{3.5pt}
\begin{tabular*}{\textwidth}{@{\extracolsep{\fill}}lrrrrr rr rr@{}}
\toprule
 & & & & & & \multicolumn{2}{c}{\textnormal{LCD-Iso}} & \multicolumn{2}{c}{\textnormal{LCD}} \\
\cmidrule(lr){7-8}\cmidrule(l){9-10}
Instance & $n$ & $T$ & Seed & Cosets & $V^\star$ & No lattice & Lattice & No lattice & Lattice \\
\midrule
Case-Pack & 2 & 24 & 0 & 4 & 1,613.98 & 1,583.88 & 1,613.48 & 1,591.83 & 1,613.98 \\
Workforce-Cohort & 2 & 24 & 0 & 6 & 1,053.33 & 1,053.33 & 1,053.33 & 1,053.33 & 1,053.33 \\
Co-Product & 2 & 24 & 0 & 5 & 412.22 & 400.89 & 412.14 & 407.17 & 411.34 \\
Bike-Sharing & 2 & 24 & 0 & 36 & 219.77 & 182.31 & 178.86 & 193.04 & 187.21 \\
Co-Product & 4 & 24 & 0 & 35 & --- & 4,431.58 & 4,419.32 & 4,442.88 & 4,431.21 \\
Bike-Sharing & 3 & 24 & 0 & 64 & --- & 482.43 & 482.16 & 482.88 & 482.68 \\
Serial-Inventory & 5 & 24 & 0 & 27 & 270.02 & 172.19 & 162.97 & 180.79 & 181.49 \\
Serial-Inventory & 5 & 24 & 1 & 27 & 344.76 & 242.73 & 186.49 & 237.49 & 216.14 \\
Case-Pack-Truck & 2 & 10 & 0 & 4 & 887.74 & 850.87 & 851.31 & 868.16 & 864.54 \\
Case-Pack-Truck & 2 & 24 & 0 & 4 & 2,058.71 & 1,816.78 & 1,818.82 & 1,896.19 & 1,877.07 \\
Case-Pack-Truck & 4 & 6 & 0 & 48 & 446.44 & 414.48 & 368.31 & 424.56 & 421.83 \\
Workforce-Cohort & 6 & 6 & 0 & 600 & 3,879.89 & 3,879.89 & 3,879.89 & 3,879.89 & 3,879.89 \\
Shift & 12 & 24 & 6101 & 36 & --- & 751.00 & 752.08 & 752.46 & 759.17 \\
Rotate & 12 & 24 & 6101 & 36 & --- & 661.45 & 676.40 & 677.00 & 692.80 \\
Shear & 12 & 24 & 6101 & 64 & --- & 484.55 & 474.71 & 684.16 & 683.88 \\
Odd-Cycle & 96 & 64 & 6101 & 4 & --- & 12,433.96 & 12,436.24 & 12,829.96 & 12,829.40 \\
 
\bottomrule
\end{tabular*}
\end{singlespace}
\figurenote{``Cosets'' is the number of cosets per stage of the lattice. ``No lattice'' is the same algorithm with $\Gamma_t=\Z^n$.}
\end{table}

With the lattice, the bound of \textnormal{LCD} is higher on 5 instances, tied on 3, and lower on 8. The five instances are the one-item Case-Pack instance, on which the lattice raises the bound from 1,591.83 to the optimum, the one-recipe Co-Product instance, one of the two Serial-Inventory instances, Shift, and Rotate. The eight instances are the three Case-Pack-Truck instances, the two Bike-Sharing instances, the two-recipe Co-Product instance, the other Serial-Inventory instance, and Shear. For \textnormal{LCD-Iso}, the corresponding counts are 7, 2, and 7. The comparison of the two algorithms does not depend on the lattice: \textnormal{LCD} outperforms \textnormal{LCD-Iso} on 13 instances, ties on 2, and loses on 1, both with and without the lattice. Therefore, at this time limit, the advantage of \textnormal{LCD} over \textnormal{LCD-Iso} comes from the components in which the two algorithms differ, namely, the directional and pooled closing cuts and the implementation of Section~\ref{sec:algorithm}, and not from the lattice. On all five instances on which the lattice raises the bound of \textnormal{LCD}, lattice cuts alone are exact or nearly exact (Table~\ref{tab:d-class}). The converse does not hold, because the representation loss of lattice cuts is also zero on the one-station Bike-Sharing instance and on the other Serial-Inventory instance.

\emph{Long-budget crossover up to 1,800 seconds.} Table~\ref{tab:x-crossover} tracks the certified bounds of \textnormal{LCD} with and without the lattice at 45, 300, and 1,800 seconds across 25 instances: four 48-stage single-item instances, all 11 multiple-item instances of Panel~A of Table~\ref{tab:scope}, all four truck fixed-charge instances, five two-item six-stage instances, and Hub-Consolidation. Across the 11 multiple-item instances, the win--tie--loss count of \textnormal{LCD} with the lattice against \textnormal{LCD} without the lattice shifts from 3--1--7 at 180 seconds to 2--6--3 at 1,800 seconds as six instances reach the same bound while generating roughly half as many closing cuts (e.g., 886 against 2,101 and 1,176 against 2,117 on the Case-Pack instances with lead time 3, and 908 against 2,458 on the two-station Bike-Sharing instance). Across the four truck fixed-charge instances, the two 10-stage runs converge to the exact optimal value at 1,800 seconds (0--2--2 overall, compared with 0--0--4 at 45 and 180 seconds). On the 48-stage single-item Case-Pack and Co-Product instances, \textnormal{LCD} with the lattice reaches the optimum with 7 and 2,739 closing cuts against 5,453 and 9,498 without the lattice.

\begin{table}[H]
\centering
\caption{Certified Lower Bounds at 45, 300, and 1,800 Seconds and Closing Cuts Generated by \textnormal{LCD} With and Without the Lattice}\label{tab:x-crossover}
\begin{singlespace}\scriptsize
\setlength{\tabcolsep}{2.2pt}
\resizebox{\textwidth}{!}{%
\begin{tabular}{@{}lrrrr rr rr rr rr@{}}
\toprule
 & & & & & \multicolumn{2}{c}{Bound at 45 s} & \multicolumn{2}{c}{Bound at 300 s} & \multicolumn{2}{c}{Bound at 1,800 s} & \multicolumn{2}{c}{Closing cuts} \\
\cmidrule(lr){6-7}\cmidrule(lr){8-9}\cmidrule(lr){10-11}\cmidrule(l){12-13}
Instance & $n$ & $T$ & Cosets & $V^\star$ & Lattice & No lat. & Lattice & No lat. & Lattice & No lat. & Lattice & No lat. \\
\midrule
Case-Pack & 2 & 48 & 4 & 3,247.7 & 3,247.0 & 2,349.2 & 3,247.6 & 3,184.1 & 3,247.6 & 3,230.7 & 7 & 5,453 \\
Workforce-Cohort & 2 & 48 & 6 & 2,113.1 & 2,113.1 & 2,113.1 & 2,113.1 & 2,113.1 & 2,113.1 & 2,113.1 & 105 & 526 \\
Co-Product & 2 & 48 & 5 & 828.7 & 809.0 & 807.2 & 828.1 & 820.4 & 828.7 & 828.4 & 2,739 & 9,498 \\
Bike-Sharing & 2 & 48 & 36 & 412.1 & 342.5 & 342.3 & 354.5 & 376.7 & 366.5 & 393.1 & 2,138 & 6,263 \\
\addlinespace[2pt]
Case-Pack, lead time 3 & 4 & 10 & 8 & --- & 709.4 & 704.2 & 710.0 & 710.0 & 710.0 & 710.0 & 886 & 2,101 \\
Case-Pack, lead time 3 & 4 & 24 & 8 & --- & 1,417.6 & 1,401.8 & 1,467.5 & 1,453.4 & 1,479.2 & 1,472.3 & 1,176 & 2,117 \\
Case-Pack, lead time 3 & 4 & 24 & 8 & 557.4 & 499.3 & 490.9 & 552.1 & 549.4 & 557.3 & 557.4 & 1,678 & 3,468 \\
Co-Product & 4 & 10 & 35 & --- & 1,836.7 & 1,840.2 & 1,839.5 & 1,840.2 & 1,840.2 & 1,840.2 & 1,373 & 1,543 \\
Co-Product & 4 & 24 & 35 & --- & 4,431.8 & 4,442.4 & 4,437.1 & 4,444.0 & 4,443.6 & 4,444.0 & 3,156 & 3,214 \\
Co-Product & 4 & 24 & 35 & --- & 3,416.3 & 3,421.1 & 3,421.1 & 3,422.0 & 3,422.0 & 3,422.1 & 2,581 & 4,974 \\
Bike-Sharing & 3 & 10 & 64 & --- & 190.7 & 190.6 & 190.8 & 190.7 & 191.2 & 191.1 & 973 & 1,659 \\
Bike-Sharing & 3 & 24 & 64 & --- & 482.9 & 482.6 & 483.3 & 484.1 & 484.4 & 484.5 & 908 & 2,458 \\
Serial-Inventory & 5 & 10 & 27 & 90.8 & 80.5 & 89.3 & 90.8 & 90.8 & 90.8 & 90.8 & 98 & 341 \\
Serial-Inventory & 5 & 24 & 27 & 270.0 & 177.1 & 176.0 & 207.4 & 243.6 & 256.1 & 269.9 & 1,083 & 2,664 \\
Serial-Inventory & 5 & 24 & 27 & 344.8 & 204.9 & 209.1 & 276.0 & 305.3 & 326.9 & 339.9 & 1,334 & 2,342 \\
\addlinespace[2pt]
Case-Pack-Truck & 2 & 10 & 4 & 887.7 & 862.1 & 863.6 & 884.8 & 886.2 & 887.7 & 887.7 & 1,099 & 1,384 \\
Case-Pack-Truck & 2 & 10 & 6 & 773.6 & 732.7 & 750.8 & 772.0 & 772.1 & 773.6 & 773.6 & 1,160 & 1,951 \\
Case-Pack-Truck & 2 & 24 & 4 & 2,058.7 & 1,868.1 & 1,878.1 & 1,950.3 & 1,976.8 & 2,031.5 & 2,034.9 & 2,026 & 3,311 \\
Case-Pack-Truck & 2 & 24 & 6 & 1,900.3 & 1,627.7 & 1,670.6 & 1,790.0 & 1,818.5 & 1,868.4 & 1,879.2 & 1,991 & 3,537 \\
\addlinespace[2pt]
Case-Pack-Truck & 4 & 6 & 48 & 446.4 & 419.2 & 413.1 & 431.3 & 439.7 & 443.2 & 446.4 & 520 & 774 \\
Case-Pack & 4 & 6 & 48 & 353.5 & 343.0 & 345.8 & 353.3 & 353.5 & 353.5 & 353.5 & 187 & 321 \\
Workforce-Cohort & 4 & 6 & 80 & 3,685.3 & 3,685.3 & 3,685.3 & 3,685.3 & 3,685.3 & 3,685.3 & 3,685.3 & 14 & 21 \\
Co-Product & 4 & 6 & 35 & 905.4 & 901.8 & 905.4 & 905.2 & 905.4 & 905.4 & 905.4 & 688 & 132 \\
Bike-Sharing & 3 & 6 & 64 & 66.5 & 62.1 & 60.7 & 64.5 & 64.3 & 65.7 & 66.2 & 779 & 1,339 \\
\addlinespace[2pt]
Hub-Consolidation & 130 & 52 & 36 & --- & 73,953.3 & 73,933.3 & 73,970.9 & 73,966.7 & 74,083.9 & 74,794.6 & 21 & 36 \\
 
\bottomrule
\end{tabular}}
\end{singlespace}
\figurenotes{Runs stop at 1,800 seconds or upon reaching a known optimal value $V^\star$ to within the tolerance of Section~\ref{sec:exp-implementation}. ``No lat.'' is \textnormal{LCD} with $\Gamma_t=\Z^n$.}
\end{table}

\emph{Three seeds at 180 seconds.} Table~\ref{tab:x-seeds} reports the mean and the standard deviation of the bound over three seeds for six instances with multiple items, recipes, or stations and a limit of 180 seconds. SDDmiP is run with a threshold of $10^{-3}$, a bundle tolerance of $10^{-3}$, and a stage gap of $10^{-6}$.

\begin{table}[H]
\centering
\caption{Mean and Standard Deviation of the Certified Lower Bound Over Three Seeds, and Mean Policy Cost (Time Limit of 180 Seconds)}\label{tab:x-seeds}
\begin{singlespace}\scriptsize
\setlength{\tabcolsep}{2.5pt}
\resizebox{\textwidth}{!}{%
\begin{tabular}{@{}lrrr rrrr rrrr@{}}
\toprule
 & & & & \multicolumn{4}{c}{Certified lower bound, mean (standard deviation)} & \multicolumn{4}{c}{Mean policy cost} \\
\cmidrule(lr){5-8}\cmidrule(l){9-12}
Instance & $n$ & $T$ & Cosets & \textnormal{LCD} & \textnormal{LCD}, no lattice & SDDmiP & SDDiP & \textnormal{LCD} & \textnormal{LCD}, no lattice & SDDmiP & SDDiP \\
\midrule
Co-Product & 4 & 24 & 35 & 4,436.1 (0.1) & 4,444.0 (0.0) & 4,437.0 (1.0) & 4,430.6 (0.1) & 4,469 & 4,463 & 4,463 & 4,475 \\
Bike-Sharing & 3 & 24 & 64 & 483.3 (0.0) & 483.1 (0.1) & 481.3 (1.5) & 19.4 (2.0) & 526 & 521 & 523 & --- \\
Serial-Inventory & 5 & 24 & 27 & 262.7 (7.1) & 280.7 (6.0) & 250.8 (1.4) & 149.2 (5.9) & 426 & 402 & 411 & 595 \\
Case-Pack-Truck & 4 & 6 & 48 & 430.2 (0.5) & 437.7 (1.3) & 416.2 (1.7) & 376.1 (2.2) & 452 & 448 & 450 & 457 \\
Case-Pack-Truck & 8 & 6 & 3,072 & 947.7 (2.4) & 941.9 (4.1) & 903.8 (1.2) & 545.2 (5.2) & 1,043 & 1,043 & 1,046 & 1,231 \\
Hub-Consolidation & 130 & 52 & 36 & 73,962.7 (1.0) & 73,980.4 (53.5) & 70,481.1 (179.2) & 53,119.2 (705.3) & --- & --- & --- & --- \\
 
\bottomrule
\end{tabular}}
\end{singlespace}
\figurenotes{The three runs of a method on an instance use the random seeds 0, 1, and 2. The policy costs are means over the three seeds of the mean cost on 32 common sampled paths. A dash means that the evaluation did not finish within 30 seconds.}
\end{table}

The mean bound of \textnormal{LCD} without the lattice exceeds that of \textnormal{LCD} on four of the six instances. On three of them, namely, the Co-Product, Serial-Inventory, and two-item Case-Pack-Truck instances, which have 35, 27, and 48 cosets per stage, the difference is more than twice the larger of the two standard deviations. \textnormal{LCD} has the higher mean on the Bike-Sharing instance and on the four-item Case-Pack-Truck instance, by less than 1\%. Compared with SDDmiP, the mean bound of \textnormal{LCD} is higher on five instances and lower on the Co-Product instance (4,436.1 against 4,437.0, with a standard deviation of 1.0 for SDDmiP), and the mean bound of \textnormal{LCD} without the lattice is higher on all six. The standard deviation over the seeds is up to 2.7\% of the mean for \textnormal{LCD} (below 0.6\% on the 12 instances of Table~\ref{tab:c1} at 300 seconds), up to 2.1\% for \textnormal{LCD} without the lattice, up to 0.6\% for SDDmiP, and up to 10\% for SDDiP.

\emph{Many cosets at 60 seconds.} Table~\ref{tab:x-highcoset} considers ten instances with 1,296 or more cosets per stage and a limit of 60 seconds. In addition to \textnormal{LCD}, it reports \textnormal{LCD} with the lattice disabled, which means that no intercepts of individual cosets are computed, so that every cut has a single intercept.

\begin{table}[H]
\centering
\caption{Certified Lower Bounds After 60 Seconds on Ten Instances With Many Cosets}\label{tab:x-highcoset}
\begin{singlespace}\scriptsize
\setlength{\tabcolsep}{3pt}
\begin{tabular*}{\textwidth}{@{\extracolsep{\fill}}lrrr rr rrrr@{}}
\toprule
 & & & & \multicolumn{2}{c}{States visited by \textnormal{LCD}} & \multicolumn{4}{c}{Certified lower bound} \\
\cmidrule(lr){5-6}\cmidrule(l){7-10}
Instance & $n$ & $T$ & Cosets & All & Coset seen before & \textnormal{LCD} & \textnormal{LCD}, no lattice & SDDmiP & SDDiP \\
\midrule
Pooled-Capacity & 24 & 8 & 1,296 & 22 & 1 & 7,939.0 & 8,383.4 & 8,638.9 & 1,625.9 \\
Serial-Inventory, truck charge & 12 & 6 & 5,832 & 11 & 2 & 271.5 & 292.9 & 307.6 & 134.5 \\
Co-Product & 12 & 12 & 21,875 & 54 & 7 & 5,642.5 & 5,655.3 & 5,648.5 & 5,575.1 \\
Workforce-Cohort & 18 & 12 & 409,600 & 54 & 9 & 12,468.0 & 12,348.9 & 12,706.8 & 3,425.4 \\
Case-Pack-Truck & 16 & 12 & 7,962,624 & 14 & 0 & 693.1 & 818.7 & 2,032.5 & 415.0 \\
Case-Pack & 16 & 6 & 7,962,624 & 9 & 0 & 671.7 & 751.8 & 1,056.3 & 369.4 \\
Serial-Inventory & 24 & 6 & $7.2\times10^{7}$ & 11 & 2 & 417.5 & 492.5 & 683.5 & 221.5 \\
Bike-Sharing & 13 & 24 & $1.2\times10^{9}$ & 61 & 0 & 2,023.1 & 2,023.1 & 1,796.9 & 59.7 \\
Lead-Time-Pipeline & 48 & 12 & $6.9\times10^{12}$ & 23 & 0 & 1,549.8 & 1,611.1 & 1,859.4 & 862.9 \\
Airline-RM & 144 & 6 & $1.7\times10^{55}$ & 11 & 2 & -212,198.2 & -212,351.7 & -213,931.1 & -254,210.0 \\
 
\bottomrule
\end{tabular*}
\end{singlespace}
\figurenote{``Coset seen before'' is the number of visited states whose coset had been visited earlier in the run.}
\end{table}

Few of the visited states lie in a coset that the run has visited before, namely, between none and 9 of up to 61. Disabling the lattice raises the bound on seven instances, leaves it unchanged on the Bike-Sharing instance, and lowers it on the Workforce-Cohort and Airline-RM instances. With or without the lattice, \textnormal{LCD} loses to SDDmiP on most of these instances (on eight with the lattice and on seven without it) and outperforms SDDiP on all ten. Within the 60 seconds, \textnormal{LCD} completes between one and three iterations on these instances, and SDDmiP between 10 and 66.

\subsection{Reuse of Intercepts}\label{app:add-reuse}

A lattice cut pays off when an intercept that has been computed for a coset is used again at another state of this coset. We instrumented \textnormal{LCD} to count this event on ten instances with between 4 and 7,962,624 cosets per stage and a limit of 45 seconds (Table~\ref{tab:x-reuse}).

\begin{table}[H]
\centering
\caption{Reuse of Coset Intercepts in Runs of \textnormal{LCD} (Time Limit of 45 Seconds)}\label{tab:x-reuse}
\begin{singlespace}\footnotesize
\setlength{\tabcolsep}{3.5pt}
\begin{tabular*}{\textwidth}{@{\extracolsep{\fill}}lrrr rrrrrr@{}}
\toprule
Instance & $n$ & $T$ & Cosets & Iter. & States & Priced & Reuses & MIPs & Closing \\
\midrule
Case-Pack & 2 & 24 & 4 & 60 & 222 & 87 & 10,652 & 14,220 & 0 \\
Co-Product & 2 & 24 & 5 & 18 & 136 & 115 & 2,579 & 5,762 & 248 \\
Workforce-Cohort & 2 & 24 & 6 & 129 & 116 & 99 & 12,558 & 12,298 & 42 \\
Serial-Inventory & 5 & 24 & 27 & 4 & 104 & 189 & 69 & 2,578 & 27 \\
Rotate & 12 & 24 & 36 & 5 & 181 & 157 & 397 & 4,787 & 0 \\
Bike-Sharing & 3 & 24 & 64 & 10 & 435 & 355 & 868 & 3,614 & 29 \\
Pooled-Capacity & 24 & 8 & 1,296 & 1 & 13 & 25 & 3 & 387 & 0 \\
Sterile-Trays & 256 & 84 & 2,304 & 0 & 64 & 75 & 0 & 1,328 & 0 \\
Case-Pack-Truck & 8 & 6 & 3,072 & 2 & 21 & 57 & 3 & 665 & 18 \\
Case-Pack-Truck & 16 & 12 & 7,962,624 & 1 & 18 & 43 & 3 & 600 & 0 \\
 
\bottomrule
\end{tabular*}
\end{singlespace}
\figurenotes{``States'' is the number of distinct states visited, summed over the stages, and ``Priced'' is the number of pairs of a stage and a coset for which an intercept was computed. ``Reuses'' is the number of times that an intercept computed at one state gave the value of a cut at another visited state of the same coset or was taken from storage. ``MIPs'' is the number of MIPs solved for intercepts, and ``Closing'' is the number of closing cuts.}
\end{table}

On the six instances with at most 64 cosets per stage, the computed intercepts are used again between 69 and 12,558 times in a run. On the four instances with 1,296 cosets or more, they are used again at most three times in a run. This is the mechanism behind the second limit of Section~\ref{sec:exp-regimes}.

\subsection{Exact Evaluation of the Policies}\label{app:add-policy}

The policy costs of Section~\ref{sec:experiments} are means over sampled paths. For eight one-item instances, we evaluate the policies exactly instead. After a training run of 30 seconds (10 stages) or 45 seconds (24 stages), the cuts are frozen, and the expected cost $J$ of the policy is computed by a recursion over all states that the policy can reach, in which the stage problem of the policy is solved once for every reachable state and outcome. This covers all $3^{9}$ or $3^{23}$ scenario paths without sampling. Table~\ref{tab:x-policy} reports the relative gap of the bound, $(V^\star-L)/V^\star$, and the relative loss of the policy, $(J-V^\star)/V^\star$, for \textnormal{LCD}, \textnormal{LCD} without the lattice, and SDDmiP with the threshold $0.1$ and with a threshold of $10^{-3}$ (bundle tolerance $10^{-4}$, stage gap $10^{-9}$).

\begin{table}[H]
\centering
\caption{Relative Gap of the Bound and Relative Loss of the Exactly Evaluated Policy on Eight Instances}\label{tab:x-policy}
\begin{singlespace}\scriptsize
\setlength{\tabcolsep}{2.5pt}
\resizebox{\textwidth}{!}{%
\begin{tabular}{@{}lrrr cccc cccc@{}}
\toprule
 & & & & \multicolumn{4}{c}{Relative gap of the bound} & \multicolumn{4}{c}{Relative loss of the policy} \\
\cmidrule(lr){5-8}\cmidrule(l){9-12}
 & & & Limit & & \textnormal{LCD}, & SDDmiP, & SDDmiP, & & \textnormal{LCD}, & SDDmiP, & SDDmiP, \\
Instance & $T$ & $V^\star$ & (s) & \textnormal{LCD} & no lattice & $10^{-1}$ & $10^{-3}$ & \textnormal{LCD} & no lattice & $10^{-1}$ & $10^{-3}$ \\
\midrule
Case-Pack & 10 & 702.70 & 30 & $5.1\times10^{-7}$ & $2.3\times10^{-3}$ & $9.4\times10^{-3}$ & $1.0\times10^{-2}$ & 0 & $6.3\times10^{-4}$ & $2.9\times10^{-4}$ & $4.8\times10^{-4}$ \\
Workforce-Cohort & 10 & 468.67 & 30 & $4.0\times10^{-7}$ & $4.2\times10^{-7}$ & $2.8\times10^{-4}$ & $5.2\times10^{-5}$ & 0 & 0 & 0 & 0 \\
Co-Product & 10 & 168.26 & 30 & $3.7\times10^{-5}$ & $1.1\times10^{-3}$ & $2.3\times10^{-2}$ & $3.8\times10^{-2}$ & 0 & $5.9\times10^{-4}$ & $8.2\times10^{-5}$ & $1.5\times10^{-2}$ \\
Bike-Sharing & 10 & 50.46 & 30 & $1.8\times10^{-1}$ & $1.2\times10^{-1}$ & $9.8\times10^{-2}$ & $1.9\times10^{-1}$ & $1.3\times10^{-1}$ & $2.1\times10^{-2}$ & $1.4\times10^{-2}$ & $1.7\times10^{-2}$ \\
Case-Pack-Truck & 10 & 887.74 & 30 & $3.3\times10^{-2}$ & $2.6\times10^{-2}$ & $6.7\times10^{-2}$ & $7.7\times10^{-2}$ & $1.2\times10^{-2}$ & $2.9\times10^{-2}$ & $2.6\times10^{-2}$ & $2.9\times10^{-2}$ \\
Serial-Inventory & 10 & 90.82 & 30 & $2.0\times10^{-1}$ & $3.9\times10^{-2}$ & $1.6\times10^{-1}$ & $1.9\times10^{-1}$ & $1.5\times10^{-1}$ & $4.7\times10^{-1}$ & $7.6\times10^{-2}$ & $4.6\times10^{-2}$ \\
Case-Pack & 24 & 1,613.98 & 45 & $1.2\times10^{-6}$ & $1.4\times10^{-2}$ & $1.7\times10^{-2}$ & $1.8\times10^{-2}$ & 0 & $8.6\times10^{-4}$ & $8.4\times10^{-4}$ & $8.0\times10^{-4}$ \\
Workforce-Cohort & 24 & 1,053.33 & 45 & $1.1\times10^{-6}$ & $1.2\times10^{-6}$ & $2.4\times10^{-4}$ & $9.2\times10^{-5}$ & 0 & 0 & 0 & 0 \\
 
\bottomrule
\end{tabular}}
\end{singlespace}
\figurenotes{A loss of 0 means that the expected cost of the policy equals $V^\star$ to within $5\times10^{-7}$ in absolute terms. All instances have seed 0, and the Serial-Inventory instance has $n=5$.}
\end{table}

The policy of \textnormal{LCD} is optimal on five of the eight instances. On two of them, the Workforce-Cohort instances, the policies of all four methods are optimal. On the other three, the two Case-Pack instances and the Co-Product instance, only the policy of \textnormal{LCD} is optimal, and the policies of the other methods lose between $8.2\times10^{-5}$ and $1.5\times10^{-2}$ of the optimal value. On the Case-Pack-Truck instance, the policy of \textnormal{LCD} loses 1.2\% and those of the other methods lose between 2.6\% and 2.9\%. On the Bike-Sharing and Serial-Inventory instances, on which the bound of \textnormal{LCD} is still 18\% and 20\% below the optimum at the time limit, its policy is worse than that of SDDmiP: it loses 13\% and 15\%, against 1.4\% and 7.6\% for SDDmiP with the threshold $0.1$. Therefore, a policy of \textnormal{LCD} is as good as its bound: where the bound is tight, the policy is optimal, and where the bound is still loose after a short run, the policy can be worse than those of the other methods.

\subsection{Instances Without Batch Structure}\label{app:add-nobatch}

Our method is designed for models with batches. For completeness, we report what it does on 12 instances without batch structure, on which Algorithm~\ref{alg:lattice} returns the trivial lattice at every stage. Four of them are capacitated lot-sizing instances with binary setups, a shared capacity, bounded inventories, and lost sales, of the kind used by \citet{zou2019}, with 2, 3, 4, and 6 products. The others are Case-Pack-Truck, Serial-Inventory, and Bike-Sharing instances in which every batch size is one. On these instances, lattice cuts are affine cuts, \textnormal{LCD} reduces to affine cuts with directional and pooled closing cuts, and \textnormal{LCD-Iso} reduces to the control without the lattice of Section~\ref{sec:exp-methods}. Table~\ref{tab:x-nobatch} reports the bounds and policy costs after 30 seconds.

\begin{table}[H]
\centering
\caption{Certified Lower Bounds and Mean Policy Costs After 30 Seconds on 12 Instances Without Batch Structure}\label{tab:x-nobatch}
\begin{singlespace}\scriptsize
\setlength{\tabcolsep}{3pt}
\begin{tabular*}{\textwidth}{@{\extracolsep{\fill}}lrr rrrr rrrr@{}}
\toprule
 & & & \multicolumn{4}{c}{Certified lower bound} & \multicolumn{4}{c}{Mean policy cost} \\
\cmidrule(lr){4-7}\cmidrule(l){8-11}
Instance & $n$ & $T$ & \textnormal{LCD} & \textnormal{LCD-Iso} & SDDmiP & SDDiP & \textnormal{LCD} & \textnormal{LCD-Iso} & SDDmiP & SDDiP \\
\midrule
Lot-Sizing & 2 & 6 & 682.5 & 533.3 & 537.8 & 524.2 & 836 & 822 & 836 & 986 \\
Lot-Sizing & 3 & 8 & 533.4 & 380.3 & 473.8 & 561.7 & 1,216 & 1,410 & 1,455 & 1,294 \\
Lot-Sizing & 4 & 12 & 1,092.0 & 828.3 & 960.0 & 947.6 & 3,038 & 3,638 & 3,354 & 3,109 \\
Lot-Sizing & 6 & 16 & 1,463.3 & 1,462.9 & 1,448.4 & 1,588.5 & 7,180 & 7,801 & 6,839 & --- \\
Case-Pack-Truck, pack size 1 & 2 & 10 & 684.8 & 675.3 & 676.6 & 407.1 & 693 & 694 & 691 & 730 \\
Case-Pack-Truck, pack size 1 & 2 & 24 & 1,603.5 & 1,504.0 & 1,610.7 & 696.0 & 1,638 & 1,714 & 1,635 & 1,881 \\
Case-Pack-Truck, pack size 1 & 8 & 12 & 2,143.1 & 1,790.1 & 2,249.0 & 609.0 & 2,604 & 2,772 & 2,712 & 2,932 \\
Case-Pack-Truck, pack size 1 & 16 & 12 & 2,359.4 & 2,562.8 & 4,676.0 & 1,163.6 & --- & --- & --- & --- \\
Serial-Inventory, batch size 1 & 4 & 12 & 236.0 & 228.9 & 226.4 & 98.9 & 272 & 248 & 268 & 298 \\
Serial-Inventory, batch size 1 & 12 & 12 & 1,355.0 & 1,325.0 & 1,682.0 & 345.0 & 2,036 & 2,152 & 1,880 & 2,233 \\
Bike-Sharing, batch size 1 & 4 & 12 & 427.2 & 427.2 & 426.7 & 22.0 & 441 & 441 & 441 & 513 \\
Bike-Sharing, batch size 1 & 7 & 16 & 728.1 & 721.7 & 715.6 & 26.7 & 727 & 728 & 733 & --- \\
 
\bottomrule
\end{tabular*}
\end{singlespace}
\figurenotes{The lattice is trivial for both of our methods. Policy costs are means over 32 common sampled paths, and a dash means that the evaluation did not finish within 15 seconds.}
\end{table}

The bound of \textnormal{LCD} is higher than that of \textnormal{LCD-Iso} on 10 instances, higher than that of SDDiP on 10, and higher than that of SDDmiP on 8, and it is lower than that of SDDmiP on 4. The losses against SDDmiP include the instance with 16 state variables, on which the bound of \textnormal{LCD} is about half that of SDDmiP (2,359.4 against 4,676.0), and the mean relative difference of the two bounds over the 12 instances is 1.8\% in favor of SDDmiP. On the 11 instances on which the policies were evaluated, the policy of \textnormal{LCD} has the lower mean cost than that of SDDmiP on five and the higher on six, where two of these differences are smaller than 0.2. We also ran the variant of SDDmiP with regularized instead of normalized cuts (Reg-A in \citealt{bansal2026}) on these instances, and the counts on the bound are the same: \textnormal{LCD} outperforms it on 8 instances and loses on 4. Therefore, without batch structure, our algorithm is comparable with SDDmiP at this time limit and not better, and these instances are outside the regime for which we propose it.

\subsection{Numerical Checks of Two Statements}\label{app:add-theory}

\emph{Loss of affine cuts when the stock carries over.} Proposition~\ref{prop:loss} assumes a direction $d$ with $Ad=0$, that is, a change of the state that does not affect the next state. In the Case-Pack model, the stock carries over to the next stage, and the proposition does not apply. We therefore examine the 10-stage one-item Case-Pack instance with pack size 4 numerically. Let $e$ denote the unit vector of the stock. For every stage $t\ge 2$, every $q\in\{1,2,3\}$, and every state $s^0$ such that $s^0$, $s^0+qe$, and $s^0+4e$ belong to $D_t$, we compute from the exact cost-to-go function the quantity
\[
 \Delta=V_t(s^0+qe)-\bigl(1-\tfrac q4\bigr)V_t(s^0)-\tfrac q4\,V_t(s^0+4e).
\]
By the argument of Proposition~\ref{prop:loss}, every affine function $u\le V_t$ satisfies $V_t(s^0+qe)-u(s^0+qe)\ge\Delta$. We find $\Delta>0$ for some states at seven of the nine stages, and the largest value is $\Delta=2.545$, which is attained at stage 2 with $q=2$, where $V_t$ takes the values $582.698$, $577.243$, and $566.698$ at the three states. Therefore, affine cuts also lose a fixed amount at states between two batch multiples in a model in which $Ad\ne 0$.

\emph{Conclusion of Theorem~\ref{thm:optimal-period} on small instances.} Theorem~\ref{thm:optimal-period} states that the periodic part of the cost-to-go function vanishes exactly on the lattice $K_t$ and is positive elsewhere. We check this conclusion on four small instances by dynamic programming: Shift and Rotate with $n=6$ and $T=3$, and two three-stage models with a single state variable and batch sizes 6 and 4, in which the stock can be changed by whole batches at no cost and by single units at a cost, demand comes in whole batches, and every unit of mismatch is penalized. In each case, we compute the cost-to-go function of the second stage, subtract its linear trend, and compare the result with the lattice returned by Algorithm~\ref{alg:lattice}. On all four instances, the periodic part is zero on one coset of this lattice and positive on all other cosets (at least 5.0 on Shift and Rotate, and at least 14 and 18 on the two other models). On Shift and Rotate, the cost-to-go function is moreover a linear function of the state plus a function of its coset, with no deviation at any of the states checked, as Proposition~\ref{prop:phasevalue} predicts. We emphasize that these are checks of conclusions on small instances. They do not verify Assumption~\ref{assump:period} for the benchmark families of Section~\ref{sec:experiments}.

\clearpage
\phantomsection
\addcontentsline{toc}{section}{References}
\begin{singlespace}
\end{singlespace}

\end{document}